\documentclass[12pt,a4paper]{article}%{amsart}
\usepackage[english]{babel}
\usepackage{amsmath}
\usepackage{amssymb}
\usepackage{amscd}
\usepackage{amsthm}
\usepackage{amsfonts}
\usepackage{amsbsy}
\usepackage{enumerate}
\usepackage{graphicx}
\usepackage{color}
\usepackage{array}
\usepackage[dvips]{geometry}
\usepackage{epsfig,afterpage}
\usepackage[dvips]{psfrag}
\usepackage{psfrag}
\usepackage{graphicx,blindtext}
\usepackage[active]{srcltx}
\usepackage{hyphenat}
\usepackage{multicol}
\usepackage{multirow}
\usepackage{pdflscape}
\usepackage{morefloats}
\usepackage{url}
\usepackage{longtable}
\usepackage{lscape}
\newcommand{\QES}{\bf{Q}{\widehat{ES}}}
\newcommand{\QESA}{\bf{Q}{{\widehat{ES}(A)}}}
\newcommand{\QESB}{\bf{Q}{{\widehat{ES}(B)}}}
\newcommand{\QESC}{\bf{Q}{{\widehat{ES}(C)}}}

\newtheorem{theorem}{Theorem}
\newtheorem{lemma}{Lemma}
\newtheorem{definition}{Definition}
\newtheorem{proposition}{Proposition}
\newtheorem{remark}{Remark}

\newtheorem{corollary}{Corollary}

\begin{document}

\title{Phase portraits for quadratic systems possessing an infinite elliptic--saddle or an infinite nilpotent saddle}

\author{\textsc{Joan C. Art\'es}\\
        Departament de Matem\`atiques, Universitat Aut\`onoma de Barcelona\\
        08193, Barcelona, Spain\\
        E-mail: joancarles.artes@uab.cat\\
 \and
        \textsc{Marcos C. Mota}\\
        Instituto de Ci\^encias Matem\'aticas e de Computa\c{c}\~{a}o, \\ Universidade de S\~{a}o Paulo\\
        13566--590, S\~{a}o Carlos, S\~{a}o Paulo, Brazil\\
        E-mail: coutinhomotam@gmail.com\\
 \and
		\textsc{Alex C. Rezende}\\
		Departamento de Matem\'atica, Universidade Federal de S\~{a}o Carlos\\
		13565-905, S\~ao Carlos, S\~ao Paulo, Brazil\\
		E-mail: alexcr@ufscar.br}

\date{}

\maketitle

\begin{abstract}
\noindent This paper presents a global study of the class $\QES$ of all real quadratic polynomial differential systems possessing exactly one elemental infinite singular point and one triple infinite singular point, which is either an infinite nilpotent elliptic--saddle or a nilpotent saddle. This class can be divided into three different families, namely, $\QESA$ of phase portraits possessing three real finite singular points, $\QESB$ of phase portraits possessing one real and two complex finite singular points, and $\QESC$ of phase portraits possessing one real triple finite singular point. Here we provide the complete study of the geometry of these three families. Modulo the action of the affine group and time homotheties, families $\QESA$ and $\QESB$ are three--dimensional and family $\QESC$ is two--dimensional. We study the respective bifurcation diagrams of their closures with respect to specific normal forms, in subsets of real Euclidean spaces. The bifurcation diagram of family $\QESA$ (respectively, $\QESB$ and $\QESC$) yields 1274 (respectively, 89 and 14) subsets with 91 (respectively, 27 and 12) topologically distinct phase portraits for systems in the closure $\overline{\QESA}$ (respectively, $\overline{\QESB}$ and $\overline{\QESC}$) within the representatives of $\QESA$ (respectively, $\QESB$ and $\QESC$) given by a specific normal form.
\end{abstract}


%\runningheads{J.C. Art\'{e}s et al.}{Phase portraits for QS with an infinite elliptic--saddle or an infinite nilpotent saddle}

%\tableofcontents

\section{Introduction, brief review of the literature and statement of the results}\label{sec:int}

\indent Here we call \textit{quadratic differential systems}, or simply \textit{quadratic systems}, 
differential systems of the form
\begin{equation} 
\begin{array}{lcccl}
\dot{x}&=& p(x,y), \\
\dot{y}&=& q(x,y), \\
\end{array} \label{eq:qs} 
\end{equation}
where $p, q\in\mathbb{R}[x,y]$ verify $\max\{\deg(p),\deg(q)\}=2$. To such systems we can 
associate the quadratic vector field
\begin{equation} 
\xi=p\frac{\partial}{\partial x}+q\frac{\partial}{\partial y} \label{eq:qvf},
\end{equation} 
as well as the differential equation
\begin{equation} 
q\,dx-p\,dy=0.
\label{eq:de}
\end{equation} 
Along this paper we shall use indistinctly the expressions \textit{quadratic systems} and 
\textit{quadratic vector fields} to refer to either \eqref{eq:qs}, or \eqref{eq:qvf}, or \eqref{eq:de}.

The class of all quadratic differential systems is denoted by \textbf{QS}.

We can also write systems \eqref{eq:qs} as
\begin{equation} 
\begin{array}{lcccl}
\dot{x}&=&p_0+p_{1}(x,y)+p_{2}(x,y)\equiv p(x,y), \\
\dot{y}&=&q_0+q_{1}(x,y)+q_{2}(x,y)\equiv q(x,y), \\
\end{array} \label{2l1} 
\end{equation}    
where $p_i$ and $q_i$ are homogeneous polynomials of degree $i$ in the variables $x$ and $y$ 
with real coefficients and $p_{2}^2+q_{2}^2 \neq 0$.

Even after hundreds of studies on the topology of real planar quadratic vector fields, it is somewhat impossible at this point to fully characterize their phase portraits and try to topologically classify them (which is very common in applications) due to the large number of parameters involved. 

The main purpose of this paper is to present the study of the bifurcation diagrams of the class of quadratic systems possessing exactly one elemental infinite singular point and one triple infinite singular point, being an infinite nilpotent elliptic--saddle (which can be of three types: $\widehat{\!{1\choose 2}\!\!}\ PHP-E$, $\widehat{\!{1\choose 2}\!\!}\ H-E$, or $\widehat{\!{1\choose 2}\!\!}\ PEP-H$) or a nilpotent saddle $\widehat{\!{1\choose 2}\!\!}\ HHH-H$ (see \cite{Artes-Llibre-Schlomiuk-Vulpe-2021a} for details on this notation). We denote this class by $\QES$. A nilpotent singularity is a point where both eigenvalues are zero but the Jacobian matrix is nonzero.

Whenever one wants to study a specific family of differential systems sharing a common property, it is necessary to select one (or several) normal form which contains all the phase portraits sharing the desired property. However, except for a few trivial cases, it is impossible that the normal form does not contain other phase portraits, normally more degenerate than the cases under study. These other phase portraits are very important for understanding the bifurcations that occur within the chosen normal form. Therefore, we always check not only the family of systems with the desired properties, but also the clousure of the normal form which contains that family. That is, we examine the entire parameter space of the chosen normal form, whether or not it leads to the desired property. However, it is possible that a different normal form could have been chosen, in which case the generic elements of the family should be the same, but the elements in the border might not be. That is, some phase portraits in the border of one normal form could be common or not, with elements in the border of the second normal form.

It is well known that quadratic systems possess at most four real simple finite singular points and at most three pairs of infinite singular points. As our aim is to study \textbf{QS} possessing an infinite singular point of multiplicity three, formed by the coalescence of one finite singular point with one double infinite singular point, a quadratic differential system from the class $\overline{\QES}$ can have at most three simple real finite singular points and, in case it has total multiplicity 3 of finite singularities, it will have two pairs of infinite singular points, being one simple and the other one triple. So, inside the class $\overline{\QES}$ we must consider the following families:
\begin{itemize}
	\item $\QESA$: quadratic systems possessing three real finite singular points, either an infinite nilpotent elliptic--saddle or an infinite nilpotent saddle, and an elemental infinite singularity;
	\item $\QESB$: quadratic systems possessing one real and two complex finite singular points, either an infinite nilpotent elliptic--saddle or an infinite nilpotent saddle, and an elemental infinite singularity;
	\item $\QESC$: quadratic systems possessing one real triple finite singular point, either an infinite nilpotent elliptic--saddle or an infinite nilpotent saddle, and an elemental infinite singularity.
\end{itemize}

For our proposed study, we followed the pattern specified in \cite{Artes-Rezende-Oliveira-2015,Artes-Mota-Rezende-2021c} and, in order to avoid repeating technical common sections, we refer to the mentioned papers for more complete information.

All the phase portraits in this paper are drawn in the Poincar\'{e} disc (for its definition we refer to \cite{Dumortier-Llibre-Artes-2006,Artes-Rezende-Oliveira-2015}). In the sequel, we give the concept of \textit{graphics}, which play an important role in obtaining limit cycles when they arise, for example, from connection of separatrices.

A \textit{(nondegenerate) graphic} as defined in \cite{Dumortier-Roussarie-Rousseau-1994} is formed by a finite 
sequence of singular points $r_1,r_2,\ldots,r_n$ (with possible repetitions) and non--trivial connecting orbits $\gamma_i$ 
for $i=1,\ldots,n$ such that  $\gamma_i$ has $r_i$ as $\alpha$--limit set and $r_{i+1}$ as $\omega$--limit set for $i<n$ 
and $\gamma_n$ has $r_n$ as  $\alpha$--limit set and $r_{1}$ as $\omega$--limit set. Also normal orientations $n_j$ of the 
non--trivial orbits must be coherent  in the sense that if $\gamma_{j-1}$ has left--hand orientation then so does $\gamma_j$. 
A \textit{polycycle} is a graphic which has  a Poincar\'{e} return map. 

A \textit{degenerate graphic} is formed by a finite sequence of singular points $r_1,r_2,\ldots,r_n$ (with possible repetitions) 
and non--trivial connecting orbits and/or segments of curves of singular points $\gamma_i$ for $i=1,\ldots,n$ such that 
$\gamma_i$ has $r_i$ as $\alpha$--limit set and $r_{i+1}$ as $\omega$--limit set for $i<n$ and $\gamma_n$ has $r_n$ as 
$\alpha$--limit set and $r_{1}$ as $\omega$--limit set. Also normal orientations $n_j$ of the non--trivial orbits must be coherent 
in the sense that if $\gamma_{j-1}$ has left--hand orientation then so does $\gamma_j$. For more details, see 
\cite{Dumortier-Roussarie-Rousseau-1994}.

In \cite{Artes-Kooij-Llibre-1998} the authors proved the existence of 44 topologically different phase portraits for the structurally stable quadratic planar differential systems modulo limit cycles, also known as the codimension--zero quadratic systems. Roughly speaking, these systems are characterized by having all singularities, finite and infinite, simple, no separatrix connection, and where any nest of limit cycles counts as a single point with the stability of the outer limit cycle. 

In addition, in \cite{Artes-Llibre-Rezende-2018} the authors classified the structurally unstable quadratic systems of codimension one modulo limit cycles which have one and only one of the simplest structurally unstable objects: a saddle--node of multiplicity two (finite or infinite), a separatrix from one saddle point to another, or a separatrix forming a loop for a saddle point with its divergence nonzero. All the phase portraits of codimension one are split into four sets according to the possession of a structurally unstable element: \label{page:groups-cod1} 
\begin{itemize}
	\item[(A)] possessing a finite semi--elemental saddle--node;
	\item[(B)] possessing an infinite semi--elemental saddle--node $\overline{\!{0\choose2}\!\!}\ SN$;
	\item[(C)] possessing an infinite semi--elemental saddle--node $\overline{\!{1\choose1}\!\!}\ SN$; and
	\item[(D)] possessing a separatrix connection.
\end{itemize}
The study of the codimension--one systems was carried out during a period of approximately 20 years, and this study yielded at least 204 (and at most 211) topologically distinct phase portraits of codimension one modulo limit cycles. Some recent research (already at preprint level) showed two mistakes in that book and reduced (and confirmed) the number of cases to 202 (and a most 209).

The next step is to study the structurally unstable quadratic systems of codimension two, modulo limit cycles. The approach is the same as used in the previous two works \cite{Artes-Kooij-Llibre-1998,Artes-Llibre-Rezende-2018}. One starts by looking for all the potential topological phase portraits of codimension two, and then tries to realize all of them or show that some of them are impossible. So, it is also very convenient to have studied a bifurcation diagram that helps us to solve the realization problem. In many publications of this last type where families of phase portraits have been studied, it is quite common that the authors have missed one or several phase portraits, as we discuss in Appendix~\ref{ap:incomp-QES}. This may happen either because they have not interpreted correctly some of the bifurcation parts, or they have missed the existence of some nonalgebraic bifurcation, or there may exist some small ``island'' as they are described in Sec.~\ref{sec:islands-QESA}, \ref{sec:islands-QESB}, and \ref{sec:islands-QESC}. However, when examining all the potential topological phase portraits and systematically compiling error--free list, then there is no possibility of missing a realizable case. It is just a problem of finding examples of realization or producing irrefutable proofs of the impossibility of realization of phase portraits.

Research on codimension--two quadratic systems is already ongoing. In \cite{Artes-Oliveira-Rezende-2020b} the authors have considered set (AA) obtained by the existence of a cusp point, or two saddle--nodes or the coalescence of three finite singular points forming a semi--elemental singularity, yielding either a triple saddle, or a triple node. They obtained all the possible topological phase portraits of set (AA) and proved their realization. In their study, they got 34 new topologically distinct phase portraits in the Poincar\'e disc modulo limit cycles. Moreover, they proved the impossibility of one phase portrait among the $204$ phase portraits presented in \cite{Artes-Llibre-Rezende-2018}.

Moreover, the bifurcation diagram for the class of the quadratic systems possessing a finite saddle--node $\overline{sn}_{(2)}$ and an infinite saddle--node $\overline{\!{0\choose2}\!\!}\ SN$ was studied in \cite{Artes-Rezende-Oliveira-2014,Artes-Rezende-Oliveira-2015}, in which all the phase portraits obtained belong to the closure of set (AB). Also, in \cite{Artes-Mota-Rezende-2021b,Artes-Mota-Rezende-2021c} the authors studied the bifurcation diagram for the class of quadratic systems possessing a finite saddle--node $\overline{sn}_{(2)}$ and an infinite saddle--node $\overline{\!{1\choose1}\!\!}\ SN$ and all the phase portraits obtained belong to the closure of set (AC).

The topological classification of sets (AB) and (AC) was done in \cite{Artes-Mota-Rezende-2021d}. In this study, the authors obtained 71 topologically distinct phase portraits modulo limit cycles for the set (AB), and for the set (AC) they got 40 ones.

Consider now the set (BC), characterized by quadratic systems possessing two types of coalescence of singular points:
\begin{itemize}
	\item coalescence of two infinite elemental singular points; and
	\item coalescence of a finite elemental singular point with an infinite one.
\end{itemize}
In a near future we will present a paper that includes the study of the bifurcation diagram of quadratic systems with infinite saddle--nodes $\overline{\!{0\choose2}\!\!}\ SN$ and $\overline{\!{1\choose1}\!\!}\ SN$. 

Since here we want to study quadratic systems with exactly one elemental infinite singular point and one triple infinite singular point (in the sense that it is the coalescence of two infinite singularities plus a finite one), families $\QESA$ and $\QESB$ can be considered as codimension--two cases from the border of set (BC) and family $\QESC$ can be seeing as a codimension--four case from the border of set (BC).

In the normal form \eqref{eq:nf-QES-A}, see page \pageref{eq:nf-QES-A}, the class $\overline{\QESA}$ is partitioned into 1274 parts: 288 three--dimen\-sional ones, 573 two--dimensional ones, 351 one--dimensional ones, and 62 points. This partition is obtained by considering all the bifurcation surfaces of singularities, and bifurcation surfaces related to the presence of invariant straight lines, the presence of invariant parabolas, and connections of separatrices, modulo ``islands'' (see Sec. \ref{sec:islands-QESA}).

Also, in the normal form \eqref{eq:nf-QES-B}, see page \pageref{eq:nf-QES-B}, the class $\overline{\QESB}$ is partitioned into 89 parts: 26 three--dimen\-sional ones, 39 two--dimensional ones, 20 one--dimensional ones, and four points. This partition is obtained by considering all the bifurcation surfaces of singularities, and bifurcation surfaces related to the presence of invariant straight lines, the presence of invariant parabolas, the presence of curves filled up with singular points, and connections of separatrices, modulo ``islands'' (see Sec. \ref{sec:islands-QESB}).

Finally, in the normal form \eqref{eq:nf-QES-C}, see page \pageref{eq:nf-QES-C}, the class $\overline{\QESC}$ is partitioned into 14 parts: four two--dimensional ones, seven one--dimensional ones, and three points. This partition is obtained by considering all the bifurcation surfaces of singularities, the presence of curves filled up with singular points, and bifurcation surfaces related to the presence of invariant straight line and invariant parabola, modulo ``islands'' (see Sec. \ref{sec:islands-QESC}).

\begin{figure}[h!]
	\centering
	\includegraphics[width=0.95\textwidth]{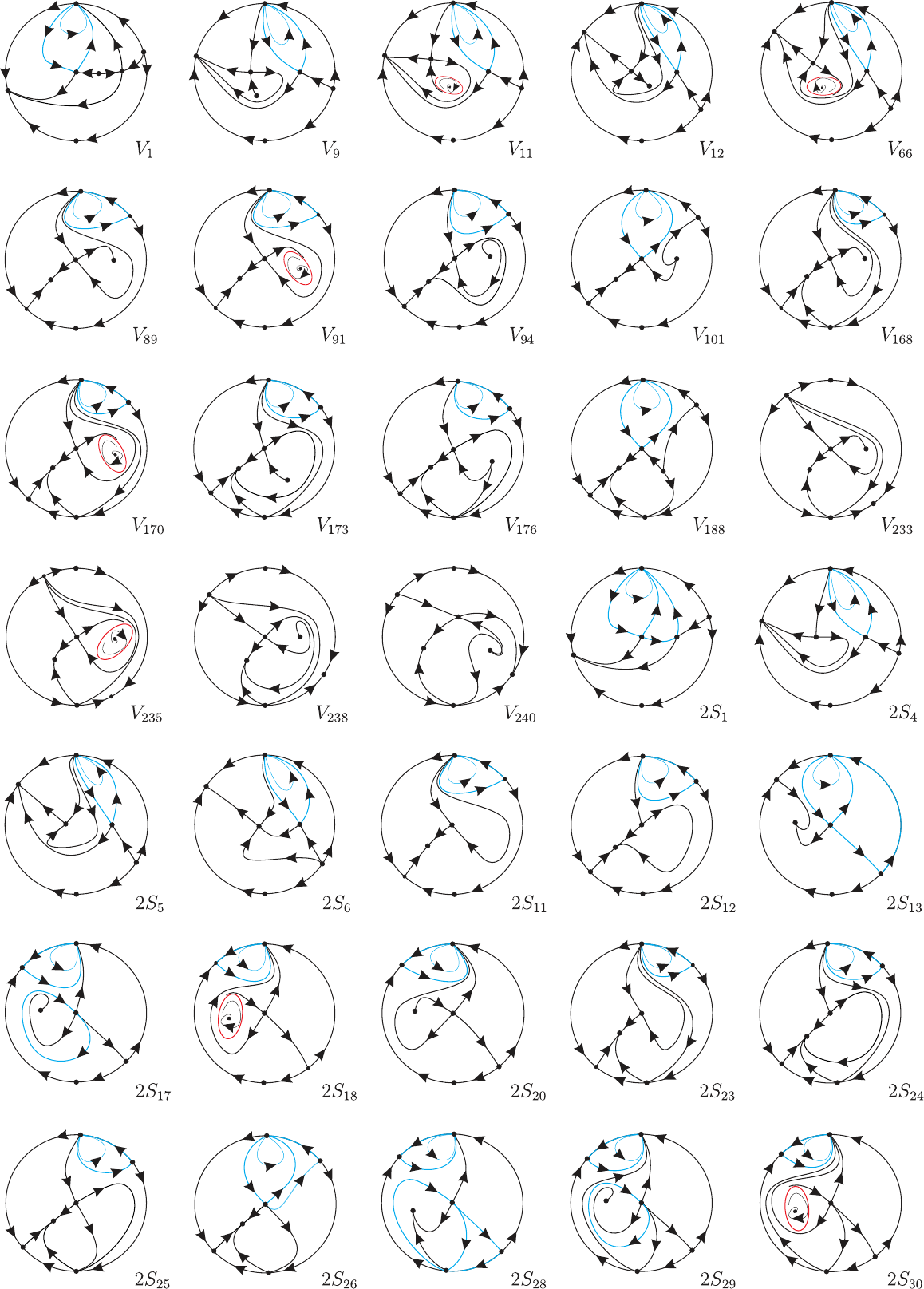}
	\caption{\small \label{fig:pp-QES-A-1} Phase portraits for quadratic vector fields possessing three real finite singular points, either an infinite nilpotent elliptic--saddle or an infinite nilpotent saddle, and an elemental infinite singularity, from class $\overline{\QESA}$}
\end{figure}

\begin{figure}[h!]
	\centering
	\includegraphics[width=0.99\textwidth]{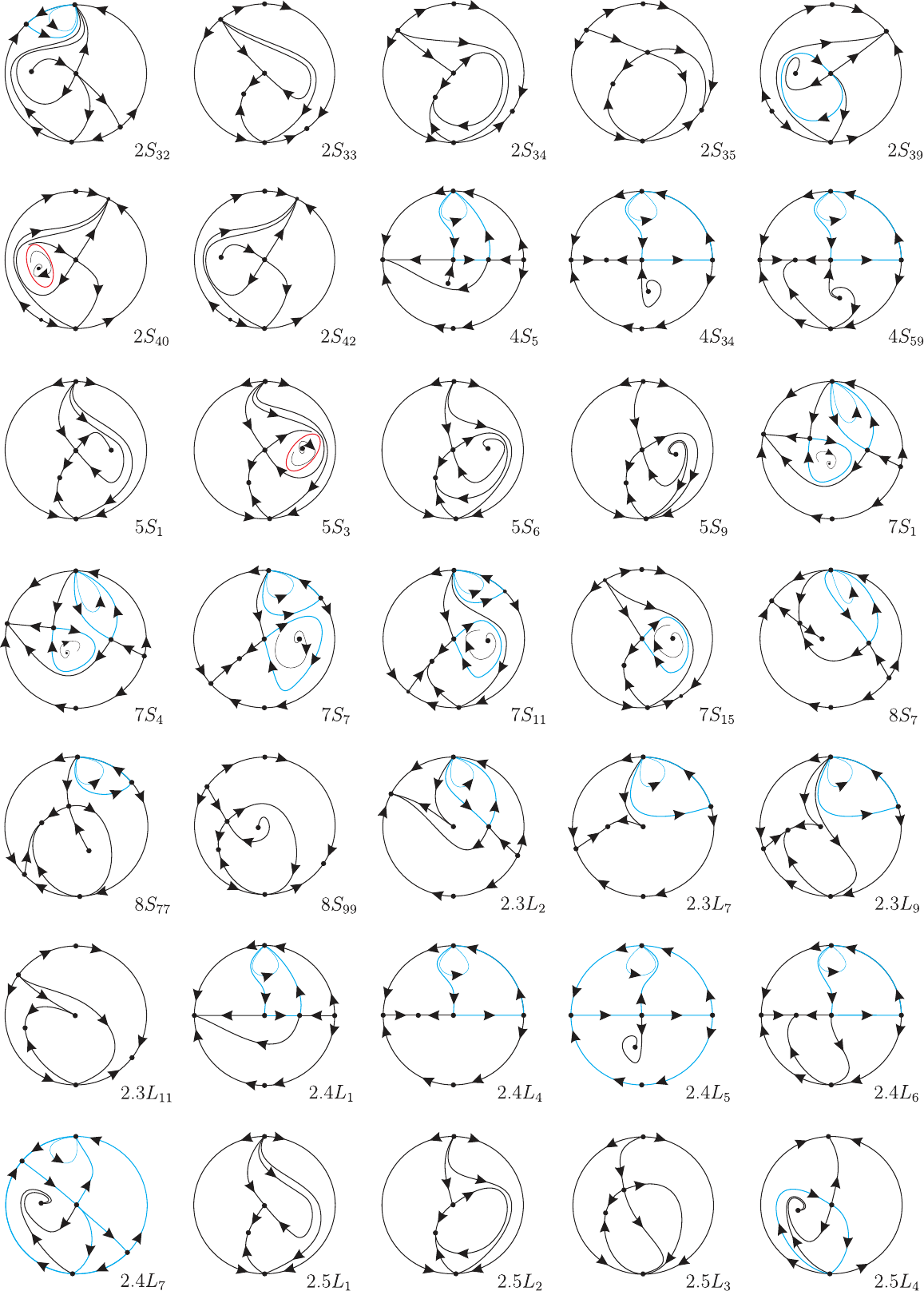}
	\caption{\small \label{fig:pp-QES-A-2} Continuation of Fig.~\ref{fig:pp-QES-A-1}}
\end{figure}

\begin{figure}[h!]
	\centering
	\includegraphics[width=0.95\textwidth]{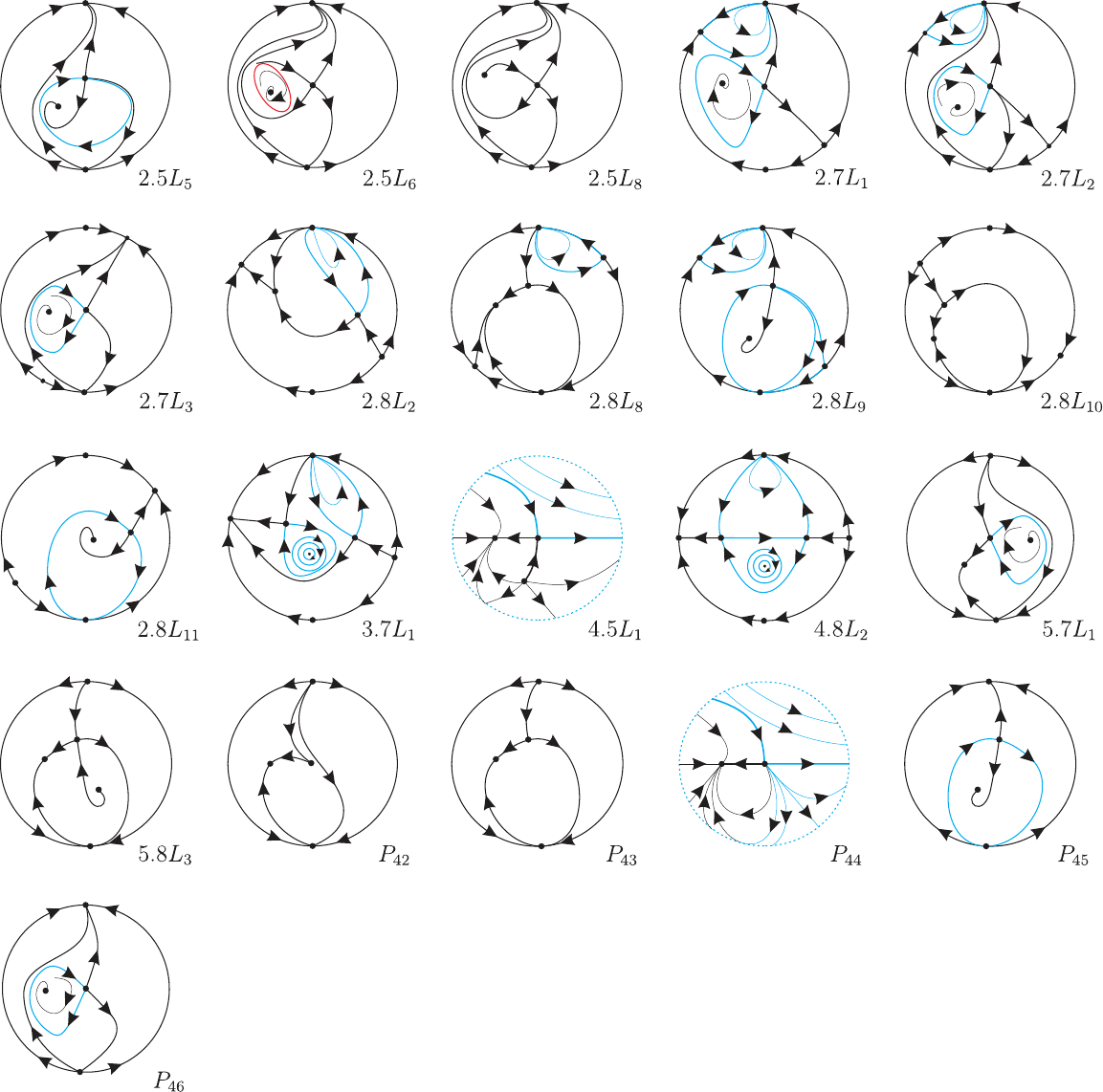}
	\caption{\small \label{fig:pp-QES-A-3} Continuation of Fig.~\ref{fig:pp-QES-A-2}}
\end{figure}

\begin{theorem} \label{th:main-thm-QES-A} There are $91$ topologically distinct phase portraits for the closure of the family of quadratic vector fields possessing three real finite singular points, either an infinite nilpotent elliptic--saddle or an infinite nilpotent saddle, and an elemental infinite singularity, and given by the normal form \eqref{eq:nf-QES-A} (class $\overline{\QESA}$). The bifurcation diagram for this class is given in the parameter space which is a subset of the real Euclidean three--dimensional space $\mathbb{R}^3$. All these phase portraits are shown in Figs. \ref{fig:pp-QES-A-1} to \ref{fig:pp-QES-A-3}. Also, for this class, the following statements hold:
	\begin{enumerate}
		\item[(a)] there are $18$ topologically distinct phase portraits in $\QESA$, namely, $V_{1}$, $V_{9}$, $V_{11}$, $V_{12}$, $V_{66}$, $V_{89}$, $V_{91}$, $V_{94}$, $V_{101}$, $V_{168}$, $V_{170}$, $V_{173}$, $V_{176}$, $V_{188}$, $V_{233}$, $V_{235}$, $V_{238}$, and $V_{240}$;	
		\item[(b)] consider the $18$ phase portraits from the previous item. Such phase portraits can be split according to the type of infinite singularities:
		\begin{itemize}
			\item phase portraits $V_{1}$, $V_{9}$, $V_{11}$, $V_{12}$, and $V_{66}$ possess an infinite nilpotent elliptic--saddle $\widehat{\!{1\choose 2}\!\!}\ PEP-H$ and also an infinite elemental node;
			\item phase portraits $V_{89}$, $V_{91}$, $V_{94}$, and $V_{101}$ possess an infinite nilpotent elliptic--saddle $\widehat{\!{1\choose 2}\!\!}\ PEP-H$ and also an infinite elemental saddle;
			\item phase portraits $V_{168}$, $V_{170}$, $V_{173}$, $V_{176}$, and $V_{188}$ possess an infinite nilpotent elliptic--saddle $\widehat{\!{1\choose 2}\!\!}\ E-PHP$ and also an infinite elemental saddle;
			\item phase portraits $V_{233}$, $V_{235}$, $V_{238}$, and $V_{240}$ possess an infinite nilpotent saddle $\widehat{\!{1\choose 2}\!\!}\ H-HHH$ and also an infinite elemental node;
		\end{itemize}
		in addition, from the study of the bifurcation diagram of class $\overline{\QESA}$ we observe the existence of $35$ two--dimensional regions (modulo islands) in which the corresponding phase portraits possess an infinite nilpotent elliptic--saddle $\widehat{\!{1\choose 2}\!\!}\ H-E$ and also an infinite elemental saddle;
		\item[(c)] there are ten phase portraits possessing exactly one simple limit cycle (or an odd number of them taking into account their multiplicity), and they are in the parts $V_{11}$, $V_{66}$, $V_{91}$, $V_{170}$, $V_{235}$, $2S_{18}$, $2S_{30}$, $2S_{40}$, $5S_{3}$,  and $2.5L_{6}$;
		\item[(d)] phase portraits $4.5L_{1}$ and $P_{44}$ possess the line at infinity filled up with singular points. Moreover, they have one infinite family of degenerate graphics;
		\item[(e)] there are nine phase portraits possessing only one nondegenerate graphic (surrounding a focus). More precisely, phase portraits $2S_{39}$, $7S_{15}$, $2.5L_{5}$, $2.7L_{3}$, $5.7L_{1}$, and $P_{46}$ have only one finite graphic and phase portraits $2.5L_{4}$, $2.8L_{11}$, and $P_{45}$ have only one infinite graphic;
		\item[(f)] there are $56$ phase portraits having only one infinite family of nondegenerate graphics (with no singularity inside), and these phase portraits are in the parts $V_{1}$, $V_{9}$, $V_{11}$, $V_{12}$, $V_{66}$, $V_{89}$, $V_{91}$, $V_{94}$, $V_{101}$, $V_{168}$, $V_{170}$, $V_{173}$, $V_{176}$, $V_{188}$, $2S_{1}$, $2S_{4}$, $2S_{5}$, $2S_{6}$, $2S_{11}$, $2S_{12}$, $2S_{13}$, $2S_{17}$, $2S_{18}$, $2S_{20}$, $2S_{23}$, $2S_{24}$, $2S_{25}$, $2S_{26}$, $2S_{28}$, $2S_{29}$, $2S_{30}$, $2S_{32}$, $4S_{5}$, $4S_{34}$, $4S_{59}$, $7S_{1}$, $7S_{4}$, $7S_{7}$, $7S_{11}$, $8S_{7}$, $8S_{77}$, $2.3L_{2}$, $2.3L_{7}$, $2.3L_{9}$, $2.4L_{1}$, $2.4L_{4}$, $2.4L_{5}$, $2.4L_{6}$, $2.4L_{7}$, $2.7L_{1}$, $2.7L_{2}$, $2.8L_{2}$, $2.8L_{8}$, $2.8L_{9}$, $3.7L_{1}$, and $4.8L_{2}$;
		\item[(g)] there are phase portraits that possess an infinite family of nondegenerate graphics (with no singularity inside) plus a finite number of nondegenerate graphics (which do not belong to the infinite family):
		\begin{itemize}
			\item phase portraits $2S_{1}$, $2S_{13}$, and $2S_{26}$ possess an infinite family of nondegenerate graphics plus one nondegenerate graphic with no singularity inside;
			\item phase portraits $2S_{17}$, $2S_{29}$, $7S_{1}$, $7S_{4}$, $7S_{7}$, $7S_{11}$, $2.7L_{1}$, and $2.7L_{2}$ possess an infinite family of nondegenerate graphics plus one nondegenerate graphic surrounding a focus;
			\item phase portraits $3.7L_{1}$ and $4.8L_{2}$ possess an infinite family of nondegenerate graphics plus one nondegenerate graphic surrounding a center;
			\item phase portraits $2S_{28}$ and $2.8L_{9}$ possess an infinite family of nondegenerate graphics plus two nondegenerate graphics surrounding the same focus;
			\item phase portrait $2.4L_{5}$ possesses an infinite family of nondegenerate graphics plus two nondegenerate graphics in which one of them surrounds a focus and the other one with no singularity inside;
			\item phase portrait $2.4L_{7}$ possesses an infinite family of nondegenerate graphics plus three nondegenerate graphics in which two of them surround the same focus and the other one with no singularity inside;
		\end{itemize}
		\item[(h)]  phase portraits $V_{11}$, $V_{66}$, $V_{91}$, $V_{170}$, $2S_{18}$, and $2S_{30}$ possess an infinite family of nondegenerate graphics plus one limit cycle.
	\end{enumerate}
\end{theorem}

\begin{figure}[h!]
	\centering
	\includegraphics[width=0.95\textwidth]{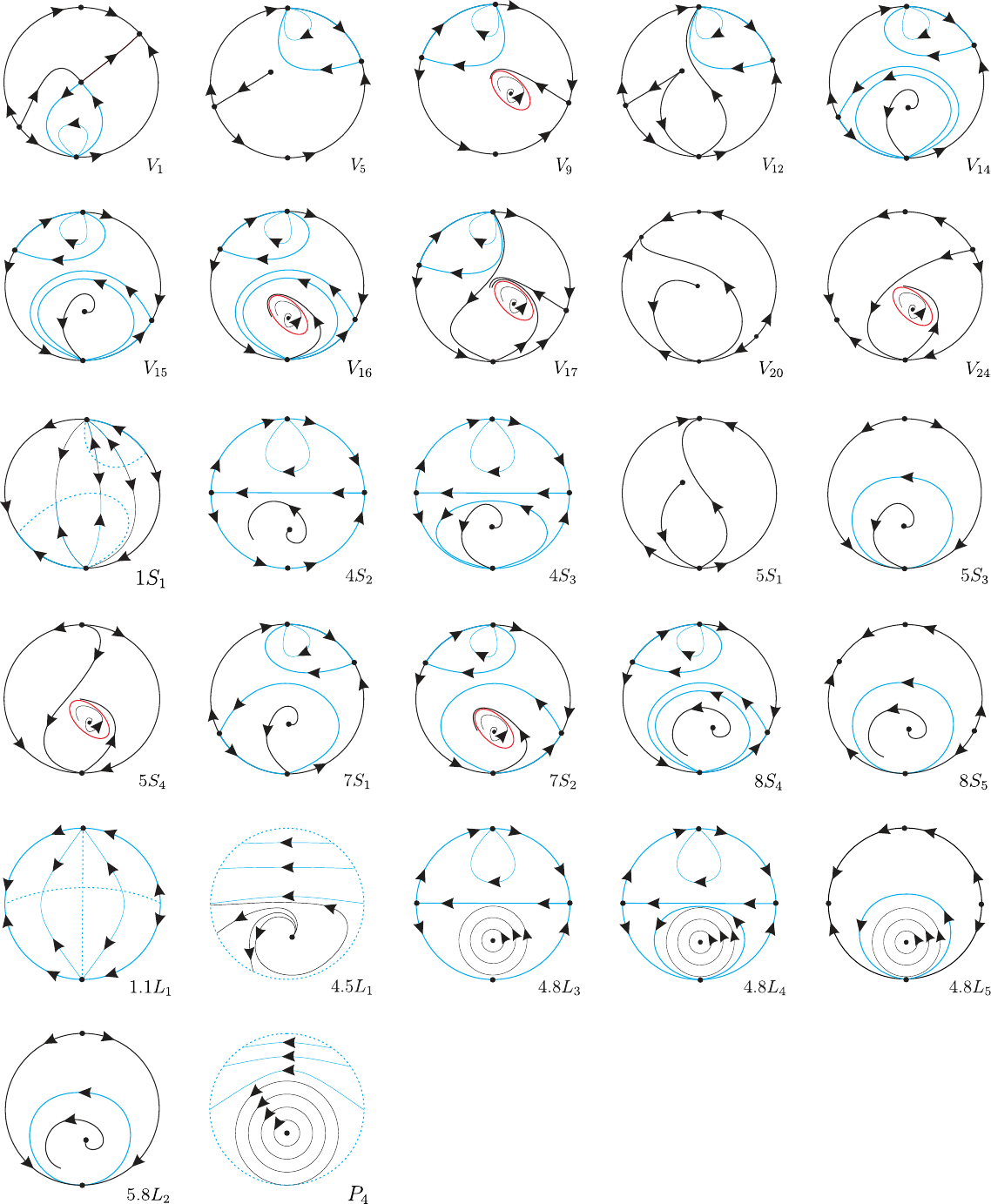}
	\caption{\small \label{fig:pp-QES-B} Phase portraits for quadratic vector fields possessing one real and two complex finite singular points, either an infinite nilpotent elliptic--saddle or an infinite nilpotent saddle, and an elemental infinite singularity, from class $\overline{\QESB}$}
\end{figure}

\begin{theorem} \label{th:main-thm-QES-B} There are $27$ topologically distinct phase portraits for the closure of the family of quadratic vector fields possessing one real and two complex finite singular points, either an infinite nilpotent elliptic--saddle or an infinite nilpotent saddle, and an elemental infinite singularity, and given by the normal form \eqref{eq:nf-QES-B} (class $\overline{\QESB}$). The bifurcation diagram for this class is given in the parameter space which is a subset of the real Euclidean three--dimensional space $\mathbb{R}^3$. All these phase portraits are shown in Fig. \ref{fig:pp-QES-B}. Also, for this class, the following statements hold:
	\begin{enumerate}
		\item[(a)] there are ten topologically distinct phase portraits in $\QESB$, namely, $V_{1}$, $V_{5}$, $V_{9}$, $V_{12}$, $V_{14}$, $V_{15}$, $V_{16}$, $V_{17}$, $V_{20}$, and $V_{24}$;	
		\item[(b)] consider the ten phase portraits from the previous item. Such phase portraits can be split according to the type of infinite singularities:
		\begin{itemize}
			\item phase portrait $V_{1}$ possesses an infinite nilpotent elliptic--saddle $\widehat{\!{1\choose 2}\!\!}\ PEP-H$ and also an infinite elemental node;
			\item phase portraits $V_{5}$ and $V_{9}$ possess an infinite nilpotent elliptic--saddle $\widehat{\!{1\choose 2}\!\!}\ PEP-H$ and also an infinite elemental saddle;
			\item phase portraits $V_{12}$, $V_{14}$, $V_{15}$, $V_{16}$, and $V_{17}$ possess an infinite nilpotent elliptic--saddle $\widehat{\!{1\choose 2}\!\!}\ E-PHP$ and also an infinite elemental saddle;
			\item phase portraits $V_{20}$ and $V_{24}$ possess an infinite nilpotent saddle $\widehat{\!{1\choose 2}\!\!}\ H-HHH$ and also an infinite elemental node;
		\end{itemize}
		in addition, from the study of the bifurcation diagram of class $\overline{\QESB}$ we observe the existence of five two--dimensional regions (modulo islands) in which the corresponding phase portraits possess an infinite nilpotent elliptic--saddle $\widehat{\!{1\choose 2}\!\!}\ H-E$ and also an infinite elemental saddle;
		\item[(c)] there are six phase portraits possessing exactly one simple limit cycle (or an odd number of them taking into account their multiplicity), and they are in the parts $V_{9}$, $V_{16}$, $V_{17}$, $V_{24}$, $5S_{4}$ and $7S_{2}$;
		\item[(d)] phase portraits $1S_{1}$ and $1.1L_{1}$ possess curves filled up with singular points. Moreover, they have one infinite family of degenerate graphics;
		\item[(e)] phase portraits $4.5L_{1}$ and $P_{4}$ possess the line at infinity filled up with singular points. Moreover, they have one infinite family of degenerate graphics;
		\item[(f)] there are three phase portraits possessing only one nondegenerate infinite graphic (surrounding a focus) and they are in the parts $5S_{3}$, $8S_{5}$ and $5.8L_{2}$. In addition, phase portrait $4.8L_{5}$ possesses only one nondegenerate infinite graphic (surrounding a center).
		\item[(g)] there are $15$ phase portraits having only one infinite family of nondegenerate graphics (with no singularity inside), and these phase portraits are in the parts $V_{1}$, $V_{5}$, $V_{9}$, $V_{12}$, $V_{14}$, $V_{15}$, $V_{16}$, $V_{17}$, $4S_{2}$, $4S_{3}$, $7S_{1}$, $7S_{2}$, $8S_{4}$, $4.8L_{3}$, and $4.8L_{4}$;
		\item[(h)] there are phase portraits that possess an infinite family of nondegenerate graphics (with no singularity inside) plus a finite number of nondegenerate graphics (which do not belong to the infinite family):
		\begin{itemize}
			\item phase portraits $4S_{2}$ and $7S_{1}$ possess an infinite family of nondegenerate graphics plus one nondegenerate graphic surrounding a focus;
			\item phase portrait $4.8L_{3}$ possesses an infinite family of nondegenerate graphics plus one nondegenerate graphic surrounding a center;
			\item phase portrait $7S_{2}$ possesses an infinite family of nondegenerate graphics plus one nondegenerate graphic surrounding a limit cycle;
			\item phase portraits $V_{14}$, $V_{15}$, $4S_{3}$, and $8S_{4}$ possess an infinite family of nondegenerate graphics plus two nondegenerate graphics surrounding the same focus;
			\item phase portrait $V_{16}$ possesses an infinite family of nondegenerate graphics plus two nondegenerate graphics surrounding the same limit cycle;
			\item phase portrait $4.8L_{4}$ possesses an infinite family of nondegenerate graphics plus two nondegenerate graphics surrounding the same center;	
		\end{itemize}
		\item[(i)]  phase portraits $V_{9}$, $V_{16}$, $V_{17}$, and $7S_{2}$ possess an infinite family of nondegenerate graphics plus one limit cycle.
	\end{enumerate}
\end{theorem}

\begin{figure}[h!]
	\centering
	\includegraphics[width=0.95\textwidth]{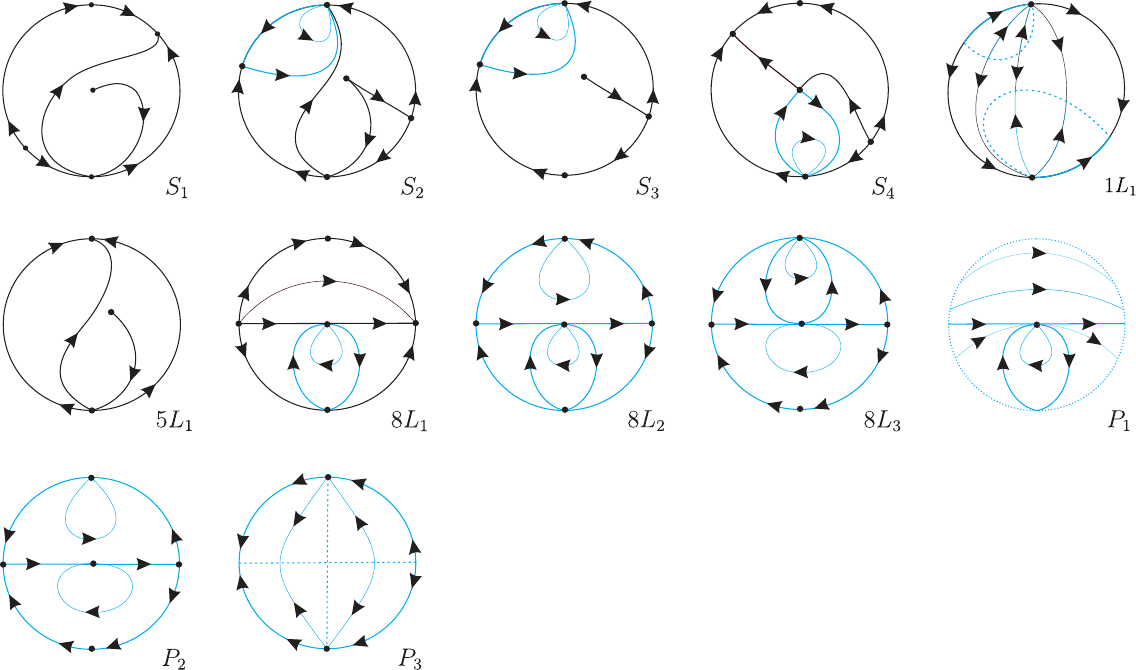}
	\caption{\small \label{fig:pp-QES-C} Phase portraits for quadratic vector fields possessing one real triple finite singular point, either an infinite nilpotent elliptic--saddle or an infinite nilpotent saddle, and an elemental infinite singularity, from class $\overline{\QESC}$}
\end{figure}

\begin{theorem} \label{th:main-thm-QES-C} There are twelve topologically distinct phase portraits for the closure of the family of quadratic vector fields possessing one real triple finite singular point, either an infinite nilpotent elliptic--saddle or an infinite nilpotent saddle, and an elemental infinite singularity, and given by the normal form \eqref{eq:nf-QES-C} (class $\overline{\QESC}$). The bifurcation diagram for this class is given in the parameter space which is a subset of the real Euclidean two--dimensional space $\mathbb{R}^2$. All these phase portraits are shown in Fig. \ref{fig:pp-QES-C}. Also, for this class, the following statements hold:
	\begin{enumerate}
		\item[(a)] there are four topologically distinct phase portraits in $\QESC$, namely, $S_{1}$, $S_{2}$, $S_{3}$, and $S_{4}$;	
		\item[(b)] the four phase portraits from the previous item can be split according to the type of infinite singularities:
		\begin{itemize}
			\item phase portrait $S_{1}$ possesses an infinite nilpotent saddle $\widehat{\!{1\choose 2}\!\!}\ H-HHH$ and also an infinite elemental node;
			\item phase portrait $S_{2}$ possesses an infinite nilpotent elliptic--saddle $\widehat{\!{1\choose 2}\!\!}\ E-PHP$ and also an infinite elemental saddle;
			\item phase portrait $S_{3}$ possesses an infinite nilpotent elliptic--saddle $\widehat{\!{1\choose 2}\!\!}\ PEP-H$ and also an infinite elemental saddle;
			\item phase portrait $S_{4}$ possesses an infinite nilpotent elliptic--saddle $\widehat{\!{1\choose 2}\!\!}\ PEP-H$ and also an infinite elemental node;
		\end{itemize}
		in addition, from the study of the bifurcation diagram of class $\overline{\QESC}$ we observe the existence of one one--dimensional region (modulo islands) in which the corresponding phase portrait possesses an infinite nilpotent elliptic--saddle $\widehat{\!{1\choose 2}\!\!}\ H-E$ and also an infinite elemental saddle;
		\item[(c)] there are no phase portraits possessing a limit cycle;
		\item[(d)] phase portraits $1L_{1}$ and $P_{3}$ possess curves filled up with singular points. Moreover, they have one infinite family of degenerate graphics;
		\item[(e)] phase portrait $P_{1}$ possesses the line at infinity filled up with singular points. Moreover, it has two infinite families of degenerate graphics;
		\item[(f)] there is no phase portraits possessing only one nondegenerate graphic;
		\item[(g)] there are five phase portraits having only one infinite family of nondegenerate graphics (with no singularity inside), and these phase portraits are in the parts $S_{2}$, $S_{3}$, $S_{4}$, $8L_{1}$, and $P_{1}$. Moreover, phase portraits $8L_{2}$, $8L_{3}$, and $P_{2}$ possess more than one infinite family of nondegenerate graphics;
		\item[(h)] there is no phase portrait possessing a finite number of nondegenerate graphics;
		\item[(i)]  there is no phase portrait possessing an infinite family of nondegenerate graphics plus one limit cycle.
	\end{enumerate}
\end{theorem}

\begin{proposition}\label{prop:number-l-c}
	There are $13$ topologically distinct phase portraits of codimension two, modulo limit cycles, in family $\QESA$ and six in family $\QESB$. The four topologically distinct phase portraits of codimension four without limit cycles in family $\QESC$ are topologically equivalent to phase portraits from family $\QESB$. So, in total we have $19$ topologically distinct phase portraits, modulo limit cycles.
\end{proposition}

\begin{corollary}\label{cor:comp-QESA-and-its-border}
	In Table~\ref{tab:comparison-QESA} (respectively, Tables~\ref{tab:comparison-QESB} and \ref{tab:comparison-QESC}) we give the numbers of phase portraits of both families $\QESA$ (respectively, $\QESB$ and $\QESC$) and its closure for special types of phase portraits.
	\begin{table}[t]\caption{\small Comparison between the set $\QESA$ and its border (the numbers represent the absolute values in each subclass)}\label{tab:comparison-QESA}
		\vspace{-0.3cm}
		\begin{center}
			%\resizebox{\columnwidth}{!}{%
				\begin{tabular}{c|cc}
					\cline{2-3}
					\multicolumn{1}{c|}{} & \multirow{2}{*}{$\QESA$} & Border of  \\
					\multicolumn{1}{c|}{} &  & $\QESA$ \\
					\hline
					Distinct phase portraits & $18$ & $73$ \\
					\hline
					Phase portraits with exactly one  & \multirow{2}{*}{$5$} & \multirow{2}{*}{$5$} \\
					simple limit cycle &  &  \\
					\hline
					Phase portraits with exactly one     & \multirow{2}{*}{$0$} & \multirow{2}{*}{$9$} \\
					nondegenerate graphic     &  &  \\
					\hline
					Phase portraits with at least       & \multirow{3}{*}{$14$} & \multirow{3}{*}{$42$} \\
					one infinite family of     &  &  \\
					nondegenerate graphics   &  &  \\
					\hline
					Phase portraits with degenerate & \multirow{2}{*}{$0$} & \multirow{2}{*}{$2$} \\
					graphics &  &  \\
					\hline
				\end{tabular}
			%}
		\end{center}
	\end{table}
	
	\begin{table}[t]\caption{\small Comparison between the set $\QESB$ and its border (the numbers represent the absolute values in each subclass)}\label{tab:comparison-QESB}
		\vspace{-0.3cm}
		\begin{center}
			%\resizebox{\columnwidth}{!}{%
				\begin{tabular}{c|cc}
					\cline{2-3}
					\multicolumn{1}{c|}{} & \multirow{2}{*}{$\QESB$} & Border of  \\
					\multicolumn{1}{c|}{} &  & $\QESB$ \\
					\hline
					Distinct phase portraits & $10$ & $15$ \\
					\hline
					Phase portraits with exactly one  & \multirow{2}{*}{$4$} & \multirow{2}{*}{$2$} \\
					simple limit cycle &  &  \\
					\hline
					Phase portraits with exactly one     & \multirow{2}{*}{$0$} & \multirow{2}{*}{$4$} \\
					nondegenerate graphic &  &  \\
					\hline
					Phase portraits with at least       & \multirow{3}{*}{$8$} & \multirow{3}{*}{$7$} \\
					one infinite family of     &  &  \\
					nondegenerate graphics   &  &  \\
					\hline
					Phase portraits with degenerate & \multirow{2}{*}{$0$} & \multirow{2}{*}{$4$} \\
					graphics &  &  \\
					\hline
				\end{tabular}
			%}
		\end{center}
	\end{table}
	
	\begin{table}[t]\caption{\small Comparison between the set $\QESC$ and its border (the numbers represent the absolute values in each subclass)}\label{tab:comparison-QESC}
		\vspace{-0.3cm}
		\begin{center}
			%\resizebox{\columnwidth}{!}{%
				\begin{tabular}{c|cc}
					\cline{2-3}
					\multicolumn{1}{c|}{} & \multirow{2}{*}{$\QESC$} & Border of  \\
					\multicolumn{1}{c|}{} &  & $\QESC$ \\
					\hline
					Distinct phase portraits & $4$ & $6$ \\
					\hline
					%Phase portraits with exactly one  & \multirow{2}{*}{$0$} & \multirow{2}{*}{$0$} \\
					%simple limit cycle &  &  \\
					%\hline
					%Phase portraits with exactly one     & \multirow{2}{*}{$0$} & \multirow{2}{*}{$0$} \\
					%nondegenerate graphic     &  &  \\
					%\hline
					Phase portraits with at least       & \multirow{3}{*}{$3$} & \multirow{3}{*}{$5$} \\
					one infinite family of     &  &  \\
					nondegenerate graphics   &  &  \\
					\hline
					Phase portraits with degenerate & \multirow{2}{*}{$0$} & \multirow{2}{*}{$3$} \\
					graphics &  &  \\
					\hline
				\end{tabular}
			%}
		\end{center}
	\end{table}
\end{corollary}

\begin{corollary}\label{cor:famBandC} There are seven topologically distinct phase portraits which appear simultaneously in both classes $\overline{\QESB}$ and $\overline{\QESC}$. The correspondences are indicated in Table~\ref{tab:equiv-pp-famb-famc} and the phase portraits in each row are topologically equivalent.
\end{corollary}

\begin{table}[h!]\caption{\small Topological equivalence between phase portraits from classes $\overline{\QESB}$ and $\overline{\QESC}$}\label{tab:equiv-pp-famb-famc}
	\begin{center}
		\begin{tabular}{cc}
			$\overline{\QESB}$ & $\overline{\QESC}$  \\
			\hline
			$V_{1}$ & $S_{4}$ \\ \hline
			$V_{7}$ & $S_{3}$ \\ \hline
			$V_{18}$ & $S_{2}$ \\ \hline
			$V_{28}$ & $S_{1}$ \\ \hline
			$1S_{1}$ & $1L_{1}$ \\ \hline
			$5S_{1}$ & $5L_{1}$ \\ \hline
			$1.1L_{1}$ & $P_{3}$ \\ \hline
		\end{tabular}
	\end{center}
\end{table}

In Figs. \ref{fig:pp-QES-A-1} to \ref{fig:pp-QES-C} we have illustrated all the singularities with a small disc. In case of degenerate systems we have also illustrated the infinite singular point belonging to the degenerate set with a small disc only if this point is an infinite singularity of the reduced system. We have drawn with thicker curves the separatrices and we have added some thinner orbits to avoid confusion in some  cases.

We have drawn all the limit cycles (and loops) possessing a convex shape (see Lemma 3.31 from  \cite{Artes-Llibre-Rezende-2018}). The limit cycles are colored in red (as in \cite{Artes-Mota-Rezende-2021c}, for instance) and all the graphics are colored in blue.%% We notice that, as the weak foci are graphics reduced to a point, the weak foci could be included in the definition of graphics and then we should have colored them in blue. However, in order to follow the same pattern as in previous similar papers and according to our definition of graphics we keep all of them in black.

\begin{remark} \label{rem:label-p-p} 
	We label the phase portraits according to the parts of the bifurcation diagram where they occur. Here we call {\it volumes} ($V$) the three--dimensional parts of the bifurcation diagram, {\it surfaces} ($S$) the two--dimensional ones, {\it curves} ($L$) the one--dimensional ones, and {\it points} ($P$) the zero--dimensional ones. These labels could be different for two topologically equivalent phase portraits occurring in distinct parts. Some of the phase portraits in three--dimensional parts also occur in some lower dimensional parts bordering these three--dimensional parts. An example occurs when a node turns into a focus. An analogous situation happens for phase portraits in two--dimensional or one--dimensional parts, coinciding with some phase portraits situated on their border. Moreover, as in \cite{Artes-Llibre-Schlomiuk-2006,Artes-Rezende-Oliveira-2015,Artes-Mota-Rezende-2021c}, we use the same pattern in order to indicate the elements ($V$), ($S$), ($L$) and ($P$) in the bifurcation diagram.
\end{remark}

This paper is organized as follows. In this section we have presented an introduction to this study, a brief review and some results already existent on the literature, and the statement of our main results.

In Sec.~\ref{sec:n-fs} we describe the normal forms that describe families $\QESA$, $\QESB$, and $\QESC$. Moreover, in such a section we present a study of invariant algebraic curves (straight lines and parabolas) for each family.

In Sec.~\ref{sec:b-d} we present the study of the three bifurcation diagrams. More precisely, in Sec.~\ref{subsec:bd-QES-A} (respectively, Sec.~\ref{subsec:bd-QES-B} and Sec.~\ref{subsec:bd-QES-C}) we present the bifurcation diagram of family $\QESA$ (respectively, $\QESB$ and $\QESC$). Related to the study of each family we present three subsections discussing, respectively, on the possible existence of ``islands'' in the corresponding bifurcation diagram, on the classification (up to topological equivalence) of the phase portraits, and on the completion of the proof of the correspondent theorem from Sec.~\ref{sec:int}.

In Appendix~\ref{ap:incomp-QES} we present some incompatibilities found in previous studies of phase portraits possessing specific properties on its singularities.

\medskip

\section{Normal forms and invariant algebraic curves from class $\QES$}\label{sec:n-fs}

In Table 6.1 from the book \cite{Artes-Llibre-Schlomiuk-Vulpe-2021a} one can obtain canonical 
forms of quadratic systems possessing different kinds of singular points. In this section we use the invariant 
theory in order to perform some affine transformations and time rescaling so that we obtain the normal 
forms that describe families $\QESA$, $\QESB$, and $\QESC$.

\begin{proposition} \label{prop:n-f-QESA}
	Every nondegenerate quadratic system possessing three real finite singular points plus either an infinite 
	nilpotent elliptic--saddle or an infinite nilpotent saddle, can be brought via affine transformations and 
	time rescaling to the following normal form
	\begin{equation}
	\begin{aligned}
	&x^\prime=cx+y-cx^2,\\
	&y^\prime=ex+\left(-1+\dfrac{e+f}{c}\right)y-ex^2+2xy,\\
	\end{aligned}
	\label{eq:nf-QES-A}
	\end{equation}
	where $c\in\mathbb{R}\setminus\{0\}$, $f\in\mathbb{R}^+\cup\{0\}$, and $e\in\mathbb{R}$ are 
	parameters, describing family $\QESA$.
\end{proposition}

\begin{proof}
	In fact, from \cite[Table 6.1]{Artes-Llibre-Schlomiuk-Vulpe-2021a} we get the so called \textit{canonical form 10}
	(see systems \eqref{eq:nf-10-book}), obtained by using affine transformations and time rescaling, 
	which describes quadratic systems possessing three real finite singular points and one infinite 
	singular point of multiplicity two (formed by the coalescence of one finite and one infinite elemental
	singular points).
	\begin{equation}
	\begin{aligned}
	&x^\prime=cx+dy-cx^2+2hxy,\\
	&y^\prime=ex+fy-ex^2+2mxy,\\
	\end{aligned}
	\label{eq:nf-10-book}
	\end{equation}
	where $c, d, h, e, f, m$ are real parameters, verifying conditions
	{\small
		$$(eh\!-\! cm)(cf\!-\! de)(dm\!-\! fh)(2(eh\!-\! cm)\!-\!(cf\!-\! de))\!\ne\!0.$$
	}
	For these systems, computations show that
	$$\begin{aligned}
	&\mu_{0}=0,\\
	&\mu_{1}=-4(eh-cm)(fh-dm)x,\\
	&\eta=4 h^2 (-8 e h + (c + 2 m)^2),\\
	&\widetilde{M}\!=\!-8(\!(\!-6eh\!+\!(\! c\!+\!2m\!)^2)x^2\!-\!2h(\! c\!+\!2m\!)xy\!+\!4h^2\! y^2),\\
	&\kappa=-128h^2 (e h - c m).\\
	\end{aligned}$$
	According to \cite[Diagram 6.3]{Artes-Llibre-Schlomiuk-Vulpe-2021a} we observe that
	in order to have three real elemental finite singularities and two singular points at infinity, being 
	one real elemental singularity and the other one a triple point formed by the coalescence of
	one finite singularity with two infinite ones, the previous invariants must verify
	$$\mu_0=0, \quad \mu_1\ne0, \quad \eta=0, \quad \widetilde{M}\ne0, \quad \kappa=0,$$
	respectively. So, by considering $h=0$ at systems~\eqref{eq:nf-10-book} we have systems
	\begin{equation}
	\begin{aligned}
	&x^\prime=cx+dy-cx^2,\\
	&y^\prime=ex+fy-ex^2+2mxy,\\
	\end{aligned}
	\label{eq:nf-10-book-h0}
	\end{equation}
	where $c, d, e, f, m$ are real parameters, verifying conditions
	\begin{equation}
	cdm(cf-de)\left[2c m+(cf-de)\right]\!\ne\!0,
	\label{eq:conditions-ne0}
	\end{equation}
	and, for systems~\eqref{eq:nf-10-book-h0},
	$$\begin{aligned}
	&\mu_{0}=0,\\
	&\mu_{1}=-4cdm^2x,\\
	&\eta=0,\\
	&\widetilde{M}=-8(c+2m)^2x^2,\\
	&\kappa=0.\\
	\end{aligned}$$
	Since $d\ne0$ and $m\ne0$ (due to \eqref{eq:conditions-ne0}), we perform the change
	$$(x,y,t)\rightarrow(x,(m/d)y,t/m),$$
	and we get systems
	\begin{equation*}
	\begin{aligned}
	&x^\prime=\frac{c}{m}x+y-\frac{c}{m}x^2,\\
	&y^\prime=\frac{de}{m^2}x+\frac{f}{m}y-\frac{de}{m^2}x^2+2xy.\\
	\end{aligned}
	\end{equation*}
	By renaming
	$$\frac{c}{m}\rightarrow c, \quad \frac{de}{m^2}\rightarrow e, \quad \frac{f}{m}\rightarrow f,$$
	we obtain systems~\eqref{eq:nf-10-book-h0} with $d=m=1$.\\
	Now we compute the following polynomial invariants:
	\begin{equation*}
	\begin{aligned}
	&\mathcal{B}_{1}=2(c-f-2)(c+f)\left[e+c(c-e+cf)\right],\\
	&\mathbb{D}=-192(cf-e)^2(2c-e+cf)^2.\\
	\end{aligned}
	\end{equation*}
	These polynomial invariants (whose meaning will be explained later) shall define bifurcation surfaces. 
	From the factors of $\mathcal{B}_{1}$ we observe that we can perform a translation
	$$f=F-1,$$
	and we obtain
	$$\mathbb{D}=-192(c+e-cF)^2(c-e+cF)^2.$$
	We rewrite the factors of $\mathbb{D}$ as a pair of horizontal parallel straight lines, i.e. we solve 
	$$-c-e+c(f+1)=F-\tilde{c} \quad \text{and} \quad c-e+c(f+1)=F+\tilde{c},$$
	which yield
	$$f=\frac{e-c+F}{c}, \quad \tilde{c}=c,$$
	and we rename $F=f$.
	%%or yet $$f=\frac{-c+e+f}{c}=-1+\frac{e+f}{c}, \quad \tilde{c}=c,$$ by renaming $F=f$. 
	Remember that $c\ne0$ due to conditions \eqref{eq:conditions-ne0}.
	Therefore we arrive at systems~\eqref{eq:nf-QES-A}. Indeed, by considering the change
	$$(x,y,t)\rightarrow(-x+1,y,-t)$$
	we obtain systems
	\begin{equation*}
	\begin{aligned}
	&x^\prime=cx+y-cx^2,\\
	&y^\prime=-ex-\left(1+\dfrac{e+f}{c}\right)y+ex^2+2xy,\\
	\end{aligned}
	\end{equation*}
	i.e. $(c,e, f)\rightarrow(c,-e, -f)$, so one can consider $f\in\mathbb{R}^+\cup\{0\}$.
\end{proof}

The next two results assure the existence of invariant straight lines and invariant parabolas, respectively, 
under certain conditions for family~\eqref{eq:nf-QES-A}.

\begin{lemma}\label{lemma:S4-inv-lines-QES-A}
	Family~\eqref{eq:nf-QES-A} possesses the following invariant straight line if and only if the corresponding 
	condition is satisfied:
	\begin{enumerate}[(i)]
		\item $\{y=0\}\Leftrightarrow e=0$;
		\item $\{c-f-(c-f)x+2y=0\}\Leftrightarrow e=(2+c)(c-f)/2$;
		\item $\{(c+f)x+2y=0\}\Leftrightarrow e=-(2+c)(c+f)/2$.
	\end{enumerate}
\end{lemma}

\begin{proof} We consider the algebraic curves
	\[
	\begin{aligned}
	f_1(x,y)&\equiv y=0, \\ f_2(x,y)&\equiv -c+f+(c-f)x-2y=0, \\ f_3(x,y)&\equiv (c+f)x+2y=0,\\
	\end{aligned}
	\]
	and we show that the polynomials
	\[
	\begin{aligned}
	K_1(x,y)&=2x+(f-c)/c, \\ K_2(x,y)&=2x, \\ K_3(x,y)&=2(x-1),
	\end{aligned}
	\]
	are the cofactors of $f_1=0$, $f_2=0$, and $f_3=0$, respectively, after restricting family~\eqref{eq:nf-QES-A} to the respective conditions.
\end{proof}

\begin{lemma}\label{lemma:S8-inv-parab-QES-A}
	Family~\eqref{eq:nf-QES-A} possesses the following invariant parabola if and only if the corresponding 
	condition is satisfied:
	\begin{enumerate}[(i)]
		\item $\left\{-\dfrac{c+c^2+e}{c}+\dfrac{2c+2c^2+e}{c}x-(1+c)x^2+y=0\right\}\Leftrightarrow f=-(2c^2+c+2e)$;
		\item $\left\{\dfrac{e}{c}x-(1+c)x^2+y=0\right\}\Leftrightarrow f=2c^2+c-2e$;
		\item $\{(1+c)x-(1+c)x^2+y=0\}\Leftrightarrow e=-f(1+c)$;
	\end{enumerate}
\end{lemma}

\begin{proof} We consider the algebraic curves
	\[
	\begin{aligned}
	g_1(x,y)&\equiv -\dfrac{c+c^2+e}{c}+\dfrac{2c+2c^2+e}{c}x-(1+c)x^2+y=0, \\ g_2(x,y)&\equiv \dfrac{e}{c}x-(1+c)x^2+y=0, \\ g_3(x,y)&\equiv (1+c)x-(1+c)x^2+y=0,\\
	\end{aligned}
	\]
	and we show that the polynomials
	\[
	\begin{aligned}
	H_1(x,y)&=-2cx, \\ H_2(x,y)&=2c(1-x), \\ H_3(x,y)&=c-f-2cx,
	\end{aligned}
	\]
	are the cofactors of $g_1=0$, $g_2=0$, and $g_3=0$, respectively, after restricting family~\eqref{eq:nf-QES-A} to the respective conditions.
\end{proof}

The study of the bifurcation diagram of family~\eqref{eq:nf-QES-A} is presented in Sec.~\ref{subsec:bd-QES-A}.

\begin{proposition} \label{prop:n-f-QESB}
	Every nondegenerate quadratic system possessing one real and two complex finite singular points plus 
	either an infinite nilpotent elliptic--saddle or an infinite nilpotent saddle, can be brought via affine transformations 
	and time rescaling to the following normal form
	\begin{equation}
	\begin{aligned}
	&x^\prime=-2gux+g(1+u^2)y+gx^2,\\
	&y^\prime=-2(\ell u-1)x+\ell(1+u^2)y+\ell x^2-2xy,\\
	\end{aligned}
	\label{eq:nf-QES-B}
	\end{equation}
	where $g\in\mathbb{R}\setminus\{0\}$, $u\in\mathbb{R}^+\cup\{0\}$, and $\ell\in\mathbb{R}$ are 
	parameters, describing family $\QESB$.
\end{proposition}

\begin{proof}
	In fact, from \cite[Table 6.1]{Artes-Llibre-Schlomiuk-Vulpe-2021a} we get the so called \textit{canonical form 11}
	(see systems \eqref{eq:nf-11-book}), obtained by using affine transformations and time rescaling, 
	which describes quadratic systems possessing one real and two complex finite singular points plus one infinite 
	singular point of multiplicity two (formed by the coalescence of one finite and one infinite elemental
	singular points).
	\begin{equation}
	\begin{aligned}
	&x^\prime= 2(h-gu)x+g(u^2+1)y+gx^2-2hxy,\\
	&y^\prime=2(m-\ell u)x+\ell(u^2+1)y+\ell x^2-2mxy,\\
	\end{aligned}
	\label{eq:nf-11-book}
	\end{equation}
	where $h, g, u, m, \ell$ are real parameters, verifying condition
	$$gm-h\ell \ne0.$$
	For these systems, computations show that
	$$\begin{aligned}
	\mu_{0}=& \ 0,\\
	\mu_{1}=& \ 4(h\ell -gm)^2(1+u^2)x,\\
	\eta=& \ 4h^2\left[(g+2m)^2-8h\ell\right],\\
	\widetilde{M}=& \ -8[\left((g+2m)^2-6h\ell\right)x^2-2h(g+2m)xy+4h^2y^2],\\
	\kappa=& \ 128h^2(gm-h\ell).\\
	\end{aligned}$$
	As in the proof of Proposition~\ref{prop:n-f-QESA}, from 
	\cite[Diagram 6.3]{Artes-Llibre-Schlomiuk-Vulpe-2021a}, the previous invariants must verify
	$$\mu_0=0, \quad \mu_1\ne0, \quad \eta=0, \quad \widetilde{M}\ne0, \quad \kappa=0,$$
	respectively. So, by considering $h=0$ at systems~\eqref{eq:nf-11-book} we have systems
	\begin{equation}
	\begin{aligned}
	&x^\prime= -2gux+g(u^2+1)y+gx^2,\\
	&y^\prime=2(m-\ell u)x+\ell(u^2+1)y+\ell x^2-2mxy,\\
	\end{aligned}
	\label{eq:nf-11-book-h0}
	\end{equation}
	where $g, u, m, \ell$ are real parameters, verifying condition
	\begin{equation}
	gm \ne0
	\label{eq:conditions-ne0-QESB}
	\end{equation}
	and, for systems~\eqref{eq:nf-11-book-h0},
	$$\begin{aligned}
	&\mu_{0}=0,\\
	&\mu_{1}=4g^2m^2(1+u^2)x,\\
	&\eta=0,\\
	&\widetilde{M}=-8(g+2m)^2x^2,\\
	&\kappa=0.\\
	\end{aligned}$$
	Since $m\ne0$ (due to \eqref{eq:conditions-ne0-QESB}), we perform the change
	$$(x,y,t)\rightarrow(x,y,t/m),$$
	and we get systems
	\begin{equation*}
	\begin{aligned}
	&x^\prime= -2\frac{g}{m}ux+\frac{g}{m}(u^2+1)y+\frac{g}{m}x^2,\\
	&y^\prime=2\left(1-\frac{\ell}{m}u\right)x+\frac{\ell}{m}(u^2+1)y+\frac{\ell}{m}x^2-2xy.\\
	\end{aligned}
	\end{equation*}
	By renaming
	$$\frac{g}{m}\rightarrow g, \quad \frac{\ell}{m}\rightarrow \ell,$$
	we obtain systems~\eqref{eq:nf-11-book-h0} with $m=1$, i.e., we arrive
	at normal form~\eqref{eq:nf-QES-B}, in which $g\ne0$ due to \eqref{eq:conditions-ne0-QESB}. Indeed, by considering the change
	$$(x,y,t)\rightarrow(-x,y,-t),$$
	we obtain systems
	\begin{equation*}
	\begin{aligned}
	&x^\prime=2gux+g(1+u^2)y+gx^2,\\
	&y^\prime=-2(-1+\ell u)x-\ell(1+u^2)y-\ell x^2-2xy,\\
	\end{aligned}
	\end{equation*}
	i.e. $(u, \ell, g)\rightarrow(-u, -\ell, g)$, so one can consider $u\in\mathbb{R}^+\cup\{0\}$.
\end{proof}

In the next result we prove the existence of invariant algebraic curves (straight lines and parabolas) under 
certain conditions for systems~\eqref{eq:nf-QES-B}.

\begin{lemma}\label{lemma:S4-inv-curves-QES-B}
	Systems~\eqref{eq:nf-QES-B} possess the following invariant algebraic curves if and only if the corresponding 
	condition is satisfied:
	\begin{enumerate}[(i)]
		\item $\{y-1=0\}\Leftrightarrow \ell=0$;
		\item $\left\{\dfrac{\ell x^2-2 \ell u x+2u}{\ell u^2+\ell-2 u}+y=0\right\}\Leftrightarrow g=\dfrac{\ell u^2+\ell-2 u}{2 u}$;
		\item $\left\{\dfrac{(g+1) x^2+1}{g}+y=0\right\}\Leftrightarrow \ell=u=0$.
	\end{enumerate}
\end{lemma}

\begin{proof} We consider the algebraic curves
	\[
	\begin{aligned}
	f_1(x,y)&\equiv y-1=0, \\ 
	f_2(x,y)&\equiv \dfrac{\ell x^2-2 \ell u x+2u}{\ell u^2+\ell-2 u}+y=0, \\ 
	f_3(x,y)&\equiv \dfrac{(g+1) x^2+1}{g}+y=0,\\
	\end{aligned}
	\]
	and we show that the polynomials
	\[
	\begin{aligned}
	K_1(x,y)&=-2x, \\ K_2(x,y)&=\frac{x \left(\ell u^2+\ell-2 u\right)}{u}, \\ K_3(x,y)&=2gx,
	\end{aligned}
	\]
	are the cofactors of $f_1=0$, $f_2=0$, and $f_3=0$, respectively, after restricting systems~\eqref{eq:nf-QES-B} to the respective conditions.
\end{proof}

The bifurcation diagram of systems~\eqref{eq:nf-QES-B} is studied in Sec.~\ref{subsec:bd-QES-B}.

\begin{proposition} \label{prop:n-f-QESC}
	Every nondegenerate quadratic system possessing one triple real finite singular point plus either 
	an infinite nilpotent elliptic--saddle or an infinite nilpotent saddle, can be brought via affine transformations 
	and time rescaling to the following normal form
	\begin{equation}
	\begin{aligned}
	&x^\prime=gy+gx^2,\\
	&y^\prime=\ell y+2xy+\ell x^2,\\
	\end{aligned}
	\label{eq:nf-QES-C}
	\end{equation}
	where $g\in\mathbb{R}\setminus\{0\}$ and $\ell\in\mathbb{R}^+\cup\{0\}$ are 
	parameters, describing family $\QESC$.
\end{proposition}

\begin{proof}
	In fact, from \cite[Table 6.1]{Artes-Llibre-Schlomiuk-Vulpe-2021a} we get the so called \textit{canonical form 13}, 
	obtained by using affine transformations and time rescaling (see systems \eqref{eq:nf-13-book}),
	which describes quadratic systems possessing one real triple finite singular point and one infinite 
	singular point of multiplicity two (formed by the coalescence of one finite and one infinite elemental
	singular points).
	\begin{equation}
	\begin{aligned}
	&x^\prime=gy+gx^2+2hxy,\\
	&y^\prime= \ell y+\ell x^2+2mxy,\\
	\end{aligned}
	\label{eq:nf-13-book}
	\end{equation}
	where $g, h, \ell, m$ are real parameters, verifying condition
	$$gm-\ell h\ne0.$$
	For these systems, computations show that
	$$\begin{aligned}
	\mu_{0}=&0,\\
	\mu_{1}=&4(h\ell -gm)^2x,\\
	\eta=&4h^2\left[8h\ell +(g-2m)^2\right],\\
	\widetilde{M}=&-8\left[\left((g-2m)^2-6h\ell\right)x^2-2h(g-2m)xy+4h^2y^2\right],\\
	\kappa=&128h^2(h\ell -gm).\\
	\end{aligned}$$
	As in the proof of Propositions~\ref{prop:n-f-QESA} and \ref{prop:n-f-QESB}, from 
	\cite[Diagram 6.3]{Artes-Llibre-Schlomiuk-Vulpe-2021a}, the previous invariants must verify
	$$\mu_0=0, \quad \mu_1\ne0, \quad \eta=0, \quad \widetilde{M}\ne0, \quad \kappa=0,$$
	respectively. So, by considering $h=0$ at systems~\eqref{eq:nf-13-book} we have systems
	\begin{equation}
	\begin{aligned}
	&x^\prime=gy+gx^2,\\
	&y^\prime= \ell y+\ell x^2+2mxy,\\
	\end{aligned}
	\label{eq:nf-13-book-h0}
	\end{equation}
	where $g, \ell, m$ are real parameters, verifying condition
	\begin{equation}
	gm \ne0
	\label{eq:conditions-ne0-QESC}
	\end{equation}
	and, for systems~\eqref{eq:nf-13-book-h0},
	$$\begin{aligned}
	&\mu_{0}=0,\\
	&\mu_{1}=4g^2m^2x,\\
	&\eta=0,\\
	&\widetilde{M}=-8(g-2m)^2x^2,\\
	&\kappa=0.\\
	\end{aligned}$$
	Since $m\ne0$ (due to \eqref{eq:conditions-ne0-QESC}), we perform the change
	$$(x,y,t)\rightarrow(x,y,t/m),$$
	and we get systems
	\begin{equation*}
	\begin{aligned}
	&x^\prime=\frac{g}{m}y+\frac{g}{m}x^2,\\
	&y^\prime=\frac{\ell}{m}y+\frac{\ell}{m}x^2+2xy.\\
	\end{aligned}
	\end{equation*}
	By renaming
	$$\frac{g}{m}\rightarrow g, \quad \frac{\ell}{m}\rightarrow \ell,$$
	we obtain systems~\eqref{eq:nf-13-book-h0} with $m=1$, i.e., we arrive
	at normal form~\eqref{eq:nf-QES-C}, in which $g\ne0$ due to \eqref{eq:conditions-ne0-QESC}. Indeed, by considering the change
	$$(x,y,t)\rightarrow(-x,y,-t),$$
	we obtain systems
	\begin{equation*}
	\begin{aligned}
	&x^\prime=gy+gx^2,\\
	&y^\prime=-\ell y+2xy-\ell x^2,\\
	\end{aligned}
	\end{equation*}
	i.e. $(g, \ell)\rightarrow(g, -\ell)$, so one can consider $\ell\in\mathbb{R}^+\cup\{0\}$.
\end{proof}

In what follows we prove the existence of invariant algebraic curves (straight lines and parabolas) under 
certain conditions for family~\eqref{eq:nf-QES-C}.

\begin{lemma}\label{lemma:S4-inv-curves-QES-C}
	Family~\eqref{eq:nf-QES-C} possesses the following invariant algebraic curves if and only if the corresponding 
	condition is satisfied:
	\begin{enumerate}[(i)]
		\item $\{y=0\}\Leftrightarrow \ell=0$;
		\item $\left\{\dfrac{(g-1) x^2}{g}+y=0\right\}\Leftrightarrow \ell=0$.
	\end{enumerate}
\end{lemma}

\begin{proof} We consider the algebraic curves
	\[
	\begin{aligned}
	f_1(x,y)&\equiv y=0, \\ f_2(x,y)&\equiv \frac{(g-1) x^2}{g}+y=0, \\
	\end{aligned}
	\]
	and we show that the polynomials
	\[
	\begin{aligned}
	K_1(x,y)&=2x, \\ K_2(x,y)&=2gx,
	\end{aligned}
	\]
	are the cofactors of $f_1=0$ and $f_2=0$, respectively, after restricting systems~\eqref{eq:nf-QES-C} 
	to the respective conditions.
\end{proof}

In Sec.~\ref{subsec:bd-QES-C} we present the study of the bifurcation diagram of normal form~\eqref{eq:nf-QES-C}.

\section{The bifurcation diagrams from class $\QES$}\label{sec:b-d}

In this paper we intend to perform the study of three bifurcation diagrams. And to achieve this goal we shall 
use algebraic and topological invariants. The algebraic invariants make results independent of specific normal 
forms. They also distinguish the phase portraits as the topological invariants also do. In this paper we use the 
concepts of algebraic invariant and T--comitant as formulated by the Sibirsky's School for differential equations. 
For a quick summary of the general theory of these polynomial invariants and their relevance in working with 
polynomial differential systems we recommend Sec. 7 of \cite{Artes-Llibre-Schlomiuk-2006}.

It is worth mentioning that from Sec.~7 of \cite{Artes-Llibre-Vulpe-2008} and \cite{Vulpe-2011} we get 
formulas which give the bifurcation algebraic sets of singularities in $\mathbb{R}^{12}$, produced by 
changes that may occur in the local nature of finite singularities. Also, from \cite{Schlomiuk-Vulpe-2005} 
we get equivalent formulas for the infinite singular points. All of these formulas were lately compiled and 
improved in the book \cite{Artes-Llibre-Schlomiuk-Vulpe-2021a}. In the next three subsections we shall
use several results of such a book.

\medskip

\subsection{The bifurcation diagram of family $\QESA$}\label{subsec:bd-QES-A}

In this section we present the study of the bifurcation diagram of family $\QESA$, given by systems
\eqref{eq:nf-QES-A}. 

Initially remember that family~\eqref{eq:nf-QES-A} is described by the parameters $c\in\mathbb{R}\setminus\{0\}$, 
$f\in\mathbb{R}^+\cup\{0\}$, and $e\in\mathbb{R}$. So we shall consider the bifurcation diagram formed by
planes $c=c_0\neq0$ and, in each plane, the Cartesian coordinates are $(e,f)$ with $f\geq0$.

Also, from \cite[Lemma 5.2]{Artes-Llibre-Schlomiuk-Vulpe-2021a} we calculate
$$\mu_{0}=0, \quad \mu_{1}=-4cx.$$
The condition $c\ne0$ implies $\mu_{1}\ne0$ and, therefore, we have nondegenerate systems.

Now we present the value of the algebraic invariants and T--comitants (with respect to family~\eqref{eq:nf-QES-A})
which are relevant in our study.

\medskip

\noindent \textbf{Bifurcation surfaces due to multiplicities of singularities}

\medskip

\noindent\textbf{$(\mathcal{S}_{2})$} This is the bifurcation surface in $\mathbb{R}^3$ due to multiplicity of finite 
singular points, formed by the coalescence of at least two finite singular points. For family~\eqref{eq:nf-QES-A}, 
according to \cite[Table 5.1]{Artes-Llibre-Schlomiuk-Vulpe-2021a} we calculate 
$$\begin{aligned}
\mathbf{D}\!=\!&-192(c-f)^2(c+f)^2,
\end{aligned}$$
and we define the surface
$$\begin{aligned}
(\mathcal{ S}_{2})\!:& \ (c-f)(c+f)=0,
\end{aligned}$$
which is clear formed by two planes in $\mathbb{R}^3$. Additionally as the comitant
$$\mathbf{P}\mathbf{R}\big{\vert}_{c=\pm f}=768f^4x^4y^2$$
is nonzero, from \cite[Table 5.1]{Artes-Llibre-Schlomiuk-Vulpe-2021a} we conclude that along surface 
$(\mathcal{S}_{2})$ we have one double and one simple real finite singular points.

\medskip

\noindent {\bf (${\cal S}_{5}$)} This is the bifurcation surface due to multiplicity of infinite singular points. 
Previously we mentioned that an infinite elliptic--saddle is a triple infinite singular point formed by the 
coalescence of one finite singular point with two infinite ones. So, for family~\eqref{eq:nf-QES-A} we 
have at most two pairs of infinite singular points. According to \cite[Lemma 5.5]{Artes-Llibre-Schlomiuk-Vulpe-2021a}, 
for this family we calculate
$$\eta=0, \quad \widetilde{M}=-8(2+c)^2x^2, \quad C_2=x^2\left[ex-(2+c)y\right],$$
and we observe that along surface (in fact a plane in $\mathbb{R}^3$)
$$\begin{aligned}
(\mathcal{ S}_{5})\!:& \ c+2=0,
\end{aligned}$$
we have a coalescence of infinite singular points. In addition, due to the mentioned result, on the plane 
$c=-2$ all the phase portraits corresponding to $e=0$ have the line at infinity filled up with singular points.

\medskip

\noindent \textbf{The surface of $C^{\infty}$ bifurcation points due to a strong saddle or a strong 
	focus changing the sign of their traces (weak saddle or weak focus)}

\medskip

\noindent {\bf (${\cal S}_{3}$)} This is the bifurcation surface due to weak finite singularities, which 
occurs when the trace of a finite singular point is zero. According to \cite{Vulpe-2011}, for 
family~\eqref{eq:nf-QES-A} we calculate
$$\begin{aligned}
\mathcal{T}_{4}=&\mathcal{T}_{3}=\mathcal{T}_{2}=\mathcal{T}_{1}=0,\\
\sigma=&\dfrac{c(c-1)+e+f-2c(c-1)x}{c},
\end{aligned}$$
then due to the results on the mentioned paper, in the case in which $\sigma$ is generically nonzero,
the family under consideration could possess one and only one weak singularity. Moreover as
$$\begin{aligned}
\mathcal{F}_{1}=&-2\left[3e+(2+c)f\right], \qquad \mathcal{H}=0,\\
\mathcal{B}_{1}=&\dfrac{2\left[c^2(c-1^2)-(e+f)^2\right](e+cf)}{c},\\
\mathcal{B}_{2}=&\dfrac{-2 (c-1)^2 \left[c^4-2 c^3+c^2-2 c f (e+f)\right]+2 (c-1)^2 (e+f) (3 e+f)}{c},
\end{aligned}$$
assuming $\mathcal{F}_{1}\ne0$, for family~\eqref{eq:nf-QES-A} we can obtain one weak 
singularity ($s^{(1)}$ or $f^{(1)}$) along the surface given by $\mathcal{B}_{1}=0$, i.e.
$$\begin{aligned}
({\cal S}_{3})\!:& \ \dfrac{\left[c^2(c-1^2)-(e+f)^2\right](e+cf)}{c}=0.
\end{aligned}$$
%We advise the reader to plot surface $({\cal S}_{3})$ in order to have a three--dimensional picture of this interesting surface.\\
We highlight that this bifurcation can produce a topological change if the weak point is a focus or just a 
$C^\infty$ change if it is a saddle, except when this bifurcation coincides with a loop bifurcation associated 
with the same saddle, in which case, the change may also be topological (see for instance 
\cite[p. 50]{Artes-Rezende-Oliveira-2015}).

\begin{remark}\label{weak-singularities} 
	\begin{enumerate}
		\item We just saw that in order to define surface {\bf (${\cal S}_{3}$)} we considered 
		$\sigma\ne0$ and $\mathcal{F}_{1}\ne0$. However, according to \cite[item $(e)$]{Vulpe-2011}, 
		when $\sigma\ne0$ and $\mathcal{F}_{1}=0$ we can have either an integrable saddle or a center.  
		Later we shall analyze when we have an integrable saddle. Now we investigate when we have a finite 
		singular point which is a center. In fact, as we already have $\mathcal{H}=0$, from the mentioned paper
		we solve $\mathcal{F}_{1}=\mathcal{B}_{1}=0$  (together with $\sigma\ne0$ and $f\geq0$), and we 
		obtain two solutions
		\begin{equation}
		\{e=0, f=0\}, \ \{e=-c(c+2), f=3c\}.
		\label{eq:sol-weak-sing}
		\end{equation}
		Also, when we compute $\mathcal{B}_{2}$ along these two solutions we obtain, in each case,
		$-8 (c-1)^4 c$, which is generically negative if $c>0$. Note that we must have $c\ne1$, because
		each one of the two solutions together with $c=1$ imply $\sigma=0$.\\
		Therefore, from \cite[item ($e_4$)--$\beta$]{Vulpe-2011}, this study shows that for $c>0$ we 
		shall always find a center type singular point when we have \eqref{eq:sol-weak-sing}.
		\item We observe that, independently of $x$, for $c\ne0$, we have $\sigma=0$ if and only if 
		$f=-e$ and $c=1$. Under these conditions, we have that $\mu_{0}=0$, $\mathbf{D}=-192 \left(f^2-1\right)^2$, 
		and $\mathbf{R}=48x^2$. So, according to \cite[item $(f_{3})$]{Vulpe-2011} we have three 
		finite singular points, being two integrable saddles and one center. In other words, when $c=1$,
		during the study of the curve $f=-e$ we shall always obtain a phase portrait containing two 
		integrable saddles and one center type singular point.
	\end{enumerate}
\end{remark}

\medskip

\noindent \textbf{The surface of $C^{\infty}$ bifurcation due to a node becoming a focus}

\medskip

\noindent {\bf (${\cal S}_{6}$)} This surface contains the points of the parameter space where a finite node 
of the systems turns into a focus. This surface is a $C^{\infty}$ but not a topological bifurcation surface. In fact, 
when we only cross the surface (${\cal S}_{6}$) in the bifurcation diagram, the topological phase portraits do 
not change. However, this surface is relevant for isolating the regions where a limit cycle surrounding an antisaddle 
cannot exist. According to \cite[Table 6.2]{Artes-Llibre-Schlomiuk-Vulpe-2021a} we calculate
$$\begin{aligned}
\mu_{0}=&0, \quad \mathbf{D}=-192(c^2-f^2)^2, \quad \mathbf{R}=48c^2x^2,\\
\widetilde{K}=&-4cx^2, \quad G_{9}=0,
\end{aligned}$$
and for the mentioned table we conclude that the invariant $W_{7}$ is responsible for describing the node--focus 
bifurcation. We compute this invariant polynomial and we define surface {\bf (${\cal S}_{6}$)} by the zero set of
$${\small
	\begin{aligned}
	& \dfrac{1}{c^4}\left[2 c^3 - e^2 - 2 c e f - c f^2(2 + c)\right]\times\\
	& \left[2 c^3 + c^4 + c^2 (1 + 2 e - 2 f) - 2 c (e + f) + (e + f)^2\right]\times\\
	& \left[2 c^3 + c^4 + c^2 (1 - 2 e + 2 f)+ 2 c (e + f) + (e + f)^2 \right]=0. \\ 
	\end{aligned}}$$

\medskip

\noindent \textbf{Bifurcation surface in $\mathbb{R}^3$ due to the presence of invariant straight lines}

\medskip

\noindent {\bf (${\cal S}_{4}$)} This surface contains the points of the parameter space in which there appear 
invariant straight lines (see Lemma~\ref{lemma:S4-inv-lines-QES-A}). This surface is split into some regions. 
Depending on these regions, the straight line may contain connections of separatrices from different points or not. 
So, in some cases, it may imply a topological bifurcation and, in others, just a $C^{\infty}$ bifurcation. 
According to \cite{Artes-Llibre-Schlomiuk-Vulpe-2021a}, the equation of this surface is given by the invariant 
$B_1$. It is worth mentioning that $B_1=0$ is only a necessary condition for the existence of an invariant straight 
line, but it is not sufficient (see Corollary 4.6 from \cite{Schlomiuk-Vulpe-2004}), i.e. we may find some 
component of $B_1=0$ that does not represent an invariant straight line. For family~\eqref{eq:nf-QES-A} 
we compute the invariant $B_1$ and we define the surface
$$\begin{aligned}
({\cal S}_{4})\!: & \ e\left[c(2+c-f)-2(e+f)\right]\left[c(2+c+f)+2(e+f)\right]=0,
\end{aligned}$$
which is the union of one plane together with two quadric surfaces.

\medskip

\noindent \textbf{Bifurcation surface in $\mathbb{R}^3$ due to the presence of invariant parabolas}

\medskip

\noindent {\bf (${\cal S}_{8}$)} This surface contains the points of the parameter space in which there appear 
invariant parabolas. As in the case of surface (${\cal S}_{4}$), this surface is split into some regions. 
Depending on these regions, the parabola may contain connections of separatrices from different points or not. 
So, in some cases, it may imply a topological bifurcation and, in others, just a $C^{\infty}$ bifurcation. 
According to the conditions stated in Lemma~\ref{lemma:S8-inv-parab-QES-A} we define this surface by
$$\begin{aligned}
({\cal S}_{8})\!: & \ -(e+f+cf)\left[(c + 2 c^2)^2 - (2 e + f)^2\right]=0.
\end{aligned}$$
We suggest the reader to plot surface $({\cal S}_{8})$ in order to visualize a three--dimensional picture. 

\medskip

\noindent \textbf{Bifurcation surface in $\mathbb{R}^3$ due to the infinite elliptic--saddle}

\medskip

\noindent {\bf (${\cal S}_{0}$)} Along the plane $c=-1$ the corresponding phase portraits possess an infinite 
singularity of the type $\widehat{\!{1\choose 2}\!\!}\ E-H$, which is the transition between the singularities 
$\widehat{\!{1\choose 2}\!\!}\ PEP-H$ and $\widehat{\!{1\choose 2}\!\!}\ E-PHP$. Such a plane is needed 
for the coherence of the bifurcation diagram. In fact, according to \cite{Artes-Llibre-Schlomiuk-Vulpe-2021a} 
we know that the comitant $\widetilde{N}$ is related to this phenomenon. Moreover, $\widetilde{N}$ 
``behaves like'' ${\cal{T}}_{4}$, in the sense that $\widetilde{N}=0$  splits the parameter space into two 
distinct canonical regions and the phase portrait over $\widetilde{N}=0$ is topologically equivalent to the 
phase portrait in one of its sides and topologically distinct to the one in the other side (see this phenomenon in 
\cite{Artes-Mota-Rezende-2021c}). In such a way we need to determine the points on the parameter space 
that verifies the equation $\widetilde{N}=0$. Calculations yield
$$\widetilde{N}=-4 (c+1)x^2.$$
It is clear that the plane $c+1=0$ verifies this equation. Therefore we define surface {\bf (${\cal S}_{0}$)} 
by the equation
$$({\cal S}_{0})\!:c+1=0.$$

\medskip

The bifurcation surfaces listed previously are all algebraic and they, except $({\cal S}_{4})$ and $({\cal S}_{8})$, 
are the bifurcation surfaces of singularities of family~\eqref{eq:nf-QES-A} in the parameter space. We shall detect 
other bifurcation surface not necessarily algebraic. In such a nonalgebraic surface the family has global connection 
of separatrices different from those given by $({\cal S}_{4})$ and $({\cal S}_{8})$. The equation of this 
bifurcation surface can only be determined approximately by means of numerical tools. Using arguments of 
continuity in the phase portraits we can prove the existence of this component not necessarily algebraic in the part 
where it appears, and we can check it numerically. We shall name it surface $({\cal S}_{7})$.

\begin{remark} Even though we can draw pictures of the algebraic bifurcation surfaces in $\mathbb{R}^3$, it is 
	pointless to see a single image of all these bifurcation surfaces together. As we shall see later, the partition of the 
	parameter space obtained from these bifurcation surfaces together with the nonalgebraic one has 1274 parts.
\end{remark}

Due to the last remark and, as we already said before, we shall foliate the three--dimensional bifurcation diagram in 
$\mathbb{R}^3$ by planes $c=c_0\ne0$, with $c_0$ constant and we shall give pictures of the resulting bifurcation 
diagram on these planar sections in which the Cartesian coordinates are $(e,f)$, where the horizontal line is the 
$e$--axis and $f\geq0$.

As the final bifurcation diagram is quite complex, it is useful to introduce colors which will be used to refer to the 
bifurcation surfaces:
\begin{enumerate}[(a)]
	\item surface (${\cal S}_{2}$) is drawn in green (coalescence of finite singular points);
	\item surface (${\cal S}_{3}$) is drawn in yellow (when the trace of a singular point becomes zero). We draw 
	it as a continuous curve if  the singular point is a focus or as a dashed curve if it is a saddle;
	\item surface (${\cal S}_{4}$) is drawn in purple (presence of at least one invariant straight line). We draw it as 
	a continuous curve if it implies a topological change or as a dashed curve otherwise;
	\item surface (${\cal S}_{6}$) is drawn in black (an antisaddle is on the edge of turning from a node to a focus 
	or vice versa). In the papers \cite{Artes-Llibre-Schlomiuk-2006,Artes-Rezende-Oliveira-2014,Artes-Rezende-Oliveira-2015,Artes-Mota-Rezende-2021c}
	the authors draw surface (${\cal S}_{6}$) as a continous curve. However, as it does not imply a topological change,
	we decided, from now on, to draw it as a dashed line.
	\item nonalgebraic surface (${\cal S}_{7}$) is also drawn in purple (connections of separatrices); and
	\item surface (${\cal S}_{8}$) is drawn in cyan (presence of an invariant parabola). We draw it as 
	a continuous curve if it implies a topological change or as a dashed curve otherwise.
\end{enumerate}

\begin{remark}\label{rmk-colors} Regarding the colors we use to draw the bifurcation surfaces, it is important to mention that:
	\begin{itemize}
		\item Here we use the same color for drawing (${\cal S}_{4}$) and (${\cal S}_{7}$), in order to follow the same
		pattern used in \cite{Artes-Rezende-Oliveira-2015,Artes-Mota-Rezende-2021c} for instance.
		\item In the mentioned papers surface (${\cal S}_{5}$) was drawn in red (when two infinite singular points coalesce). 
		However, for family~\eqref{eq:nf-QES-A} we are considering we saw that surface (${\cal S}_{5}$) defines the entire 
		plane $c=-2$. So, in order to avoid the utilization of several colors in the same plane, here we decided to follow the 
		pattern used in \cite{Artes-Mota-Rezende-2021b} and not to draw this entire plane in red color.
		\item In \cite{Artes-Mota-Rezende-2021c} the bifurcation line related to a presence of an infinite singular point of type 
		$\widehat{\!{1\choose 2}\!\!}\ E-H$ was drawing using brown color. However, for family~\eqref{eq:nf-QES-A} in the 
		current paper we saw that surface (${\cal S}_{0}$) defines the entire plane $c=-1$. Then, by the same reason 
		explained in the previous item, here we decided not to draw this plane in brown color.
	\end{itemize}			
\end{remark}

Having defined the bifurcation surfaces related to the study of the bifurcation diagram of family~\eqref{eq:nf-QES-A} 
we are now interested in studying the geometrical behavior of all of these algebraic surfaces for $c\ne0$, that is, 
their singularities, their intersection points and their extrema (maxima and minima) with respect to the coordinate $c$
(in other words, we have the ``tangencies'' with planes $c=c_0\ne0)$. 
Since this study requires a lot of computations which would take a very large number of pages to present all the details
(as in \cite{Artes-Rezende-Oliveira-2015,Artes-Mota-Rezende-2021c} for instance), in order to be more succinct 
here we are using the same algorithm (written in software Mathematica) already used in \cite{Artes-Mota-Rezende-2021b}.
Such an algorithm, applied to family~\eqref{eq:nf-QES-A}, is available for free download through the link 
\url{http://mat.uab.cat/~artes/articles/qvfES/qvfES-A.nb}
(some previous knowledge of Mathematica is recommended for using this algorithm).
In order to avoid repetitions, we recommend paper \cite{Artes-Mota-Rezende-2021b} for more details on the notation
used in this study and on the description and meaning of the so--called \textit{lists of objects}. 

\begin{remark}\label{param-space} In the papers \cite{Artes-Llibre-Schlomiuk-2006,Artes-Rezende-Oliveira-2014,Artes-Rezende-Oliveira-2015,Artes-Mota-Rezende-2021b,Artes-Mota-Rezende-2021c} in which families of quadratic systems were studied, the corresponding bifurcation diagram was done in an appropriate projective space, in which it was possible to analyze the slice at infinity and also to verify coherence in continuity (modulo islands) between the phase portraits on the infinite slice and phase portraits on the ``highest'' slice in the affine part. In those studies, with this approach the authors had the guarantee that they did not loose any phase portrait when one goes from the affine part towards the infinity. Due to the nature of normal form~\eqref{eq:nf-QES-A}, it is not possible to perform an analogous study for family $\QESA$. The next result presents all the algebraic values of the parameter $c$ corresponding to singular slices (or planes) in the bifurcation diagram and the greatest singular value of $c$ is $c=2$. In addition, in Proposition~\ref{prop:algebraic-values-of-c} we have that the first algebraic slice is given by $c=5$. So, taking into consideration the approach used in previous papers, we may say that in our study there is a possibility of finding a phase portrait in an slice corresponding to a value $c>5$, which would be topologically distinct to those ones obtained in the study of slices $c\leq5$. However, we believe that in case there exists such a different phase portrait in an slice $c>5$, this phase portrait would belong to a region bordered by nonalgebraic bifurcations due to connections of separatrices, since our study of the geometrical behavior of all algebraic surfaces showed that we do not have to consider any slice $c>5$.
\end{remark}

Its proof follows from the study done with the help of the mentioned algorithm.

\begin{lemma}\label{lemma:values-geom-study}
	Consider the algebraic bifurcation surfaces defined before. The study of their singularities, their intersection points, and their 
	tangencies with planes $c=c_0\ne0$ provides the following set of 12 singular values of the parameter $c$:
	{\small$$
		\left\{2,\sqrt{3},1,\frac{1}{\sqrt{3}},\frac{1}{2},-\frac{1}{2},-\frac{1}{\sqrt{3}},-\frac{2}{3},-1,-\frac{3}{2},-\sqrt{3},-2
		\right\}\!.
		$$}
\end{lemma}

Note that, when we obtained the differential equations that define family~\eqref{eq:nf-QES-A}, we proved that due to the 
symmetry on the bifurcation diagram, it is enough to consider the parameter $f\geq0$. So, apart from the previous study 
we have also to consider all the possible intersections of algebraic bifurcation curves that occur along $f=0$, since from 
such intersection points some open regions on $f>0$ could arise (or disappear). In the following result we present 
the values of the parameter $c$ in which there exist intersection of bifurcation surfaces along $f=0$. The proof is done on 
the mentioned Mathematica file and, as it is quite trivial, it is not presented here.

\begin{lemma}\label{lemma:values-intersec-f0}
	Consider the algebraic bifurcation surfaces defined before. When restricted to $f=0$, such surfaces have intersection points 
	on the planes corresponding to the following 13 values of the parameter $c$:
	{$$
		\left\{4,\sqrt{3}+2,2,1,\frac{1}{2},2-\sqrt{3},\frac{1}{4},-\frac{1}{4},-\frac{9}{16},-1,-\frac{16}{9},-2,-4
		\right\}\!.
		$$}
\end{lemma}

We shall consider the planes corresponding to these intersection points also as singular planes (in fact, the previous two lists 
have nonempty intersection). So we collect the values of the parameter $c$ obtained from Lemmas~\ref{lemma:values-geom-study} 
and \ref{lemma:values-intersec-f0} and, in the next result we present the complete list of algebraic singular planes corresponding 
to values of the parameter $c$.

\begin{proposition}\label{prop:algebraic-values-of-c}
	The full set of needed algebraic singular slices in the bifurcation diagram of family~\eqref{eq:nf-QES-A} is formed by 20 
	elements which correspond to the values of $c$ in \eqref{eq:algebraic-values-of-c}.
\end{proposition}
\begin{equation}\label{eq:algebraic-values-of-c}
\begin{gathered}
c_{1}=4, \ c_{3}=\sqrt{3}+2, \ c_{5}=2, \ c_{7}=\sqrt{3},\ 
c_{9}=1, \ c_{11}=\frac{1}{\sqrt{3}}, \ c_{13}=\frac{1}{2}, \ c_{15}=2-\sqrt{3}, \\
c_{17}=\frac{1}{4}, \ c_{21}=-\frac{1}{4}, \ c_{23}=-\frac{1}{2}, \ c_{25}=-\frac{9}{16}, \ 
c_{27}=-\frac{1}{\sqrt{3}}, \ c_{29}=-\frac{2}{3}, \ c_{31}=-1, \\ c_{33}=-\frac{3}{2},\
c_{35}=-\sqrt{3}, \ c_{37}=-\frac{16}{9}, \ c_{39}=-2, \ c_{41}=-4.\\
\end{gathered}
\end{equation}

The numeration in \eqref{eq:algebraic-values-of-c} is not consecutive  since we reserve numbers for generic slices. 
We point out that we have not found nonalgebraic slices, as in \cite{Artes-Mota-Rezende-2021c}, for instance.

In order to determine all the parts generated by the bifurcation surfaces from $({\cal S}_{0})$ to $({\cal S}_{8})$, 
we first draw the horizontal slices of the three--dimensional parameter space which correspond to the explicit values 
of $c$ obtained in Proposition~\ref{prop:algebraic-values-of-c}. However, as it will be discussed later, the presence 
of nonalgebraic bifurcation surfaces will be detected and their behavior as we move from slice to slice will be approximately 
determined. We add to each interval of singular values of $c$ an intermediate value for which we represent the bifurcation 
diagram of singularities. The diagram will remain essentially unchanged in these open intervals except the parts affected by 
the bifurcation. All the 42 sufficient values of $c$ are shown in \eqref{eq:values-of-c-QES-A}.

The values indexed by positive odd indices in \eqref{eq:values-of-c-QES-A} correspond to explicit values of $c$ for 
which there is a bifurcation in the behavior of the systems on the slices. Those indexed by even values are just intermediate 
points which are necessary to the coherence of the bifurcation diagram. Note that we skip index $c_{19}$ since such
an index would correspond to $c=0$, in which family~\eqref{eq:nf-QES-A} is not defined.

We now begin the analysis of the bifurcation diagram by studying completely one generic slice and after by moving
from slice to slice and explaining all the changes that occur. As an exact drawing of the curves produced by intersecting 
the surfaces with the slices gives us very small parts which are difficult to distinguish, and points of tangency are almost 
impossible to recognize, we have produced topologically equivalent figures where parts are enlarged and tangencies 
are easy to observe. From this reason, pictures corresponding to entire planes $(e,f)$ are split into two parts, see 
for instance Fig.~\ref{fig:slice-QES-A-01a} and \ref{fig:slice-QES-A-01b}.

\begin{equation}\label{eq:values-of-c-QES-A}
\hspace{-2mm}\arraycolsep=0.3cm\begin{array}{llllll}
c_{0}=5	          & c_{9}=1           & c_{18}=1/10 &  c_{27}=-1/\sqrt{3} & c_{36}=-175/100	 \vspace{0.1cm} \\ 
c_{1}=4 	      & c_{10}=3/4        & c_{19}=??? &  c_{28}=-62/100 & c_{37}=-16/9 \vspace{0.1cm} \\
c_{2}=385/100 	  & c_{11}=1/\sqrt{3} & c_{20}=-1/10 &   c_{29}=-2/3 &  c_{38}=-19/10       \vspace{0.1cm}\\
c_{3}=\sqrt{3}+2  & c_{12}=55/100 & c_{21}=-1/4 &   c_{30}=-85/100  & c_{39}=-2  \vspace{0.1cm} \\
c_{4}=3  	      & c_{13}=1/2 & c_{22}=-35/100 &    c_{31}=-1 &   c_{40}=-3     \vspace{0.1cm} \\
c_{5}=2 	      & c_{14}=38/100 & c_{23}=-1/2 &  c_{32}=-125/100 &   c_{41}=-4    \vspace{0.1cm} \\
c_{6}=185/100  	  & c_{15}=2-\sqrt{3} & c_{24}=-53/100 &   c_{33}=-3/2 &    c_{42}=-5     \vspace{0.1cm}	\\
c_{7}=\sqrt{3}    & c_{16}=26/100 & c_{25}=-9/16 &  c_{34}=-16/10 &                 \vspace{0.1cm}	\\
c_{8}=14/10   	  & c_{17}=1/4 & c_{26}=-57/100 &   c_{35}=-\sqrt{3} &        \vspace{0.1cm}	\\
\end{array}
%\arraycolsep=0.6cm\begin{array}{ll}
%c_{0}=5	 & c_{22}=-35/100 	 \vspace{0.1cm} \\ 
%c_{1}=4 	& c_{23}=-1/2  \vspace{0.1cm} \\
%c_{2}=385/100 	& c_{24}=-53/100           \vspace{0.1cm}\\
%c_{3}=\sqrt{3}+2  & c_{25}=-9/16    \vspace{0.1cm} \\
%c_{4}=3  	& c_{26}=-57/100           \vspace{0.1cm} \\
%c_{5}=2 	& c_{27}=-1/\sqrt{3}       \vspace{0.1cm} \\
%c_{6}=185/100  	& c_{28}=-62/100           \vspace{0.1cm}	\\
%c_{7}=\sqrt{3}   & c_{29}=-2/3                 \vspace{0.1cm}	\\
%c_{8}=14/10   	& c_{30}=-85/100           \vspace{0.1cm}	\\
%c_{9}=1  	& c_{31}=-1                    \vspace{0.1cm}	\\
%c_{10}=3/4   	& c_{32}=-125/100         	\vspace{0.1cm}	\\
%c_{11}=1/\sqrt{3}   & c_{33}=-3/2                 	\vspace{0.1cm}	\\
%c_{12}=55/100  	& c_{34}=-16/10             \vspace{0.1cm}	\\
%c_{13}=1/2   & c_{35}=-\sqrt{3}         \vspace{0.1cm} \\
%c_{14}=38/100  	& c_{36}=-175/100        \vspace{0.1cm} \\
%c_{15}=2-\sqrt{3}      			& c_{37}=-16/9              \vspace{0.1cm} \\
%c_{16}=26/100  								& c_{38}=-19/10 						\vspace{0.1cm}	\\
%c_{17}=1/4     									& c_{39}=-2                   \vspace{0.1cm}	\\
%c_{18}=1/10    									& c_{40}=-3         					\vspace{0.1cm} \\
%c_{20}=-1/10 										& c_{41}=-4                   \vspace{0.1cm} \\
%c_{21}=-1/4 											& c_{42}=-5            												           \\
%\end{array}
\end{equation}

The reader may find the exact pictures of the 20 singular slices (containing only the algebraic surfaces) described in 
\eqref{eq:algebraic-values-of-c} in a PDF file available at the link \url{http://mat.uab.es/~artes/articles/qvfES/qvfES-A.pdf}. 

We now describe the labels used for each part of the bifurcation space. As we have mentioned in Remark~\ref{rem:label-p-p}, 
the subsets of dimensions 3, 2, 1 and 0, of the partition of the parameter space will be denoted respectively by $V$, $S$, $L$, 
and $P$ for Volume, Surface, Line and Point, respectively. The surfaces are named using a number which corresponds to each 
bifurcation surface which is placed on the left side of the letter $S$. To describe the portion of the surface we place an index. 
The curves that are intersection of surfaces are named by using their corresponding numbers on the left side of the letter $L$, 
separated by a point. To describe the segment of the curve we place an index. Volumes and Points are simply indexed 
(since three or more surfaces may be involved in such an intersection).

We consider an example: surface $({\cal S}_{2})$ splits into 42 different two--dimensional parts labeled from $2S_{1}$ to 
$2S_{42}$, plus some one--dimensional arcs labeled as $2.iL_{j}$ (where $i$ denotes the other surface intersected by 
$({\cal S}_{2})$ and $j$ is a number), and some zero--dimensional parts. In order to simplify the labels in all figures we 
see \textbf{V1} which stands for the {\TeX} notation $V_{1}$. Analogously, \textbf{2S1} (respectively, \textbf{2.3L1}) 
stands for $2S_{1}$ (respectively, $2.3L_{1}$), see Fig.~\ref{fig:slice-QES-A-01a} and \ref{fig:slice-QES-A-01b}, for 
example.

\begin{figure}[h!]
	\centering\includegraphics[width=1\textwidth]{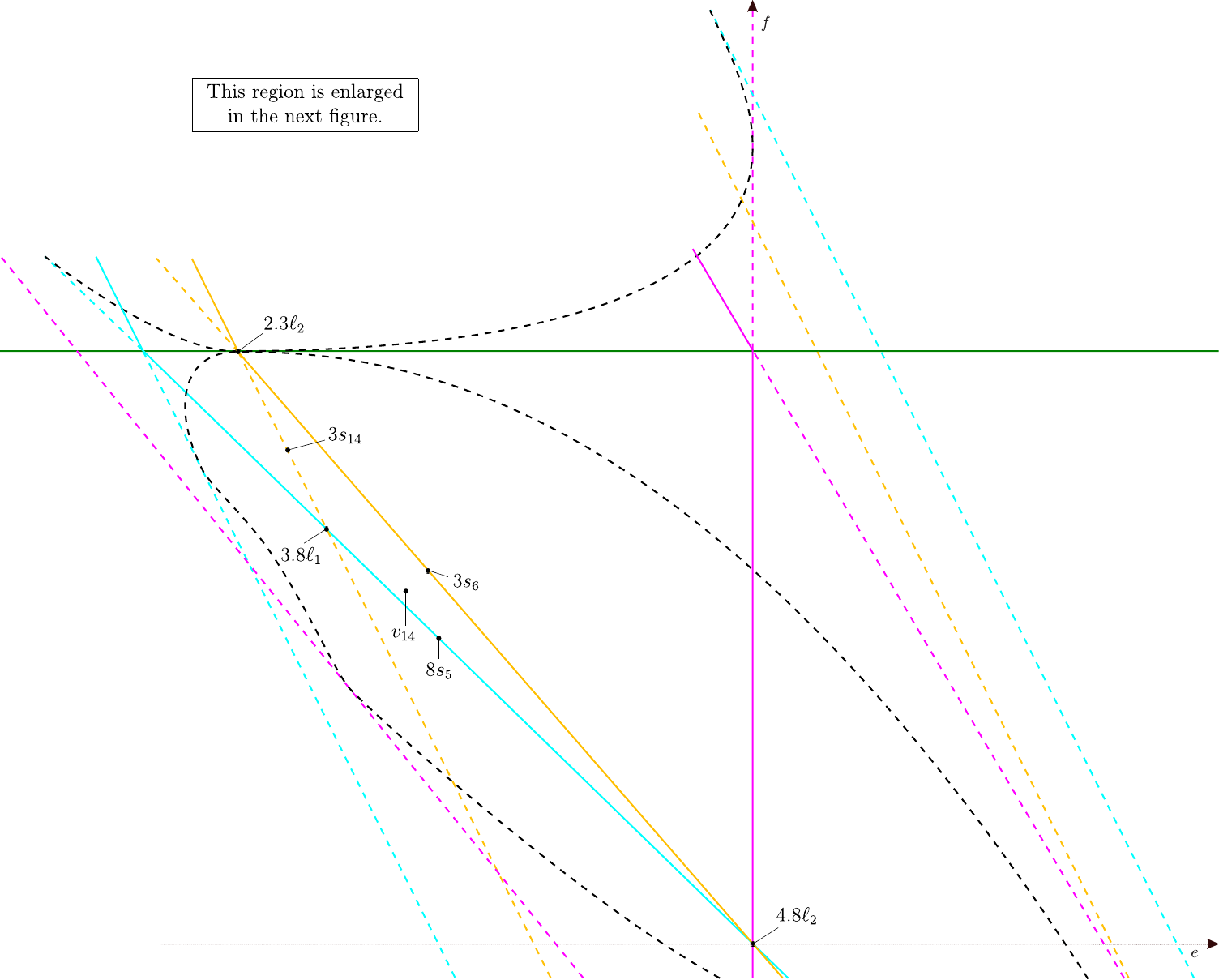} 
	\caption{\small \label{fig:slice-QES-A-01-alg1} Piece of generic slice of the parameter space when $c=5$ (only algebraic 
		surfaces), see also Fig.~\ref{fig:slice-QES-A-01-alg2}}
\end{figure}

\begin{figure}[h!]
	\centering\includegraphics[width=1\textwidth]{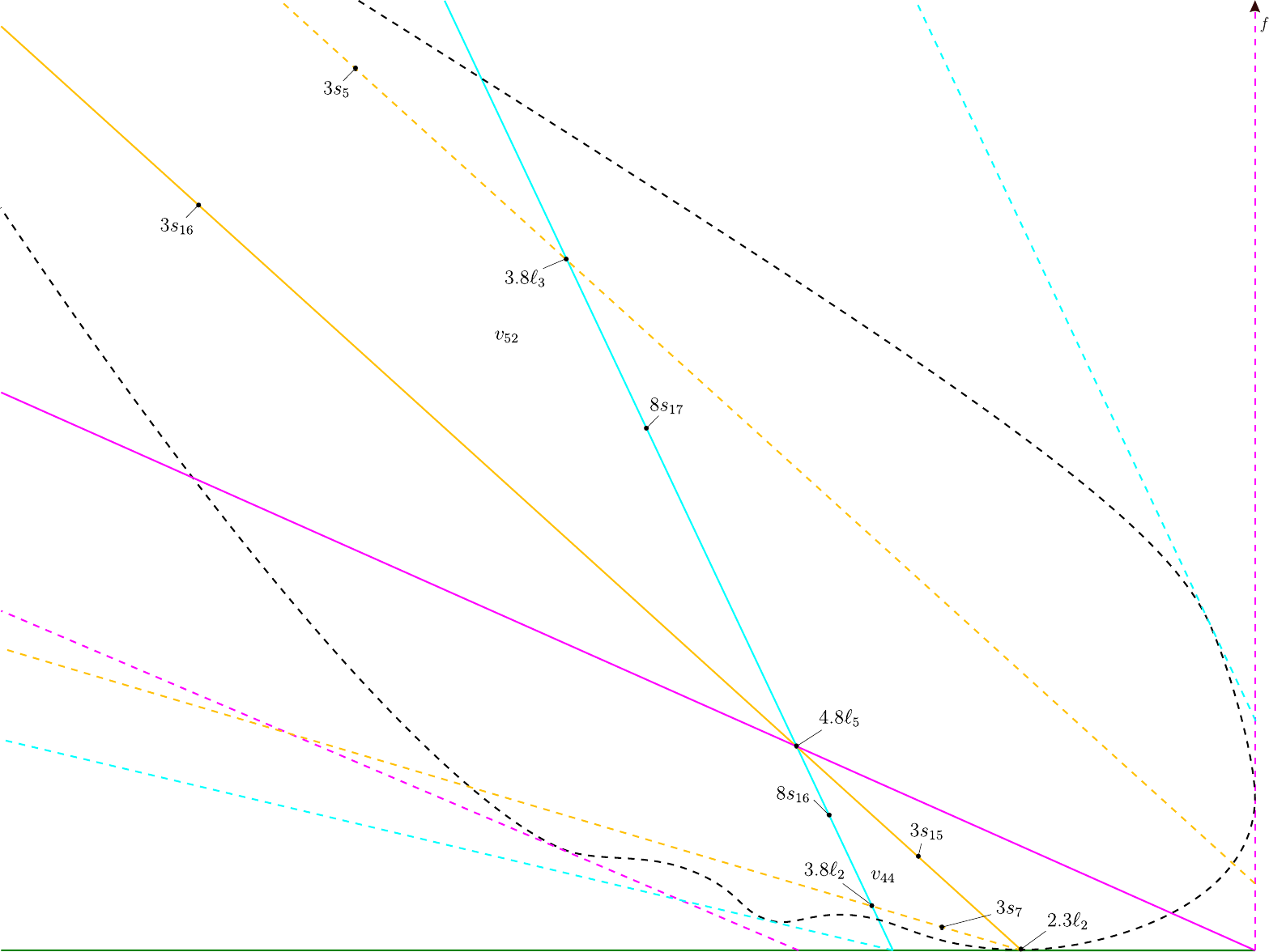} 
	\caption{\small \label{fig:slice-QES-A-01-alg2} Continuation of Fig.~\ref{fig:slice-QES-A-01-alg1}}
\end{figure}

In Fig.~\ref{fig:slice-QES-A-01-alg1} and \ref{fig:slice-QES-A-01-alg2} we represent the generic slice of the 
parameter space when $c=c_{0}=5$, showing only the algebraic surfaces. We note that there are some dashed 
branches of surface (${\cal{S}}_{3}$) (in yellow),  (${\cal{S}}_{4}$) (in purple), and (${\cal{S}}_{8}$) (in blue). 
This means the existence of a weak saddle, in the case of surface (${\cal{S}}_{3}$), the existence of an invariant 
straight line without separatrix connection, in the case of surface (${\cal{S}}_{4}$), and the existence of an invariant 
parabola without separatrix connection, in the case of surface (${\cal{S}}_{8}$); they do not mean a topological 
change in the phase portraits but a $C^\infty$ change. In the next figures we shall use the same representation for 
these characteristics of these three surfaces.

With the purpose to explain all the changes in the bifurcation diagram, we would have to present two versions of the 
picture of each slice: one of them without labels and the other with labels in each new part (as it was done, for instance, 
in \cite{Artes-Rezende-Oliveira-2013b} and \cite{Artes-Rezende-Oliveira-2014}).

However, as the number of slices is considerably large (see equation~\eqref{eq:values-of-c-QES-A} -- 42 slices 
to be more precise) we would have to present 84 pictures, which would occupy a large number of pages. Then, we 
shall present only the labeled drawings (just the ``important part'' in each slice) containing the algebraic and nonalgebraic 
bifurcation surfaces. Along this study we prove the existence of such nonalgebraic surfaces and their necessity for the 
coherence of the bifurcation diagram.

\begin{remark}\label{rem:f-n}
	Wherever two parts of equal dimension $d$ are separated only by a part of dimension $d-1$ of the black bifurcation 
	surface $({\cal S}_{6})$, their respective phase portraits are topologically equivalent since the only difference between 
	them is that a finite antisaddle has turned into a focus without change of stability and without appearance of limit cycles. 
	We denote such parts with different labels, but we do not give specific phase portraits in pictures attached to 
	Theorems~\ref{th:main-thm-QES-A}, \ref{th:main-thm-QES-B}, and \ref{th:main-thm-QES-C},
	for the parts with the focus. We only give portraits for the parts with nodes, except in the case of existence of a limit cycle 
	or a graphic where the singular point inside them is portrayed as a focus. Neither do we give specific invariant description 
	in Sec.~\ref{sec:invariants-QESA} distinguishing between these nodes and foci.
\end{remark}

Now we explain the generic slice when $c=5$ presented in Fig.~\ref{fig:slice-QES-A-01-alg1} and \ref{fig:slice-QES-A-01-alg2}. 
In this slice we shall make a complete study of all its parts, whereas in the next slices we only describe the changes. 
Some singular slices will produce only few changes which are easy to describe, but others can produce simultaneously 
many changes, even a complete change of all parts and these will require a more detailed description.

As we said before, in the mentioned figures we present the slice when $c=5$ with only the algebraic surfaces. We now 
place for each set of the partition on this slice the local behavior of the flow around the singular points. For a specific 
value of the parameters of each one of the sets in this partition we compute the global phase portrait with the 
numerical program P4 \cite{progP4,Dumortier-Llibre-Artes-2006}.

In this slice we have a partition in two--dimensional parts bordered by curved polygons, some of them bounded, 
others bordered by infinity. From now on, we use lower--case letters provisionally to describe the sets found 
algebraically in order to do not interfere with the final partition described with capital letters. 

For each two--dimensional part we obtain a phase portrait which is coherent with those of all their borders. 
Except for three parts, which are shown in Fig.~\ref{fig:slice-QES-A-01-alg1} and \ref{fig:slice-QES-A-01-alg2}
and named as follows:
\begin{itemize}\label{page:regions-c-5}
	\item $v_{14}$: the triangle bordered by yellow and blue curves (in Fig.~\ref{fig:slice-QES-A-01-alg1});
	\item $v_{44}$: the triangle bordered by yellow and blue curves (in Fig.~\ref{fig:slice-QES-A-01-alg2});
	\item $v_{52}$: the quadrilateral bordered by yellow and blue curves and infinity (in Fig.~\ref{fig:slice-QES-A-01-alg2}).
\end{itemize}
The study of these parts is important for the coherence of the bifurcation diagram. That is why we have 
decided to present only these parts in the mentioned figures.

We begin with the analysis of part $v_{14}$. We consider the segment $3s_{6}$ in Fig.~\ref{fig:slice-QES-A-01-alg1}, 
which is one of the borders of part $v_{14}$. On this segment, the corresponding phase portrait possesses 
a weak focus (of order one) and, consequently, this branch of surface (${\cal{S}}_{3}$) corresponds to a 
Hopf bifurcation. This means that the phase portrait corresponding to one of the sides of this segment must 
have a limit cycle; in fact it is in the triangle $v_{14}$.

\begin{figure}[h!]
	\centering\includegraphics[width=0.7\textwidth]{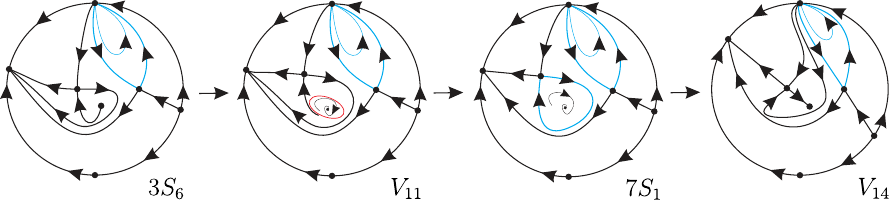} 
	\caption{\small \label{fig:bif-7S1} Sequence of phase portraits in parts $V_{11}$ and $V_{14}$ of slice 
		$c=5$ (the labels are according to Fig.~\ref{fig:slice-QES-A-01a})}
\end{figure}

\begin{figure}[h!]
	\centering\includegraphics[width=1\textwidth]{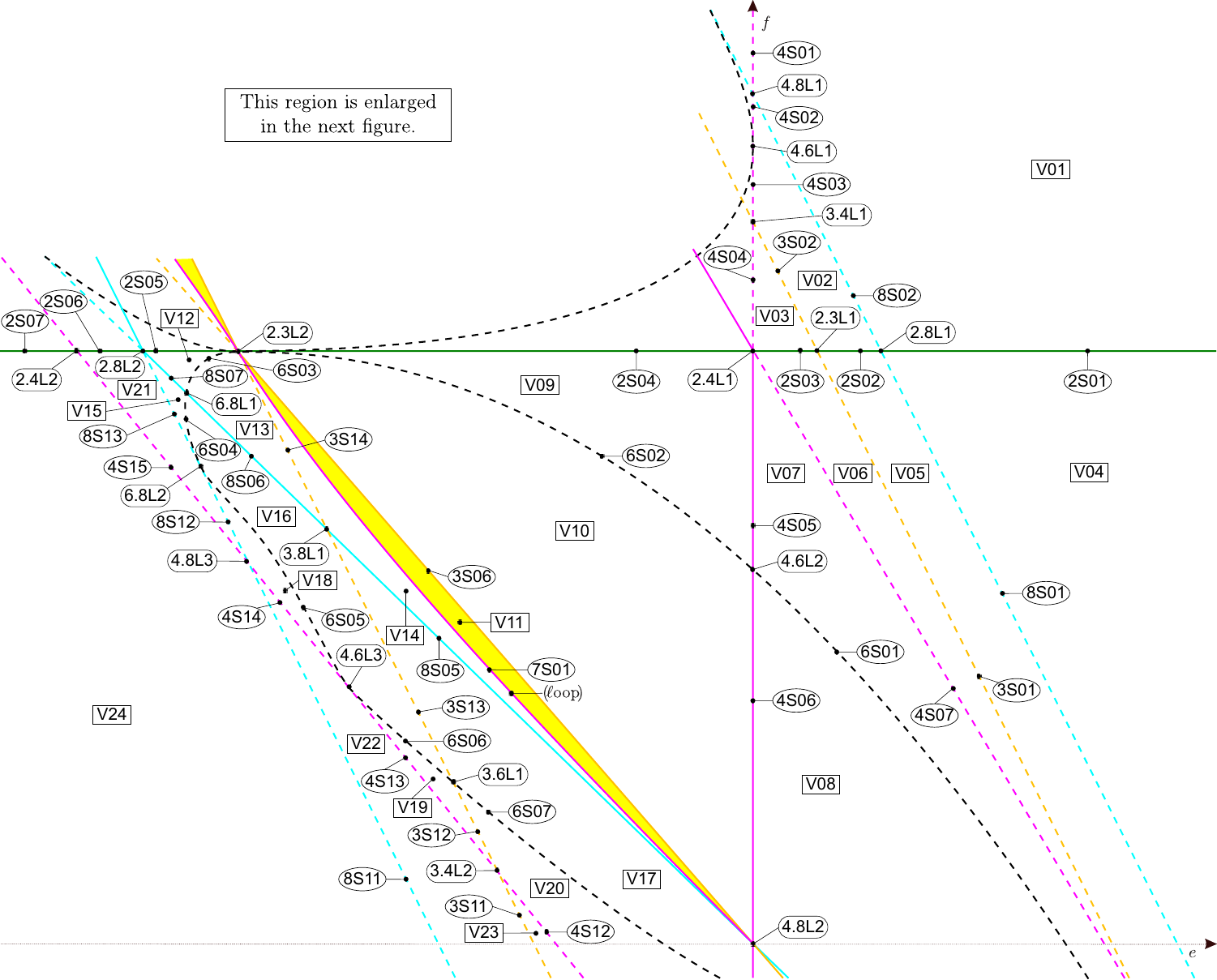} 
	\caption{\small \label{fig:slice-QES-A-01a} Piece of generic slice of the parameter space when $c=5$, 
		see also Fig.~\ref{fig:slice-QES-A-01b}}
\end{figure}

\begin{figure}[h!]
	\centering\includegraphics[width=1\textwidth]{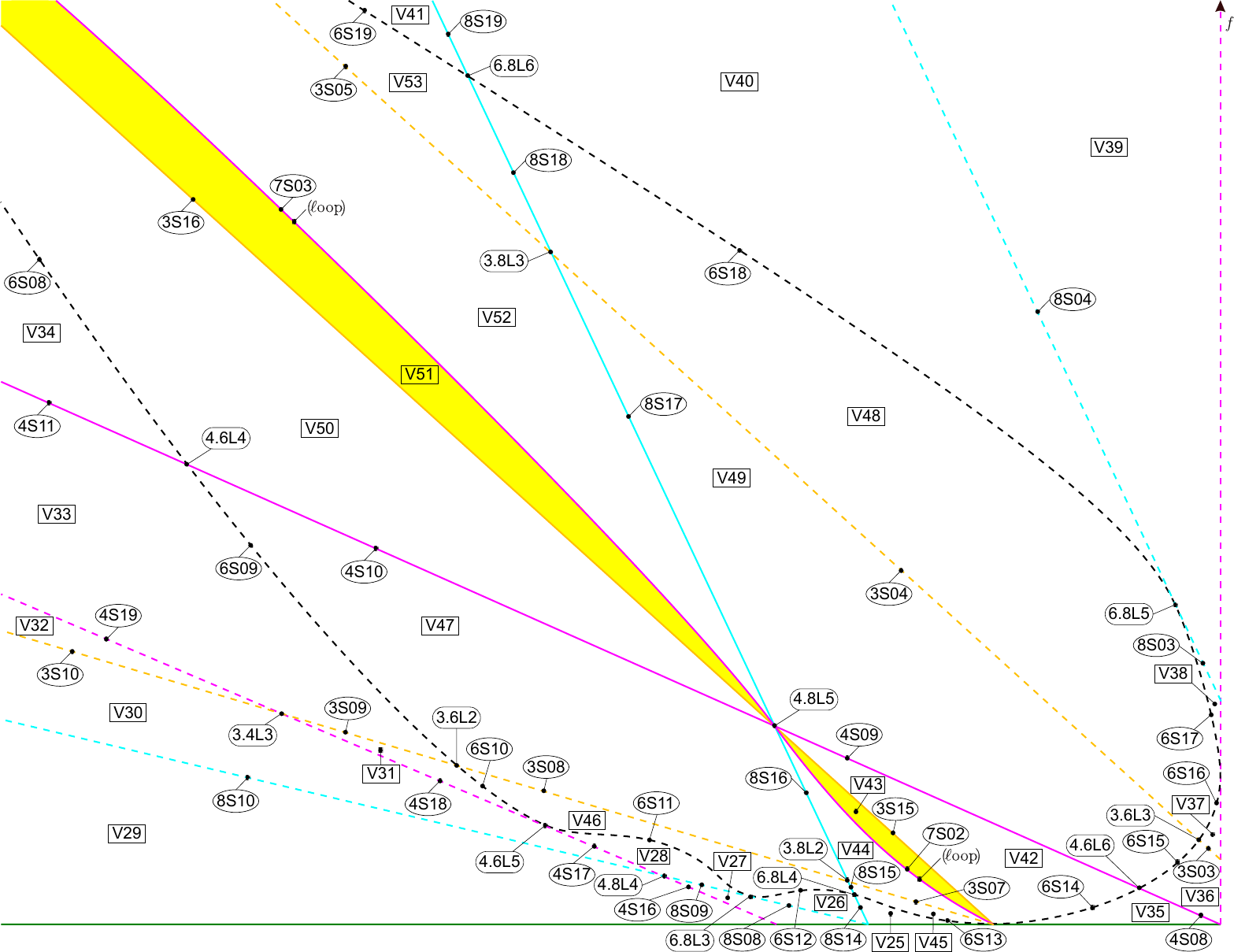} 
	\caption{\small \label{fig:slice-QES-A-01b} Continuation of Fig.~\ref{fig:slice-QES-A-01a}}
\end{figure}

However, when we get close to $8s_{5}$ and $3s_{14}$, the limit cycle has been lost, which implies the 
existence of at least one element of surface (${\cal{S}}_{7}$) (see $7S_{1}$ in Fig.~\ref{fig:slice-QES-A-01a}), 
in a neighborhood of $3s_{6}$, due to a connection of separatrices from a saddle to itself (i.e. a loop--type connection). 
In Lemma \ref{lem:endpoints-7S1} we show that $7S_{1}$ is bounded and it has its endpoints at the curves 
$4.8\ell_{2}$ and $2.3\ell_{2}$. We draw the sequence of phase portraits along these subsets (using the notation
from Fig.~\ref{fig:slice-QES-A-01a}) in Fig.~\ref{fig:bif-7S1} and we plot the complete bifurcation diagram 
for this part in Fig.~\ref{fig:slice-QES-A-01a}.

\begin{lemma}\label{lem:endpoints-7S1}
	The nonalgebraic curve $7S_{1}$  is bounded and it has its endpoints at the curves 
	$4.8\ell_{2}$ and $2.3\ell_{2}$.
\end{lemma}

\begin{proof} Numerical analysis indicate the veracity of the result. Indeed, note that if one of the 
	endpoints of this surface is any point of $3s_{6}$, then a portion of this subset must not refer to a 
	Hopf bifurcation, which contradicts the fact that on $3s_{6}$ we have a weak focus of order one.
	Also, observe that it is not possible that the starting point of these surfaces is on $3s_{14}$, since 
	on this portion of the yellow surface we have only a $C^\infty$ bifurcation (weak saddle).
	Finally, the endpoints cannot be on $8s_{5}$ because, in order to have this, first we need to break 
	the invariant parabola. Then, the only possible endpoints of surface $7S_{1}$ are $4.8\ell_{2}$ and $2.3\ell_{2}$.
\end{proof}

Now we consider parts $v_{44}$ and $v_{52}$ in Fig.~\ref{fig:slice-QES-A-01-alg2}. When 
are very close to the yellow curves $3s_{15}$ and $3s_{16}$ we have the existence of a 
limit cycle in the phase portraits corresponding to parts $v_{44}$ and $v_{52}$, respectively. 
However, when we move away from these yellow curves we observe that the limit cycles disappear.
So there exist at least one element of surface (${\cal{S}}_{7}$) (see $7S_{2}$ and $7S_{3}$ in 
Fig.~\ref{fig:slice-QES-A-01b}), in a neighborhood of $3s_{15}$ and $3s_{16}$, respectively, 
due to a loop--type connection. In fact, numerical verification shows the existence of such nonalgebraic
surfaces. Moreover, as we have that:
\begin{itemize}
	\item $3s_{6}$, $3s_{15}$, and $3s_{16}$ provide topologically equivalent phase portraits,
	\item $3s_{14}$, $3s_{7}$, and $3s_{5}$ provide topologically equivalent phase portraits,
	\item $8s_{5}$, $8s_{16}$, and $8s_{17}$ provide topologically equivalent phase portraits, and
	\item $4.8\ell_{2}$ and $4.8\ell_{5}$ provide topologically equivalent phase portraits, 
\end{itemize}
from the analysis we made from region $v_{14}$ it is easy to conclude the following result.

\begin{lemma}\label{lem:endpoints-7S2-and-7S3}
	The nonalgebraic curves $7S_{2}$ and $7S_{3}$ are continuation of $7S_{1}$. Moreover,
	$7S_{2}$ is bounded and it has its endpoints at $4.8\ell_{5}$ and $2.3\ell_{2}$, and $7S_{3}$ 
	is not bounded and starts at $4.8\ell_{5}$.
\end{lemma}

The complete bifurcation diagram for this part can be seeing in Fig.~\ref{fig:slice-QES-A-01b}.

Regarding Remark~\ref{weak-singularities}, item 1, in equation \eqref{eq:sol-weak-sing} we
obtained regions of the parameter space in which the corresponding phase portrait possesses 
center type singular point. The regions we obtained in that equation correspond to the curves 
$4.8L_{2}$ (see Fig.~\ref{fig:slice-QES-A-01a}) and $4.8L_{5}$ (see Fig.~\ref{fig:slice-QES-A-01b}),
respectively.

We have added in the bifurcation diagram a label associated to each part of the bifurcation 
(${\cal{S}}_7$) indicating the type of connection produced by this bifurcation. More precisely,
in the pictures where it appears ``($\l$oop)'' we are indicating this type of separatrix connection.

Having analyzed all the parts pointed out on page \pageref{page:regions-c-5} and explained the 
existence of all possible nonalgebraic surfaces in there (modulo islands), we have finished the 
study of the generic slice $c=5$. However, we cannot be sure that these are all the additional 
bifurcation surfaces in this slice. There could exist others which are closed surfaces small enough 
to escape our numerical research. For all other two--dimensional parts of the partition of this slice, 
whenever we join two points which are close to different borders of the part, the two phase portraits 
are topologically equivalent. So, we do not encounter more situations than the ones mentioned before. 
In short, it is expected that the complete bifurcation diagram for $c=5$ is the one shown in 
Fig.~\ref{fig:slice-QES-A-01a} and \ref{fig:slice-QES-A-01b}. In these and in the next figures 
we have colored in light yellow the open regions with one limit cycle, in black the labels referring 
to new parts which are created in a slice and in red the labels corresponding to parts which has already 
appeared in previous slices. 

Due to the computation we mentioned before, we already know that there are no more singular 
slices for $c>5$. Moreover, as we discussed in Remark~\ref{param-space}, 
because normal form~\eqref{eq:nf-QES-A}
does not allow the study of the slice at infinity, we cannot guarantee that for $c>5$ it does not exist 
a nonalgebraic singular slice. So, having finished the complete study of slice $c=5$, the next step is to decrease 
the values of $c$, according to equation~\eqref{eq:values-of-c-QES-A}, and make an analogous 
study for each one of the slices that we need to consider and also search for changes when going 
from one slice to the next one.

We now start decreasing the values of the parameter $c$ in order to explain as much as we can 
the bifurcations in the parameter space.

\begin{figure}[h!]
	\centering\includegraphics[width=0.6\textwidth]{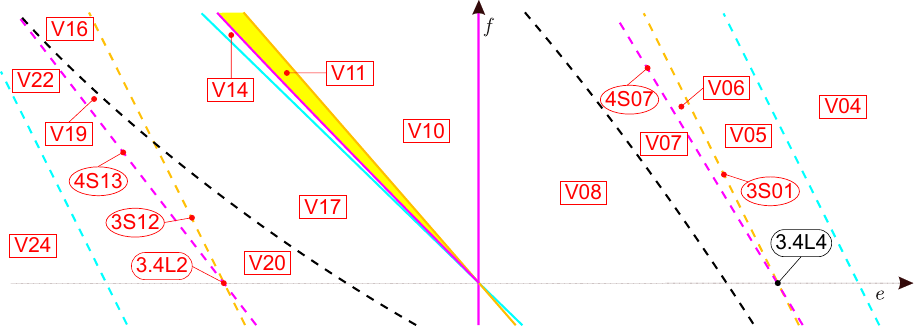} 
	\caption{\small \label{fig:slice-QES-A-02} Piece of singular slice of the parameter space when $c=4$}
\end{figure}

Consider Fig.~\ref{fig:slice-QES-A-01a}. When we move down from $c=5$ to $c=4$ (a singular slice) the 
curve $3.4L_{2}$ goes to $f=0$ and the bifurcation curves $3S_{1}$ and $4S_{7}$ intersect themselves
on $f=0$, more precisely, at $3.4L_{4}$, see Fig.~\ref{fig:slice-QES-A-02}.

\begin{figure}[h!]
	\centering\includegraphics[width=0.6\textwidth]{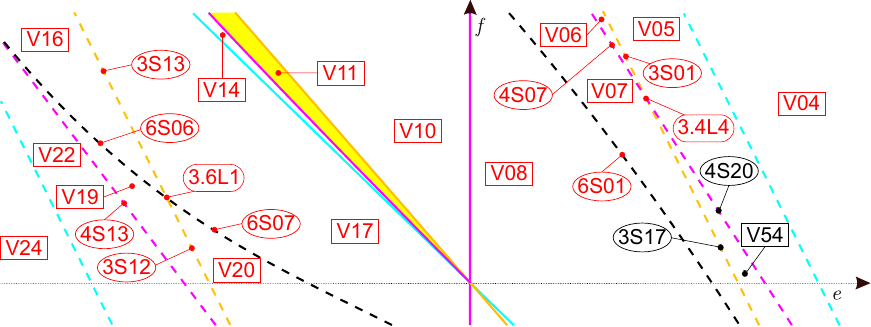} 
	\caption{\small \label{fig:slice-QES-A-03} Piece of generic slice of the parameter space when $c=385/100$}
\end{figure}

Taking $c=385/100$ we observe that $3.4L_{2}$ goes to $f<0$ and from $3.4L_{4}$ it arises the volume
region $V_{54}$, see Fig.~\ref{fig:slice-QES-A-03}.

\begin{figure}[h!]
	\centering\includegraphics[width=0.6\textwidth]{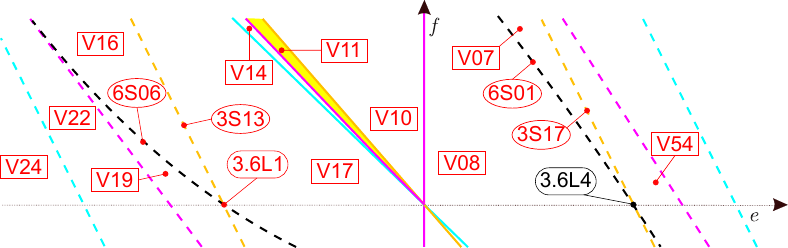} 
	\caption{\small \label{fig:slice-QES-A-04} Piece of singular slice of the parameter space when $c=2+\sqrt{3}$}
\end{figure}

When $c=2+\sqrt{3}$ we have that $3.6L_{1}$ goes to $f=0$ and the bifurcation curves $3S_{17}$ 
and $6S_{1}$ intersect themselves on $f=0$, more precisely, at $3.6L_{4}$, see Fig.~\ref{fig:slice-QES-A-04}.

\begin{figure}[h!]
	\centering\includegraphics[width=1\textwidth]{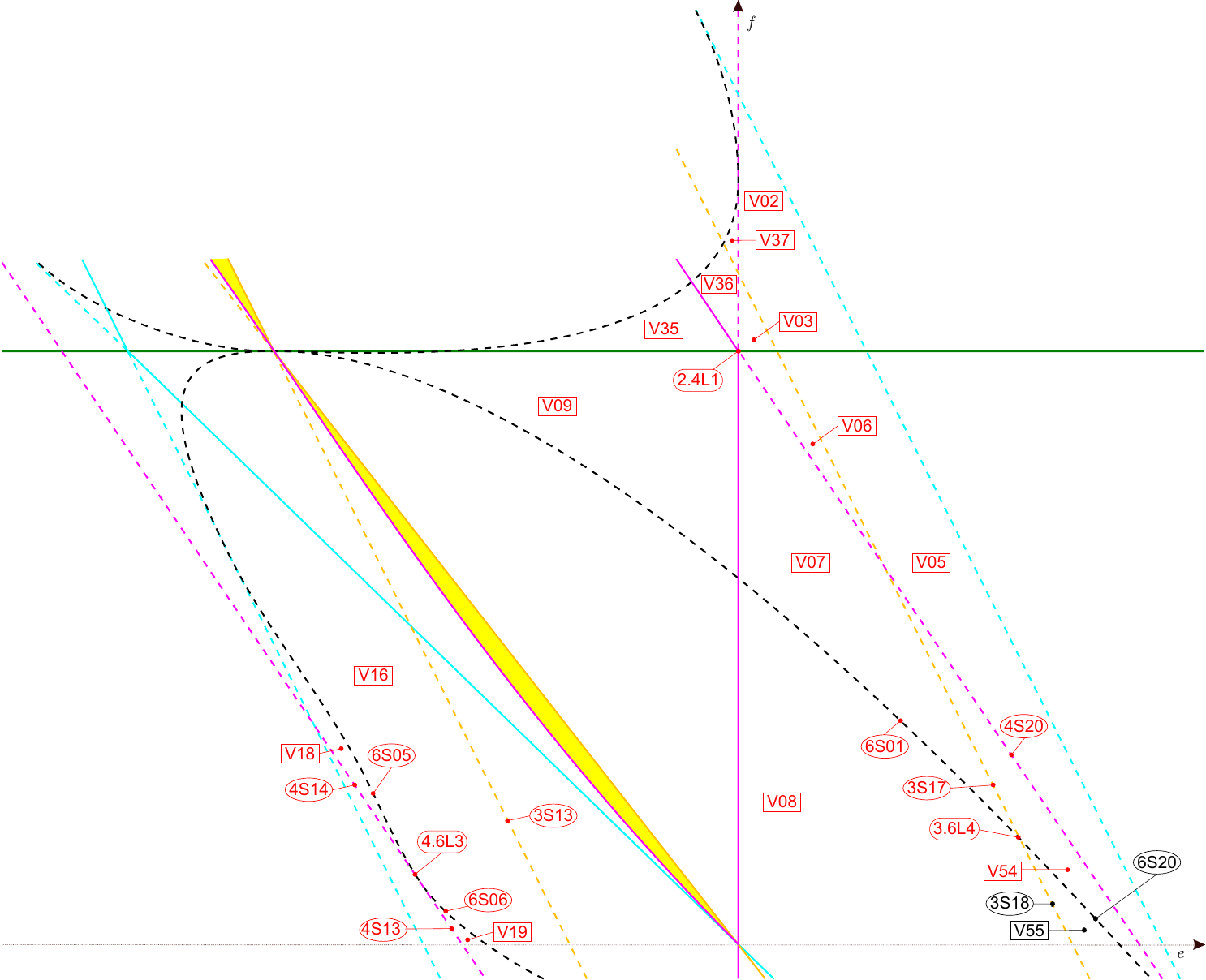} 
	\caption{\small \label{fig:slice-QES-A-05} Piece of generic slice of the parameter space when $c=3$}
\end{figure}

When we consider $c=3$ we notice that $3.6L_{1}$ goes to $f<0$ and from $3.6L_{4}$ it arises the volume
region $V_{55}$. In Fig.~\ref{fig:slice-QES-A-05} we present a piece of this generic slice, where we label
the mentioned regions and also another regions that appear in the sequence.

\begin{figure}[h!]
	\centering\includegraphics[width=1\textwidth]{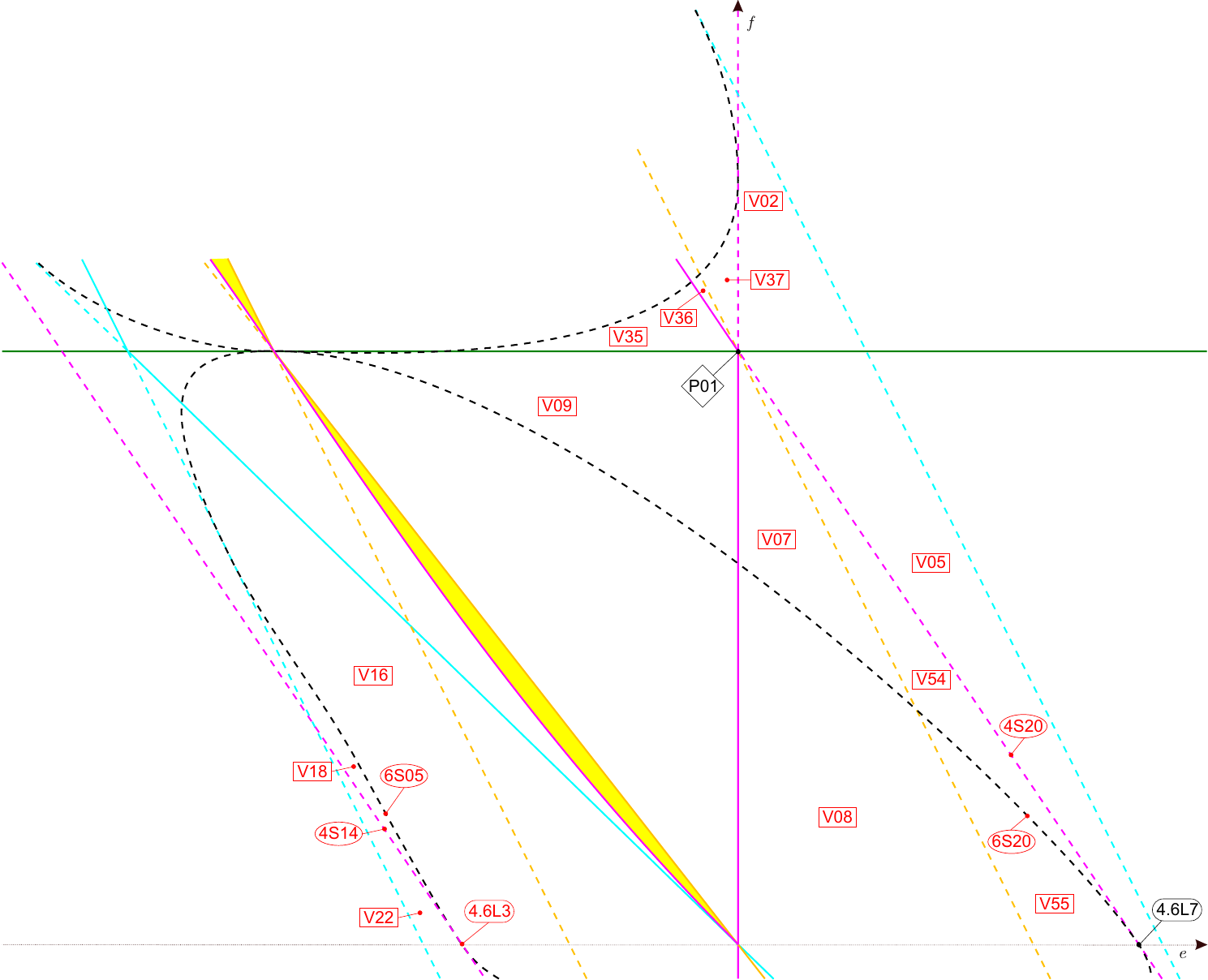} 
	\caption{\small \label{fig:slice-QES-A-06a} Piece of singular slice of the parameter space when $c=2$, 
		see also Fig.~\ref{fig:slice-QES-A-06b}}
\end{figure}

\begin{figure}[h!]
	\centering\includegraphics[width=0.35\textwidth]{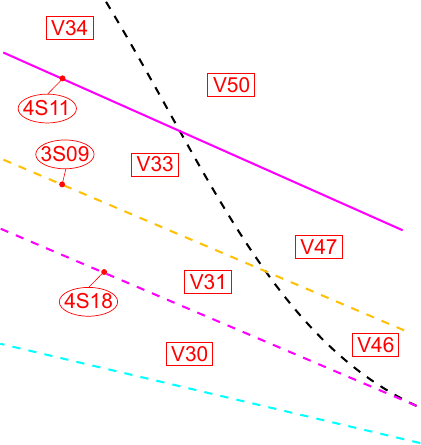} 
	\caption{\small \label{fig:slice-QES-A-06b} Another piece of singular slice of the parameter space when $c=2$, 
		compare this region with Fig.~\ref{fig:slice-QES-A-01b}}
\end{figure}

Consider Fig.~\ref{fig:slice-QES-A-05}. When we study the singular slice $c=2$ we observe that:
\begin{itemize}
	\item the triangles $V_{3}$ and $V_{6}$ coalesce at $2.4L_{1}$, generating point $P_{1}$;
	\item bifurcation curve $6S_{20}$ intercepts $4S_{20}$ at $4.6L_{7}$ (on $f=0$); and
	\item $4.6L_{3}$ goes to $f=0$, making $V_{19}$ go to $f<0$.
\end{itemize}
Also, by considering Fig.~\ref{fig:slice-QES-A-01b}, we note that when $c=2$ the bifurcation
straight lines $3S_{9}$, $4S_{11}$, and $4S_{18}$ are parallel, making both $3.4L_{3}$ and $V_{32}$ 
go to infinity. The singular slice under consideration is presented in Fig.~\ref{fig:slice-QES-A-06a} and 
\ref{fig:slice-QES-A-06b}, in which we label only the regions that are relevant in this slice.

\begin{figure}[h!]
	\centering\includegraphics[width=1\textwidth]{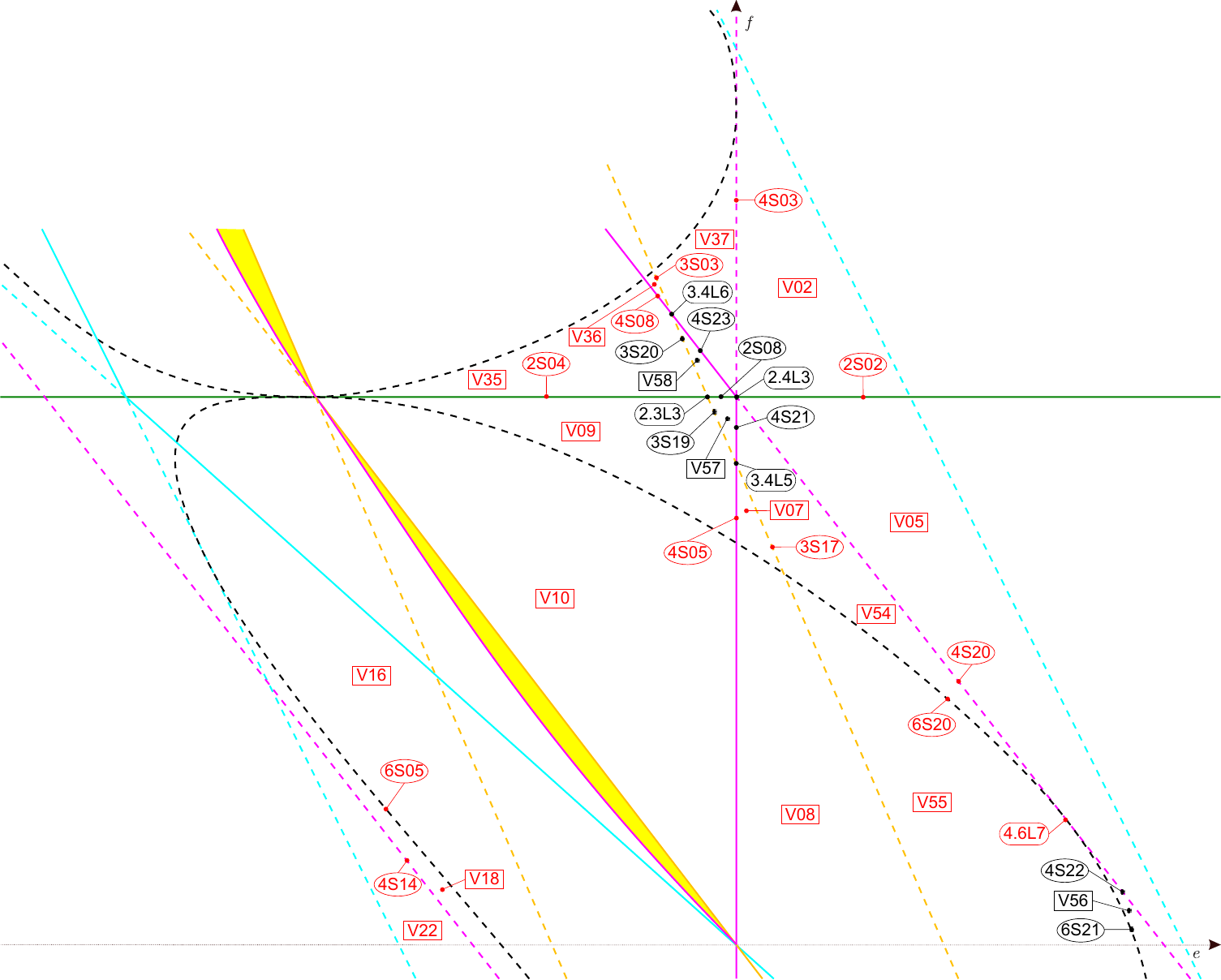} 
	\caption{\small \label{fig:slice-QES-A-07a} Piece of generic slice of the parameter space when $c=185/100$, 
		see also Fig.~\ref{fig:slice-QES-A-07b}}
\end{figure}

\begin{figure}[h!]
	\centering\includegraphics[width=0.35\textwidth]{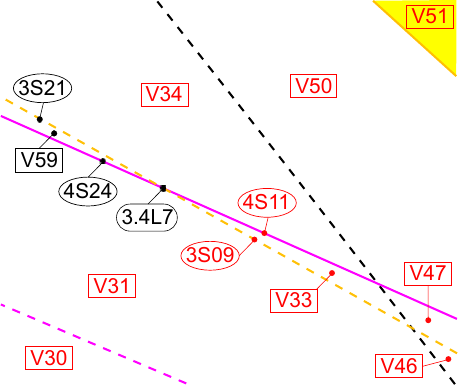} 
	\caption{\small \label{fig:slice-QES-A-07b} Another piece of generic slice of the parameter space when $c=185/100$, 
		compare this region with Fig.~\ref{fig:slice-QES-A-06b}}
\end{figure}

Now we consider the generic slice $c=185/100$. By studying completely this slice we observe that:
\begin{itemize}
	\item $4.6L_{3}$ goes to $f<0$;
	\item $4.6L_{7}$ goes to $f>0$ and it arises volume region $V_{56}$;
	\item from point $P_{1}$ we get two new volume regions, namely, $V_{57}$ and $V_{58}$;
\end{itemize}
see Fig.~\ref{fig:slice-QES-A-07a}. Moreover, we have that the yellow straight line $3S_{9}$ now
intercepts $4S_{11}$ at $3.4L_{7}$ and it arises volume region $V_{59}$, see Fig.~\ref{fig:slice-QES-A-07b}.

\begin{figure}[h!]
	\centering\includegraphics[width=0.25\textwidth]{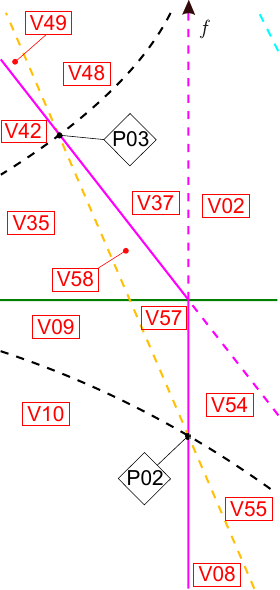} 
	\caption{\small \label{fig:slice-QES-A-08a} Piece of singular slice of the parameter space when $c=\sqrt{3}$, 
		compare this region with Fig.~\ref{fig:slice-QES-A-07a} and see also Fig.~\ref{fig:slice-QES-A-08b}}
\end{figure}

\begin{figure}[h!]
	\centering\includegraphics[width=0.35\textwidth]{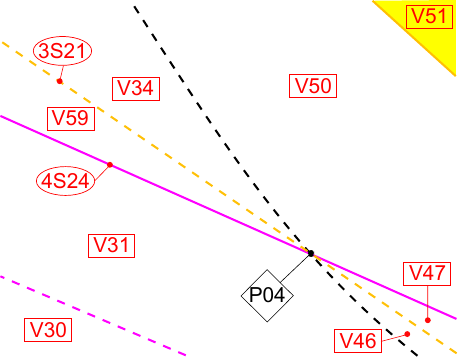} 
	\caption{\small \label{fig:slice-QES-A-08b} Another piece of singular slice of the parameter space when $c=\sqrt{3}$, 
		compare this region with Fig.~\ref{fig:slice-QES-A-07b}}
\end{figure}

When we move down and consider the singular slice $c=\sqrt{3}$ we note that the volume regions $V_{7}$ and 
$V_{36}$ are reduced to the points $P_{2}$ and $P_{3}$, respectively (see Fig.~\ref{fig:slice-QES-A-08a}).
We also have that at this value of the parameter $c$ the volume region $V_{33}$ is reduced to the point $P_{4}$,
which can be seeing in Fig.~\ref{fig:slice-QES-A-08b}.

\begin{figure}[h!]
	\centering\includegraphics[width=0.35\textwidth]{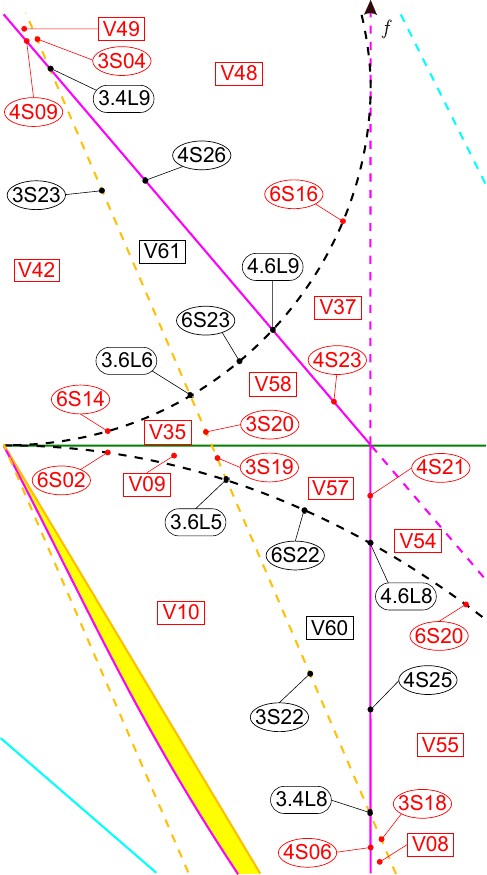} 
	\caption{\small \label{fig:slice-QES-A-09a} Piece of generic slice of the parameter space when $c=14/10$, 
		compare this region with Fig.~\ref{fig:slice-QES-A-08a} and see also Fig.~\ref{fig:slice-QES-A-09b}}
\end{figure}

\begin{figure}[h!]
	\centering\includegraphics[width=0.35\textwidth]{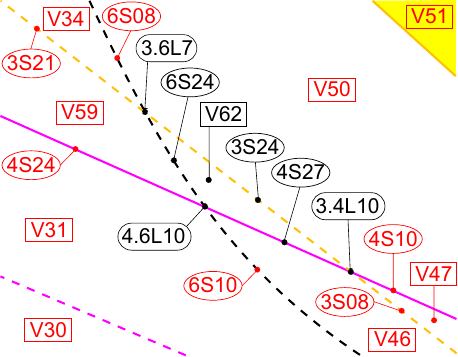} 
	\caption{\small \label{fig:slice-QES-A-09b} Another piece of generic slice of the parameter space when $c=14/10$, 
		compare this region with Fig.~\ref{fig:slice-QES-A-08b}}
\end{figure}

During the study of the generic slice $c=14/10$ we observe that from the points $P_{2}$ and $P_{3}$ arise
the volume regions $V_{60}$ and $V_{61}$, respectively (see Fig.~\ref{fig:slice-QES-A-09a}), and we 
also have that from the point $P_{4}$ it arises the volume region $V_{62}$, as it can be seeing in Fig.~\ref{fig:slice-QES-A-09b}.

\begin{figure}[h!]
	\centering\includegraphics[width=1\textwidth]{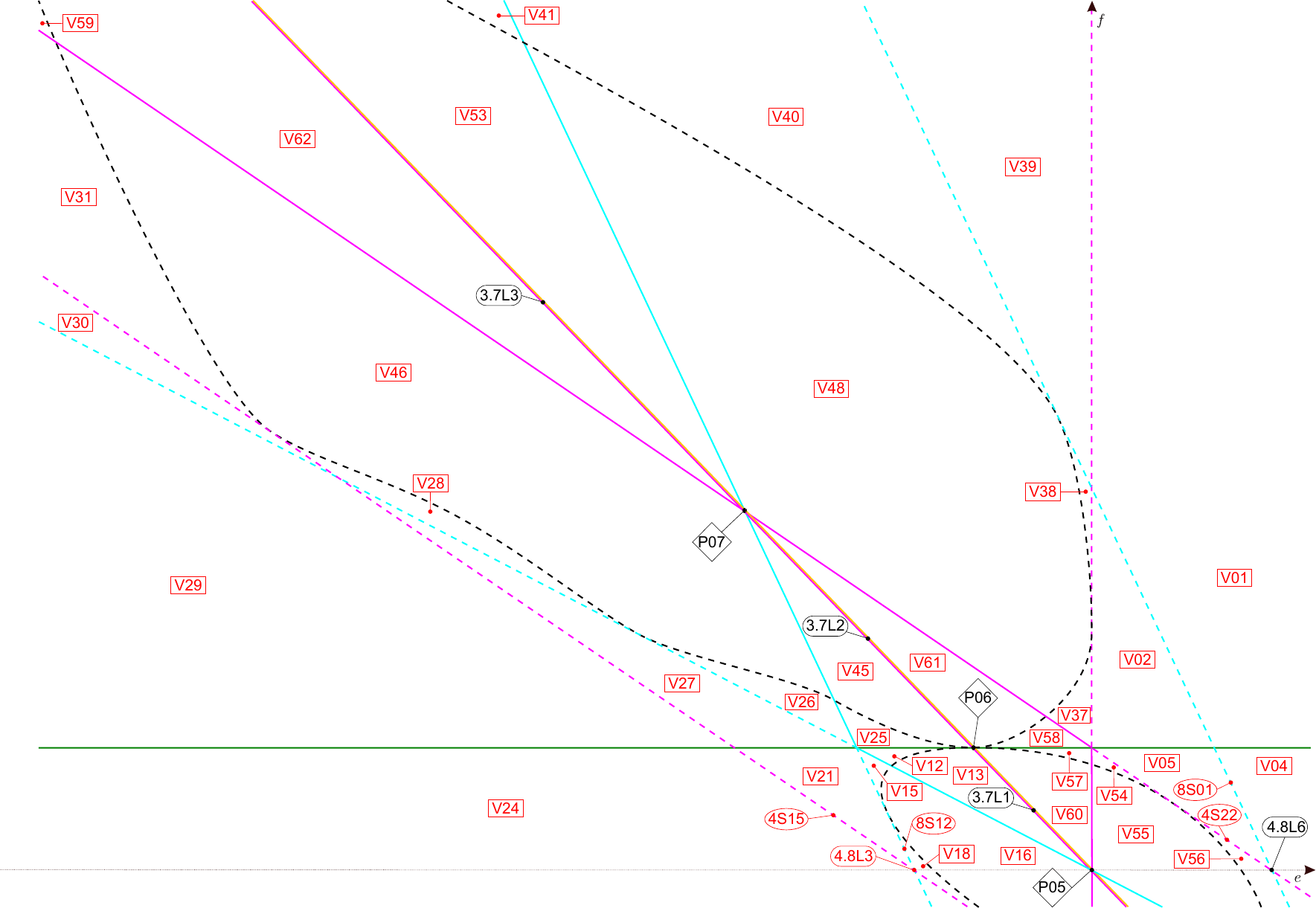} 
	\caption{\small \label{fig:slice-QES-A-10} Singular slice of the parameter space when $c=1$}
\end{figure}

Now we sum up the study of the singular slice $c=1$. At this slice there are several phenomena happening
simultaneously.
\begin{enumerate}
	\item Line $4.8L_{3}$ goes to $f=0$ and $V_{22}$ goes to $f<0$;
	\item the bifurcation curves $4S_{22}$ and $8S_{1}$ intercept themselves along $f=0$, more precisely,
	at $4.8L_{6}$;
	\item remember that, up to here we had, in each plane, the existence of three yellow straight lines and 
	one nonalgebraic curve.	However, at $c=1$ all of these bifurcation curves coalesce along the straight 
	line $f=-e$ (in fact, $({\cal S}_{3})\vert_{c=1}=-(e+f)^3$). And from this coalescence we have that:
	\begin{enumerate}
		\item the following 15 volume regions disappear along $f=-e$: $V_{8}$, $V_{9}$, $V_{10}$, 
		$V_{11}$, $V_{14}$, $V_{17}$, $V_{35}$, $V_{42}$, $V_{43}$, $V_{44}$, $V_{47}$, 
		$V_{49}$, $V_{50}$, $V_{51}$, $V_{52}$; and
		\item volume region $V_{34}$ goes to infinity.
		\item In addition, remember Remark~\ref{weak-singularities}, item 2, in which we verified that, for $c=1$
		and $f=-e$ the corresponding phase portrait possesses one center type singular point.
	\end{enumerate}
\end{enumerate}
In Fig.~\ref{fig:slice-QES-A-10} we present the entire singular slice $c=1$ properly labeled.

\begin{figure}[h!]
	\centering\includegraphics[width=1\textwidth]{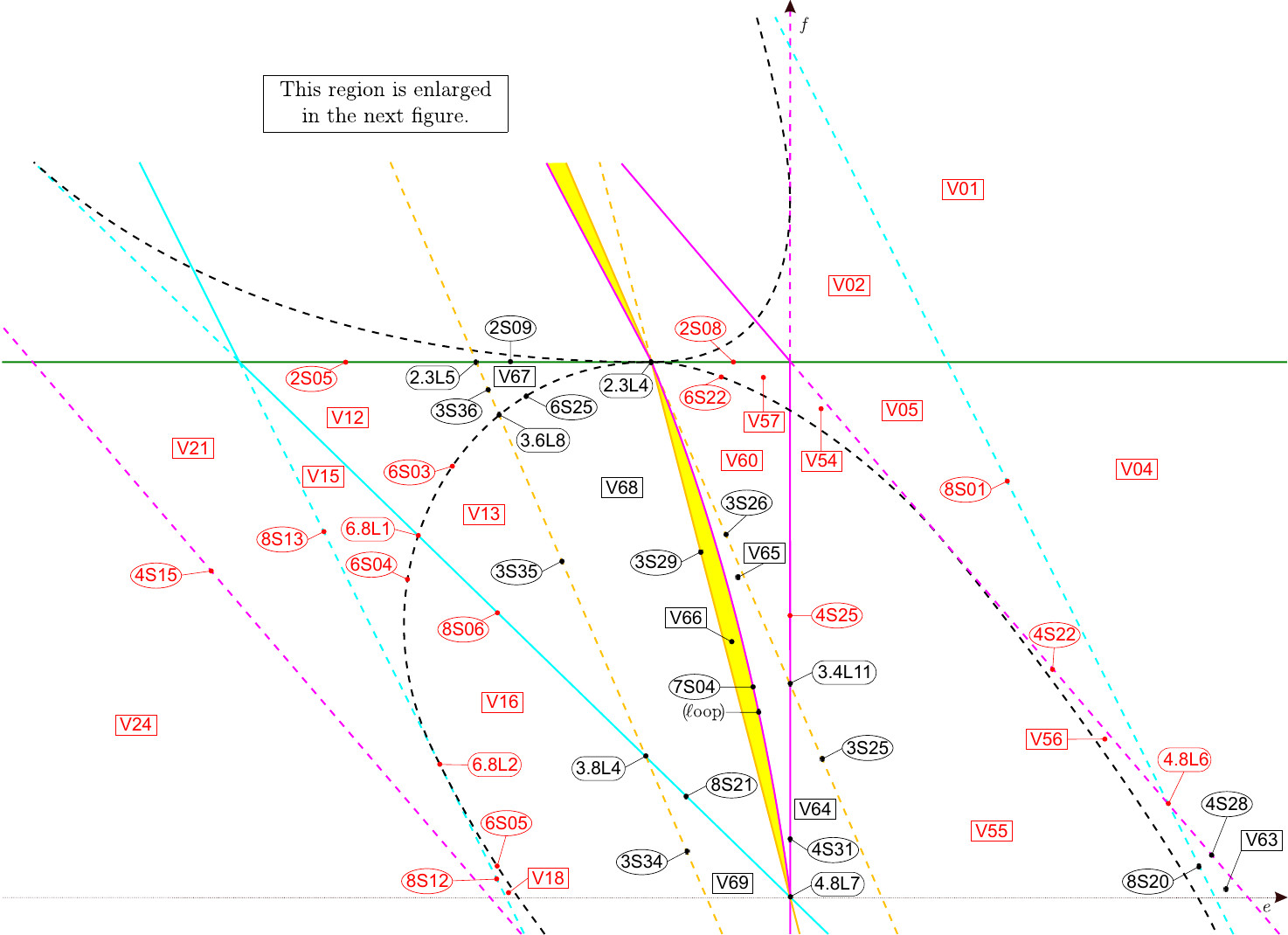} 
	\caption{\small \label{fig:slice-QES-A-11a} Piece of generic slice of the parameter space when $c=3/4$, 
		compare this region with Fig.~\ref{fig:slice-QES-A-10} and see also Fig.~\ref{fig:slice-QES-A-11b}}
\end{figure}

\begin{figure}[h!]
	\centering\includegraphics[width=1\textwidth]{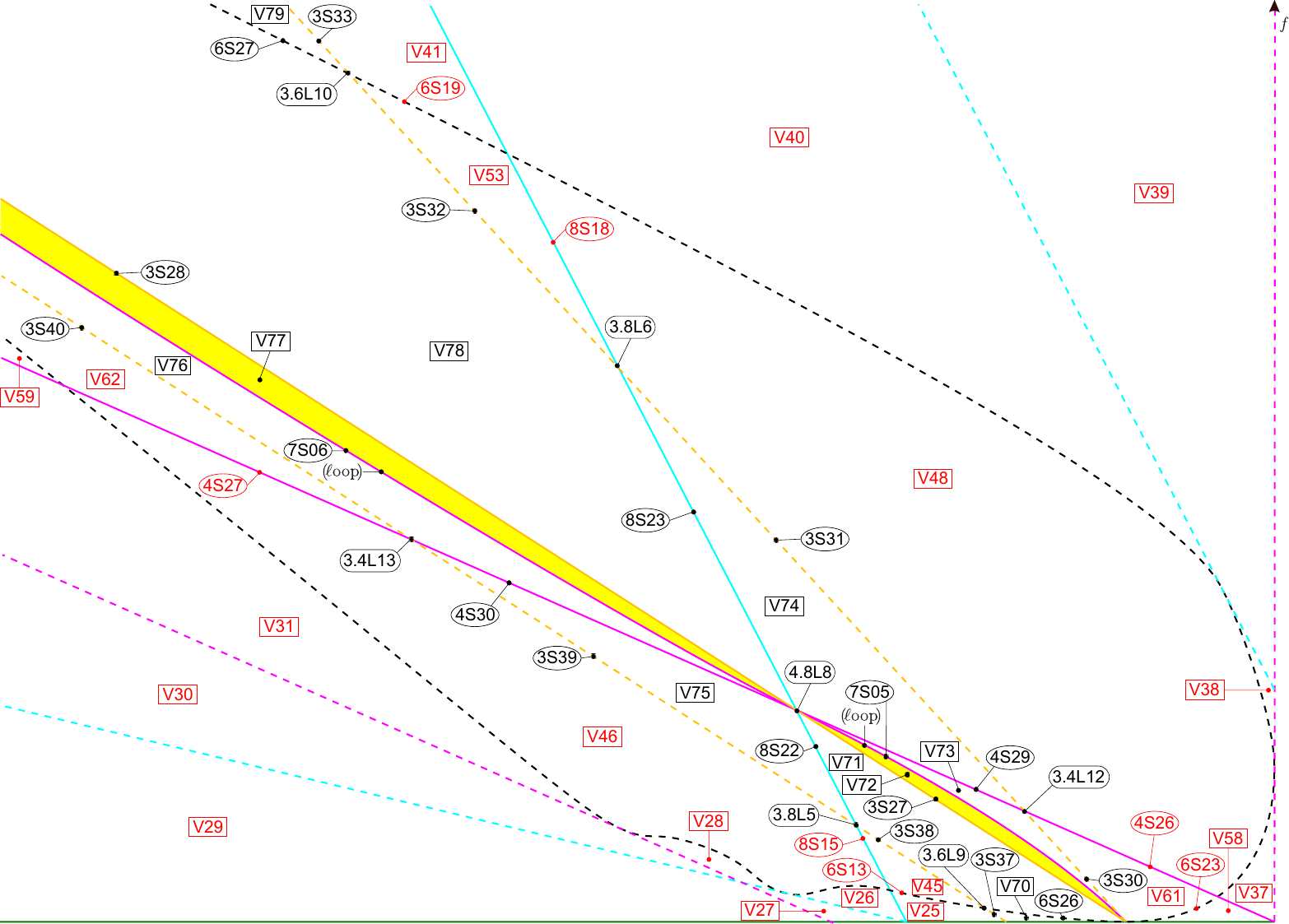} 
	\caption{\small \label{fig:slice-QES-A-11b} Another piece of generic slice of the parameter space when $c=3/4$, 
		compare this region with Fig.~\ref{fig:slice-QES-A-10}}
\end{figure}

Now, as it was expected, the generic slice $c=3/4$ brings several new information, as we describe
in the sequence.
\begin{enumerate}
	\item Line $4.8L_{3}$ goes to $f<0$;
	\item $4.8L_{6}$ goes to $f>0$ and it arises the volume region $V_{63}$;
	\item consider the bifurcation straight line $f=-e$ presented at slice $c=1$. For $c=3/4$ this 
	straight line splits itself into three yellow straight lines together with one nonalgebraic bifurcation curve.
	As a consequence, it arise the following 16 volume regions: $V_{64}$ up to $V_{79}$.
\end{enumerate}
We present this slice in Fig.~\ref{fig:slice-QES-A-11a} and \ref{fig:slice-QES-A-11b}.

Regarding the nonalgebraic curves $7S_{4}$ up to $7S_{6}$ that there appear in the mentioned
figures, we point out that their existence can be proved by using numerical tools and, by analogous
arguments as the ones we presented before, the following result can be easily proved. 

\begin{lemma}\label{lem:bif-7S4-till-7S6}
	In the generic slice $c=3/4$ there exist three pieces of nonalgebraic surfaces, denoted by
	$7S_{4}$, $7S_{5}$, and $7S_{6}$. These curves are displayed as in
	Fig.~\ref{fig:slice-QES-A-11a} and \ref{fig:slice-QES-A-11b}. Moreover, 
	$7S_{5}$ and $7S_{6}$ are continuation of $7S_{4}$.
\end{lemma}

\begin{figure}[h!]
	\centering\includegraphics[width=0.25\textwidth]{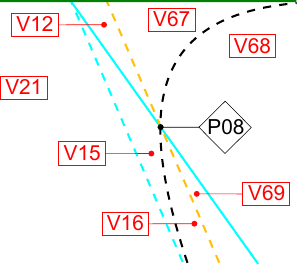} 
	\caption{\small \label{fig:slice-QES-A-12a} Piece of singular slice of the parameter space when $c=1/\sqrt{3}$, 
		compare this region with Fig.~\ref{fig:slice-QES-A-11a} and see also Fig.~\ref{fig:slice-QES-A-12b} and \ref{fig:slice-QES-A-12c}}
\end{figure}

\begin{figure}[h!]
	\centering\includegraphics[width=0.25\textwidth]{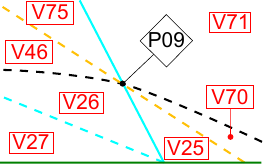} 
	\caption{\small \label{fig:slice-QES-A-12b} Piece of singular slice of the parameter space when $c=1/\sqrt{3}$, 
		compare this region with Fig.~\ref{fig:slice-QES-A-11b} and see also \ref{fig:slice-QES-A-12c}}
\end{figure}

\begin{figure}[h!]
	\centering\includegraphics[width=0.25\textwidth]{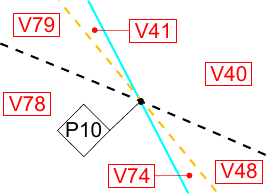} 
	\caption{\small \label{fig:slice-QES-A-12c} Piece of singular slice of the parameter space when $c=1/\sqrt{3}$, 
		compare this region with Fig.~\ref{fig:slice-QES-A-11b}}
\end{figure}

Now, for the singular slice $c=1/\sqrt{3}$ we observe that volume region $V_{13}$ coalesces at $P_{8}$
(see Fig.~\ref{fig:slice-QES-A-12a}), $V_{45}$ coalesces at $P_{9}$ (see Fig.~\ref{fig:slice-QES-A-12b}), 
and $V_{53}$ coalesces at $P_{10}$ (see Fig.~\ref{fig:slice-QES-A-12c}).

\begin{figure}[h!]
	\centering\includegraphics[width=0.30\textwidth]{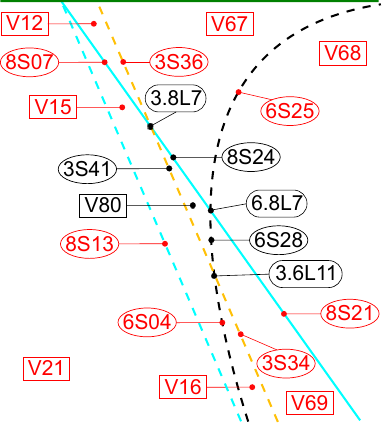} 
	\caption{\small \label{fig:slice-QES-A-13a} Piece of generic slice of the parameter space when $c=55/100$, 
		compare this region with Fig.~\ref{fig:slice-QES-A-12a} and see also Fig.~\ref{fig:slice-QES-A-13b} and \ref{fig:slice-QES-A-13c}}
\end{figure}

\begin{figure}[h!]
	\centering\includegraphics[width=0.35\textwidth]{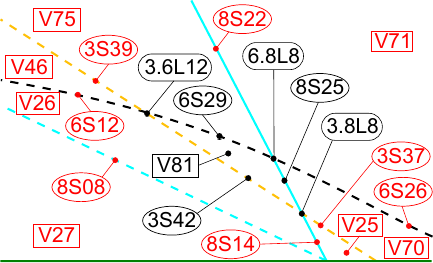} 
	\caption{\small \label{fig:slice-QES-A-13b} Piece of generic slice of the parameter space when $c=55/100$, 
		compare this region with Fig.~\ref{fig:slice-QES-A-12b} and see also Fig.~\ref{fig:slice-QES-A-13c}}
\end{figure}

\begin{figure}[h!]
	\centering\includegraphics[width=0.3\textwidth]{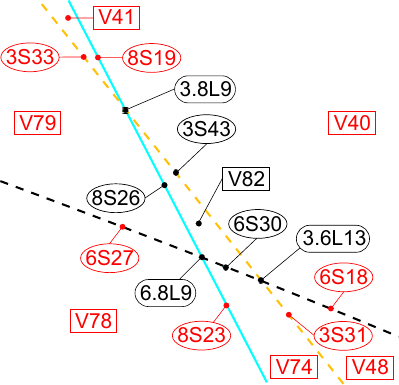} 
	\caption{\small \label{fig:slice-QES-A-13c} Piece of generic slice of the parameter space when $c=55/100$, 
		compare this region with Fig.~\ref{fig:slice-QES-A-12c}}
\end{figure}

In the generic slice $c=55/100$ we observe that from $P_{8}$ it arises the volume region $V_{80}$
(see Fig.~\ref{fig:slice-QES-A-13a}), from $P_{9}$ we get $V_{81}$ (see Fig.~\ref{fig:slice-QES-A-13b}), 
and  from $P_{10}$ we have $V_{82}$ (see Fig.~\ref{fig:slice-QES-A-13c}).

\begin{figure}[h!]
	\centering\includegraphics[width=0.55\textwidth]{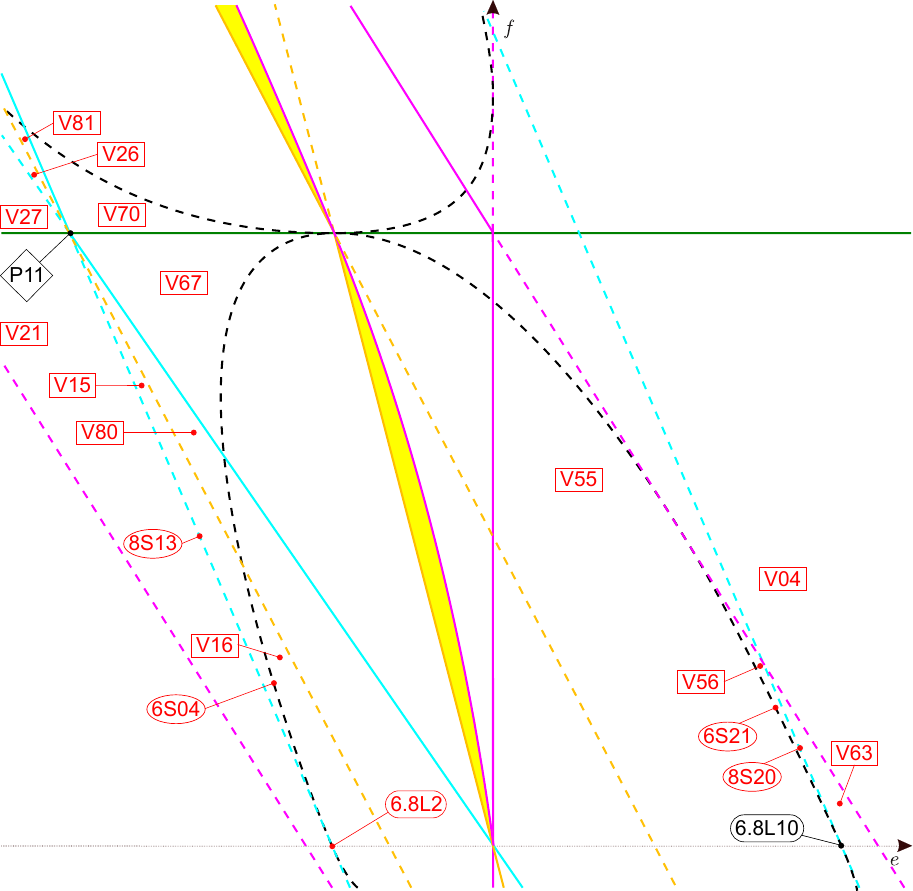} 
	\caption{\small \label{fig:slice-QES-A-14a} Piece of singular slice of the parameter space when $c=1/2$, 
		compare this region with Fig.~\ref{fig:slice-QES-A-13a} and Fig.~\ref{fig:slice-QES-A-13b} and see also 
		Fig.~\ref{fig:slice-QES-A-14b}}
\end{figure}

We now pass to describe the result of the study of the singular slice $c=1/2$.
\begin{itemize}
	\item Consider volume regions $V_{12}$ (Fig.~\ref{fig:slice-QES-A-13a}) and $V_{25}$ (Fig.~\ref{fig:slice-QES-A-13b}).
	By studying the singular slice $c=1/2$ we observe that these two volume regions coalesce at $P_{11}$.
	\item We also have that $6.8L_{2}$ goes to $f=0$; and
	\item $6S_{21}$ intercepts $8S_{20}$ on $f=0$, more precisely, at $6.8L_{10}$.
\end{itemize}
In Fig.~\ref{fig:slice-QES-A-14a} one can see these movements of the algebraic bifurcation surfaces.

\begin{figure}[h!]
	\centering\includegraphics[width=0.3\textwidth]{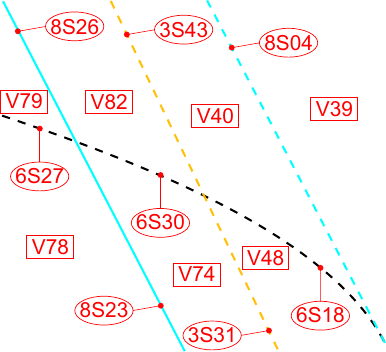} 
	\caption{\small \label{fig:slice-QES-A-14b} Another piece of singular slice of the parameter space when $c=1/2$, 
		compare this region with Fig.~\ref{fig:slice-QES-A-13c}}
\end{figure}

In addition to the previous description, when we have $c=1/2$, curve $3.8L_{9}$ together with $V_{41}$ 
(see Fig.~\ref{fig:slice-QES-A-13c}) go to infinity and the straight lines $3S_{43}$, $8S_{4}$, and $8S_{26}$
are now parallel (see Fig.~\ref{fig:slice-QES-A-14b}).

\begin{figure}[h!]
	\centering\includegraphics[width=0.3\textwidth]{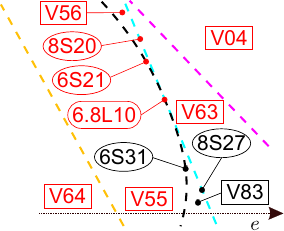} 
	\caption{\small \label{fig:slice-QES-A-15a} Piece of generic slice of the parameter space when $c=38/100$, 
		compare this region with Fig.~\ref{fig:slice-QES-A-14a} and see also Fig.~\ref{fig:slice-QES-A-15b} and \ref{fig:slice-QES-A-15c}}
\end{figure}

\begin{figure}[h!]
	\centering\includegraphics[width=0.25\textwidth]{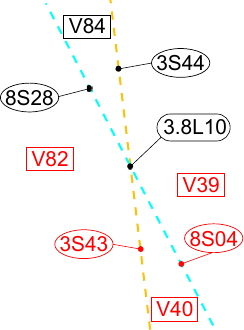} 
	\caption{\small \label{fig:slice-QES-A-15b} Piece of generic slice of the parameter space when $c=38/100$, 
		compare this region with Fig.~\ref{fig:slice-QES-A-14b} and see also Fig.~\ref{fig:slice-QES-A-15c}}
\end{figure}

\begin{figure}[h!]
	\centering\includegraphics[width=0.35\textwidth]{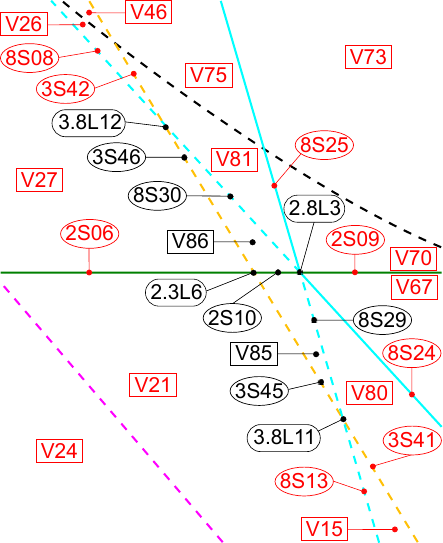} 
	\caption{\small \label{fig:slice-QES-A-15c} Piece of generic slice of the parameter space when $c=38/100$, 
		compare this region with Fig.~\ref{fig:slice-QES-A-14a}}
\end{figure}

After studying the singular slice $c=1/2$, if we consider $c=38/100$ as a generic value of the parameter $c$, 
we observe that:
\begin{itemize}
	\item $6.8L_{2}$ from Fig.~\ref{fig:slice-QES-A-14a} goes to $f<0$;
	\item $6.8L_{10}$ goes to $f>0$ and it arises volume region $V_{83}$ (see Fig.~\ref{fig:slice-QES-A-15a});
	\item $3S_{43}$ intercepts $8S_{4}$ at $3.8L_{10}$, generating volume region $V_{84}$ (see Fig.~\ref{fig:slice-QES-A-15b});
	\item from $P_{11}$ arise volume regions $V_{85}$ and $V_{86}$ (see Fig.~\ref{fig:slice-QES-A-15c}).
\end{itemize}

\begin{figure}[h!]
	\centering\includegraphics[width=0.4\textwidth]{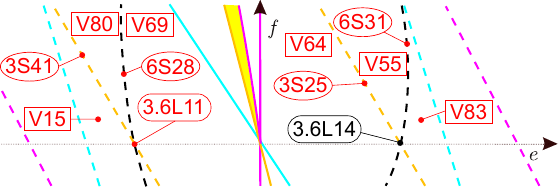} 
	\caption{\small \label{fig:slice-QES-A-16} Piece of singular slice of the parameter space when $c=2-\sqrt{3}$, 
		compare this region with Fig.~\ref{fig:slice-QES-A-14a} and \ref{fig:slice-QES-A-15a}}
\end{figure}

Now, when we consider the singular value $c=2-\sqrt{3}$ we observe that $3.6L_{11}$ goes to $f=0$
and $3S_{25}$ intercepts $6S_{31}$ at $3.6L_{14}$ (also on $f=0$), see Fig.~\ref{fig:slice-QES-A-16}.

\begin{figure}[h!]
	\centering\includegraphics[width=0.4\textwidth]{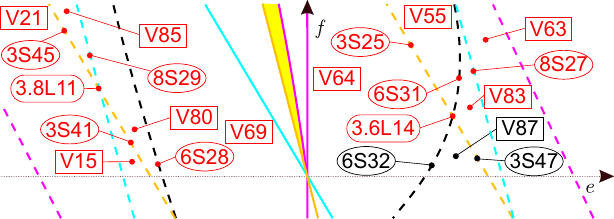} 
	\caption{\small \label{fig:slice-QES-A-17} Piece of generic slice of the parameter space when $c=26/100$, 
		compare this region with Fig.~\ref{fig:slice-QES-A-16}}
\end{figure}

In Fig.~\ref{fig:slice-QES-A-17} we present piece of generic slice $c=26/100$. For this value of the parameter
$c$ we observe that $3.6L_{11}$ goes to $f<0$ and $3.6L_{14}$ goes to $f>0$ and this provokes the 
appearance of volume region $V_{87}$.

\begin{figure}[h!]
	\centering\includegraphics[width=0.4\textwidth]{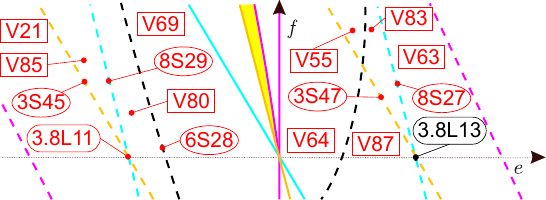} 
	\caption{\small \label{fig:slice-QES-A-18} Piece of singular slice of the parameter space when $c=1/4$, 
		compare this region with Fig.~\ref{fig:slice-QES-A-17}}
\end{figure}

Consider Fig.~\ref{fig:slice-QES-A-17}. During the study of the singular slice $c=1/4$ we notice that $3.8L_{11}$  
goes to $f=0$ and then $V_{15}$ goes to $f<0$. Moreover, $3S_{47}$ intercepts $8S_{27}$ on $f=0$, more
precisely, at $3.8L_{13}$. In Fig.~\ref{fig:slice-QES-A-18} one can see a piece of the parameter space corresponding
to this singular slice.

\begin{figure}[h!]
	\centering\includegraphics[width=0.4\textwidth]{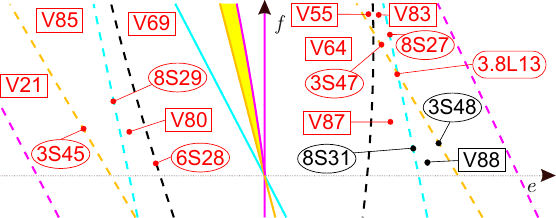} 
	\caption{\small \label{fig:slice-QES-A-19} Piece of generic slice of the parameter space when $c=1/10$, 
		compare this region with Fig.~\ref{fig:slice-QES-A-18}}
\end{figure}

Now we consider the last generic slice corresponding to $c>0$. In fact, for $c=1/10$ we see 
that $3.8L_{11}$ goes to $f<0$ and $3.8L_{13}$ goes to $f>0$ which allows the appearance 
of volume region $V_{88}$. Moreover, we point out that numerical verification shows that the
nonalgebraic curves maintain their position (with respect to the algebraic curves) as it was verified
in slice $c=3/4$. In Fig.~\ref{fig:slice-QES-A-19} we present the corresponding piece of the 
generic slice under consideration.

According to \eqref{eq:values-of-c-QES-A} now we start the study of the regions of the bifurcation diagram 
corresponding to negative values of the parameter $c$. The first generic slice to be considered is given by 
$c=-1/10$.

\begin{figure}[h!]
	\centering\includegraphics[width=1\textwidth]{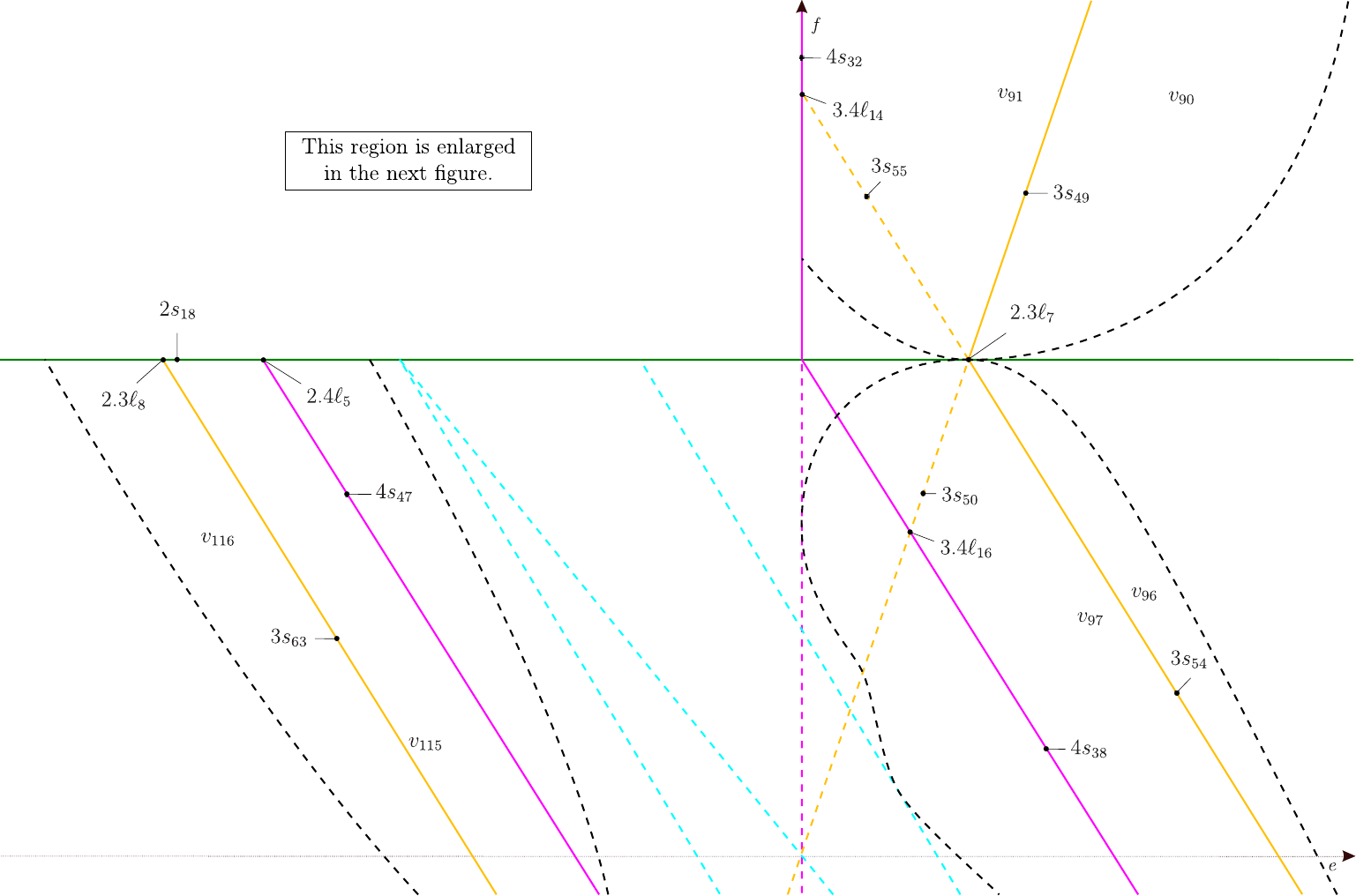} 
	\caption{\small \label{fig:slice-QES-A-20-alg1} Piece of generic slice of the parameter space when $c=5$ (only algebraic 
		surfaces), see also Fig.~\ref{fig:slice-QES-A-20-alg2}}
\end{figure}

\begin{figure}[h!]
	\centering\includegraphics[width=1\textwidth]{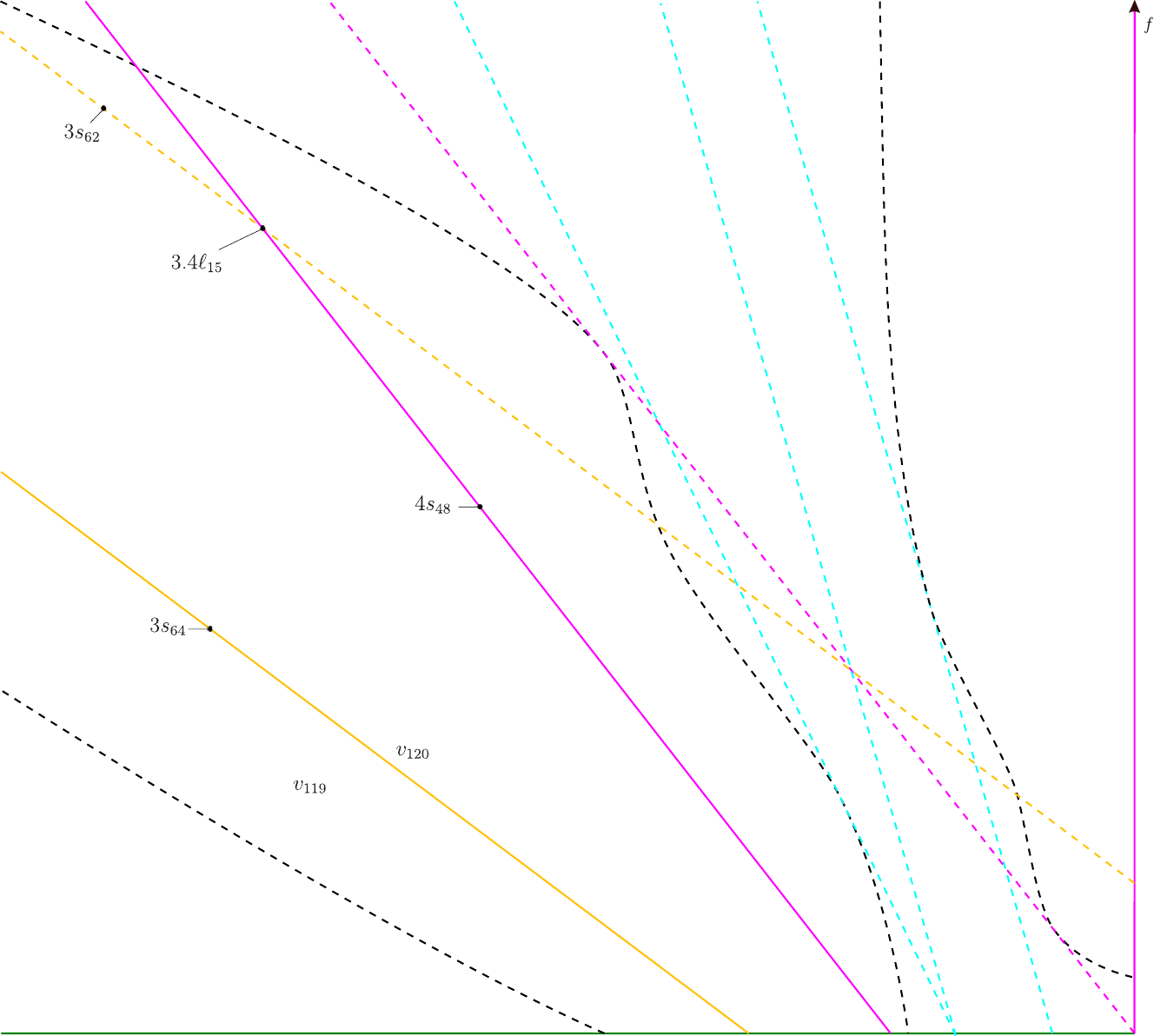} 
	\caption{\small \label{fig:slice-QES-A-20-alg2} Continuation of Fig.~\ref{fig:slice-QES-A-20-alg1}}
\end{figure}

As in the case of slice $c=5$, here we have a partition in two--dimensional parts bordered by curved polygons, 
some of them bounded and others bordered by infinity. And we use lower--case letters provisionally to describe 
the sets found algebraically in order to do not interfere with the final partition described with capital letters, see
the algebraic slice in Fig.~\ref{fig:slice-QES-A-20-alg1} and \ref{fig:slice-QES-A-20-alg2}.

For each two--dimensional part we obtain a phase portrait which is coherent with those of all their borders. 
Except for four parts, which are shown in Fig.~\ref{fig:slice-QES-A-20-alg1} and \ref{fig:slice-QES-A-20-alg2}
and named as follows:
\begin{itemize}\label{page:regions-c-minus1over10}
	\item $v_{91}$: the quadrilateral bordered by yellow and purple curves and also by the line at infinity (in Fig.~\ref{fig:slice-QES-A-20-alg1});
	\item $v_{97}$: the quadrilateral bordered by yellow, purple, and (due to the symmetry) green curves (in Fig.~\ref{fig:slice-QES-A-20-alg1});
	\item $v_{115}$: the quadrilateral bordered by green, purple, and (due to symmetry) yellow curves (in Fig.~\ref{fig:slice-QES-A-20-alg1});
	\item $v_{120}$: the quadrilateral bordered by yellow and purple curves and infinity (in Fig.~\ref{fig:slice-QES-A-20-alg2}).
\end{itemize}
The study of these parts is important for the coherence of the bifurcation diagram. That is why we have 
decided to present only these parts in the mentioned figures (in Fig.~\ref{fig:slice-QES-A-20a} and 
Fig.~\ref{fig:slice-QES-A-20b} one can see the complete bifurcation diagram for this slice). 

We start the study of part $v_{91}$. Segment $3s_{49}$ in Fig.~\ref{fig:slice-QES-A-20-alg1} is one of 
the borders of this part and, the phase portrait corresponding to this segment possesses a weak focus (of order one), 
so this branch of surface (${\cal{S}}_{3}$) corresponds to a Hopf bifurcation. This means that the phase portrait 
corresponding to one of the sides of this segment must have a limit cycle; in fact it is in region $v_{91}$.

\begin{figure}[h!]
	\centering\includegraphics[width=0.7\textwidth]{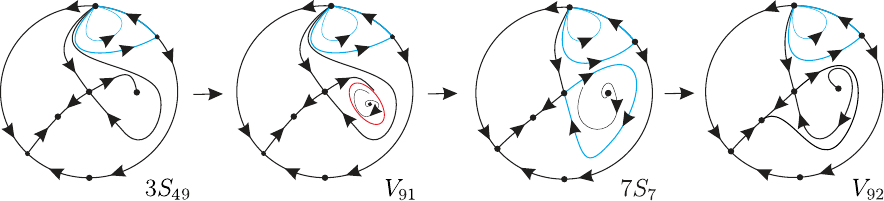} 
	\caption{\small \label{fig:bif-7S7} Sequence of phase portraits in parts $V_{91}$ and $V_{92}$ of slice 
		$c=-1/10$ (the labels are according to Fig.~\ref{fig:slice-QES-A-20a})}
\end{figure}

\begin{figure}[h!]
	\centering\includegraphics[width=1\textwidth]{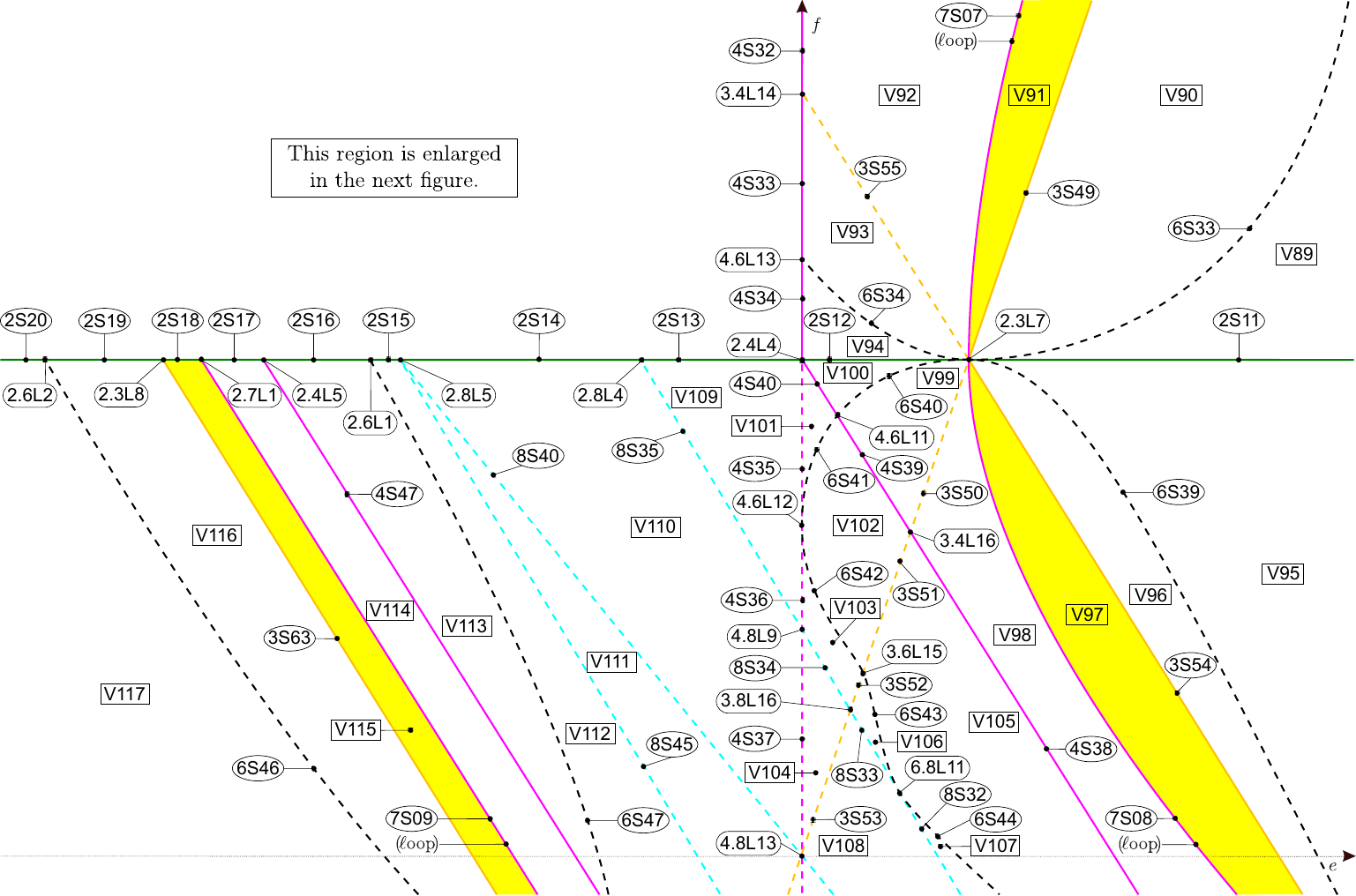} 
	\caption{\small \label{fig:slice-QES-A-20a} Piece of generic slice of the parameter space when 
		$c=-1/10$, see also Fig.~\ref{fig:slice-QES-A-20b}}
\end{figure}

\begin{figure}[h!]
	\centering\includegraphics[width=0.8\textwidth]{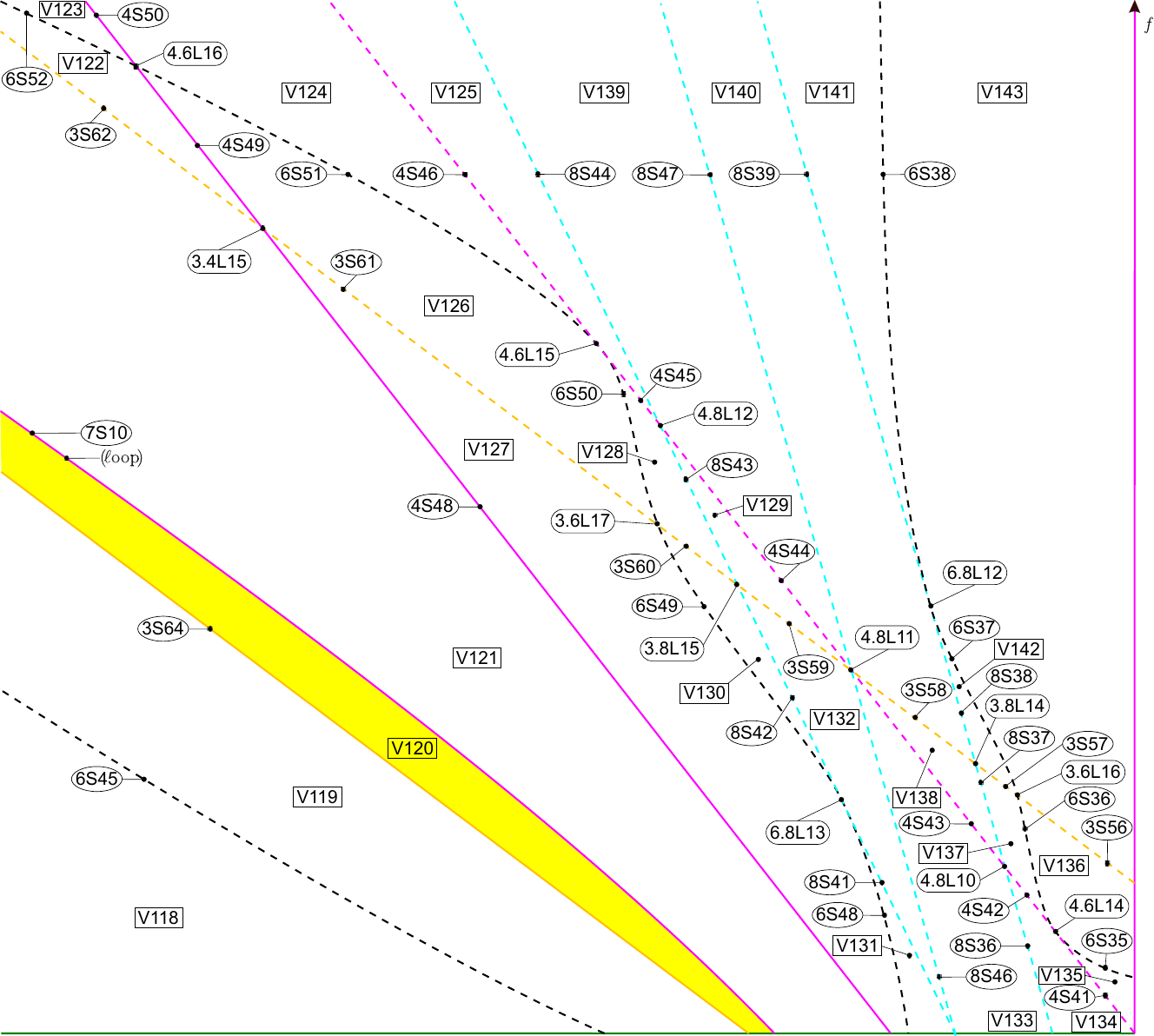} 
	\caption{\small \label{fig:slice-QES-A-20b} Continuation of Fig.~\ref{fig:slice-QES-A-20a}}
\end{figure}

However, when we approach $4s_{32}$ and $3s_{55}$, the limit cycle has been lost, which implies the 
existence of at least one element of surface (${\cal{S}}_{7}$) (see $7S_{7}$ in Fig.~\ref{fig:slice-QES-A-20a}), 
in a neighborhood of $3s_{49}$, due to a connection of separatrices from a saddle to itself (i.e. a loop--type connection). 
In Lemma \ref{lem:endpoints-7S7} we show that $7S_{7}$ is unbounded and it has one of its endpoints 
at the curve $2.3\ell_{7}$. We draw the sequence of phase portraits along these subsets (using the notation
from Fig.~\ref{fig:slice-QES-A-20a}) in Fig.~\ref{fig:bif-7S7} and we plot the complete bifurcation diagram 
for this part in Fig.~\ref{fig:slice-QES-A-20a}.

\begin{lemma}\label{lem:endpoints-7S7}
	The nonalgebraic curve $7S_{7}$  is unbounded and it has one of its endpoints at the curve 
	$2.3\ell_{7}$.
\end{lemma}

\begin{proof} Numerical analysis suggest that this result is true. In fact, note that if one of the 
	endpoints of this surface is any point of $3s_{49}$, then a portion of this subset must not refer to a 
	Hopf bifurcation, which contradicts the fact that on $3s_{49}$ we have a weak focus of order one.
	Also, observe that it is not possible that the starting point of this surface is on $3s_{55}$, since 
	on this portion of the yellow surface we have only a $C^\infty$ bifurcation (weak saddle).
	Finally, the endpoints cannot be on $4s_{32}$ because, in order to have this, first we need to break 
	the invariant straight line. Then, the only possible endpoint of surface $7S_{7}$ is $2.3\ell_{7}$.
\end{proof}

Consider Fig.~\ref{fig:slice-QES-A-20a}. Note that here we have an interesting situation. 
On one hand, $2.3L_{7}$ is a transition between $2S_{11}$ and $2S_{12}$, i.e. one can 
see a cusp point being a transition of different types of saddle--nodes. On the other hand, being 
$2.3L_{7}$ an endpoint of $7S_{7}$ we observe a cusp point formed by the coalescence 
of a focus with a saddle.

Now we consider parts $v_{97}$, $v_{115}$, and $v_{120}$ in Fig.~\ref{fig:slice-QES-A-20-alg1}
and \ref{fig:slice-QES-A-20-alg2}. As we have that:
\begin{itemize}
	\item $3s_{49}$ produces a phase portrait that is topologically equivalent to the ones in $3s_{54}$, $3s_{63}$,
	and $3s_{64}$;
	\item $3s_{55}$ produces a phase portrait that is topologically equivalent to the ones in $3s_{50}$ and $3s_{62}$;
	\item $4s_{32}$ produces a phase portrait that is topologically equivalent to the ones in $4s_{38}$, $4s_{47}$,
	and $4s_{48}$;
\end{itemize}
by the same arguments used in the study of part $v_{91}$ we conclude the existence of nonalgebraic surfaces 
$7S_{8}$, $7S_{9}$, and $7S_{10}$ in Fig.~\ref{fig:slice-QES-A-20a} and \ref{fig:slice-QES-A-20b}.
Moreover we also have that $7S_{10}$ is not bounded, $7S_{8}$ and $7S_{9}$ are bounded (due to the 
symmetry on the bifurcation diagram), $7S_{8}$ is a continuation of $7S_{7}$, and $7S_{10}$ is a 
continuation of $7S_{9}$.

\begin{figure}[h!]
	\centering\includegraphics[width=0.4\textwidth]{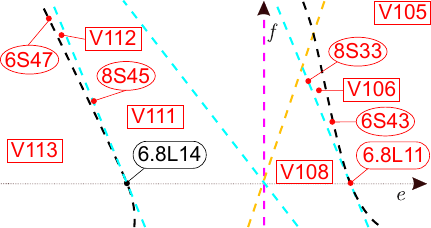} 
	\caption{\small \label{fig:slice-QES-A-21} Piece of singular slice of the parameter space when $c=-1/4$, 
		compare this region with Fig.~\ref{fig:slice-QES-A-20a}}
\end{figure}

Now we take the singular value $c=-1/4$. For this value of the parameter $c$ we notice that
$6.8L_{11}$ goes to $f=0$ and $6S_{47}$ intercepts $8S_{45}$ at $6.8L_{14}$, see 
these phenomena along $f=0$ in Fig.~\ref{fig:slice-QES-A-21}.

\begin{figure}[h!]
	\centering\includegraphics[width=0.4\textwidth]{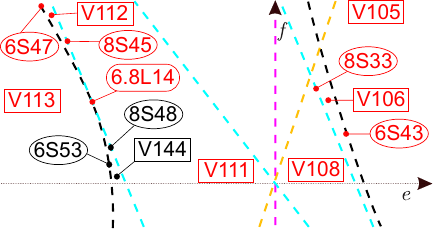} 
	\caption{\small \label{fig:slice-QES-A-22} Piece of generic slice of the parameter space when $c=-35/100$, 
		compare this region with Fig.~\ref{fig:slice-QES-A-21}}
\end{figure}

By considering $c=-35/100$ as a generic slice, two expected situations are detected, namely, $6.8L_{11}$
goes to $f<0$ and $6.8L_{14}$ goes to $f>0$ giving place to the appearance of volume region $V_{144}$,
see Fig.~\ref{fig:slice-QES-A-22}.

\begin{figure}[h!]
	\centering\includegraphics[width=1\textwidth]{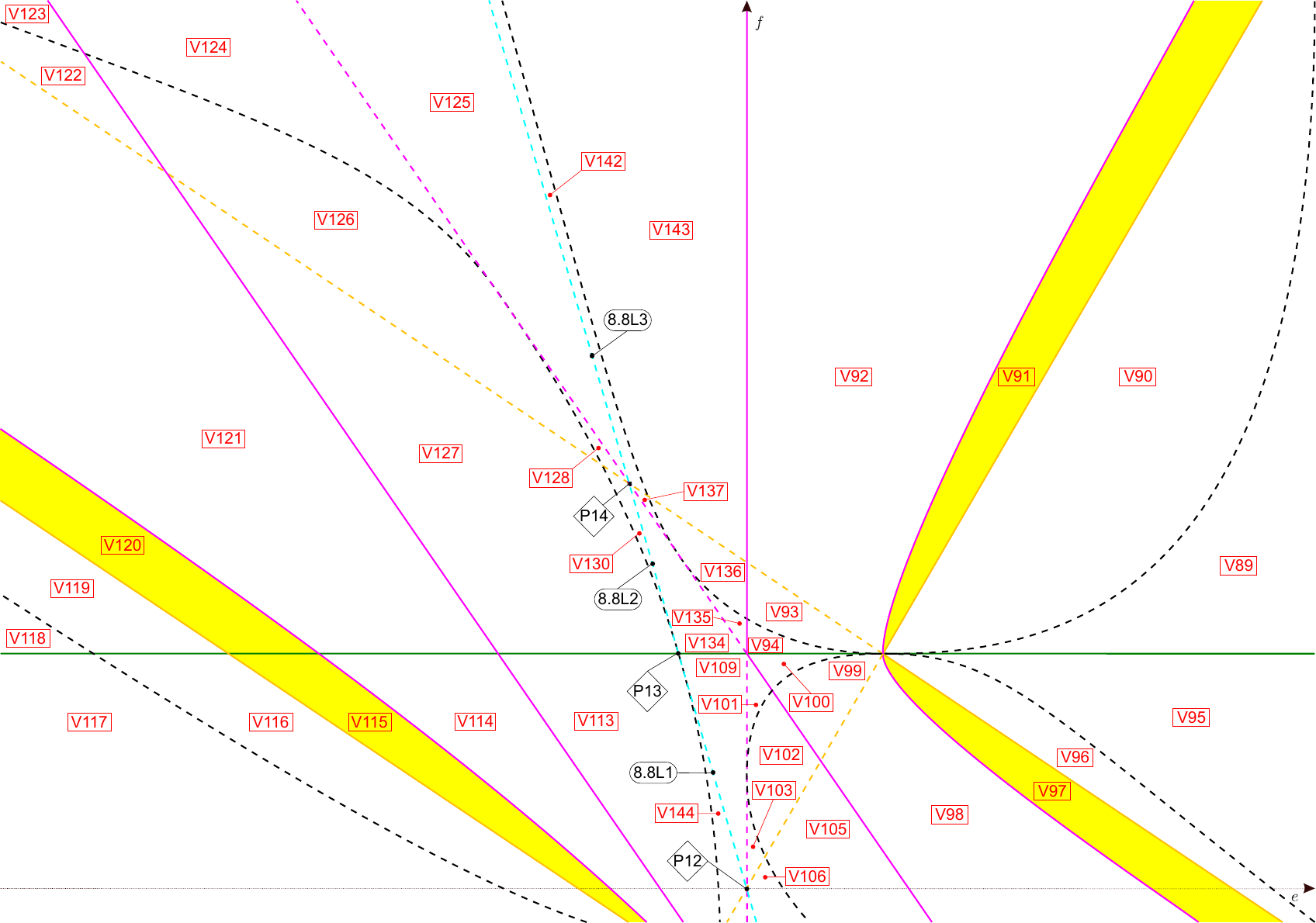} 
	\caption{\small \label{fig:slice-QES-A-23} Piece of singular slice of the parameter space when 
		$c=-1/2$}
\end{figure}

Now we consider the singular slice $c=-1/2$. Up to here we had, in each plane, the existence of three 
cyan straight lines. However, at $c=-1/2$ these bifurcation curves coalesce along the straight line 
$f=-2e$ (indeed, $({\cal S}_{8})\vert_{c=-1/2}=(2e+f)^3/2$). And from this coalescence we have that:
\begin{enumerate}
	\item Volume regions $V_{112}$ (Fig.~\ref{fig:slice-QES-A-22}) and $V_{131}$ 
	(Fig.~\ref{fig:slice-QES-A-20b}) coalesce at $P_{13}$;
	\item $6.8L_{12}$ together with $V_{141}$ (Fig.~\ref{fig:slice-QES-A-20b}) go to infinity; and
	\item the following ten volume regions disappear along $f=-2e$: $V_{104}$, $V_{108}$, $V_{110}$, 
	$V_{111}$, $V_{129}$, $V_{132}$, $V_{133}$, $V_{138}$, $V_{139}$, and $V_{140}$.
	We advise the reader to remember their location in Fig.~\ref{fig:slice-QES-A-20a} and \ref{fig:slice-QES-A-20b}.
\end{enumerate}
In Fig.~\ref{fig:slice-QES-A-23} we present the entire singular slice $c=-1/2$ completely labeled.

\begin{figure}[h!]
	\centering\includegraphics[width=1\textwidth]{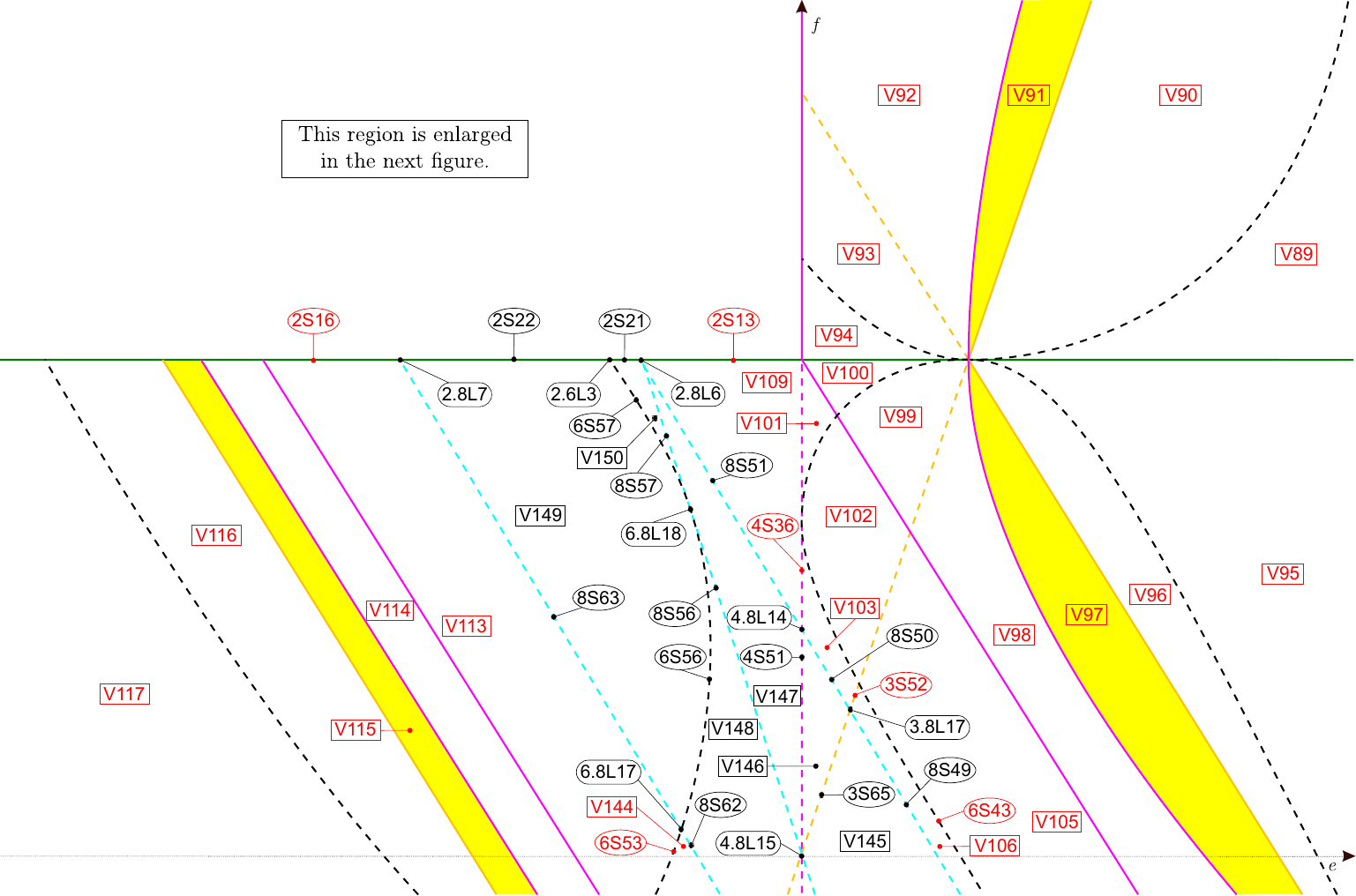} 
	\caption{\small \label{fig:slice-QES-A-24a} Piece of generic slice of the parameter space when 
		$c=-53/100$, see also Fig.~\ref{fig:slice-QES-A-24b}}
\end{figure}

\begin{figure}[h!]
	\centering\includegraphics[width=0.52\textwidth]{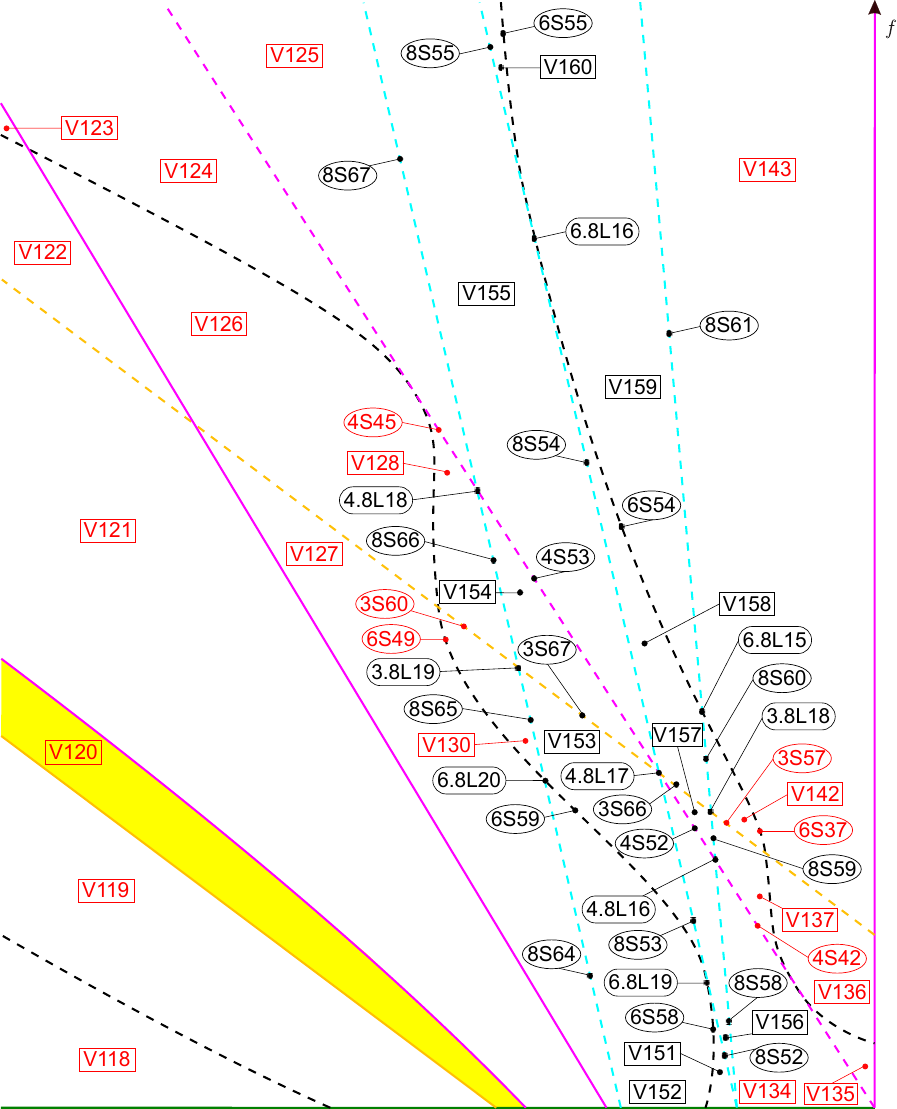} 
	\caption{\small \label{fig:slice-QES-A-24b} Continuation of Fig.~\ref{fig:slice-QES-A-24a}}
\end{figure}

When we consider the generic slice $c=-53/100$ we observe that the triple cyan bifurcation straight line 
(obtained in the previous slice) splits itself into 16 new volume regions, namely, $V_{145}$ up to $V_{160}$. 
These volume regions are displayed as in Fig.~\ref{fig:slice-QES-A-24a} and \ref{fig:slice-QES-A-24b}.

\begin{figure}[h!]
	\centering\includegraphics[width=0.4\textwidth]{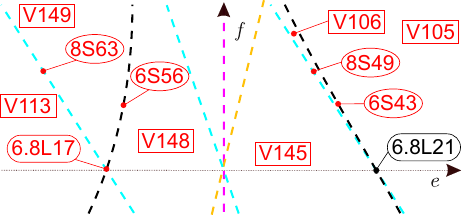} 
	\caption{\small \label{fig:slice-QES-A-25} Piece of singular slice of the parameter space when 
		$c=-9/16$}
\end{figure}

After the analysis of the generic slice $c=-53/100$ we study the singular slice $c=-9/16=-0.5625$. 
Consider Fig.~\ref{fig:slice-QES-A-24a}. For this singular value of the parameter $c$ we observe 
that $6.8L_{17}$ goes to $f=0$ and then $V_{144}$ goes to $f<0$. Also, we have that $6S_{43}$ 
intercepts $8S_{49}$ at $6.8L_{21}$. In Fig.~\ref{fig:slice-QES-A-25} we present the piece of
slice of the parameter space corresponding to these regions.

\begin{figure}[h!]
	\centering\includegraphics[width=0.4\textwidth]{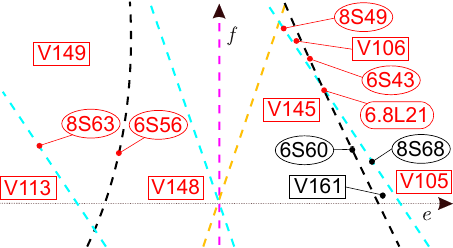} 
	\caption{\small \label{fig:slice-QES-A-26} Piece of generic slice of the parameter space when 
		$c=-57/100$}
\end{figure}

Now if we consider $c=-57/100$ as a generic slice one can detect the expected phenomena: $6.8L_{17}$
goes to $f<0$ and $6.8L_{21}$ goes to $f>0$ (from which it arises $V_{161}$), see Fig.~\ref{fig:slice-QES-A-26}. 

\begin{figure}[h!]
	\centering\includegraphics[width=0.2\textwidth]{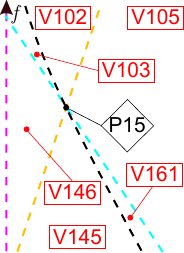} 
	\caption{\small \label{fig:slice-QES-A-27a} Piece of singular slice of the parameter space when 
		$c=-1/\sqrt{3}$, see also Fig.~\ref{fig:slice-QES-A-27b}}
\end{figure}

\begin{figure}[h!]
	\centering\includegraphics[width=0.35\textwidth]{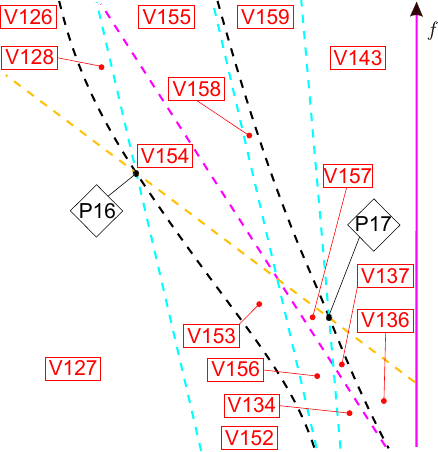} 
	\caption{\small \label{fig:slice-QES-A-27b} Another piece of singular slice of the parameter space when 
		$c=-1/\sqrt{3}$, see also Fig.~\ref{fig:slice-QES-A-27a}}
\end{figure}

Moving on with the study of the list of slices presented in \eqref{eq:values-of-c-QES-A}, now we consider
the singular slice $c=-1/\sqrt{3}$. During the study of this slice we observe that volume regions 
$V_{106}$ (see Fig.~\ref{fig:slice-QES-A-26}), $V_{130}$, and $V_{142}$ 
(see Fig.~\ref{fig:slice-QES-A-24b}) are reduced to the points $P_{15}$, $P_{16}$, and
$P_{17}$, respectively, as we illustrate in Fig.~\ref{fig:slice-QES-A-27a} and \ref{fig:slice-QES-A-27b}.

\begin{figure}[h!]
	\centering\includegraphics[width=0.25\textwidth]{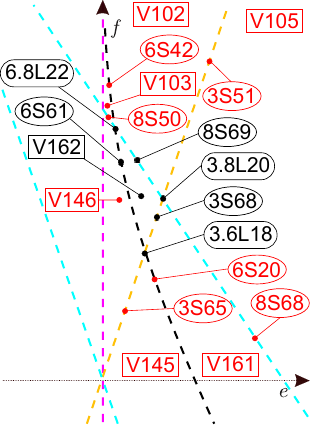} 
	\caption{\small \label{fig:slice-QES-A-28a} Piece of generic slice of the parameter space when 
		$c=-62/100$, see also Fig.~\ref{fig:slice-QES-A-28b}}
\end{figure}

\begin{figure}[h!]
	\centering\includegraphics[width=0.4\textwidth]{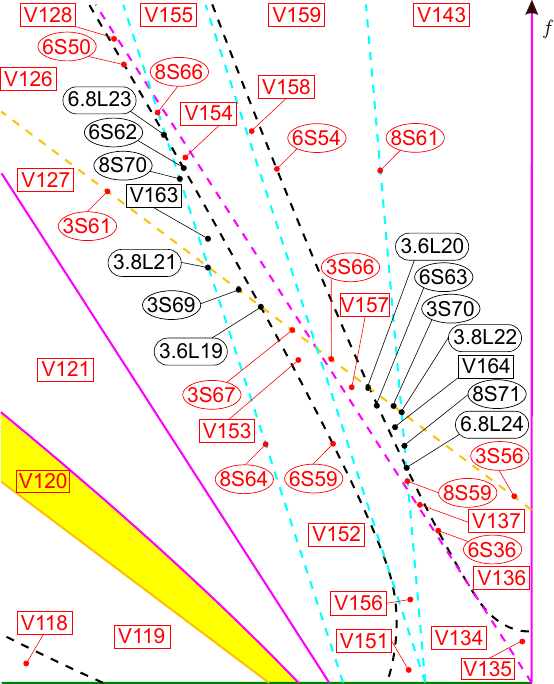} 
	\caption{\small \label{fig:slice-QES-A-28b} Another piece of generic slice of the parameter space when 
		$c=-62/100$, see also Fig.~\ref{fig:slice-QES-A-28a}}
\end{figure}

Taking $c=-62/100$ as a generic slice, we observe that from the points $P_{15}$, $P_{16}$, and
$P_{17}$ arise the volume regions $V_{162}$, $V_{163}$, and $P_{164}$, respectively. A 
draw of these regions can be seeing in Fig.~\ref{fig:slice-QES-A-28a} and \ref{fig:slice-QES-A-28b}.

\begin{figure}[h!]
	\centering\includegraphics[width=0.18\textwidth]{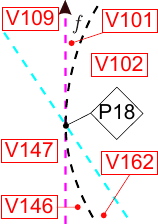} 
	\caption{\small \label{fig:slice-QES-A-29a} Piece of singular slice of the parameter space when 
		$c=-2/3$, see also Fig.~\ref{fig:slice-QES-A-29b}}
\end{figure}

\begin{figure}[h!]
	\centering\includegraphics[width=0.35\textwidth]{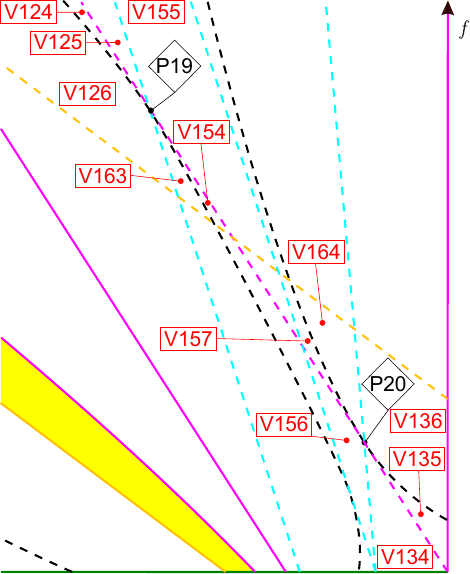} 
	\caption{\small \label{fig:slice-QES-A-29b} Another piece of singular slice of the parameter space when 
		$c=-2/3$, see also Fig.~\ref{fig:slice-QES-A-29a}}
\end{figure}

Now when we perform the study of singular slice $c=-2/3$ we observe that volume regions $V_{103}$ 
(Fig.~\ref{fig:slice-QES-A-28a}), $V_{128}$, and $V_{137}$ (Fig.~\ref{fig:slice-QES-A-28b}) are
reduced to the points $P_{18}$, $P_{19}$, and $P_{20}$, respectively. These points are drawn 
in Fig.~\ref{fig:slice-QES-A-29a} and \ref{fig:slice-QES-A-29b}.

\begin{figure}[h!]
	\centering\includegraphics[width=0.9\textwidth]{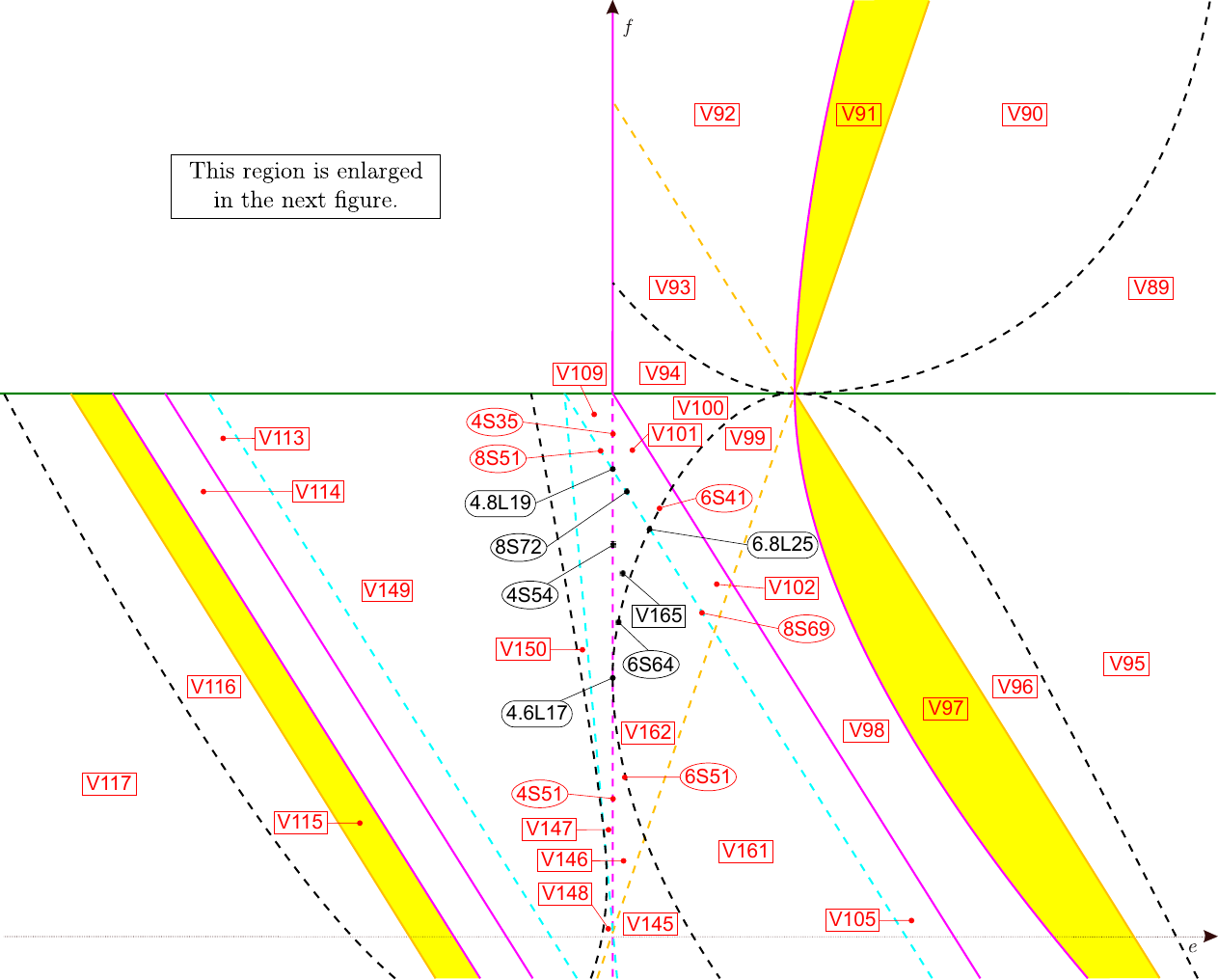} 
	\caption{\small \label{fig:slice-QES-A-30a} Piece of generic slice of the parameter space when 
		$c=-85/100$, see also Fig.~\ref{fig:slice-QES-A-30b}}
\end{figure}

\begin{figure}[h!]
	\centering\includegraphics[width=0.5\textwidth]{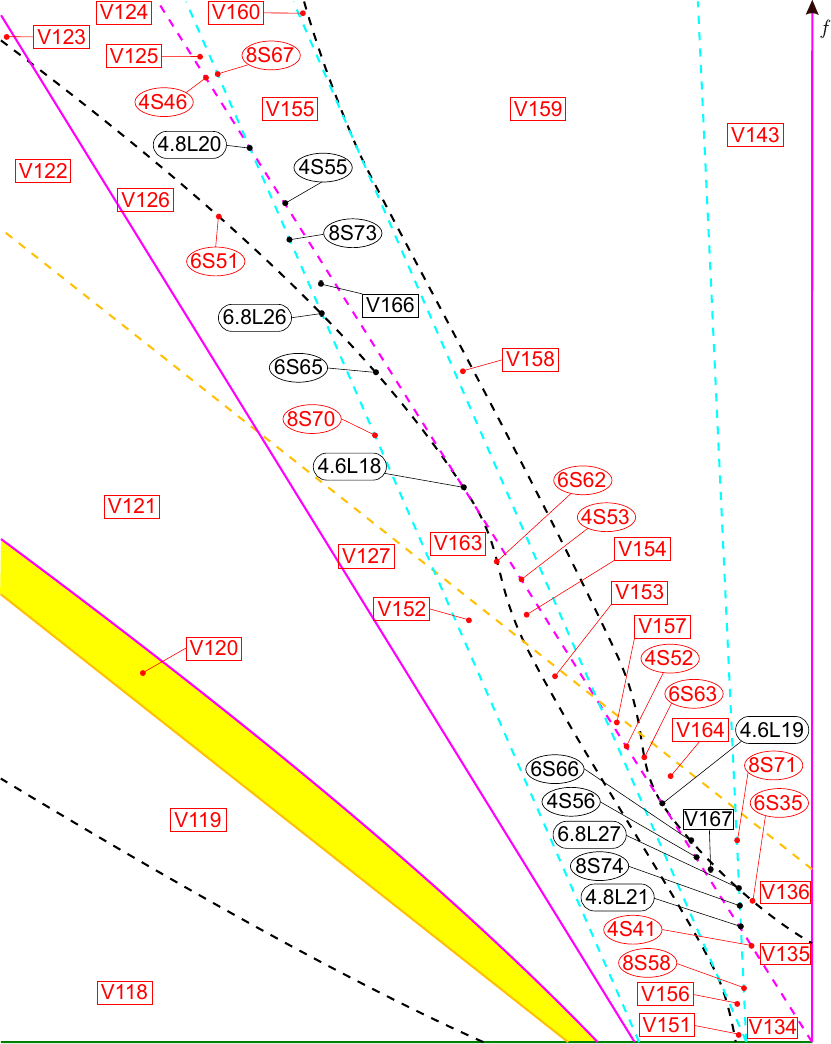} 
	\caption{\small \label{fig:slice-QES-A-30b} Continuation of Fig.~\ref{fig:slice-QES-A-30a}}
\end{figure}

Now we consider the generic slice $c=-85/100$. From the points $P_{15}$, $P_{16}$, and
$P_{17}$ arise the volume regions $V_{165}$, $V_{166}$, and $V_{167}$, respectively, which 
can be seeing in Fig.~\ref{fig:slice-QES-A-30a} and \ref{fig:slice-QES-A-30b}.

\begin{figure}[h!]
	\centering\includegraphics[width=1\textwidth]{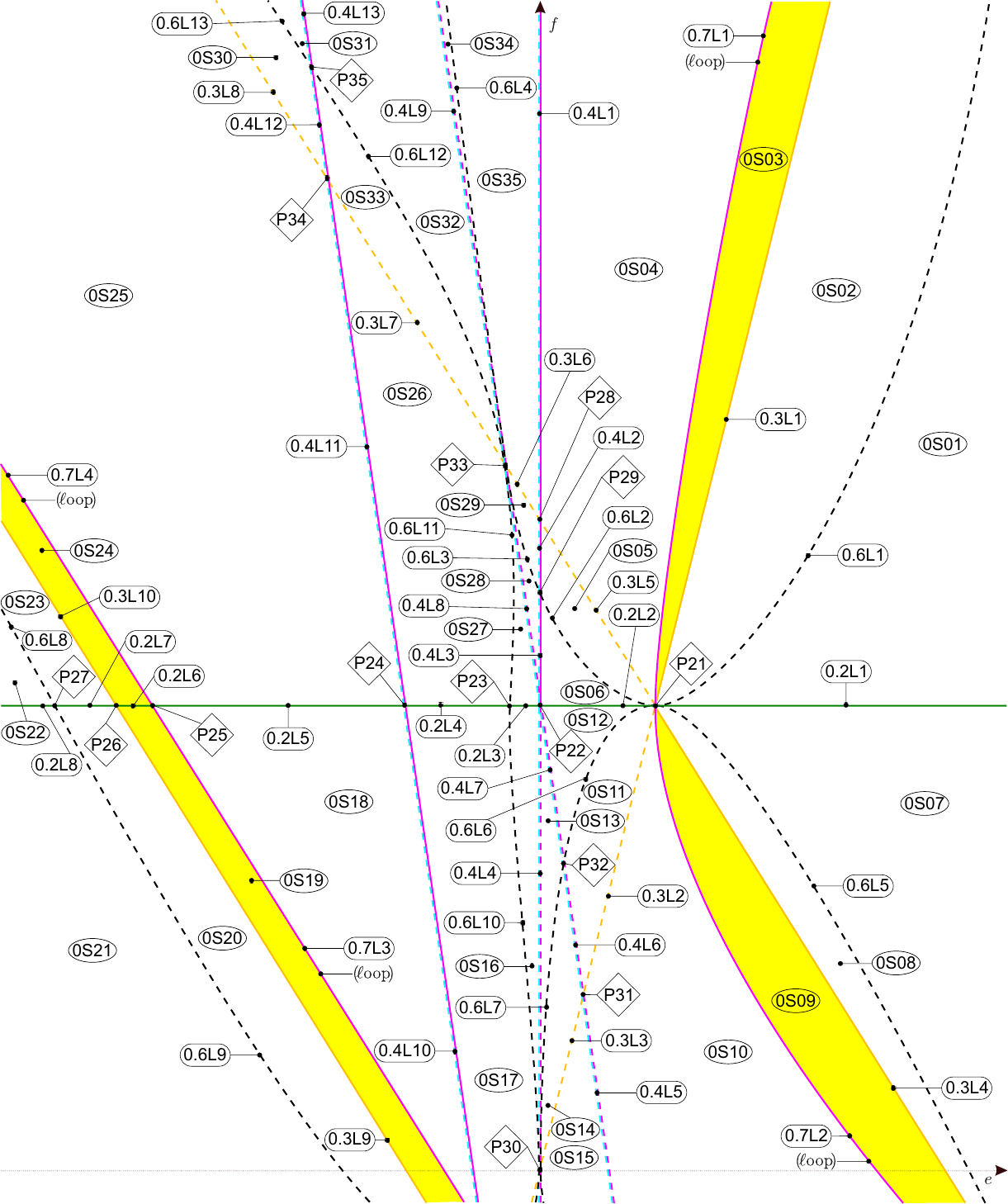} 
	\caption{\small \label{fig:slice-QES-A-31} Singular slice of the parameter space when $c=-1$}
\end{figure}

Now we consider the singular slice $c=-1$. One may say that this is a quite interesting singular slice, because:
\begin{itemize}
	\item Previously we mentioned that surface (${\cal S}_{0}$), related to a presence of an infinite elliptic--saddle of type 
	$\widehat{\!{1\choose 2}\!\!}\ E-H$, defines the entire plane $c=-1$. As it was pointed out in \cite{Artes-Mota-Rezende-2021c}
	each phase portrait obtained in the study of this slice is topologically equivalent to a phase portrait obtained in a neighborhood of
	this plane. However, in order to have a coherent bifurcation diagram, this plane must be studied. Here we follow the pattern established
	in Remark~\ref{rmk-colors} and we shall not draw this plane in brown color.
	\item So far we had the existence of three purple bifurcation straight lines and three cyan bifurcation straight lines. 
	For this value of the parameter $c$ we observe a coalescence among pairs of these straight lines. In fact, 
	calculation show that
	$$({\cal S}_{4})\vert_{c=-1}=({\cal S}_{8})\vert_{c=-1}=e(-2e-f-1)(2e+f-1),$$
	so the bifurcation straight lines $e=0$, $f=-2e-1$, and $f=-2e+1$ have multiplicity two.
\end{itemize}
In  Fig.~\ref{fig:slice-QES-A-31} we present the entire slice $c=-1$ completely labeled. In such a 
figure we use the pattern set out in \cite{Artes-Mota-Rezende-2021b} in order to draw the bifurcation 
straight lines which are double, and in order to present a label for each region (in this case the open
regions are labeled as pieces of surface (${\cal S}_{0}$), a bifurcation curve $X$ is labeled as 
$0.XL_{j}$, $j\in\mathbb{N}$, and each intersection of two or more bifurcation curves is indicated as a point.)

\begin{figure}[h!]
	\centering\includegraphics[width=1\textwidth]{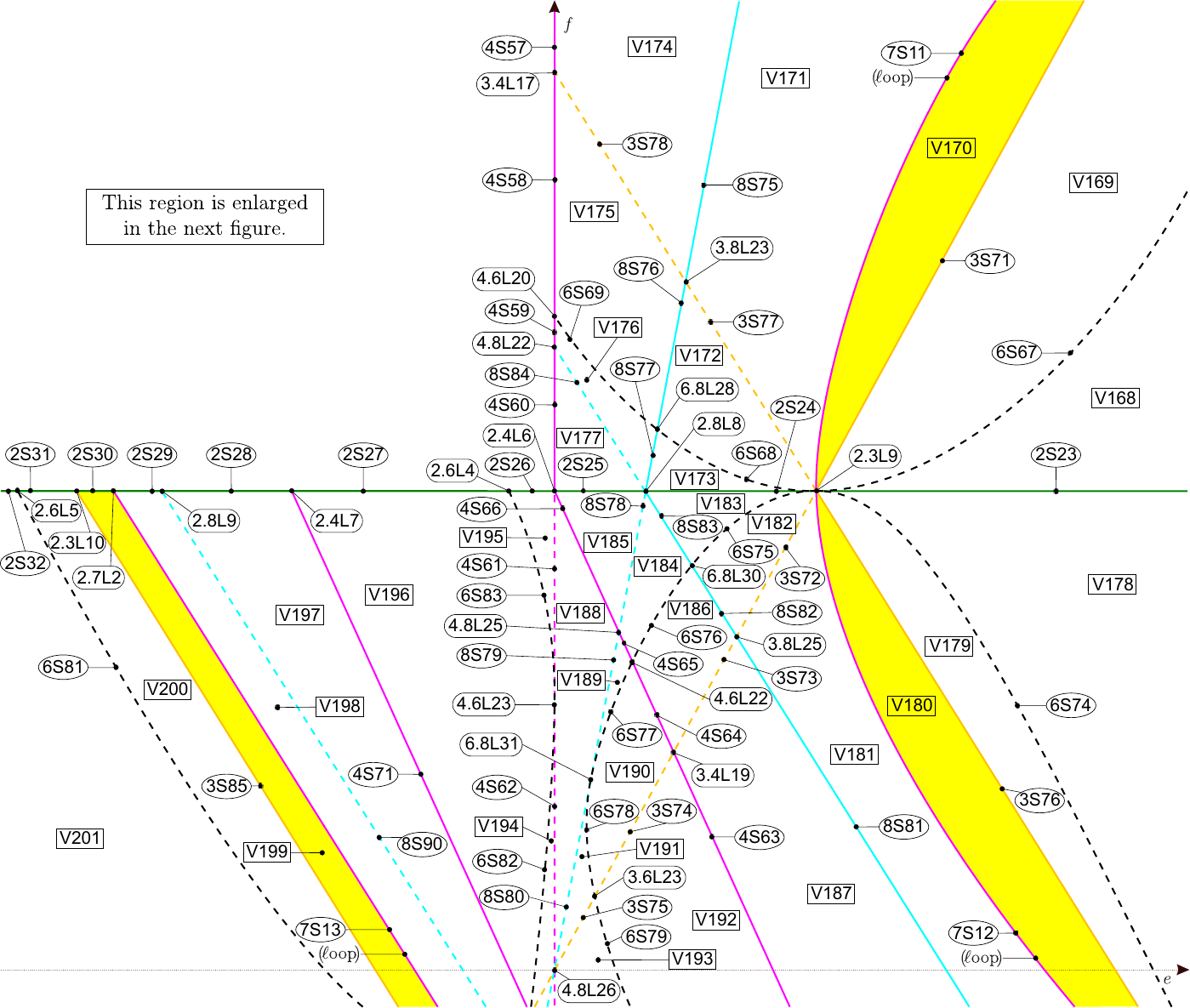} 
	\caption{\small \label{fig:slice-QES-A-32a} Piece of generic slice of the parameter space when 
		$c=-125/100$, see also Fig.~\ref{fig:slice-QES-A-32b}}
\end{figure}

\begin{figure}[h!]
	\centering\includegraphics[width=0.7\textwidth]{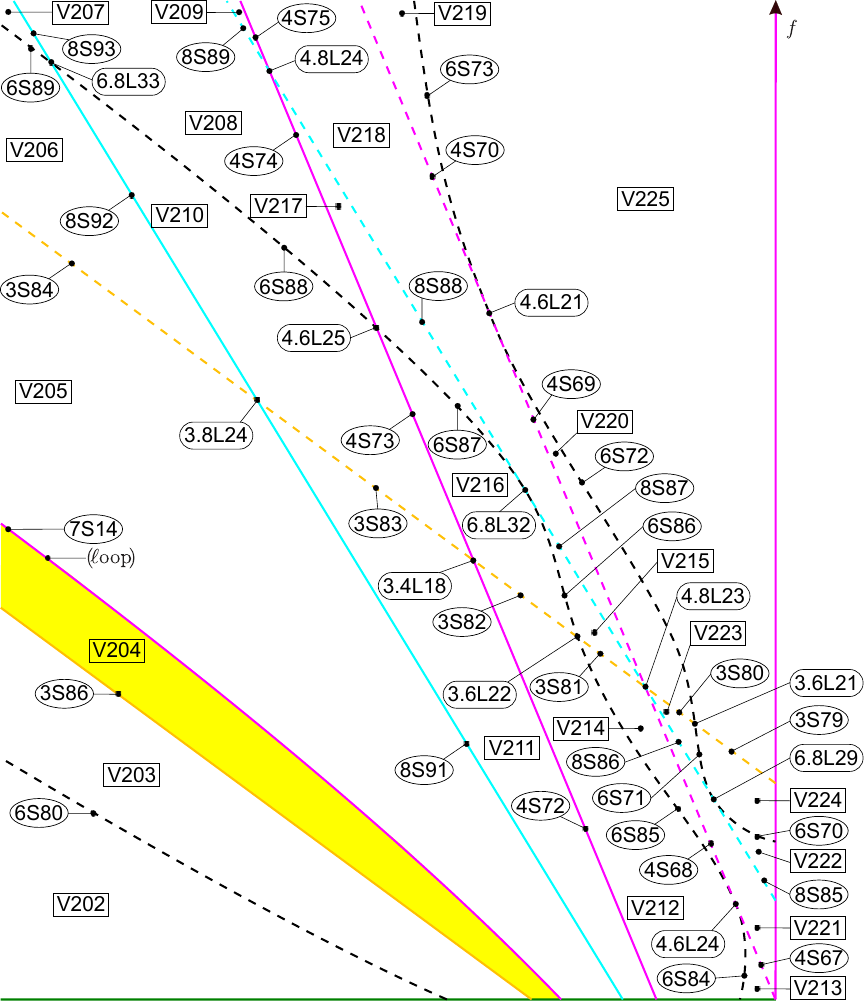} 
	\caption{\small \label{fig:slice-QES-A-32b} Continuation of Fig.~\ref{fig:slice-QES-A-32a}}
\end{figure}

From the list of slices presented in \eqref{eq:values-of-c-QES-A} we observe that the generic slice
to be considered now is $c=-125/100$. Doing the study of this entire slice we observe that the purple
and cyan bifurcation straight lines split themselves into three purple and three cyan bifurcation straight lines.
Also, it is clear that in this case we no longer have  (${\cal S}_{0}$)$\equiv0$. This generic slice is 
presented in Figs.~\ref{fig:slice-QES-A-32a} and \ref{fig:slice-QES-A-32b}.

\begin{figure}[h!]
	\centering\includegraphics[width=0.22\textwidth]{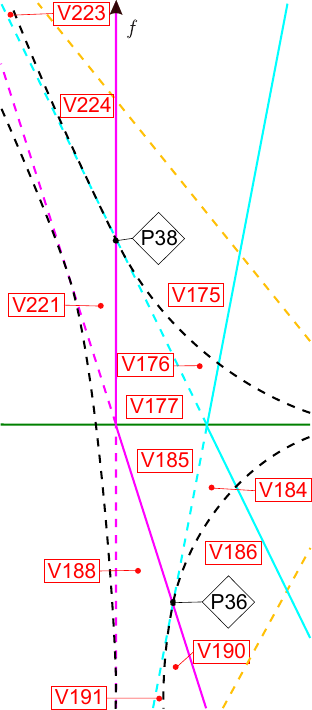} 
	\caption{\small \label{fig:slice-QES-A-33a} Piece of singular slice of the parameter space when 
		$c=-3/2$, see also Fig.~\ref{fig:slice-QES-A-33b}}
\end{figure}

\begin{figure}[h!]
	\centering\includegraphics[width=0.2\textwidth]{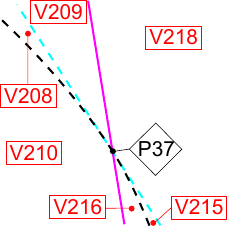} 
	\caption{\small \label{fig:slice-QES-A-33b} Another piece of singular slice of the parameter space when 
		$c=-3/2$, see also Fig.~\ref{fig:slice-QES-A-33a} and compare this region with Fig.~\ref{fig:slice-QES-A-32b}}
\end{figure}

Now we consider the singular slice $c=-3/2$. At this value, the volume regions $V_{189}$ 
(see Fig.~\ref{fig:slice-QES-A-32a}), $V_{217}$, and $V_{222}$ (see Fig.~\ref{fig:slice-QES-A-32b}) 
are reduced to the points $P_{36}$, $P_{37}$, and $P_{38}$, respectively, and these points are
presented in Fig.~\ref{fig:slice-QES-A-33a} and \ref{fig:slice-QES-A-33b}.

\begin{figure}[h!]
	\centering\includegraphics[width=0.25\textwidth]{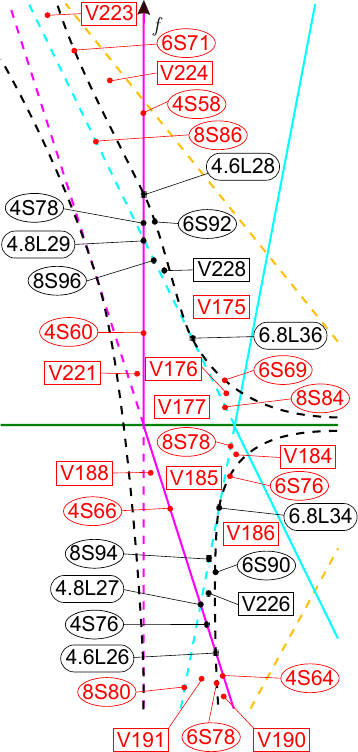} 
	\caption{\small \label{fig:slice-QES-A-34a} Piece of generic slice of the parameter space when 
		$c=-16/10$, see also Fig.~\ref{fig:slice-QES-A-34b}}
\end{figure}

\begin{figure}[h!]
	\centering\includegraphics[width=0.24\textwidth]{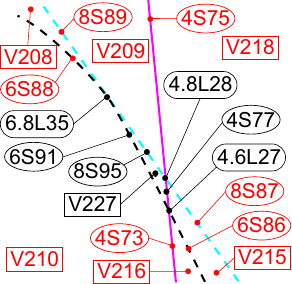} 
	\caption{\small \label{fig:slice-QES-A-34b} Another piece of generic slice of the parameter space when 
		$c=-16/10$, see also Fig.~\ref{fig:slice-QES-A-34a} and compare this region with Fig.~\ref{fig:slice-QES-A-33b}}
\end{figure}

Now, as it was expected, if we consider $c=-16/10$ as a generic slice, from the points $P_{36}$, $P_{37}$, 
and $P_{38}$ we get volume regions $V_{226}$, $V_{227}$, and $V_{228}$, respectively, which 
can be seeing in Fig.~\ref{fig:slice-QES-A-34a} and \ref{fig:slice-QES-A-34b}.

\begin{figure}[h!]
	\centering\includegraphics[width=0.16\textwidth]{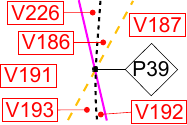} 
	\caption{\small \label{fig:slice-QES-A-35a} Piece of singular slice of the parameter space when 
		$c=-\sqrt{3}$, see also Fig.~\ref{fig:slice-QES-A-35b} and compare with Fig.~\ref{fig:slice-QES-A-32a}}
\end{figure}

\begin{figure}[h!]
	\centering\includegraphics[width=0.35\textwidth]{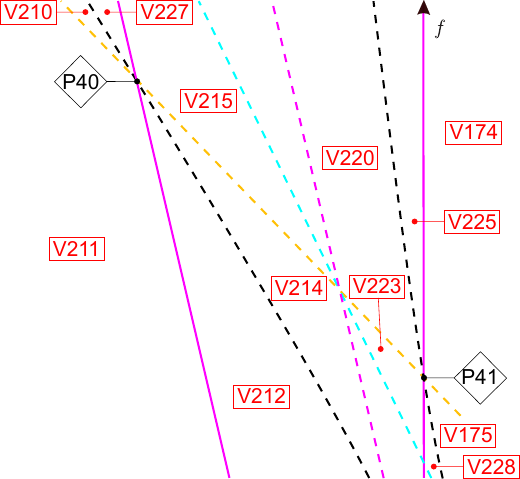} 
	\caption{\small \label{fig:slice-QES-A-35b} Another piece of singular slice of the parameter space when 
		$c=-\sqrt{3}$, see also Fig.~\ref{fig:slice-QES-A-35a} and compare with Fig.~\ref{fig:slice-QES-A-32a} 
		and \ref{fig:slice-QES-A-32b}}
\end{figure}

For the singular slice $c=-\sqrt{3}$, the volume regions $V_{190}$ (see Fig.~\ref{fig:slice-QES-A-32a}) 
$V_{216}$, and $V_{224}$ (see Fig.~\ref{fig:slice-QES-A-32b}) are reduced to the points $P_{39}$, 
$P_{40}$, and $P_{41}$, respectively, see Fig.~\ref{fig:slice-QES-A-35a} and \ref{fig:slice-QES-A-35b}.

\begin{figure}[h!]
	\centering
	\includegraphics[width=0.2\textwidth]{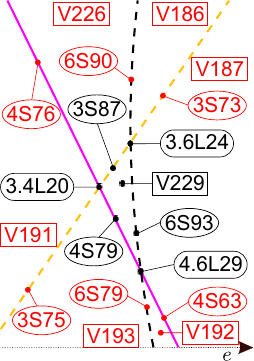}
	\caption{\small \label{fig:slice-QES-A-36a} Piece of generic slice of the parameter space when $c=-175/100$, see also Fig.~\ref{fig:slice-QES-A-36b} and \ref{fig:slice-QES-A-36c}}
\end{figure}

\begin{figure}[h!]
	\centering
	\includegraphics[width=0.4\textwidth]{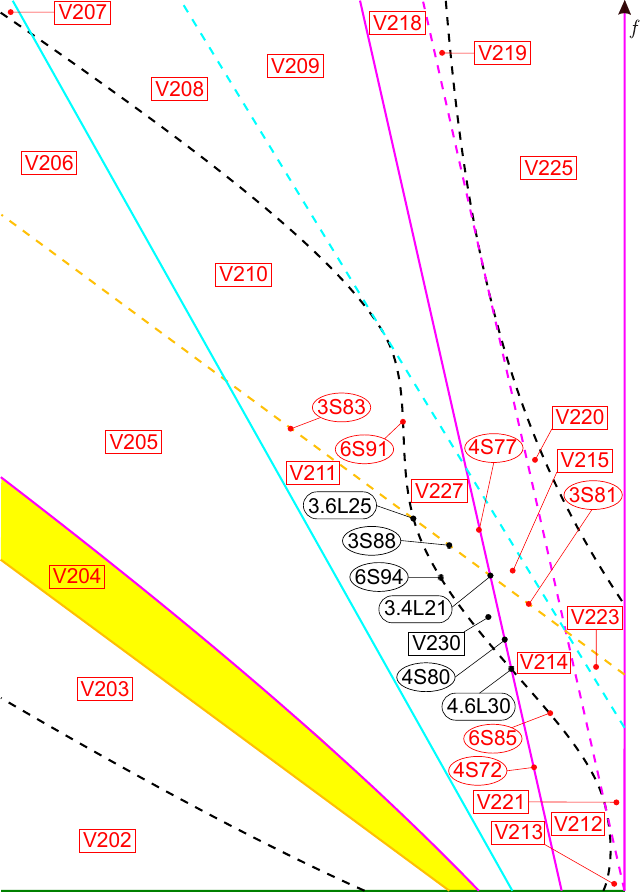}
	\caption{\small \label{fig:slice-QES-A-36b} Another piece of generic slice of the parameter space when $c=-175/100$, see also Fig.~\ref{fig:slice-QES-A-36a} and \ref{fig:slice-QES-A-36c}}
\end{figure}

\begin{figure}[h!]
	\centering\includegraphics[width=0.25\textwidth]{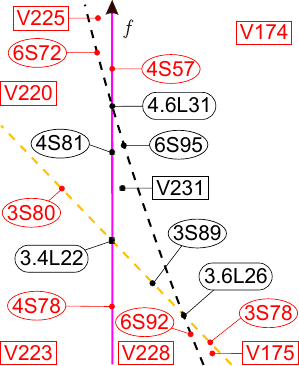}
	\caption{\small \label{fig:slice-QES-A-36c} Another piece of generic slice of the parameter space when 
		$c=-175/100$, see also Fig.~\ref{fig:slice-QES-A-36a} and \ref{fig:slice-QES-A-36b}}
\end{figure}

By considering the generic slice $c=-175/100$, from the points $P_{39}$, $P_{40}$, and $P_{41}$ 
we obtain volume regions $V_{229}$, $V_{230}$, and $V_{231}$, respectively, see Fig.~\ref{fig:slice-QES-A-36a}, \ref{fig:slice-QES-A-36b}, and \ref{fig:slice-QES-A-36c}.

\begin{figure}[h!]
	\centering\includegraphics[width=0.25\textwidth]{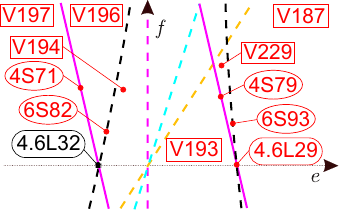} 
	\caption{\small \label{fig:slice-QES-A-37} Piece of singular slice of the parameter space when 
		$c=-16/9$, compare with Fig.~\ref{fig:slice-QES-A-32a} and \ref{fig:slice-QES-A-36a}}
\end{figure}

For the singular slice $c=-16/9$, we have that $4.6L_{29}$ (Fig.~\ref{fig:slice-QES-A-36a}) goes to $f=0$ 
and $V_{192}$ goes to $f<0$. Also, we have an intersection between $4S_{71}$ and $6S_{82}$
(Fig.~\ref{fig:slice-QES-A-32a}) at $4.6L_{32}$. See these phenomena along $f=0$ in Fig.~\ref{fig:slice-QES-A-37}.

\begin{figure}[h!]
	\centering\includegraphics[width=0.98\textwidth]{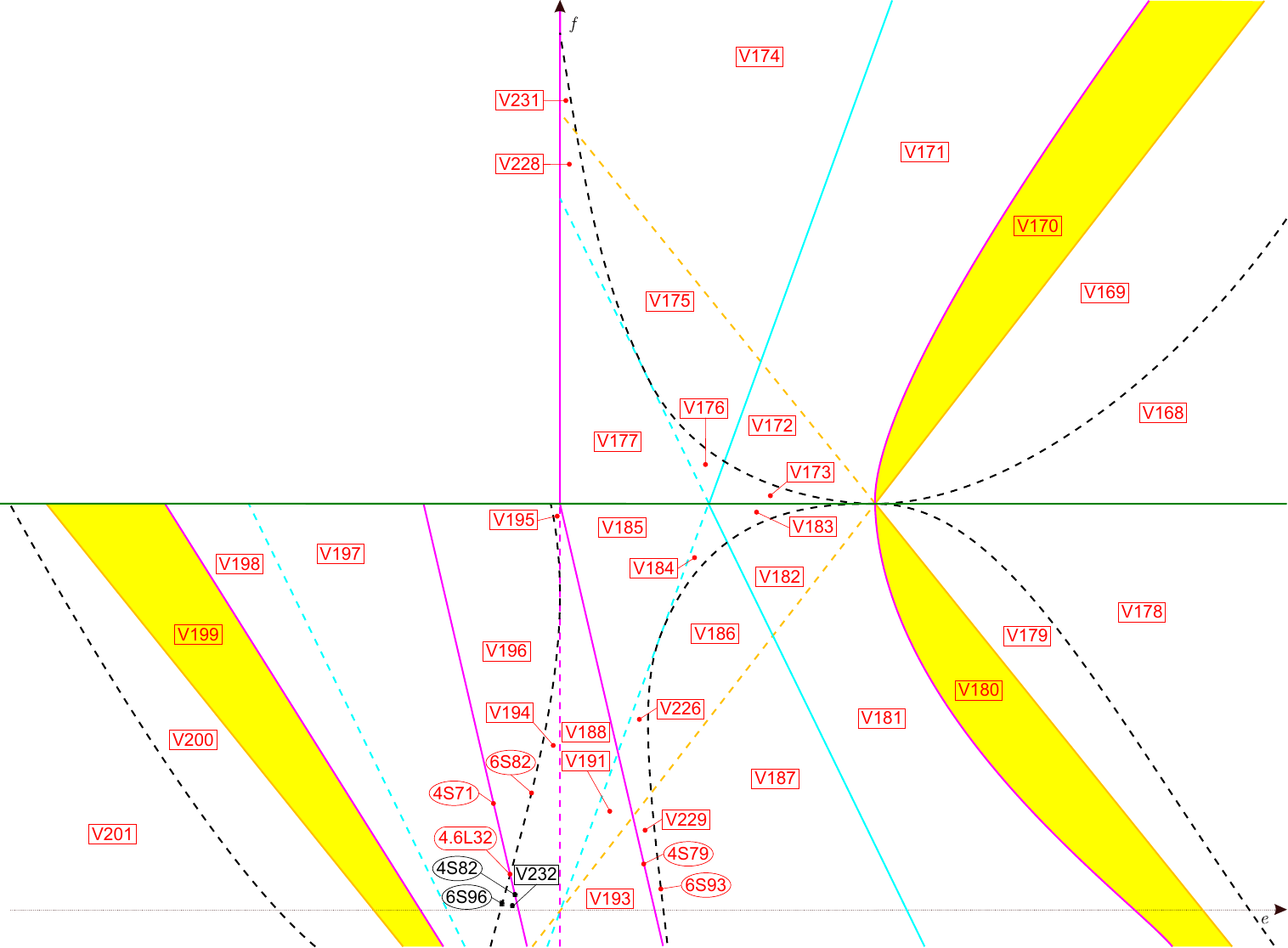} 
	\caption{\small \label{fig:slice-QES-A-38} Piece of generic slice of the parameter space when 
		$c=-19/10$, see again Fig.~\ref{fig:slice-QES-A-36b}}
\end{figure}

Taking into consideration Fig.~\ref{fig:slice-QES-A-37}, when we perform the study of the 
generic slice $c=-19/10$ we observe that $4.6L_{29}$ goes to $f<0$ e $4.6L_{32}$ goes 
to $f>0$ and it arises volume region $V_{232}$, see Fig.~\ref{fig:slice-QES-A-38}. We
point out that Fig.~\ref{fig:slice-QES-A-36b} can be considered as a continuation of 
Fig.~\ref{fig:slice-QES-A-38} since we did not detect any change in that region.

\begin{figure}[h!]
	\includegraphics[width=0.9\textwidth]{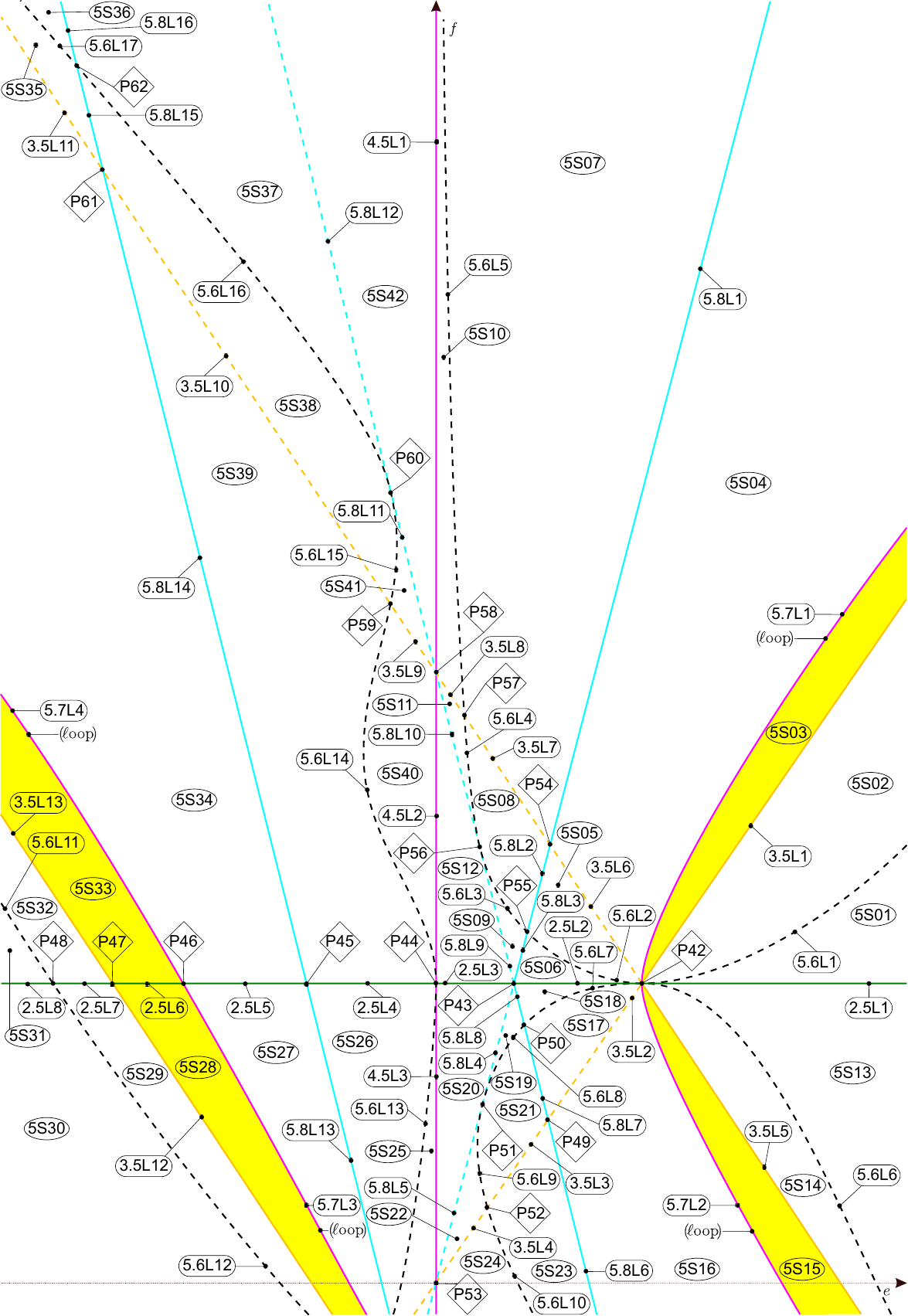}
	\caption{\small \label{fig:slice-QES-A-39} Singular slice of the parameter space when $c=-2$}
\end{figure}

Now we consider the singular slice $c=-2$. This is another interesting and important singular slice.
\begin{itemize}
	\item Surface ($\mathcal{S}_{5}$)$= c+2$ is related to a coalescence of infinite singular points. 
	Remember that if $e\ne0$ the phase portraits obtained in the study of this slice possess at most 
	one pair of infinite singular points and, if $e=0$ the corresponding phase portraits have the line 
	at infinity filled up with singularities. Here we follow Remark~\ref{rmk-colors} and we shall not 
	draw the slice $c=-2$ in red color.
	\item So far we had the existence of three purple bifurcation straight lines. For this value of the parameter 
	$c$ we observe that they coalesce along $e=0$. In fact, calculation show that
	$$({\cal S}_{4})\vert_{c=-2}=-4e^3,$$
	so the bifurcation straight line $e=0$ has multiplicity three.
\end{itemize}
In  Fig.~\ref{fig:slice-QES-A-39} we present the entire slice $c=-2$ completely labeled. In such a 
figure we use the same pattern as the one used in the slice $c=-1$ in order to present a label for 
each region.

\begin{remark}\label{type-nilp-inf-sing} It is important to mention that the infinite nilpotent
	singularity is always an elliptic--saddle of type:
	\begin{itemize}
		\item $\widehat{\!{1\choose 2}\!\!}\ PEP-H$, for all $c>-1$;
		\item $\widehat{\!{1\choose 2}\!\!}\ E-H$, for $c=-1$; and 
		\item $\widehat{\!{1\choose 2}\!\!}\ E-PHP$, for $-2<c<-1$.
	\end{itemize}
	In addition, when $c=-2$ we had an infinite nilpotent saddle--node and, for all $c<-2$ we shall have infinite nilpotent saddles $\widehat{\!{1\choose 2}\!\!}\ HHH-H$.
\end{remark}

\begin{figure}[h!]
	\centering\includegraphics[width=0.9\textwidth]{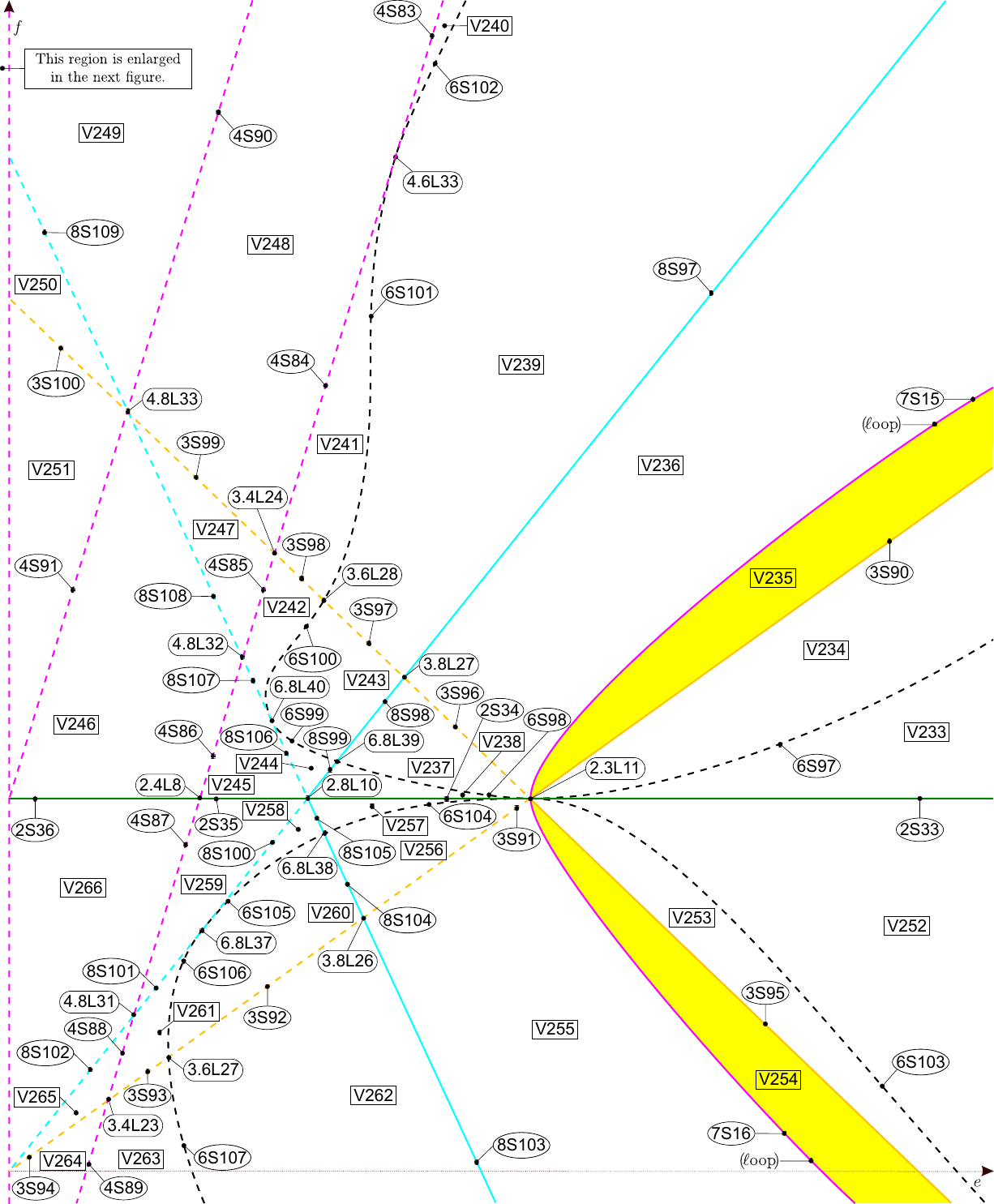} 
	\caption{\small \label{fig:slice-QES-A-40a} Piece of generic slice of the parameter space when 
		$c=-3$, see also Fig.~\ref{fig:slice-QES-A-40b}}
\end{figure}

\begin{figure}[h!]
	\centering\includegraphics[width=0.9\textwidth]{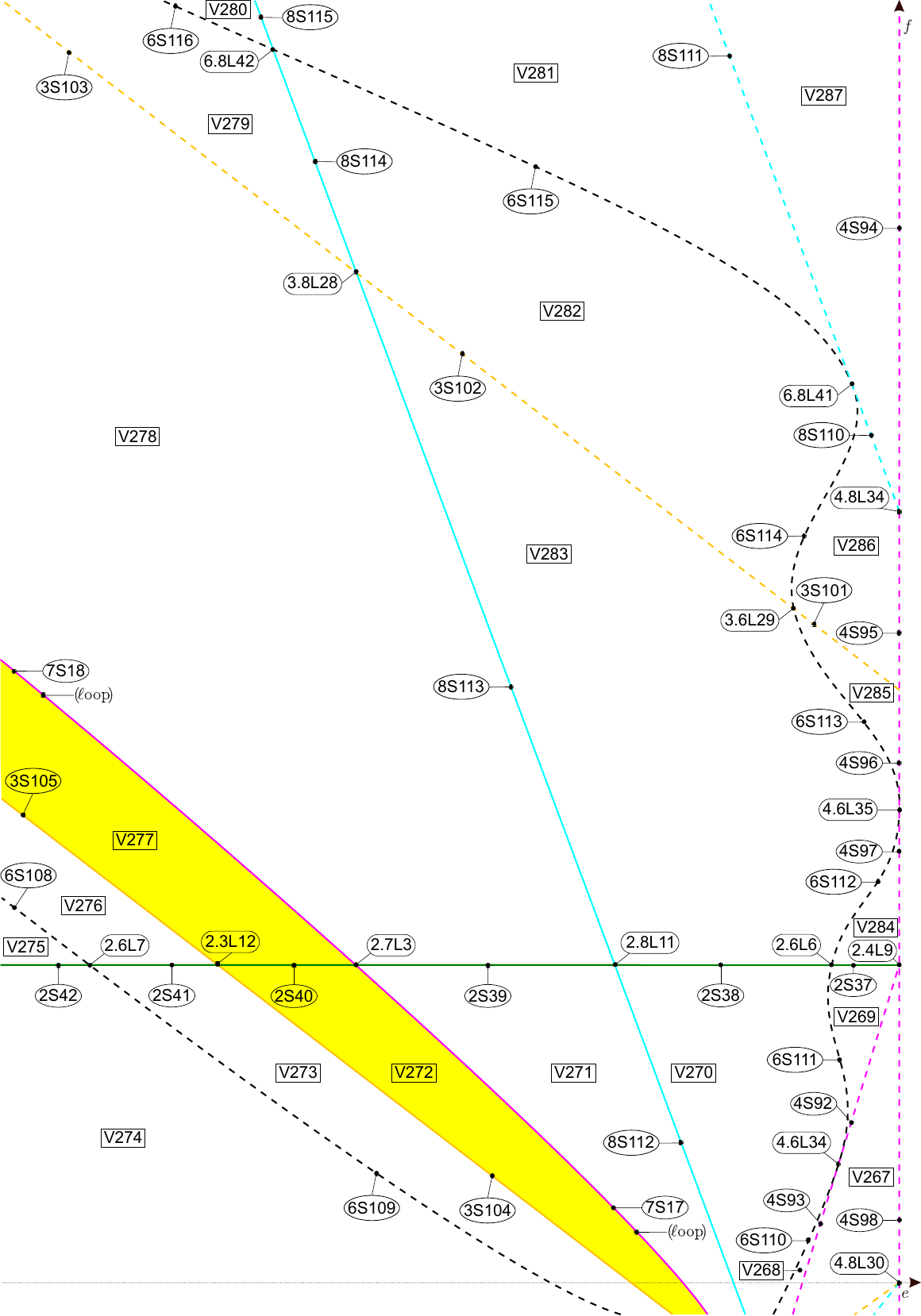} 
	\caption{\small \label{fig:slice-QES-A-40b} Continuation of Fig.~\ref{fig:slice-QES-A-40a}}
\end{figure}

Now we present the study of the generic slice $c=-3$. In this case, the triple purple bifurcation straight line
from $c=-2$ splits itself into three bifurcation straight lines. Moreover, here we no longer have a coalescence
of infinite singular points, given by surface ($\mathcal{S}_{5}$). This generic slice is presented in 
Fig.~\ref{fig:slice-QES-A-40a} and \ref{fig:slice-QES-A-40b}.

\begin{figure}[h!]
	\centering\includegraphics[width=0.3\textwidth]{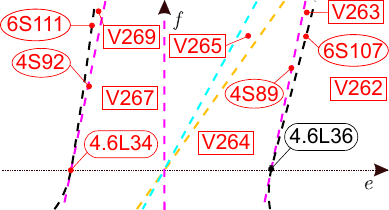} 
	\caption{\small \label{fig:slice-QES-A-41} Piece of singular slice of the parameter space when 
		$c=-4$, compare with Fig.~\ref{fig:slice-QES-A-40a} and \ref{fig:slice-QES-A-40b}}
\end{figure}

Consider Fig.~\ref{fig:slice-QES-A-40a} and \ref{fig:slice-QES-A-40b}. When we perform the study of
the singular slice $c=-4$ we notice that $4.6L_{34}$ goes to $f=0$ (carrying $V_{268}$ to $f<0$) and 
we also have that $4S_{89}$ intercepts $6S_{107}$ at $4.6L_{36}$, see Fig.~\ref{fig:slice-QES-A-41}.

\begin{figure}[h!]
	\centering\includegraphics[width=0.95\textwidth]{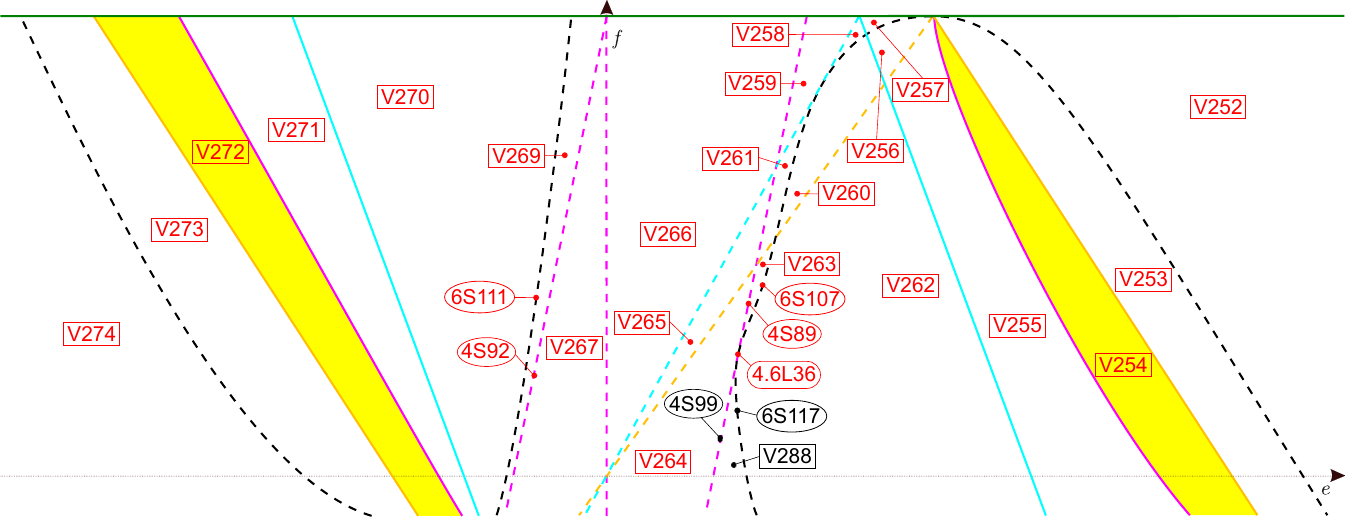} 
	\caption{\small \label{fig:slice-QES-A-42} Piece of generic slice of the parameter space when 
		$c=-5$, compare with Fig.~\ref{fig:slice-QES-A-41}}
\end{figure}

Finally we consider the last generic slice from the list presented in \eqref{eq:values-of-c-QES-A},
namely, $c=-5$. In this slice we observe that $4.6L_{34}$ goes to $f<0$ and $4.6L_{36}$ goes
to $f>0$, giving place to the appearance of volume region $V_{288}$, see Fig.~\ref{fig:slice-QES-A-42}.

Since there is coherence among the generic and singular slices presented before, no more slices are needed for the complete coherence of the bifurcation diagram. So, all the values of the parameter $c$ in \eqref{eq:values-of-c-QES-A} are sufficient for the coherence of the bifurcation diagram. Thus, we can affirm that we have described a complete bifurcation diagram for class $\overline{\QESA}$ modulo islands and modulo any other nonalgebraic slice (above or below, or very close to $c=0$), as we discuss in Sec.~\ref{sec:islands-QESA}.

\subsubsection{Other relevant facts about the bifurcation diagram} \label{sec:islands-QESA}

The bifurcation diagram we have obtained for the class $\overline{\QESA}$ is completely coherent, i.e. in this family, by taking any two points in the parameter space and joining them by a continuous curve, along this curve the changes in phase portraits that occur when crossing the different bifurcation surfaces we mention can be completely explained.

Nevertheless, we cannot be sure that this bifurcation diagram is the complete bifurcation diagram for $\overline{\QESA}$ due to the possibility of the existence of ``islands'' inside the parts bordered by unmentioned bifurcation surfaces. In case they exist, these ``islands'' would not mean any modification of the nature of the singular points. So, on the border of these ``islands'' we could only have bifurcations due to saddle connections or multiple limit cycles.

In case there were more bifurcation surfaces, we should still be able to join two representatives of any two parts of the 1274 parts of $\overline{\QESA}$ found until now with a continuous curve either without crossing such a bifurcation surface or, in case the curve crosses it, it must do it an even number of times without tangencies, otherwise one must take into account the multiplicity of the tangency, so the total number must be even. This is why we call these potential bifurcation surfaces ``\textit{islands}''.

However, we have not found a different phase portrait which could fit in such an island. A potential ``island'' would be the set of parameters for which the phase portraits possess a double limit cycle and this ``island'' would be inside the parts where $W_{4}<0$ since we have the presence of a focus.

\subsubsection{Completion of the proof of Theorem~\ref{th:main-thm-QES-A}} \label{sec:invariants-QESA} 

In the bifurcation diagram we may have topologically equivalent phase portraits belonging to distinct parts of the parameter space. As here we have 1274 distinct parts of the parameter space, to help us to identify or to distinguish phase portraits, we need to introduce some invariants and we actually choose integer valued, character and symbol invariants. Some of them were already used in \cite{Artes-Rezende-Oliveira-2013b} and \cite{Artes-Mota-Rezende-2021c}, but we recall them and introduce some needed ones. These invariants yield a classification which is easier to grasp.

\begin{definition}\label{def:invariant-I1-QESA}
	We denote by $I_{1}(S)$ the number of real finite singular points.
\end{definition}

\begin{definition}\label{def:invariant-I2-QESA}
	We denote by $I_{2}(S)$ the sum of the indices of the isolated real finite singular points.
\end{definition}

\begin{definition}\label{def:invariant-I3-QESA}
	We denote by $I_{3}(S)$ the number of real infinite singular points. We note that this number can also be infinite, which is represented by $\infty$.
\end{definition}

\begin{definition}\label{def:invariant-I4-QESA}
	For a given infinite singularity $s$ of a system $S$, let $\l_s$ be the number of global or local separatrices beginning or ending at $s$ and which do not lie on the line at infinity. We have $0\leq\l_s\leq4$. We denote by $I_{4}(S)$ the sequence of all such $\l_s$ when $s$ moves in the set of infinite singular points of the system $S$. We start the sequence at the infinite singular point which receives (or sends) the greatest number of separatrices and take the direction which yields the greatest absolute value, e.g. the values $2110$ and $2011$ for this invariant are symmetrical (and, therefore, they are the same), so we consider $2110$.
\end{definition}

\begin{definition}\label{def:invariant-I5-QESA}
	We denote by $I_{5}(S)$ the number of graphics different from the orbits of the elliptic sector (including the border of the elliptic sector).
\end{definition}

\begin{definition}\label{def:invariant-I6-QESA} 
	We denote by $I_{6}(S)$ a character from the set $\{\emptyset,\overline{sn}_{(2)},\widehat{cp}_{(2)}\}$ which indicate the following types of finite multiple singularities, respectively: none (in this case the system does not contain a finite multiple singularity), saddle--node, and cusp.
\end{definition}

\begin{definition}\label{def:invariant-I7-QESA} 
	We denote by $I_{7}(S)$ a character from the set $\{\emptyset,\ell,{f\!-\!i}\}$ which indicate the following types of separatrix connection, respectively: none (in this case the system does not contain a separatrix connection), $\ell$oop, and \textbf{f}inite--\textbf{i}nfinite.
\end{definition}

\begin{definition}\label{def:invariant-I8-QESA}
	We denote by $I_{8}(S)$ the number of limit cycles around a foci.
\end{definition}

\begin{definition}\label{def:invariant-I9-QESA}
	We denote by $I_{9}(S)$ the number of separatrices arriving or leaving one real finite antisaddle. In case we have two real finite antisaddles this invariant is given by a pair $(A,B)$ where $A$ and $B$ indicate the corresponding numbers of separatrices arriving or leaving each antisaddle.
\end{definition}

\begin{definition}\label{def:invariant-I10-QESA}
	We denote by $I_{10}(S)$ an element from the set $\{c, f(s), f(u)\}$, indicating the type of the real finite singularity located inside the region bordered by the graphic, which can be of the following types, respectively: center, stable focus, and unstable focus.
\end{definition}

As we have noted previously in Remark~\ref{rem:f-n}, we do not distinguish between phase portraits whose only difference is that in one we have a finite node and in the other a focus. Both phase portraits are topologically equivalent and they can only be distinguished within the $C^1$ class. In case we may want to distinguish between them, a new invariant may easily be introduced.

\begin{theorem} \label{th:QESA-inv}
	Consider the class $\overline{\QESA}$ and all the phase portraits that we have obtained for this family. The values of the affine invariant ${\cal I} = (I_{1}, I_{2}, I_{3}, I_{4}, I_{5}, I_{6}, I_{7}, I_{8}, I_{9}, I_{10})$ given in the diagram from Tables \ref{tab:geom-classif-QESA-pt1} to \ref{tab:geom-classif-QESA-pt4} yield a partition of these phase portraits of the class $\overline{\QESA}$.
	
	Furthermore, for each value of $\cal I$ in this diagram there corresponds a single phase portrait; i.e. $S$ and $S'$ are such that ${\cal I}(S)={\cal I}(S')$, if and only if $S$ and $S'$ are topologically equivalent.
\end{theorem}

The bifurcation diagram for $\overline{\QESA}$ has 1274 parts which produce 91 topologically different phase portraits as described in Tables~\ref{tab:geom-classif-QESA-pt1} to \ref{tab:top-equiv-QESA-pt11}. The remaining 1183 parts do not produce any new phase portrait which was not included in the 91 previous ones. The difference is basically the presence of a strong focus instead of a node and vice versa, weak points, and a presence of invariant algebraic curves (lines or parabolas) which do not represent a separatrix connection.

The phase portraits having neither limit cycle nor graphic have been denoted surrounded by parenthesis, for example $(V_{233})$; the phase portraits having one limit cycle have been denoted surrounded by brackets, for example $[V_{235}]$; the phase portraits having one graphic have been denoted surrounded by $\{\ast\}$ and those ones having two or more graphics have been denoted surrounded by $\{\!\{\ast\}\!\}$, for example $\{2S_{39}\}$ and $\{\!\{4S_{59}\}\!\}$, respectively. Moreover, the phase portraits having one limit cycle and more than one graphic have been denoted surrounded by $[\{\{\ast\}\}]$, for example $[\{\{2S_{18}\}\}]$.

\begin{proof}[Proof of Theorem~\ref{th:QESA-inv}]
	The above result follows from the results in the previous sections and a careful analysis of the bifurcation diagrams given in Sec.~\ref{subsec:bd-QES-A}, in Figs.~\ref{fig:slice-QES-A-01-alg1} and \ref{fig:slice-QES-A-01-alg2} to Fig.~\ref{fig:slice-QES-A-42}, the definition of the invariants $I_{j}$ and their explicit values for the corresponding phase portraits.
\end{proof}

We recall some observations regarding the equivalence relations used in this study: the affine and time rescaling, $C^1$ and topological equivalences.

The coarsest one among these three is the topological equivalence and the finest is the affine equivalence. We can have two systems which are topologically equivalent but not $C^1-$equivalent. For example, we could have a system with a finite antisaddle which is a structurally stable node and in another system with a focus, the two systems being topologically equivalent but belonging to distinct $C^1-$equivalence classes, separated by the surface $({\cal S}_{6})$ on which the node turns into a focus.

In Tables~\ref{tab:top-equiv-QESA-pt1} to \ref{tab:top-equiv-QESA-pt11} we list in the first column 91 parts with all the distinct phase portraits of Figs.~\ref{fig:pp-QES-A-1} to \ref{fig:pp-QES-A-3}. Corresponding to each part listed in column one we have in each row all parts whose phase portraits are topologically equivalent to the phase portrait appearing in column 1 of the same row.

In the second column we set all the parts whose systems yield topologically equivalent phase portraits to those in the first column, but which may have some algebro--geometric features related to the position of the orbits. In the third column we present all the parts which are topologically equivalent to the ones from the first column having a focus instead of a node.

In the fourth (respectively, fifth; and sixth) column we list all parts whose phase portraits have a node which is at a bifurcation point producing foci close to the node in perturbations, a node--focus to shorten (respectively, a finite weak singular point; and possess an invariant curve (straight line and/or parabola) not yielding a connection of separatrices).

The last column refers to other reasons associated to different geometrical aspects and they are described as follows:
\begin{enumerate}[(1)]
	\vspace{-2mm}
	\item  The phase portraits correspond to symmetric parts of the bifurcation diagram;
	\vspace{-2mm}
	\item  the phase portrait possesses a singularity of type $\widehat{\!{1\choose2}\!\!}\ E-H$ at infinity.
\end{enumerate}

Whenever phase portraits appear in a row in a specific column, the listing is done according to the decreasing dimension of the parts where they appear, always placing the lower dimensions on lower lines.

\subsubsection{Proof of Theorem~\ref{th:main-thm-QES-A}}

The bifurcation diagram described in Sec. \ref{subsec:bd-QES-A}, plus Tables~\ref{tab:geom-classif-QESA-pt1} to \ref{tab:geom-classif-QESA-pt4} of the geometrical invariants distinguishing the 91 phase portraits, plus Tables~\ref{tab:top-equiv-QESA-pt1} to \ref{tab:top-equiv-QESA-pt11} giving the equivalences with the remaining phase portraits lead to the proof of Theorem~\ref{th:main-thm-QES-A}.

\begin{landscape}
	
	\begin{table}\caption{\small Geometric classification for the family $\QESA$}\label{tab:geom-classif-QESA-pt1} 
		\begin{center}
			\[
			I_{1}\!=\!
			\left\{
			\begin{array}{ll}
			2 \,\, \& \,\, I_{2}\!=\!
			\left\{
			\begin{array}{ll}
			-1 \,\, \& \,\, I_{3}\!=\! 2 \,\, \& \,\, I_{4}\!=\!
			\left\{
			\begin{array}{ll}
			2210 \,\, \left\{\left\{2.4L_{1}\right\}\right\}, \\  			
			3101 \,\, \left\{\left\{2.8L_{2}\right\}\right\}, \\  			
			3201 \,\, \& \,\, I_{5}\!=\!
			\left\{
			\begin{array}{ll}
			1 \,\, \& \,\, I_{6}\!=\!
			\left\{
			\begin{array}{ll}
			\widehat{cp}_{(2)} \,\, \left\{\left\{2.3L_{2}\right\}\right\}, \\  
			\overline{sn}_{(2)} \,\, \left\{\left\{2S_{6}\right\}\right\}, \\  
			\end{array}
			\right. \\
			2 \,\, \left\{\left\{2S_{1}\right\}\right\}, \\  
			\end{array}
			\right. \\
			3310 \,\, \left\{\left\{2S_{4}\right\}\right\}, \\  
			4201 \,\, \left\{\left\{2S_{5}\right\}\right\}, \\  
			\end{array}
			\right. \\
			1 \,\, \& \,\, I_{3}\!=\!
			\left\{
			\begin{array}{ll}											
			1 \,\, \& \,\, I_{4}\!=\!
			\left\{
			\begin{array}{ll}
			21 \,\, \& \,\, I_{5}\!=\!
			\left\{
			\begin{array}{ll}
			0 \,\, \left(P_{43}\right), \\  													
			1 \,\, \left\{P_{45}\right\}, \\  													
			\end{array}
			\right. \\
			22 \,\, \& \,\, I_{5}\!=\!
			\left\{
			\begin{array}{ll}
			0 \,\, \& \,\, I_{6}\!=\!
			\left\{
			\begin{array}{ll}
			\widehat{cp}_{(2)} \,\, \left(P_{42}\right), \\  
			\overline{sn}_{(2)} \,\, \left(2.5L_{2}\right), \\  
			\end{array}
			\right. \\
			1 \,\, \& \,\, I_{6}\!=\!\overline{sn}_{(2)} \,\, \& \,\, I_{7}\!=\!
			\left\{
			\begin{array}{ll}
			\emptyset \,\, \left\{2.5L_{5}\right\}, \\  
			\ell \,\, \left\{P_{46}\right\}, \\  
			\end{array}
			\right. \\
			\end{array}
			\right. \\
			31 \,\, \& \,\, I_{5}\!=\!
			\left\{
			\begin{array}{ll}
			0 \,\, \left(2.5L_{3}\right), \\  													
			1 \,\, \left\{2.5L_{4}\right\}, \\  													
			\end{array}
			\right. \\
			32 \,\, \& \,\, I_{5}\!=\!0 \,\, \& \,\, I_{6}\!=\!\overline{sn}_{(2)} \,\, \& \,\, I_{7}\!=\!\emptyset \,\, \& \,\, I_{8}\!=\!
			\left\{
			\begin{array}{ll}
			0  \,\, \& \,\, I_{9}\!=\!
			\left\{
			\begin{array}{ll}
			1 \,\, \left(2.5L_{8}\right), \\  													
			2 \,\, \left(2.5L_{1}\right), \\  													
			\end{array}
			\right. \\
			1 \,\, \left[2.5L_{6}\right], \\  
			\end{array}
			\right. \\
			\end{array}
			\right. \\
			2 \,\, \& \,\, I_{4}\!=\!
			\left\{
			\begin{array}{ll}
			1110 \,\, \& \,\, I_{5}\!=\!
			\left\{
			\begin{array}{ll}
			1 \,\, \left\{\left\{2.4L_{4}\right\}\right\}, \\  													
			3 \,\, \left\{\left\{2.4L_{5}\right\}\right\}, \\  													
			\end{array}
			\right. \\											
			2100 \,\, \& \,\, I_{5}\!=\!
			\left\{
			\begin{array}{ll}
			0 \,\, \left(2.8L_{10}\right), \\  													
			1 \,\, \left\{2.8L_{11}\right\}, \\  													
			\end{array}
			\right. \\
			2101 \,\, \& \,\, I_{5}\!=\!
			\left\{
			\begin{array}{ll}
			0 \,\, \left(2S_{35}\right), \\  													
			1 \,\, \& \,\, I_{6}\!=\!
			\left\{
			\begin{array}{ll}
			\widehat{cp}_{(2)} \,\, \left\{\left\{2.3L_{7}\right\}\right\}, \\  
			\overline{sn}_{(2)} \,\, \left\{\left\{2S_{12}\right\}\right\}, \\  
			\end{array}
			\right. \\
			2 \,\, \& \,\, I_{6}\!=\!\overline{sn}_{(2)} \,\, \& \,\, I_{7}\!=\!
			\left\{
			\begin{array}{ll}
			\emptyset \,\, \& \,\, I_{8}\!=\!0 \,\, \& \,\, I_{9}\!=\!
			\left\{
			\begin{array}{ll}
			1 \,\, \left\{\left\{2S_{17}\right\}\right\}, \\  													
			2 \,\, \left\{\left\{2S_{13}\right\}\right\}, \\  													
			\end{array}
			\right. \\
			\ell \,\, \left\{\left\{2.7L_{1}\right\}\right\}, \\  
			\end{array}
			\right. \\												
			\end{array}
			\right. \\
			\mathcal{A}_{1} \,\, \hbox{\it (next page)}, \\
			\end{array}
			\right. \\				
			\infty \,\, \left\{\left\{P_{44}\right\}\right\}, \\  
			\end{array}
			\right. \\
			\end{array}
			\right. \\
			\mathcal{A}_{2} \,\,  \hbox{\it (next page)}, \\
			\end{array} 
			\right.
			\]
		\end{center}
	\end{table}
	
	\begin{table}\caption{\small Geometric classification for the family $\QESA$ \textit{(cont.)}}\label{tab:geom-classif-QESA-pt2}
		\begin{center}
			\[
			%I_{1}\!=\!
			%\left\{
			\begin{array}{ll}
			\begin{matrix}\mathcal{A}_{1} \\ \begin{bmatrix}I_{1}\!=\!2,\\ I_{2}\!=\!1,\\ I_{3}\!=\!2\end{bmatrix}\end{matrix} \,\, \& \,\, I_{4}\!=\!
			\left\{
			\begin{array}{ll}
			2111 \,\, \& \,\, I_{5}\!=\!
			\left\{
			\begin{array}{ll}
			1 \,\, \left\{\left\{2.4L_{6}\right\}\right\}, \\  													
			4 \,\, \left\{\left\{2.4L_{7}\right\}\right\}, \\  													
			\end{array}
			\right. \\											
			2121 \,\, \& \,\, I_{5}\!=\!
			\left\{
			\begin{array}{ll}
			1 \,\, \left\{\left\{2.8L_{8}\right\}\right\}, \\  													
			2 \,\, \left\{\left\{2S_{26}\right\}\right\}, \\  													
			3 \,\, \left\{\left\{2.8L_{9}\right\}\right\}, \\  													
			\end{array}
			\right. \\
			2200 \,\, \& \,\, I_{5}\!=\!
			\left\{
			\begin{array}{ll}									
			0 \,\, \& \,\, I_{6}\!=\!
			\left\{
			\begin{array}{ll}
			\widehat{cp}_{(2)} \,\, \left(2.3L_{11}\right), \\  
			\overline{sn}_{(2)} \,\, \left(2S_{34}\right), \\  
			\end{array}
			\right. \\
			1 \,\, \& \,\, I_{6}\!=\!\overline{sn}_{(2)} \,\, \& \,\, I_{7}\!=\!
			\left\{
			\begin{array}{ll}
			\emptyset \,\, \left\{2S_{39}\right\}, \\  
			\ell \,\, \left\{2.7L_{3}\right\}, \\  
			\end{array}
			\right. \\												
			\end{array}
			\right. \\
			3101 \,\, \& \,\, I_{5}\!=\!1 \,\, \& \,\, I_{6}\!=\!\overline{sn}_{(2)} \,\, \& \,\, I_{7}\!=\!\emptyset \,\, \& \,\, I_{8}\!=\!
			\left\{
			\begin{array}{ll}
			0  \,\, \& \,\, I_{9}\!=\!
			\left\{
			\begin{array}{ll}
			1 \,\, \left\{\left\{2S_{20}\right\}\right\}, \\  													
			2 \,\, \left\{\left\{2S_{11}\right\}\right\}, \\  													
			\end{array}
			\right. \\
			1 \,\, [\left\{\left\{2S_{18}\right\}\right\}], \\  
			\end{array}
			\right. \\
			3121 \,\, \& \,\, I_{5}\!=\!
			\left\{
			\begin{array}{ll}									
			1 \,\, \& \,\, I_{6}\!=\!
			\left\{
			\begin{array}{ll}
			\widehat{cp}_{(2)} \,\, \left\{\left\{2.3L_{9}\right\}\right\}, \\  
			\overline{sn}_{(2)}\,\, \& \,\, I_{7}\!=\!\emptyset \,\, \& \,\, I_{8}\!=\!0\,\, \& \,\, I_{9}\!=\!
			\left\{
			\begin{array}{ll}
			3 \,\, \left\{\left\{2S_{25}\right\}\right\}, \\  													
			4 \,\, \left\{\left\{2S_{24}\right\}\right\}, \\  													
			\end{array}
			\right. \\
			\end{array}
			\right. \\
			2 \,\, \& \,\, I_{6}\!=\!\overline{sn}_{(2)} \,\, \& \,\, I_{7}\!=\!
			\left\{
			\begin{array}{ll}
			\emptyset \,\, \left\{\left\{2S_{29}\right\}\right\}, \\  
			\ell \,\, \left\{\left\{2.7L_{2}\right\}\right\}, \\  
			\end{array}
			\right. \\					
			3 \,\, \left\{\left\{2S_{28}\right\}\right\}, \\  		
			\end{array}
			\right. \\			
			3200 \,\, \& \,\, I_{5}\!=\!0 \,\, \& \,\, I_{6}\!=\!\overline{sn}_{(2)} \,\, \& \,\, I_{7}\!=\!\emptyset \,\, \& \,\, I_{8}\!=\!
			\left\{
			\begin{array}{ll}
			0  \,\, \& \,\, I_{9}\!=\!
			\left\{
			\begin{array}{ll}
			1 \,\, \left(2S_{42}\right), \\  													
			2 \,\, \left(2S_{33}\right), \\  													
			\end{array}
			\right. \\
			1 \,\, \left[2S_{40}\right], \\  
			\end{array}
			\right. \\							
			4121 \,\, \& \,\, I_{5}\!=\!1 \,\, \& \,\, I_{6}\!=\!\overline{sn}_{(2)} \,\, \& \,\, I_{7}\!=\!\emptyset \,\, \& \,\, I_{8}\!=\!
			\left\{
			\begin{array}{ll}
			0  \,\, \& \,\, I_{9}\!=\!
			\left\{
			\begin{array}{ll}
			1 \,\, \left\{\left\{2S_{32}\right\}\right\}, \\  													
			3 \,\, \left\{\left\{2S_{23}\right\}\right\}, \\  													
			\end{array}
			\right. \\
			1 \,\, [\left\{\left\{2S_{30}\right\}\right\}], \\  
			\end{array}
			\right. \\																																			
			\end{array}
			\right. \\
			%3 \,\, \& \,\, I_{2}\!=\!\mathcal{A}_{2} \,\, \hbox{\it (next page)}, \\
			\end{array} 
			%\right.
			\]
		\end{center}
	\end{table}
	
	\begin{table}\caption{\small Geometric classification for the family $\QESA$ \textit{(cont.)}}\label{tab:geom-classif-QESA-pt3}
		\begin{center}
			\[
			%I_{1}\!=\!
			%   \left\{
			\begin{array}{ll}
			\begin{matrix}\mathcal{A}_{2} \\ \begin{bmatrix}I_{1}\!=\!3\end{bmatrix}\end{matrix} \,\, \& \,\, I_{2}\!=\!
			\left\{
			\begin{array}{ll}
			-1 \,\, \& \,\, I_{3}\!=\! 2 \,\, \& \,\, I_{4}\!=\!
			\left\{
			\begin{array}{ll}
			2101 \,\, \left\{\left\{4.8L_{2}\right\}\right\}, \\  			
			2210 \,\, \left\{\left\{4S_{5}\right\}\right\}, \\  			
			3101 \,\, \left\{\left\{8S_{7}\right\}\right\}, \\  			
			3201 \,\, \& \,\, I_{5}\!=\!
			\left\{
			\begin{array}{ll}
			1 \,\, \left\{\left\{V_{1}\right\}\right\}, \\   
			2 \,\, \& \,\, I_{6}\!=\!\emptyset\,\, \& \,\, I_{7}\!=\!\ell\,\, \& \,\, I_{8}\!=\!0\,\, \& \,\, I_{9}\!=0\!\,\, \& \,\, I_{10}\!=\!
			\left\{
			\begin{array}{ll}
			c \,\, \left\{\left\{3.7L_{1}\right\}\right\}, \\  
			f(s) \,\, \left\{\left\{7S_{1}\right\}\right\}, \\  
			f(u) \,\, \left\{\left\{7S_{4}\right\}\right\}, \\  
			\end{array}
			\right. \\
			\end{array}
			\right. \\
			3310 \,\, \& \,\, I_{5}\!=\!1 \,\, \& \,\, I_{6}\!=\!\emptyset\,\, \& \,\, I_{7}\!=\!\emptyset \,\, \& \,\, I_{8}\!=\!
			\left\{
			\begin{array}{ll}
			0 \,\, \left\{\left\{V_{9}\right\}\right\}, \\  
			1 \,\, [\left\{\left\{V_{11}\right\}\right\}], \\  
			\end{array}
			\right. \\
			4201 \,\, \& \,\, I_{5}\!=\!1 \,\, \& \,\, I_{6}\!=\!\emptyset\,\, \& \,\, I_{7}\!=\!\emptyset \,\, \& \,\, I_{8}\!=\!
			\left\{
			\begin{array}{ll}
			0 \,\, \left\{\left\{V_{12}\right\}\right\}, \\  
			1 \,\, [\left\{\left\{V_{66}\right\}\right\}], \\  
			\end{array}
			\right. \\
			\end{array}
			\right. \\
			1 \,\, \& \,\, I_{3}\!=\!
			\left\{
			\begin{array}{ll}											
			1 \,\, \& \,\, I_{4}\!=\!
			\left\{
			\begin{array}{ll}
			21  \,\, \left(5.8L_{3}\right), \\  
			22 \,\, \& \,\, I_{5}\!=\!
			\left\{
			\begin{array}{ll}
			0 \,\, \left(5S_{6}\right), \\  													
			1 \,\, \left\{5.7L_{1}\right\}, \\  													
			\end{array}
			\right. \\
			31  \,\, \left(5S_{9}\right), \\  
			32 \,\, \& \,\, I_{5}\!=\!0 \,\, \& \,\, I_{6}\!=\!\emptyset\,\, \& \,\, I_{7}\!=\!\emptyset \,\, \& \,\, I_{8}\!=\!
			\left\{
			\begin{array}{ll}
			0 \,\, \left(5S_{1}\right), \\  
			1 \,\, \left[5S_{3}\right], \\  
			\end{array}
			\right. \\
			\end{array}
			\right. \\
			2 \,\, \& \,\, I_{4}\!=\!
			\left\{
			\begin{array}{ll}
			1110 \,\, \left\{\left\{4S_{34}\right\}\right\}, \\  													
			2100 \,\, \left(8S_{99}\right), \\  													
			2101 \,\, \& \,\, I_{5}\!=\!
			\left\{
			\begin{array}{ll}
			0 \,\, \left(V_{240}\right), \\  													
			1\,\, \& \,\, I_{6}\!=\!\emptyset\,\, \& \,\, I_{7}\!=\!\emptyset \,\, \& \,\, I_{8}\!=\!0\,\, \& \,\, I_{9}\!=\!
			\left\{
			\begin{array}{ll}
			(1,3) \,\, \left\{\left\{V_{94}\right\}\right\}, \\  													
			(2,2) \,\, \left\{\left\{V_{101}\right\}\right\}, \\  													
			\end{array}
			\right. \\								
			2 \,\, \left\{\left\{7S_{7}\right\}\right\}, \\  													
			\end{array}
			\right. \\
			2111 \,\, \left\{\left\{4S_{59}\right\}\right\}, \\  			
			2121\,\, \& \,\, I_{5}\!=\!1 \,\, \& \,\, I_{6}\!=\!\emptyset \,\, \& \,\, I_{7}\!=\!
			\left\{
			\begin{array}{ll}
			\emptyset \,\, \left\{\left\{V_{188}\right\}\right\}, \\  													
			f\!-\!{i}\,\, \left\{\left\{8S_{77}\right\}\right\}, \\  													
			\end{array}
			\right. \\													
			2200 \,\, \& \,\, I_{5}\!=\!
			\left\{
			\begin{array}{ll}
			0 \,\, \left(V_{238}\right), \\  													
			1 \,\, \left\{7S_{15}\right\}, \\  													
			\end{array}
			\right. \\
			\mathcal{A}_{3} \,\, \hbox{\it (next page)}, \\
			\end{array}
			\right. \\				
			\infty \,\, \left\{\left\{4.5L_{1}\right\}\right\}, \\  
			\end{array}
			\right. \\
			\end{array}
			\right. \\
			\end{array} 
			%\right.
			\]
		\end{center}
	\end{table}
	
	\begin{table}\caption{\small Geometric classification for the family $\QESA$ \textit{(cont.)}}\label{tab:geom-classif-QESA-pt4}
		\begin{center}
			\[
			%I_{1}\!=\!
			%   \left\{
			\begin{array}{ll}
			\begin{matrix}\mathcal{A}_{3} \\ \begin{bmatrix}I_{1}\!=\!3,\\ I_{2}\!=\!1,\\ I_{3}\!=\!2\end{bmatrix}\end{matrix} \,\, \& \,\, I_{4}\!=\!
			\left\{
			\begin{array}{ll}
			3101 \,\, \& \,\, I_{5}\!=\!1 \,\, \& \,\, I_{6}\!=\!\emptyset\,\, \& \,\, I_{7}\!=\!\emptyset \,\, \& \,\, I_{8}\!=\!
			\left\{
			\begin{array}{ll}
			0 \,\, \left\{\left\{V_{89}\right\}\right\}, \\  
			1 \,\, [\left\{\left\{V_{91}\right\}\right\}], \\  
			\end{array}
			\right. \\
			3121\,\, \& \,\, I_{5}\!=\!
			\left\{
			\begin{array}{ll}
			1 \,\, \& \,\, I_{6}\!=\!\emptyset \,\, \& \,\, I_{7}\!=\!\emptyset\,\, \& \,\, I_{8}\!=\!0 \,\, \& \,\, I_{9}\!=\!
			\left\{
			\begin{array}{ll}
			(1,3) \,\, \left\{\left\{V_{173}\right\}\right\}, \\  
			(2,2) \,\, \left\{\left\{V_{176}\right\}\right\}, \\  
			\end{array}
			\right. \\
			2 \,\, \left\{\left\{7S_{11}\right\}\right\}, \\  													
			\end{array}
			\right. \\				
			3200 \,\, \& \,\, I_{5}\!=\!0 \,\, \& \,\, I_{6}\!=\!\emptyset\,\, \& \,\, I_{7}\!=\!\emptyset \,\, \& \,\, I_{8}\!=\!
			\left\{
			\begin{array}{ll}
			0 \,\, \left(V_{233}\right), \\  
			1 \,\, \left[V_{235}\right], \\  
			\end{array}
			\right. \\
			4121 \,\, \& \,\, I_{5}\!=\!1 \,\, \& \,\, I_{6}\!=\!\emptyset\,\, \& \,\, I_{7}\!=\!\emptyset \,\, \& \,\, I_{8}\!=\!
			\left\{
			\begin{array}{ll}
			0 \,\, \left\{\left\{V_{168}\right\}\right\}, \\  
			1 \,\, [\left\{\left\{V_{170}\right\}\right\}], \\  
			\end{array}
			\right. \\
			\end{array}
			\right. \\
			\end{array} 
			%\right.
			\]
		\end{center}
	\end{table}
	
	\begin{table}\caption{\small Topological equivalences for the family $\QESA$}\label{tab:top-equiv-QESA-pt1}
		\begin{center}
			\begin{tabular}{ccccccc}
				\hline
				Presented & Identical       & Finite      & Finite       & Finite  &  Possessing       &              \\
				phase     & under           & antisaddle  & antisaddle   & weak    &  invariant curve   &  Other reasons\\
				portrait  & perturbations   & focus       & node--focus  & point   &           (no separatrix)        &              \\
				\hline
				\multirow{5}{*}{$V_{1}$} & $V_{2}$,  $V_{3}$,  $V_{4}$,  $V_{5}$,  $V_{6}$ & $V_{8}$, $V_{48}$ &   &   & & $V_{15}^{(1)}$,  $V_{16}^{(1)}$,  $V_{17}^{(1)}$,  $V_{18}^{(1)}$, $V_{19}^{(1)}$, $V_{20}^{(1)}$ \\  
				& $V_{7}$, $V_{36}$,  $V_{37}$,  $V_{38}$ &   $V_{49}$,  $V_{55}$    & &  &   &  $V_{21}^{(1)}$,  $V_{22}^{(1)}$, $V_{23}^{(1)}$, $V_{24}^{(1)}$, $V_{26}^{(1)}$,  $V_{27}^{(1)}$ \\
				& $V_{39}$,  $V_{40}$, $V_{54}$, $V_{56}$ &  $V_{64}$,  $V_{74}$    & & &   &  $V_{28}^{(1)}$, $V_{29}^{(1)}$,  $V_{30}^{(1)}$, $V_{31}^{(1)}$, $V_{32}^{(1)}$, $V_{33}^{(1)}$\\
				& $V_{63}$, $V_{82}$,  $V_{83}$, $V_{84}$ &          & & &   &   $V_{46}^{(1)}$, $V_{47}^{(1)}$, $V_{69}^{(1)}$, $V_{75}^{(1)}$, $V_{80}^{(1)}$, $V_{81}^{(1)}$\\
				& $V_{87}$,  $V_{88}$ &         & & &   & $V_{85}^{(1)}$, $V_{86}^{(1)}$\\
				&   &   $3S_{4}$, $3S_{18}$     & $6S_{1}$, $6S_{15}$,  $6S_{16}$ &$3S_{1}$,  $3S_{2}$ & $4S_{1}$,  $4S_{2}$,  $4S_{3}$  & $3S_{8}^{(1)}$,  $3S_{9}^{(1)}$,  $3S_{10}^{(1)}$, $3S_{11}^{(1)}$, $3S_{12}^{(1)}$\\
				& &  $3S_{25}$,  $3S_{31}$       & $6S_{17}$,  $6S_{18}$,  $6S_{20}$ & $3S_{3}$, $3S_{17}$  &  $4S_{4}$, $4S_{7}$, $4S_{20}$ & $3S_{13}^{(1)}$, $3S_{34}^{(1)}$,  $3S_{39}^{(1)}$, $3S_{41}^{(1)}$,  $3S_{42}^{(1)}$  \\
				& &         & $6S_{21}$, $6S_{30}$,  $6S_{31}$ & $3S_{43}$, $3S_{44}$ &  $4S_{22}$,  $4S_{28}$, $8S_{1}$  & $3S_{45}^{(1)}$,  $3S_{46}^{(1)}$, $4S_{12}^{(1)}$,  $4S_{13}^{(1)}$,  $4S_{14}^{(1)}$  \\
				& &         &$6S_{32}$ &  $3S_{47}$, $3S_{48}$ & $8S_{2}$,  $8S_{3}$,  $8S_{4}$  &  $4S_{15}^{(1)}$, $4S_{16}^{(1)}$,  $4S_{17}^{(1)}$,  $4S_{18}^{(1)}$,  $4S_{19}^{(1)}$ \\
				& &         & &   &  $8S_{20}$,  $8S_{27}$ & $6S_{4}^{(1)}$,  $6S_{5}^{(1)}$,  $6S_{6}^{(1)}$,  $6S_{7}^{(1)}$,  $6S_{9}^{(1)}$  \\
				& &         & &   & $8S_{28}$, $8S_{31}$  & $6S_{10}^{(1)}$,  $6S_{11}^{(1)}$,  $6S_{12}^{(1)}$, $6S_{28}^{(1)}$,  $6S_{29}^{(1)}$  \\
				& &         & &   &   & $8S_{8}^{(1)}$,  $8S_{9}^{(1)}$,  $8S_{10}^{(1)}$,  $8S_{11}^{(1)}$,  $8S_{12}^{(1)}$  \\
				& &         & &   &   & $8S_{13}^{(1)}$, $8S_{29}^{(1)}$,  $8S_{30}^{(1)}$  \\
				& &         & $3.6L_{3}$,  $3.6L_{4}$ & $3.4L_{1}$,  $3.4L_{4}$  &  $4.8L_{1}$, $4.8L_{6}$ & $3.4L_{2}^{(1)}$,  $3.4L_{3}^{(1)}$, $3.6L_{1}^{(1)}$, $3.6L_{2}^{(1)}$ \\
				& &         & $3.6L_{13}$,  $3.6L_{14}$ & $3.8L_{10}$, $3.8L_{13}$ &   & $3.6L_{11}^{(1)}$,  $3.6L_{12}^{(1)}$, $3.8L_{11}^{(1)}$,  $3.8L_{12}^{(1)}$  \\
				& &         & $4.6L_{1}$, $4.6L_{7}$  &  &   & $4.6L_{3}^{(1)}$,  $4.6L_{5}^{(1)}$,  $4.8L_{3}^{(1)}$,  $4.8L_{4}^{(1)}$ \\
				& &         &  $6.8L_{5}$,  $6.8L_{10}$ &  &   &$6.8L_{2}^{(1)}$,  $6.8L_{3}^{(1)}$   \\
				\hline												
				\multirow{2}{*}{$V_{9}$} & $V_{35}$, $V_{57}$,  $V_{58}$ & $V_{10}$, $V_{42}$, $V_{60}$ &            & & & $V_{34}^{(1)}$, $V_{50}^{(1)}$, $V_{59}^{(1)}$, $V_{62}^{(1)}$, $V_{76}^{(1)}$ \\  
				& & $V_{61}$,  $V_{65}$,  $V_{73}$         &  &   && \\
				& &  $3S_{6}$,  $3S_{15}$, $3S_{22}$       &  $6S_{2}$, $6S_{14}$,  $6S_{22}$ & $3S_{19}$,  $3S_{20}$  &&  $3S_{16}^{(1)}$, $3S_{21}^{(1)}$, $3S_{24}^{(1)}$, $3S_{40}^{(1)}$ \\
				& &  $3S_{23}$, $3S_{26}$,  $3S_{30}$       &$6S_{23}$  &   && $6S_{8}^{(1)}$, $6S_{24}^{(1)}$ \\
				& &       &$3.6L_{5}$,  $3.6L_{6}$   &   &&  $3.6L_{7}^{(1)}$ \\
				\hline												
			\end{tabular}
		\end{center}
	\end{table}
	
	\begin{table}\caption{\small Topological equivalences for the family $\QESA$ \textit{(cont.)}}\label{tab:top-equiv-QESA-pt2}
		\begin{center}
			\begin{tabular}{ccccccc}
				\hline
				Presented & Identical       & Finite      & Finite       & Finite  &  Possessing       &              \\
				phase     & under           & antisaddle  & antisaddle   & weak    &  invariant curve   &  Other reasons\\
				portrait  & perturbations   & focus       & node--focus  & point   &           (no separatrix)        &              \\
				\hline
				\multirow{1}{*}{$V_{11}$}   &      $V_{43}$           &                  &               &           & & $V_{51}^{(1)}$\\  
				\hline
				\multirow{2}{*}{$V_{12}$} & $V_{25}$, $V_{67}$, $V_{70}$ & $V_{13}$,  $V_{14}$, $V_{44}$ &          &            & & $V_{41}^{(1)}$, $V_{52}^{(1)}$,  $V_{53}^{(1)}$, $V_{78}^{(1)}$,  $V_{79}^{(1)}$ \\  
				& &  $V_{45}$, $V_{68}$, $V_{71}$       &  &   & &  \\
				& &    $3S_{7}$,  $3S_{14}$,  $3S_{27}$     &$6S_{3}$,  $6S_{13}$ & $3S_{36}$,  $3S_{37}$  & & $3S_{5}^{(1)}$, $3S_{28}^{(1)}$,   $3S_{32}^{(1)}$,  $3S_{33}^{(1)}$ \\
				& &    $3S_{29}$, $3S_{35}$, $3S_{38}$    &  $6S_{25}$,  $6S_{26}$ &   & & $6S_{19}^{(1)}$, $6S_{27}^{(1)}$  \\
				& &       &  $3.6L_{8}$,  $3.6L_{9}$ &   & &  $3.6L_{10}^{(1)}$  \\
				\hline
				\multirow{1}{*}{$V_{66}$} &$V_{72}$ &  &            & & & $V_{77}^{(1)}$ \\  
				\hline
				\multirow{1}{*}{$V_{89}$} &$V_{95}$ &    $V_{90}$,  $V_{96}$       & & & & $V_{116}^{(1)}$,  $V_{117}^{(1)}$,  $V_{118}^{(1)}$,  $V_{119}^{(1)}$ \\  
				& &  $0S_{2}$, $0S_{8}$, $3S_{49}$      & $6S_{33}$,  $6S_{39}$ &  &    & $0S_{1}^{(2)}$, $0S_{7}^{(2)}$, $0S_{20}^{(1)}$,  $0S_{21}^{(1)}$,  $0S_{22}^{(1)}$ \\
				& &   $3S_{54}$    &  &  &    & $0S_{23}^{(1)}$, $3S_{63}^{(1)}$,  $3S_{64}^{(1)}$, $6S_{45}^{(1)}$,  $6S_{46}^{(1)}$ \\
				& &  $0.3L_{1}$,  $0.3L_{4}$    & $0.6L_{1}$,  $0.6L_{5}$ &  &    & $0.3L_{9}^{(1)}$,  $0.3L_{10}^{(1)}$, $0.6L_{8}^{(1)}$,  $0.6L_{9}^{(1)}$ \\
				\hline
				\multirow{1}{*}{$V_{91}$} & $V_{97}$ &            & & & & $V_{115}^{(1)}$,  $V_{120}^{(1)}$ \\  
				& &         &  &  &    & $0S_{3}^{(2)}$,  $0S_{9}^{(2)}$, $0S_{19}^{(1)}$,  $0S_{24}^{(1)}$    \\
				\hline
				\multirow{1}{*}{$V_{94}$} &    $V_{100}$    &    $V_{92}$,  $V_{93}$,  $V_{98}$,  $V_{99}$      & & & & $V_{114}^{(1)}$,  $V_{121}^{(1)}$,  $V_{122}^{(1)}$,  $V_{123}^{(1)}$\\  
				& &   $0S_{4}$,  $0S_{5}$, $0S_{10}$,  $0S_{11}$     &  $6S_{34}$,  $6S_{40}$&  &    & $0S_{6}^{(2)}$, $0S_{12}^{(2)}$, $0S_{18}^{(1)}$,  $0S_{25}^{(1)}$,  $0S_{30}^{(1)}$ \\
				& &   $3S_{50}$,  $3S_{55}$     &  &  &    &   $0S_{31}^{(1)}$, $3S_{62}^{(1)}$, $6S_{52}^{(1)}$ \\
				& &   $0.3L_{2}$,  $0.3L_{5}$     & $0.6L_{2}$,  $0.6L_{6}$ &  &    & $0.3L_{8}^{(1)}$, $0.6L_{13}^{(1)}$    \\
				\hline
			\end{tabular}
		\end{center}
	\end{table}
	
	\begin{table}\caption{\small Topological equivalences for the family $\QESA$ \textit{(cont.)}}\label{tab:top-equiv-QESA-pt3}
		\begin{center}
			\begin{tabular}{ccccccc}
				\hline
				Presented & Identical       & Finite      & Finite       & Finite  &  Possessing       &   Other           \\
				phase     & under           & antisaddle  & antisaddle   & weak    &  invariant curve   &  reasons  \\
				portrait  & perturbations   & focus       & node--focus  & point   &           (no separatrix)        &              \\
				\hline
				\multirow{7}{*}{$V_{101}$} & $V_{103}$,  $V_{104}$, $V_{106}$,  $V_{107}$, $V_{108}$, $V_{109}$ &  $V_{102}$, $V_{105}$ &            & & &\\  
				& $V_{110}$,  $V_{111}$, $V_{112}$, $V_{124}$, $V_{125}$, $V_{128}$  &  $V_{113}$, $V_{126}$        &  & &   &\\
				& $V_{129}$,  $V_{130}$,  $V_{131}$, $V_{132}$, $V_{133}$,  $V_{134}$  &   $V_{127}$, $V_{136}$      &  & &   &\\
				& $V_{135}$, $V_{137}$, $V_{138}$,  $V_{139}$,  $V_{140}$,  $V_{141}$ & $V_{143}$, $V_{149}$         &  & &   &\\
				& $V_{142}$, $V_{144}$,  $V_{145}$,  $V_{146}$,  $V_{147}$,  $V_{148}$  &  $V_{152}$, $V_{159}$       &  & &   &\\
				& $V_{150}$,  $V_{151}$, $V_{153}$,  $V_{154}$,  $V_{155}$, $V_{156}$ &   $V_{161}$,  $V_{162}$      &  & &   &\\
				&  $V_{157}$,  $V_{158}$, $V_{160}$, $V_{165}$, $V_{166}$,  $V_{167}$&    $V_{163}$,  $V_{164}$     &  & &   &\\
				& &  $0S_{14}$,  $0S_{15}$     &  $6S_{35}$,  $6S_{36}$,  $6S_{37}$,  $6S_{38}$, $6S_{41}$ & $3S_{52}$ &  $4S_{35}$,  $4S_{36}$,  $4S_{37}$, $4S_{41}$, $4S_{42}$  & $0S_{13}^{(2)}$\\
				& &  $0S_{17}$,  $0S_{26}$       & $6S_{42}$,  $6S_{43}$,  $6S_{44}$,  $6S_{47}$,  $6S_{48}$ & $3S_{53}$  & $4S_{43}$,  $4S_{44}$,  $4S_{45}$, $4S_{46}$,  $4S_{51}$   & $0S_{16}^{(2)}$\\
				& &  $0S_{29}$, $0S_{33}$       & $6S_{49}$,  $6S_{50}$,  $6S_{51}$,  $6S_{53}$,  $6S_{54}$  & $3S_{57}$ & $4S_{52}$,  $4S_{53}$, $4S_{54}$,  $4S_{55}$,  $4S_{56}$  &  $0S_{27}^{(2)}$ \\
				& &  $0S_{35}$, $3S_{51}$       & $6S_{55}$,  $6S_{56}$,  $6S_{57}$,  $6S_{58}$,  $6S_{59}$ &  $3S_{58}$ &  $8S_{32}$,  $8S_{33}$,  $8S_{34}$,  $8S_{35}$,  $8S_{36}$ & $0S_{28}^{(2)}$ \\
				& &  $3S_{56}$, $3S_{61}$       & $6S_{60}$,  $6S_{61}$,  $6S_{62}$,  $6S_{63}$,  $6S_{64}$  &   $3S_{59}$ &  $8S_{37}$,  $8S_{38}$,  $8S_{39}$,  $8S_{40}$,  $8S_{41}$  & $0S_{32}^{(2)}$ \\
				& &  $3S_{68}$,  $3S_{69}$      &  $6S_{65}$,  $6S_{66}$ & $3S_{60}$ & $8S_{42}$,  $8S_{43}$,  $8S_{44}$,  $8S_{45}$,  $8S_{46}$   & $0S_{34}^{(2)}$ \\
				& &  $3S_{70}$, $8S_{61}$      &  & $3S_{65}$ &  $8S_{47}$,  $8S_{48}$,  $8S_{49}$,  $8S_{50}$,  $8S_{51}$ &\\
				& &  $8S_{63}$,  $8S_{64}$       &  & $3S_{66}$ & $8S_{52}$,  $8S_{53}$,  $8S_{54}$,  $8S_{55}$,  $8S_{56}$  &\\
				& &  $8S_{68}$,  $8S_{69}$      &  &$3S_{67}$ &  $8S_{57}$,  $8S_{58}$,  $8S_{59}$,  $8S_{60}$, $8S_{62}$ &\\
				& &  $8S_{71}$        &  & & $8S_{65}$,  $8S_{66}$,  $8S_{67}$, $8S_{70}$, $8S_{72}$  &\\
				& &         &  & &  $8S_{73}$,  $8S_{74}$ &\\
				& &  $0.3L_{3}$       &  $0.6L_{3}$,  $0.6L_{4}$,  $0.6L_{7}$,  $0.6L_{10}$ &  $3.8L_{14}$ & $0.4L_{4}$,  $0.4L_{8}$,  $0.4L_{9}$  &\\
				& &    $0.3L_{6}$    & $0.6L_{11}$, $0.6L_{12}$,  $3.6L_{15}$,  $3.6L_{16}$   & $3.8L_{15}$ & $4.8L_{9}$, $4.8L_{10}$,  $4.8L_{11}$  &\\
				& &   $0.3L_{7}$       &  $3.6L_{17}$,  $3.6L_{18}$,  $3.6L_{19}$,  $3.6L_{20}$ & $3.8L_{16}$ &  $4.8L_{12}$, $4.8L_{13}$, $4.8L_{14}$   &\\
				& &    $3.8L_{20}$     & $4.6L_{12}$,  $4.6L_{14}$,  $4.6L_{15}$,  $4.6L_{17}$  & $3.8L_{17}$ & $4.8L_{15}$, $4.8L_{16}$,  $4.8L_{17}$  &\\
				& &    $3.8L_{21}$     &  $4.6L_{18}$,  $4.6L_{19}$, $6.8L_{11}$,  $6.8L_{12}$ & $3.8L_{18}$ & $4.8L_{18}$, $4.8L_{19}$,  $4.8L_{20}$  &\\
				& &   $3.8L_{22}$      & $6.8L_{13}$,  $6.8L_{14}$,  $6.8L_{15}$,  $6.8L_{16}$ & $3.8L_{19}$  &  $4.8L_{21}$, $8.8L_{1}$,  $8.8L_{2}$   &\\
				& &         &$6.8L_{17}$,  $6.8L_{18}$,  $6.8L_{19}$,  $6.8L_{20}$  & &   $8.8L_{3}$ &\\
				& &         &  $6.8L_{21}$,  $6.8L_{22}$,  $6.8L_{23}$,  $6.8L_{24}$ & &   &\\
				& &         &  $6.8L_{25}$,  $6.8L_{26}$,  $6.8L_{27}$ & &   &\\
				& &         & $P_{15}$,  $P_{16}$,  $P_{17}$,  $P_{18}$, $P_{19}$,  $P_{20}$ & &  $P_{12}$,  $P_{14}$ &\\
				& &         &  $P_{30}$,  $P_{33}$ & &   &\\
				\hline
			\end{tabular}
		\end{center}
	\end{table}
	
	\begin{table}\caption{\small Topological equivalences for the family $\QESA$ \textit{(cont.)}}\label{tab:top-equiv-QESA-pt4}
		\begin{center}
			\begin{tabular}{ccccccc}
				\hline
				Presented & Identical       & Finite      & Finite       & Finite  &  Possessing       &              \\
				phase     & under           & antisaddle  & antisaddle   & weak    &  invariant curve   &  Other reasons\\
				portrait  & perturbations   & focus       & node--focus  & point   &           (no separatrix)        &              \\
				\hline
				\multirow{1}{*}{$V_{168}$} &  $V_{178}$ &  $V_{169}$, $V_{179}$  & & & & $V_{200}^{(1)}$,  $V_{201}^{(1)}$,  $V_{202}^{(1)}$,  $V_{203}^{(1)}$ \\  
				& &   $3S_{71}$,  $3S_{76}$      & $6S_{67}$,  $6S_{74}$ &  &   & $3S_{85}^{(1)}$,  $3S_{86}^{(1)}$, $6S_{80}^{(1)}$,  $6S_{81}^{(1)}$ \\
				\hline
				\multirow{1}{*}{$V_{170}$} &  $V_{180}$  &   & & & &  $V_{199}^{(1)}$,  $V_{204}^{(1)}$  \\  
				\hline
				\multirow{2}{*}{$V_{173}$} &  $V_{183}$  &  $V_{171}$,  $V_{172}$,  $V_{181}$ & & & &  $V_{198}^{(1)}$,  $V_{205}^{(1)}$,  $V_{206}^{(1)}$,  $V_{207}^{(1)}$ \\  
				& &   $V_{182}$      &  &  &   &\\
				& &     $3S_{72}$,  $3S_{77}$  &  $6S_{68}$,  $6S_{75}$  &  &   & $3S_{84}^{(1)}$, $6S_{89}^{(1)}$  \\
				\hline
				\multirow{2}{*}{$V_{176}$} &  $V_{177}$, $V_{185}$, $V_{226}$  & $V_{174}$,  $V_{175}$, $V_{184}$  & & & & $V_{197}^{(1)}$,  $V_{208}^{(1)}$,  $V_{209}^{(1)}$,  $V_{210}^{(1)}$   \\  
				& $V_{228}$, $V_{229}$, $V_{231}$  & $V_{186}$, $V_{187}$     &   &  &   &  $V_{211}^{(1)}$, $V_{227}^{(1)}$, $V_{230}^{(1)}$, $V_{232}^{(1)}$  \\
				& &  $3S_{73}$,  $3S_{78}$       & $6S_{69}$,  $6S_{76}$, $6S_{90}$ &  $3S_{87}$, $3S_{89}$ & $8S_{78}$,  $8S_{84}$   & $3S_{83}^{(1)}$, $3S_{88}^{(1)}$, $6S_{88}^{(1)}$, $6S_{91}^{(1)}$ \\
				& &         &  $6S_{92}$,  $6S_{93}$, $6S_{95}$  &  & $8S_{94}$, $8S_{96}$  & $6S_{94}^{(1)}$, $6S_{96}^{(1)}$, $8S_{89}^{(1)}$, $8S_{95}^{(1)}$   \\
				& &         &  $3.6L_{24}$, $3.6L_{26}$  &  &   & $3.6L_{25}^{(1)}$, $6.8L_{35}^{(1)}$   \\
				& &         &  $6.8L_{34}$, $6.8L_{36}$  &  &   &\\
				\hline
				\multirow{5}{*}{$V_{188}$} & $V_{189}$, $V_{191}$, $V_{193}$   & $V_{190}$, $V_{192}$  & & & &    \\  
				& $V_{194}$,  $V_{195}$, $V_{213}$ &  $V_{196}$,  $V_{212}$     &   &  &   &   \\
				&$V_{214}$,  $V_{215}$, $V_{217}$ &    $V_{216}$, $V_{224}$     &  &  &   &\\
				& $V_{218}$,  $V_{219}$,  $V_{220}$ &  $V_{225}$       &  &  &   &\\
				& $V_{221}$,  $V_{222}$,  $V_{223}$ &         &  &  &   &\\
				& &  $3S_{74}$, $3S_{79}$        & $6S_{70}$,  $6S_{71}$,  $6S_{72}$, $6S_{73}$  & $3S_{75}$, $3S_{80}$ &  $4S_{61}$,  $4S_{62}$,  $4S_{67}$  &\\
				& &    $3S_{82}$     & $6S_{77}$,   $6S_{78}$, $6S_{79}$,  $6S_{82}$ &  $3S_{81}$ &   $4S_{68}$,  $4S_{69}$,  $4S_{70}$ &\\
				& &         &  $6S_{83}$, $6S_{84}$,  $6S_{85}$,  $6S_{86}$ &  & $8S_{79}$,  $8S_{80}$,  $8S_{85}$  &\\
				& &         &  $6S_{87}$ &  & $8S_{86}$,  $8S_{87}$,  $8S_{88}$  &\\
				& &         & $3.6L_{21}$,  $3.6L_{22}$,  $3.6L_{23}$  &  &  $4.8L_{23}$,  $4.8L_{26}$  &\\
				& &         & $4.6L_{21}$,  $4.6L_{23}$,  $4.6L_{24}$  &  &   &\\
				& &         & $6.8L_{29}$,  $6.8L_{31}$,  $6.8L_{32}$ &  &   &\\
				\hline
				\multirow{1}{*}{$V_{233}$} & $V_{252}$ &   $V_{234}$, $V_{253}$        & & & & $V_{273}^{(1)}$,  $V_{274}^{(1)}$,  $V_{275}^{(1)}$,  $V_{276}^{(1)}$ \\  
				& &    $3S_{90}$,  $3S_{95}$     & $6S_{97}$, $6S_{103}$ & & &$3S_{104}^{(1)}$,  $3S_{105}^{(1)}$, $6S_{108}^{(1)}$,  $6S_{109}^{(1)}$  \\
				\hline
			\end{tabular}
		\end{center}
	\end{table}
	
	\begin{table}\caption{\small Topological equivalences for the family $\QESA$ \textit{(cont.)}}\label{tab:top-equiv-QESA-pt5}
		\begin{center}
			\begin{tabular}{ccccccc}
				\hline
				Presented & Identical       & Finite      & Finite       & Finite  &  Possessing       &              \\
				phase     & under           & antisaddle  & antisaddle   & weak    &  invariant curve   &  Other reasons\\
				portrait  & perturbations   & focus       & node--focus  & point   &           (no separatrix)        &              \\
				\hline
				$V_{235}$ &  $V_{254}$ &&            & & & $V_{272}^{(1)}$,  $V_{277}^{(1)}$   \\  
				\hline
				\multirow{2}{*}{ $V_{238}$} &  $V_{257}$ &     $V_{236}$,  $V_{237}$      & & & &  $V_{271}^{(1)}$,  $V_{278}^{(1)}$ \\  
				& &    $V_{255}$,  $V_{256}$    &  & & & $V_{279}^{(1)}$,  $V_{280}^{(1)}$  \\
				& &    $3S_{91}$,  $3S_{96}$     &  $6S_{98}$,  $6S_{104}$ & & &  $3S_{103}^{(1)}$, $6S_{116}^{(1)}$  \\
				\hline
				\multirow{7}{*}{$V_{240}$} & $V_{241}$,  $V_{242}$, $V_{244}$,  $V_{245}$ &  $V_{239}$, $V_{243}$       & & & &    \\  
				&$V_{246}$,  $V_{247}$,  $V_{248}$,  $V_{249}$ &  $V_{260}$, $V_{262}$     &  & & &    \\
				& $V_{250}$,  $V_{251}$,  $V_{258}$,  $V_{259}$ &   $V_{270}$, $V_{282}$    &  & & &    \\
				& $V_{261}$, $V_{263}$,  $V_{264}$,  $V_{265}$ &  $V_{283}$     &  & & &    \\
				&$V_{266}$,  $V_{267}$,  $V_{268}$,  $V_{269}$  &       &  & & &    \\
				& $V_{281}$, $V_{284}$,  $V_{285}$,  $V_{286}$  &       &  & & &    \\
				& $V_{287}$,  $V_{288}$ &       &  & & &    \\
				& &   $3S_{92}$, $3S_{97}$    & $6S_{99}$,  $6S_{100}$,  $6S_{101}$,  $6S_{102}$  & $3S_{93}$,  $3S_{94}$ &$4S_{83}$,  $4S_{84}$,  $4S_{85}$,  $4S_{86}$ &    \\
				& &       & $6S_{105}$,  $6S_{106}$,  $6S_{107}$,  $6S_{110}$  &$3S_{98}$,  $3S_{99}$  & $4S_{87}$,  $4S_{88}$,  $4S_{89}$,  $4S_{90}$ &    \\
				& &       & $6S_{111}$,  $6S_{112}$,  $6S_{113}$,  $6S_{114}$  &$3S_{100}$,  $3S_{101}$ & $4S_{91}$,  $4S_{92}$,  $4S_{93}$,  $4S_{94}$ &    \\
				& &       & $6S_{115}$,  $6S_{117}$  & $3S_{102}$ & $4S_{95}$,  $4S_{96}$,  $4S_{97}$,  $4S_{98}$ &    \\
				& &       &  & &$4S_{99}$, $8S_{100}$,  $8S_{101}$,  $8S_{102}$  &    \\
				& &       &  & &$8S_{106}$,  $8S_{107}$,  $8S_{108}$,  $8S_{109}$  &    \\
				& &       &  & & $8S_{110}$,  $8S_{111}$ &    \\
				& &       & $3.6L_{27}$,  $3.6L_{28}$,  $3.6L_{29}$  &$3.4L_{23}$,  $3.4L_{24}$  & $4.8L_{30}$,  $4.8L_{31}$,  $4.8L_{32}$ &    \\
				& &       & $4.6L_{33}$, $4.6L_{34}$,  $4.6L_{35}$   & $3.4L_{25}$ & $4.8L_{33}$, $4.8L_{34}$  &    \\
				& &       & $4.6L_{36}$, $6.8L_{37}$,  $6.8L_{40}$  & & &    \\
				& &       & $6.8L_{41}$  & & &    \\
				& &       &  & & &    \\
				\hline
			\end{tabular}
		\end{center}
	\end{table}
	
	\begin{table}\caption{\small Topological equivalences for the family $\QESA$ \textit{(cont.)}}\label{tab:top-equiv-QESA-pt6}
		\begin{center}
			\begin{tabular}{ccccccc}
				\hline
				Presented & Identical       & Finite      & Finite       & Finite  &  Possessing       &              \\
				phase     & under           & antisaddle  & antisaddle   & weak    &  invariant curve   &  Other reasons\\
				portrait  & perturbations   & focus       & node--focus  & point   &           (no separatrix)        &              \\
				\hline
				\multirow{1}{*}{$2S_{1}$} & $2S_{2}$,  $2S_{3}$ & &          & & &\\  
				& & &         &$2.3L_{1}$ &$2.8L_{1}$ &\\  
				\hline
				\multirow{1}{*}{$2S_{4}$} & $2S_{8}$ & &          & & &\\  
				& & &         & $2.3L_{3}$ & &\\  
				\hline
				\multirow{1}{*}{$2S_{5}$} & $2S_{9}$  & &          & & &\\  
				& & &         & $2.3L_{5}$ & &\\  
				\hline
				\multirow{1}{*}{$2S_{6}$} &$2S_{7}$,  $2S_{10}$  & &          & & &\\  
				& & &         &$2.3L_{6}$ & $2.4L_{2}$ &\\  
				\hline
				\multirow{1}{*}{$2S_{11}$} &  & &          & & &\\  
				& & &         & & & $0.2L_{1}^{(2)}$ \\  
				\hline
				\multirow{1}{*}{$2S_{12}$} &  & &          & & &\\  
				& & &         & & & $0.2L_{2}^{(2)}$ \\  
				\hline
				\multirow{1}{*}{ $2S_{13}$} &  $2S_{14}$,  $2S_{15}$, $2S_{21}$ & $2S_{16}$, $2S_{22}$  &          & & & \\  
				& &   $0.2L_{4}$, $2.8L_{7}$    &  $2.6L_{1}$,  $2.6L_{3}$ & &$2.8L_{4}$,  $2.8L_{5}$,  $2.8L_{6}$ & $0.2L_{3}^{(2)}$\\  
				& &      &   $P_{23}$ & &$P_{13}$ &  \\  
				\hline
				\multirow{1}{*}{$2S_{17}$} &  & &          & & &\\  
				& & &         & & & $0.2L_{5}^{(2)}$\\  
				\hline
				\multirow{1}{*}{$2S_{18}$} &  & &          & & &\\  
				& & &         & & & $0.2L_{6}^{(2)}$\\  
				\hline
				\multirow{1}{*}{$2S_{20}$} & $2S_{19}$ & &          & & &\\  
				& & $0.2L_{7}$, $2.3L_{8}$ &   $2.6L_{2}$      & & & $0.2L_{8}^{(2)}$\\  
				& & $P_{26}$ &   $P_{27}$     & & &  \\  
				\hline
				\multirow{1}{*}{$2S_{23}$} &  & &          & & &\\  
				\hline
				\multirow{1}{*}{$2S_{24}$} &  & &          & & &\\  
				\hline
				\multirow{1}{*}{$2S_{25}$} &  & &          & & &\\  
				\hline
				\multirow{1}{*}{$2S_{26}$} &  &$2S_{27}$ &          & & &\\  
				& & &     $2.6L_{4}$    & & &\\  
				\hline
			\end{tabular}
		\end{center}
	\end{table}
	
	\begin{table}\caption{\small Topological equivalences for the family $\QESA$ \textit{(cont.)}}\label{tab:top-equiv-QESA-pt7}
		\begin{center}
			\begin{tabular}{ccccccc}
				\hline
				Presented & Identical       & Finite      & Finite       & Finite  &  Possessing       &              \\
				phase     & under           & antisaddle  & antisaddle   & weak    &  invariant curve   &  Other reasons\\
				portrait  & perturbations   & focus       & node--focus  & point   &           (no separatrix)        &              \\
				\hline
				\multirow{1}{*}{$2S_{28}$} &  & &          & & &\\  
				\hline
				\multirow{1}{*}{$2S_{29}$} &  & &          & & &\\  
				\hline
				\multirow{1}{*}{$2S_{30}$} &  & &          & & &\\  
				\hline
				\multirow{1}{*}{$2S_{32}$} &  &$2S_{31}$ &          & & &\\  
				& & $2.3L_{10}$ &   $2.6L_{5}$      & & &\\  
				\hline
				\multirow{1}{*}{$2S_{33}$} &  & &          & & &\\  
				\hline
				\multirow{1}{*}{$2S_{34}$} &  & &          & & &\\  
				\hline
				\multirow{1}{*}{$2S_{35}$} & $2S_{36}$,  $2S_{37}$  &  $2S_{38}$ &          & & &\\  
				& & &     $2.6L_{6}$    & &$2.4L_{8}$,  $2.4L_{9}$ &\\  
				\hline
				\multirow{1}{*}{$2S_{39}$} &  & &          & & &\\  
				\hline
				\multirow{1}{*}{$2S_{40}$} &  & &          & & &\\  
				\hline
				\multirow{1}{*}{$2S_{42}$} &  & $2S_{41}$ &          & & &\\  
				& & $2.3L_{12}$ &    $2.6L_{7}$      & & &\\  
				\hline
				\multirow{2}{*}{$4S_{5}$} & $4S_{8}$, $4S_{21}$  & $4S_{6}$, $4S_{9}$,  $4S_{25}$ &          & & & $4S_{10}^{(1)}$,  $4S_{11}^{(1)}$, $4S_{24}^{(1)}$, $4S_{27}^{(1)}$ \\  
				&$4S_{23}$ & $4S_{26}$, $4S_{29}$, $4S_{31}$&         & & & $4S_{30}^{(1)}$ \\  
				& &  $3.4L_{8}$,  $3.4L_{9}$ &   $4.6L_{2}$, $4.6L_{6}$      & $3.4L_{5}$, $3.4L_{6}$ & & $3.4L_{7}^{(1)}$, $3.4L_{10}^{(1)}$, $3.4L_{13}^{(1)}$ \\  
				& &$3.4L_{11}$,  $3.4L_{12}$  &    $4.6L_{8}$,  $4.6L_{9}$     & & & $4.6L_{4}^{(1)}$, $4.6L_{10}^{(1)}$ \\  
				& & &    $P_{2}$,  $P_{3}$     & & &  $P_{4}^{(1)}$ \\  
				\hline
				\multirow{1}{*}{$4S_{34}$} & $4S_{40}$  &  $4S_{32}$,  $4S_{33}$,  $4S_{38}$, $4S_{39}$ &          & & & $4S_{47}^{(1)}$,  $4S_{48}^{(1)}$,  $4S_{49}^{(1)}$,  $4S_{50}^{(1)}$ \\  
				& &$0.4L_{1}$,  $0.4L_{2}$, $0.4L_{5}$  &    $4.6L_{11}$     & & & $0.4L_{3}^{(2)}$, $0.4L_{7}^{(2)}$, $0.4L_{10}^{(1)}$ \\  
				& & $0.4L_{6}$, $3.4L_{14}$, $3.4L_{16}$  &         & & & $0.4L_{11}^{(1)}$, $0.4L_{12}^{(1)}$,  $0.4L_{13}^{(1)}$ \\  
				& &$4.6L_{13}$  &         & & & $3.4L_{15}^{(1)}$, $4.6L_{16}^{(1)}$ \\  
				& & $P_{28}$, $P_{31}$  &     $P_{29}$, $P_{32}$    & & & $P_{34}^{(1)}$,  $P_{35}^{(1)}$ \\  
				\hline
			\end{tabular}
		\end{center}
	\end{table}
	
	\begin{table}\caption{\small Topological equivalences for the family $\QESA$ \textit{(cont.)}}\label{tab:top-equiv-QESA-pt8}
		\begin{center}
			\begin{tabular}{ccccccc}
				\hline
				Presented & Identical       & Finite      & Finite       & Finite  &  Possessing       &              \\
				phase     & under           & antisaddle  & antisaddle   & weak    &  invariant curve   &  Other reasons\\
				portrait  & perturbations   & focus       & node--focus  & point   &           (no separatrix)        &              \\
				\hline
				\multirow{2}{*}{ $4S_{59}$} & $4S_{60}$, $4S_{65}$,  $4S_{66}$, $4S_{76}$  & $4S_{57}$,  $4S_{58}$ & & & & $4S_{71}^{(1)}$,  $4S_{72}^{(1)}$,  $4S_{73}^{(1)}$,  $4S_{74}^{(1)}$ \\  
				& $4S_{78}$,  $4S_{79}$, $4S_{81}$ &$4S_{63}$,  $4S_{64}$  & & & &  $4S_{75}^{(1)}$, $4S_{77}^{(1)}$, $4S_{80}^{(1)}$, $4S_{82}^{(1)}$ \\  
				& & $3.4L_{17}$, $3.4L_{19}$ &$4.6L_{20}$,  $4.6L_{22}$  & $3.4L_{20}$, $3.4L_{22}$  &  $4.8L_{22}$, $4.8L_{25}$  & $3.4L_{18}^{(1)}$, $3.4L_{21}^{(1)}$, $4.6L_{25}^{(1)}$  \\  
				& & & $4.6L_{26}$, $4.6L_{28}$ & & $4.8L_{27}$, $4.8L_{29}$  & $4.6L_{27}^{(1)}$, $4.6L_{30}^{(1)}$, $4.6L_{32}^{(1)}$ \\  
				& & & $4.6L_{29}$, $4.6L_{31}$ & & & $4.8L_{24}^{(1)}$, $4.8L_{28}^{(1)}$  \\  
				& & & $P_{36}$, $P_{38}$,  $P_{39}$, $P_{41}$ & & &  $P_{37}^{(1)}$, $P_{40}^{(1)}$ \\  																								
				\hline
				\multirow{1}{*}{$5S_{1}$} &$5S_{13}$ & $5S_{2}$,  $5S_{14}$ & & & &  $5S_{29}^{(1)}$,  $5S_{30}^{(1)}$,  $5S_{31}^{(1)}$,  $5S_{32}^{(1)}$ \\  
				& & $3.5L_{1}$,  $3.5L_{5}$ & $5.6L_{1}$,  $5.6L_{6}$ & & &  $3.5L_{12}^{(1)}$,  $3.5L_{13}^{(1)}$, $5.6L_{11}^{(1)}$ \\  
				& & & & & & $5.6L_{12}^{(1)}$ \\  
				\hline
				\multirow{1}{*}{$5S_{3}$} & $5S_{15}$ & & & & &  $5S_{28}^{(1)}$,  $5S_{33}^{(1)}$ \\  
				& & & & & &  \\  
				\hline
				\multirow{2}{*}{$5S_{6}$} & $5S_{18}$ & $5S_{4}$,  $5S_{5}$ & & & & $5S_{27}^{(1)}$,  $5S_{34}^{(1)}$,  $5S_{35}^{(1)}$,  $5S_{36}^{(1)}$  \\  
				& &  $5S_{16}$,  $5S_{17}$ & & & &  \\  
				& & $3.5L_{2}$,  $3.5L_{6}$ & & & & $3.5L_{11}^{(1)}$, $5.6L_{17}^{(1)}$ \\  
				& & $5.6L_{2}$,  $5.6L_{7}$ & & & &  \\  
				& & & & & &  \\  
				\hline
				\multirow{2}{*}{$5S_{9}$} &$5S_{10}$,  $5S_{11}$,  $5S_{12}$,  $5S_{19}$ & $5S_{7}$,  $5S_{8}$ & & & &  $5S_{26}^{(1)}$,  $5S_{37}^{(1)}$,  $5S_{38}^{(1)}$,  $5S_{39}^{(1)}$ \\  
				&$5S_{20}$, $5S_{22}$, $5S_{24}$,  $5S_{25}$ & $5S_{21}$, $5S_{23}$ & & & &  $5S_{40}^{(1)}$,  $5S_{41}^{(1)}$,  $5S_{42}^{(1)}$ \\  
				& & $3.5L_{3}$, $3.5L_{7}$  & $5.6L_{3}$,  $5.6L_{4}$ &$3.5L_{4}$, $3.5L_{8}$ & $5.8L_{4}$,  $5.8L_{5}$  & $3.5L_{9}^{(1)}$,  $3.5L_{10}^{(1)}$, $5.6L_{13}^{(1)}$ \\  
				& & & $5.6L_{5}$, $5.6L_{8}$ & & $5.8L_{9}$,  $5.8L_{10}$ & $5.6L_{14}^{(1)}$,  $5.6L_{15}^{(1)}$,  $5.6L_{16}^{(1)}$ \\  
				& & & $5.6L_{9}$,  $5.6L_{10}$  & & & $5.8L_{11}^{(1)}$,  $5.8L_{12}^{(1)}$  \\  
				& & & $P_{51}$,  $P_{52}$ & $P_{56}$,  $P_{57}$ & &  $P_{59}^{(1)}$,  $P_{60}^{(1)}$ \\  
				\hline
				\multirow{1}{*}{$7S_{1}$} & $7S_{2}$ & & & & & $7S_{3}^{(1)}$ \\  
				\hline
				\multirow{1}{*}{$7S_{4}$} &$7S_{5}$ & & & & & $7S_{6}^{(1)}$ \\  
				\hline
			\end{tabular}
		\end{center}
	\end{table}
	
	\begin{table}\caption{\small Topological equivalences for the family $\QESA$ \textit{(cont.)}}\label{tab:top-equiv-QESA-pt9}
		\begin{center}
			\begin{tabular}{ccccccc}
				\hline
				Presented & Identical       & Finite      & Finite       & Finite  &  Possessing       &              \\
				phase     & under           & antisaddle  & antisaddle   & weak    &  invariant curve   &  Other reasons\\
				portrait  & perturbations   & focus       & node--focus  & point   &           (no separatrix)        &              \\
				\hline
				\multirow{1}{*}{$7S_{7}$} & $7S_{8}$ & & & & & $7S_{9}^{(1)}$,  $7S_{10}^{(1)}$ \\  
				& & & & & & $0.7L_{1}^{(2)}$,  $0.7L_{2}^{(2)}$ \\
				& & & & & & $0.7L_{3}^{(1)}$,  $0.7L_{4}^{(1)}$  \\
				\hline
				\multirow{1}{*}{$7S_{11}$} &$7S_{12}$ & & & & & $7S_{13}^{(1)}$,  $7S_{14}^{(1)}$ \\  
				\hline
				\multirow{1}{*}{$7S_{15}$} &$7S_{16}$ & & & & & $7S_{17}^{(1)}$,  $7S_{18}^{(1)}$ \\  
				\hline
				\multirow{2}{*}{$8S_{7}$} & $8S_{14}$,  $8S_{24}$  & $8S_{5}$,  $8S_{6}$, $8S_{15}$  & & & &  $8S_{17}^{(1)}$,  $8S_{18}^{(1)}$,  $8S_{19}^{(1)}$ \\  
				&$8S_{25}$  & $8S_{16}$, $8S_{21}$,  $8S_{22}$ & & & & $8S_{23}^{(1)}$, $8S_{26}^{(1)}$ \\
				& & $3.8L_{1}$,  $3.8L_{2}$ &$6.8L_{1}$,  $6.8L_{4}$  & $3.8L_{7}$,  $3.8L_{8}$ & &  $3.8L_{3}^{(1)}$, $3.8L_{6}^{(1)}$, $3.8L_{9}^{(1)}$  \\
				& &$3.8L_{4}$,  $3.8L_{5}$ & $6.8L_{7}$,  $6.8L_{8}$ & & &  $6.8L_{6}^{(1)}$, $6.8L_{9}^{(1)}$\\
				& & &$P_{8}$,  $P_{9}$ & & &$P_{10}^{(1)}$  \\
				\hline
				\multirow{2}{*}{$8S_{77}$} &$8S_{83}$ & $8S_{75}$,  $8S_{76}$ & & & & $8S_{90}^{(1)}$,  $8S_{91}^{(1)}$ \\  
				& &  $8S_{81}$,  $8S_{82}$ & & & & $8S_{92}^{(1)}$,  $8S_{93}^{(1)}$\\
				& & $3.8L_{23}$, $3.8L_{25}$ &$6.8L_{28}$  & & & $3.8L_{24}^{(1)}$, $6.8L_{33}^{(1)}$ \\
				& & $6.8L_{30}$ & & & &  \\
				\hline
				\multirow{2}{*}{$8S_{99}$} &$8S_{105}$  &$8S_{97}$,  $8S_{98}$ & & & & $8S_{103}^{(1)}$, $8S_{112}^{(1)}$, $8S_{113}^{(1)}$ \\  
				& &$8S_{104}$ & & & & $8S_{114}^{(1)}$, $8S_{115}^{(1)}$ \\
				& &$3.8L_{26}$,  $3.8L_{27}$ & $6.8L_{38}$,  $6.8L_{39}$ & & & $3.8L_{28}^{(1)}$, $6.8L_{42}^{(1)}$ \\
				\hline
				\multirow{1}{*}{$2.3L_{2}$} &  $2.3L_{4}$ & & & & &\\  
				& & & & $P_{6}$& &\\
				\hline
				\multirow{1}{*}{$2.3L_{7}$} & & & & & &\\  
				& & &$P_{21}$ & & &\\
				\hline
				\multirow{1}{*}{$2.3L_{9}$} & & & & & &\\  
				\hline
				\multirow{1}{*}{$2.3L_{11}$} & & & & & &\\  
				\hline
				\multirow{1}{*}{$2.4L_{1}$} & $2.4L_{3}$& & & & &\\  
				& & & & $P_{1}$& &\\
				\hline
				\multirow{1}{*}{$2.4L_{4}$} & & & & & &\\  
				& & & & &$P_{22}$ &\\
				\hline
			\end{tabular}
		\end{center}
	\end{table}
	
	\begin{table}\caption{\small Topological equivalences for the family $\QESA$ \textit{(cont.)}}\label{tab:top-equiv-QESA-pt10}
		\begin{center}
			\begin{tabular}{ccccccc}
				\hline
				Presented & Identical       & Finite      & Finite       & Finite  &  Possessing       &              \\
				phase     & under           & antisaddle  & antisaddle   & weak    &  invariant curve   &  Other reasons\\
				portrait  & perturbations   & focus       & node--focus  & point   &           (no separatrix)        &              \\
				\hline
				\multirow{1}{*}{$2.4L_{5}$} & & & & & &\\  
				& & & & &$P_{24}$ &\\
				\hline
				\multirow{1}{*}{$2.4L_{6}$} & & & & & &\\  
				\hline
				\multirow{1}{*}{$2.4L_{7}$} & & & & & &\\  
				\hline
				\multirow{1}{*}{$2.5L_{1}$} & & & & & &\\  
				\hline
				\multirow{1}{*}{$2.5L_{2}$} & & & & & &\\  
				\hline
				\multirow{1}{*}{$2.5L_{3}$} & & & & & &\\  
				\hline
				\multirow{1}{*}{$2.5L_{4}$} & & & & & &\\  
				\hline
				\multirow{1}{*}{$2.5L_{5}$} & & & & & &\\  
				\hline
				\multirow{1}{*}{$2.5L_{6}$} & & & & & &\\  
				\hline
				\multirow{1}{*}{$2.5L_{8}$} & &$2.5L_{7}$ & & & &\\  
				& &$P_{47}$ &$P_{48}$ & & &\\
				\hline
				\multirow{1}{*}{$2.7L_{1}$} & & & & & &\\  
				& & $P_{25}$ & & & &\\
				\hline
				\multirow{1}{*}{$2.7L_{2}$} & & & & & &\\  
				\hline
				\multirow{1}{*}{$2.7L_{3}$} & & & & & &\\  
				\hline
				\multirow{1}{*}{$2.8L_{2}$} &$2.8L_{3}$ & & & & &\\  
				& & & & &$P_{11}$  &\\
				\hline
				\multirow{1}{*}{$2.8L_{8}$} & & & & & &\\  
				\hline
				\multirow{1}{*}{$2.8L_{9}$} & & & & & &\\  
				\hline
				\multirow{1}{*}{$2.8L_{10}$} & & & & & &\\  
				\hline
				\multirow{1}{*}{$2.8L_{11}$} & & & & & &\\  
				\hline
				\multirow{1}{*}{$3.7L_{1}$} &$3.7L_{2}$ & & & & &  $3.7L_{3}^{(1)}$\\  
				\hline
				\multirow{1}{*}{$4.5L_{1}$} & $4.5L_{2}$,  $4.5L_{3}$ & & & & &\\  
				& & & &  $P_{53}$,  $P_{58}$ & &\\
				\hline
				\multirow{1}{*}{$4.8L_{2}$} &$4.8L_{5}$,  $4.8L_{7}$,  $4.8L_{8}$ & & & & &\\  
				& & & & $P_{5}$,  $P_{7}$ & &\\
				\hline
				\multirow{1}{*}{$5.7L_{1}$} &$5.7L_{2}$  & & & & &  $5.7L_{3}^{(1)}$,  $5.7L_{4}^{(1)}$ \\ 
				\hline
			\end{tabular}
		\end{center}
	\end{table}
	
	\begin{table}\caption{\small Topological equivalences for the family $\QESA$ \textit{(cont.)}}\label{tab:top-equiv-QESA-pt11}
		\begin{center}
			\begin{tabular}{ccccccc}
				\hline
				Presented & Identical       & Finite      & Finite       & Finite  &  Possessing       &              \\
				phase     & under           & antisaddle  & antisaddle   & weak    &  invariant curve   &  Other reasons\\
				portrait  & perturbations   & focus       & node--focus  & point   &           (no separatrix)        &              \\
				\hline
				\multirow{2}{*}{$5.8L_{3}$} & $5.8L_{8}$ & $5.8L_{1}$,  $5.8L_{2}$& & & &  $5.8L_{13}^{(1)}$,  $5.8L_{14}^{(1)}$ \\ 
				& & $5.8L_{6}$,  $5.8L_{7}$ & & & & $5.8L_{15}^{(1)}$,  $5.8L_{16}^{(1)}$ \\ 
				& & $P_{49}$, $P_{54}$ & $P_{50}$, $P_{55}$ & & &$P_{61}^{(1)}$,  $P_{62}^{(1)}$  \\ 
				\hline
				\multirow{1}{*}{$P_{42}$} & & & & & &  \\ 
				\hline
				\multirow{1}{*}{$P_{43}$} & & & & & &  \\ 
				\hline
				\multirow{1}{*}{$P_{44}$} & & & & & &  \\ 
				\hline
				\multirow{1}{*}{$P_{45}$} & & & & & &  \\ 
				\hline
				\multirow{1}{*}{$P_{46}$} & & & & & &  \\ 
				\hline
			\end{tabular}
		\end{center}
	\end{table}
	
\end{landscape}

\medskip

\subsection{The bifurcation diagram of family $\QESB$}\label{subsec:bd-QES-B}

In this section we present the study of the bifurcation diagram of family $\QESB$, described by 
systems~\eqref{eq:nf-QES-B}. 

From normal form~\eqref{eq:nf-QES-B} we observe that the family under consideration depends on the
parameters $g\in\mathbb{R}\setminus\{0\}$ (in order to have nondegenerate systems), 
$u\in\mathbb{R}^+\cup\{0\}$ (due to the symmetry we proved before), and $\ell\in\mathbb{R}$.
Here we shall consider the bifurcation diagram formed by planes $g=g_0$ in which the 
Cartesian coordinates are $(u,\ell)$ with $u\geq0$.

For systems~\eqref{eq:nf-QES-B}, computations show that 
$$\mathbf{D}=12288g^6(1+u^2)^4, \quad \mathbf{R}=48g^4(1+u^2)^2x^2,$$
therefore by \cite[Table 5.1]{Artes-Llibre-Schlomiuk-Vulpe-2021a}, for $g\ne0$ systems~\eqref{eq:nf-QES-B} 
possess exactly one real simple finite singular point and two complex ones. 

\begin{remark} In order to avoid unnecessary repetitions, along this section we shall omit most of the explanations 
	similar to the ones already presented previously along the study of family~\eqref{eq:nf-QES-A}.
\end{remark}

Now we present the value of the algebraic invariants and T--comitants (with respect to systems~\eqref{eq:nf-QES-B})
which are relevant in our study. 

\medskip

\noindent \textbf{Bifurcation surface in $\mathbb{R}^3$ due to degeneracy of the system}

For family $\QESB$ we calculate
$$\mu_0=0 \quad \text{ and } \quad \mu_1=4g^2(1+u^2)x,$$
and it is clear that the comitant $\mu_1$ vanishes if and only if $g=0$. Moreover, computation show that
$$\mu_2\vert_{g=0}=\mu_3\vert_{g=0}=\mu_4\vert_{g=0}=0,$$
i.e., along the surface
$$({\cal S}_{1})\!:g=0,$$
in fact, a plane, we have degenerate systems.

\begin{remark}\label{study-deg-famB} Family $\QESB$ restricted to surface $({\cal S}_{1})$ is given by
	$$\begin{aligned}
	&x^\prime=0,\\
	&y^\prime=2(1-\ell u)x+\ell(1+u^2)y+\ell x^2-2xy,\\
	\end{aligned}$$
	and, as we mentioned before, this two--parametric family has curves filled up with singular points. 
	According to \cite[Diagram 12.1]{Artes-Llibre-Schlomiuk-Vulpe-2021a}, for these systems we calculate
	$$\eta=0, \quad \widetilde{M}=-32x^2,  \quad \kappa=\widetilde{K}=\widetilde{L}=\kappa_1={K}_1=0,$$
	and 
	$${L}_2=-6\ell(1+u^2)\left[4+\ell(\ell-4u+\ell u^2)\right]x^4.$$
	Since the discriminant of $4+\ell(\ell-4u+\ell u^2)$ is negative, we point 
	out that ${L}_2=0$ is equivalent to $\ell=0$. So, according to the mentioned reference, 
	for $\ell\neq0$ we have a hyperbola filled up with singular points,
	and for $\ell=0$ we have two real straight lines (filled up with singular points) intersecting at a finite point.
	Therefore, in the plane $g=0$ the straight line $\ell=0$ yields a bifurcation of curves filled up with singular points.
\end{remark}

\medskip

\noindent \textbf{The surface of $C^{\infty}$ bifurcation points due to weak singularities}

\medskip

\noindent {\bf (${\cal S}_{3}$)} This is the bifurcation surface due to weak finite singularities. 
According to \cite{Vulpe-2011}, for systems~\eqref{eq:nf-QES-B} we calculate
$$\begin{aligned}
\mathcal{T}_{4}=& \ \mathcal{T}_{3}=\mathcal{T}_{2}=\mathcal{T}_{1}=0,\\
\sigma=& \ \ell - 2 g u + \ell u^2 + 2 (g-1) x,
\end{aligned}$$
then due to the results on the mentioned paper, in the case in which $\sigma$ is generically nonzero,
the family under consideration could possess one and only one weak singularity. Moreover as
$$\begin{aligned}
\mathcal{F}_{1} =& \ 2 g^2 (1 + u^2) \left[2 (2 + g) u - 3 \ell (1 + u^2)\right],\qquad
\mathcal{H} =0,\\
\mathcal{B}_{1} =& \ 2 g^2 (1 + u^2) \left[2 g u - \ell (1 + u^2)\right]\left[4g(g-2) + (1 + u^2) \left(4 + \ell (\ell - 4 u + \ell u^2)\right)\right],\\
\mathcal{B}_{2}\!=& \ 2g^3 (g-1)^2  (1 + u^2)^2 \left[4 g^2 + (1 + u^2) \left(4 - 8 \ell u + 3 \ell^2 (1 + u^2)\right) \right.\\ 
&\left. -4 g \left(2 - 2 u^2 + \ell u(1 + u^2)\right)\right],
\end{aligned}$$
assuming $\mathcal{F}_{1}\ne0$, for family~\eqref{eq:nf-QES-B} we can obtain one weak 
singularity ($s^{(1)}$ or $f^{(1)}$, depending on the sign of $\mathcal{B}_{2}$) along the 
surface given by $\mathcal{B}_{1}=0$, i.e.
{$$\begin{aligned}
	({\cal S}_{3})\!\!:& \ 2 g^2 (1 + u^2) \left[2 g u - \ell (1 + u^2)\right] \left[4 g(g-2) + (1 + u^2) \left(4 + \ell (\ell - 4 u + \ell u^2)\right)\right]=0.
	\end{aligned}$$}

\begin{remark}\label{weak-singularities2} 
	\begin{enumerate}
		\item  We observe that, independently of $x$, we have $\sigma=0$ if and only if 
		$$\{g=1, \, \ell=2u/(1+u^2)\}.$$ 
		Under these conditions, we have that $\mu_{0}=0$, $\mathbf{D}=12288 \left(1+u^2\right)^4$, 
		and $\mathbf{R}=48(1+u^2)^2x^2$. So, according to \cite[item ($f_{6}$)--$\beta$]{Vulpe-2011} 
		we have one 	finite singular point, which is an integrable saddle. In other words, when $g=1$,
		during the study of the curve $\ell=2u/(1+u^2)$ we shall always obtain a phase portrait containing one 
		integrable saddle.
		\item	We just saw that in order to define surface {\bf (${\cal S}_{3}$)}
		we considered $\sigma\ne0$ and $\mathcal{F}_{1}\ne0$. However, according to \cite[item $(e)$]{Vulpe-2011}, 
		when $\sigma\ne0$ and $\mathcal{F}_{1}=0$ we can have either an integrable saddle or a center. As 
		we already have obtained conditions in order to have an integrable saddle, now we 
		analyze when we have a center. In fact, as we already have $\mathcal{H}=0$, from the mentioned paper
		we solve $\mathcal{F}_{1}=\mathcal{B}_{1}=0$  (together with $\sigma\ne0$ and $g\neq0$), and we 
		obtain the solution
		$$\{u=0, \, \ell=0\}.$$
		Also, when we compute $\mathcal{B}_{2}$ along this solution we obtain
		$8 (g-1)^4 g^3$, which is generically negative if $g<0$. Note that we must have $g\ne1$, because
		$\sigma|_{\{u=0, \, \ell=0, \, g=1\}}=0$.\\
		Therefore, from \cite[item ($e_4$)--$\beta$]{Vulpe-2011}, this study shows that for $g<0$ we 
		shall always find a center type singular point when we have $\{u=0, \, \ell=0\}$.
	\end{enumerate}
\end{remark}

\medskip

\noindent \textbf{Bifurcation surfaces in $\mathbb{R}^3$ due to the presence of invariant algebraic curves}

\medskip

\noindent {\bf (${\cal S}_{4}$)} This surface contains the points of the parameter space in which there appear 
invariant straight lines (see Lemma~\ref{lemma:S4-inv-curves-QES-B}). For systems~\eqref{eq:nf-QES-B} 
we compute the polynomial invariant $B_1$ and we define surface
$$\begin{aligned}
({\cal S}_{4})\!: & \ -8 g^6 \ell (1 + u^2)^5 \left[\ell^2 + (2 + g - \ell u)^2\right]=0.
\end{aligned}$$

\medskip

\noindent {\bf (${\cal S}_{8}$)} This surface contains the points of the parameter space in which there appear 
invariant parabolas. According to the conditions stated in Lemma~\ref{lemma:S4-inv-curves-QES-B} we define 
this surface by 
$$\begin{aligned}
({\cal S}_{8})\!: & \ \ell - 2 u - 2 g u + \ell u^2=0.
\end{aligned}$$

\medskip

\noindent \textbf{Bifurcation surface due to multiplicities of infinite singularities}

\medskip

\noindent {\bf (${\cal S}_{5}$)} This is the bifurcation surface due to multiplicity of infinite singular points.  
According to \cite[Lemma 5.5]{Artes-Llibre-Schlomiuk-Vulpe-2021a}, for this family we calculate
$$\eta=0, \quad \widetilde{M}=-8(2+g)^2x^2,   \quad   C_2=-x^2\left[\ell x-(2+g)y\right],$$
and we observe that along
$$\begin{aligned}
(\mathcal{ S}_{5})\!:& \ g+2=0,
\end{aligned}$$
we have a coalescence of infinite singular points. In addition, due to the mentioned result, on the plane 
$g=-2$ all the phase portraits corresponding to $\ell=0$ have the line at infinity filled up with singular points.

\medskip

\noindent \textbf{The surface of $C^{\infty}$ bifurcation due to a node becoming a focus}

\medskip

\noindent {\bf (${\cal S}_{6}$)} This surface contains the points of the parameter space where a finite node 
of the systems turns into a focus. According to \cite[Table 6.2]{Artes-Llibre-Schlomiuk-Vulpe-2021a} we calculate
$\mu_{0}=0, \mathbf{D}=12288 g^6 (1 + u^2)^4, \mathbf{R}=48 g^4 (1 + u^2)^2 x^2, \widetilde{K}=-4gx^2, G_{9}=0,$
and for the mentioned table we conclude that the invariant $W_{7}$ is responsible for describing the node--focus 
bifurcation. We compute this invariant polynomial and we define surface {\bf (${\cal S}_{6}$)} by the zero set of
$${
	\begin{aligned}
	12 g^6 &(1 + u^2)^4 \left[4 g^2 u^2 - 4 g (\ell u-2) (1 + u^2) + \ell^2 (1 + u^2)^2\right] \times\\
	\times&\left[16 g (1 + u^2) \left(4 + 4 \ell u - 3 \ell^2 (1 + u^2)\right) +64 g^3 + 16 g^4+(1 + u^2)^2 \left(4 + \ell (\ell - 4 u + \ell u^2)\right)^2 + \right.\\
	&\left.+8 g^2 \left(\ell^2 (1 + u^2)^2 + 4 (3 + u^2) + 12 \ell u(1 + u^2)\right)\right]=0.\\
	\end{aligned}}$$

\medskip

\noindent \textbf{Bifurcation surface in $\mathbb{R}^3$ due to the infinite elliptic--saddle}

\medskip

\noindent {\bf (${\cal S}_{0}$)} Along the plane $g=-1$ the corresponding phase portraits possess an infinite 
singularity of the type $\widehat{\!{1\choose 2}\!\!}\ E-H$. Due to results on \cite{Artes-Llibre-Schlomiuk-Vulpe-2021a} 
we compute the comitant
$$\widetilde{N}=-4 (g+1)x^2,$$
and we define surface 
$$({\cal S}_{0})\!:g+1=0.$$

\medskip

The bifurcation surfaces listed previously are all algebraic and they, except $({\cal S}_{4})$ and $({\cal S}_{8})$, 
are the bifurcation surfaces of singularities of systems~\eqref{eq:nf-QES-B} in the parameter space. We shall detect 
other bifurcation surface not necessarily algebraic in which the family has global connection of separatrices different from 
those given by $({\cal S}_{4})$ and $({\cal S}_{8})$. We shall name it surface $({\cal S}_{7})$.

As in the previous sections, here we shall foliate the three--dimensional bifurcation diagram in $\mathbb{R}^3$ by 
planes $g=g_0$, with $g_0$ constant and we shall give pictures of the resulting bifurcation 
diagram on these planar sections in which the Cartesian coordinates are $(u,\ell)$, where the horizontal line is the 
$u$--axis and $u\geq0$.

Here we also use colors to refer to the bifurcation surfaces:
\begin{enumerate}[(a)]
	\item surface (${\cal S}_{3}$) is drawn in yellow (weak singularities). We draw 
	it as a continuous curve if  the singular point is a focus or as a dashed curve if it is a saddle;
	\item surface (${\cal S}_{4}$) is drawn in purple (presence of at least one invariant straight line). We draw it as 
	a continuous curve if it implies a topological change or as a dashed curve otherwise;
	\item surface (${\cal S}_{6}$) is drawn in black and dashed (an antisaddle is on the edge of turning from a node to a focus 
	or vice versa);
	\item nonalgebraic surface (${\cal S}_{7}$) is also drawn in purple (connections of separatrices); 
	\item surface (${\cal S}_{8}$) is drawn in cyan (presence of an invariant parabola). We draw it as 
	a continuous curve if it implies a topological change or as a dashed curve otherwise.
	\item Here we follow the pattern established on Remark~\ref{rmk-colors} for surfaces (${\cal S}_{0}$) and (${\cal S}_{5}$). 
	\item As surface \label{stat-color-deg} (${\cal S}_{1}$) is the whole plane $g=0$, due to the same reason presented on Remark~\ref{rmk-colors}, we shall not use a color for describing this entire bifurcation surface. However, for indicating the bifurcation straight line $\ell=0$ (belonging to surface (${\cal S}_{1}$)) we shall use green color and draw it as a continuous line.
\end{enumerate}

As in the previous section, in order to obtain the singular slices needed for the study of the bifurcation 
diagram of systems~\eqref{eq:nf-QES-B}, here we also perform all the computations in an algorithm written 
in software Mathematica. The reader may find the computations in the file available for free download through 
the link \url{http://mat.uab.cat/~artes/articles/qvfES/qvfES-B.nb}. 

The next result presents all the algebraic values of $g$ corresponding to singular slices (or planes) in the bifurcation diagram. 
Its proof follows from the study done with the help of the mentioned algorithm.

\begin{lemma}\label{lemma:values-geom-study-QESB}
	Consider the algebraic bifurcation surfaces defined before. The study of their singularities, their intersection points, and their 
	tangencies with planes $g=g_0$ provides the following set of four singular values of the parameter $g$:
	{$$
		\left\{1,0,-1,-2\right\}\!.
		$$}
\end{lemma}

\begin{remark}\label{rmk:another-sv-g-QESB} It is easy to conclude that surfaces (${\cal S}_{6}$) and (${\cal S}_{8}$) intercept 
	themselves along 
	$$\left\{g=-\dfrac{u^2}{2(1+u^2)}, \ \ell=\dfrac{u(2+u^2)}{(1+u^2)^2}
	\right\}\!.$$
	We notice that, when $u\rightarrow\infty$, such an intersection goes to 
	$$\left\{g=-\dfrac{1}{2}, \ \ell=0
	\right\}\!.$$
	So, $g=-1/2$ can be also considered as a singular value of the parameter $g$. And
	at this singular value, surfaces (${\cal S}_{6}$) and (${\cal S}_{8}$) intercept themselves 
	at infinity (at the endpoint of straight line $\ell=0$). 
\end{remark}

We collect the values of the parameter $g$ obtained from Lemma~\ref{lemma:values-geom-study-QESB} and 
Remark~\ref{rmk:another-sv-g-QESB} and, in the next result we present the complete list of algebraic singular 
planes corresponding to values of the parameter $g$.

\begin{proposition}\label{prop:algebraic-values-of-g}
	The full set of needed algebraic singular slices in the bifurcation diagram of family~\eqref{eq:nf-QES-B} is formed by 
	five elements which correspond to the values of $g$ in \eqref{eq:algebraic-values-of-g}.
\end{proposition}
\begin{equation}\label{eq:algebraic-values-of-g}
\begin{aligned}
g_{1}&=1, \ g_{3}=0, \ g_{5}=-\frac{1}{2}, \ g_{7}=-1, \ g_{9}=-2.\\ 
\end{aligned}
\end{equation}

The numeration in \eqref{eq:algebraic-values-of-g} is not consecutive  since we reserve numbers for generic slices. 
We point out that we have not found nonalgebraic slices, as in \cite{Artes-Mota-Rezende-2021c}, for instance.

In order to determine all the parts generated by the bifurcation surfaces from $({\cal S}_{0})$ to $({\cal S}_{8})$, 
we first draw the horizontal slices of the three--dimensional parameter space which correspond to the explicit values 
of $g$ obtained in Proposition~\ref{prop:algebraic-values-of-g}. However, as it will be discussed later, the presence 
of nonalgebraic bifurcation surfaces will be detected and their behavior as we move from slice to slice will be approximately 
determined. We add to each interval of singular values of $g$ an intermediate value for which we represent the bifurcation 
diagram of singularities. The diagram will remain essentially unchanged in these open intervals except the parts affected by 
the bifurcation. All the eleven sufficient values of $g$ are shown in \eqref{eq:values-of-g-QES-B}.

\begin{equation}\label{eq:values-of-g-QES-B}
\arraycolsep=0.6cm\begin{array}{ll}
g_{0}=2		            						& g_{6}=-3/4 						\vspace{0.1cm} \\ 
g_{1}=1       											& g_{7}=-1   	           \vspace{0.1cm} \\
g_{2}=1/2         								& g_{8}=-3/2             \vspace{0.1cm}\\
g_{3}=0														& g_{9}=-2									 \vspace{0.1cm} \\
g_{4}=-1/4                  			& g_{10}=-3     						\vspace{0.1cm} \\           
g_{5}=-1/2												&              			
\end{array}
\end{equation}

The values indexed by positive odd indices in \eqref{eq:values-of-g-QES-B} correspond to explicit values of $g$ for 
which there is a bifurcation in the behavior of the systems on the slices. Those indexed by even values are just intermediate 
points which are necessary to the coherence of the bifurcation diagram. 

We now begin the analysis of the bifurcation diagram by studying completely one generic slice and after by moving
from slice to slice and explaining all the changes that occur. As an exact drawing of the curves produced by intersecting 
the surfaces with the slices gives us very small parts which are difficult to distinguish, and points of tangency are almost 
impossible to recognize, we have produced topologically equivalent figures where parts are enlarged and tangencies 
are easy to observe. 

The reader may find the exact pictures of the five singular slices (containing only the algebraic surfaces) described in 
\eqref{eq:algebraic-values-of-g} in a PDF file available at the link \url{http://mat.uab.es/~artes/articles/qvfES/qvfES-B.pdf}. 

As in the previous section we use the same pattern in order to describe each part of the bifurcation diagram 
(labels and colors) and we also use continuous and dashed (bifurcation) curves, as explained before.

\begin{figure}[h!]
	\centering
	\includegraphics[width=0.5\textwidth]{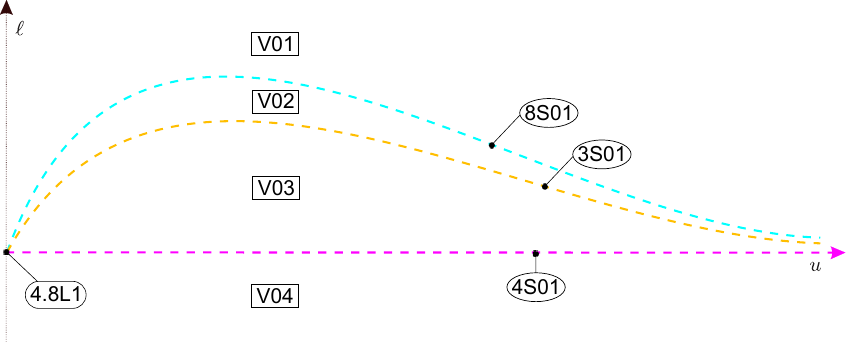}
	\caption{\small \label{fig:slice-QES-B-01} Generic slice of the parameter space when $g=2$}
\end{figure}

In Fig.~\ref{fig:slice-QES-B-01} we represent the entire generic slice of the parameter space when $g=g_{0}=2$
(remember that we proved that it is enough to consider $u\geq0$). In this figure (and in the next ones) we denote 
the $\ell$--axis with a dashed and thin black straight line.

\begin{figure}[h!]
	\centering
	\includegraphics[width=0.5\textwidth]{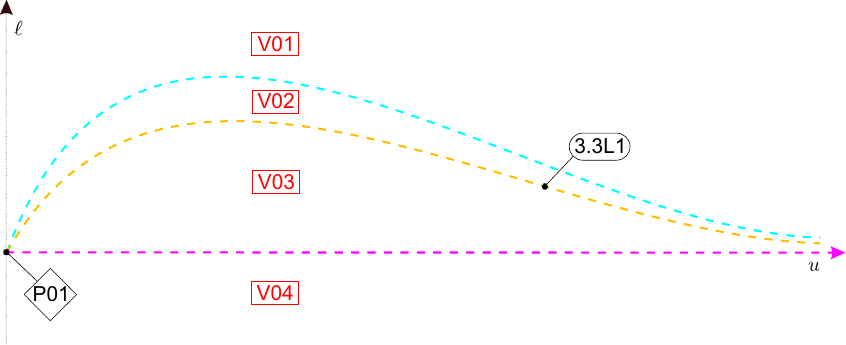}
	\caption{\small \label{fig:slice-QES-B-02} Singular slice of the parameter space when $g=1$}
\end{figure}

When we consider the singular value $g=g_{1}=1$ of the parameter $g$ we observe that surface 
(${\cal S}_{3}$) reduces to 
$$-2 (1 + u^2) (\ell - 2 u + \ell u^2)^3.$$
By discarding the factor $-2(1+u^2)$ (which does not have real roots) we observe that such a 
surface has multiplicity three. On the other hand, by item 1 of Remark~\ref{weak-singularities2} 
this change of multiplicity is related to the presence of an integrable saddle. For this case, the bifurcation 
diagram can be seeing in Fig.~\ref{fig:slice-QES-B-02}.

\begin{figure}[h!]
	\centering\includegraphics[width=0.5\textwidth]{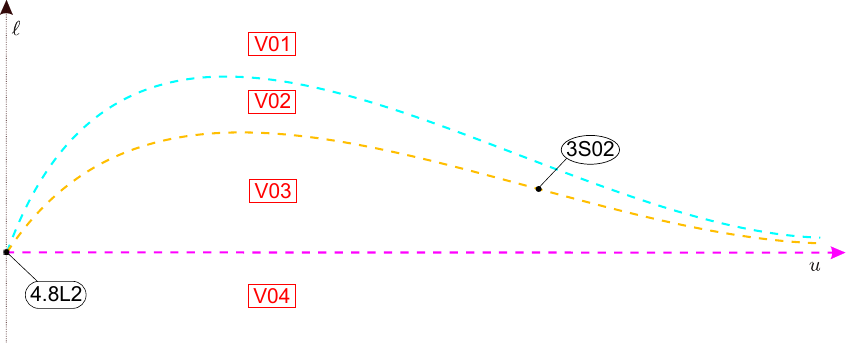} 
	\caption{\small \label{fig:slice-QES-B-03} Generic slice of the parameter space when $g=1/2$}
\end{figure}

Now, for the generic value $g=g_2=1/2$, the yellow curve is simple again (i.e. it has multiplicity one), 
see Fig.~\ref{fig:slice-QES-B-03}.

\begin{figure}[h!]
	\centering\includegraphics[width=0.5\textwidth]{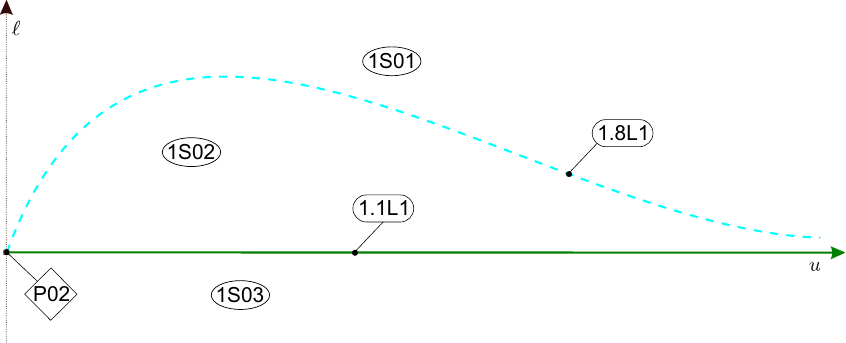} 
	\caption{\small \label{fig:slice-QES-B-g0} Singular slice of the parameter space when $g=0$}
\end{figure}

As we said before, for $g=g_3=0$ systems~\eqref{eq:nf-QES-B} are degenerate. In fact, for this value
of the parameter $g$ we have that bifurcation surfaces (${\cal S}_{1}$), (${\cal S}_{3}$),
(${\cal S}_{4}$), and (${\cal S}_{6}$) vanish and, in addition, (${\cal S}_{0}$)$\vert_{g=0}=1$,
(${\cal S}_{5}$)$\vert_{g=0}=2$, and (${\cal S}_{8}$)$\vert_{g=0}=\ell-2u+\ell u^2$. 
Moreover, Remark~\ref{study-deg-famB} provides the type of the curve filled up with singular points,
according to the value of the parameter $\ell$. In Fig.~\ref{fig:slice-QES-B-g0} we present the singular
slice $g=g_3=0$ in which we are using the colors and pattern we mentioned in page~\pageref{stat-color-deg}. 

\begin{figure}[h!]
	\centering\includegraphics[width=0.5\textwidth]{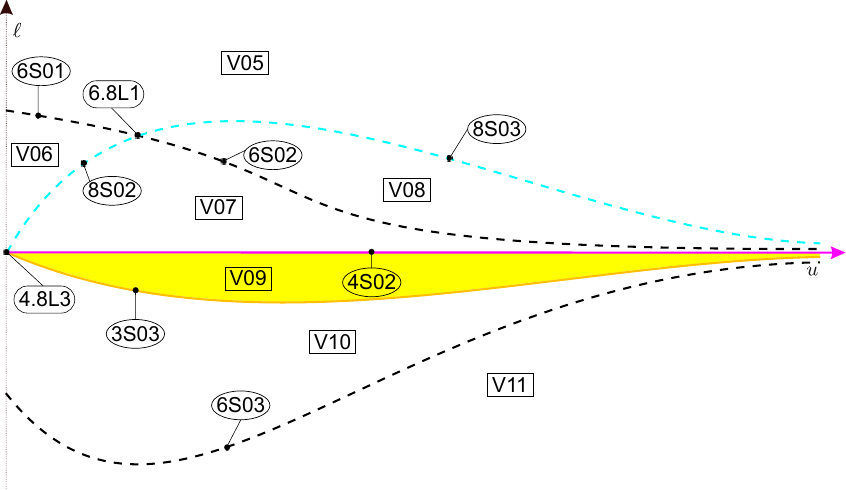} 
	\caption{\small \label{fig:slice-QES-B-04} Generic slice of the parameter space when $g=-1/4$}
\end{figure}

We start the study of the negative values of the parameter $g$ (so according to item 2 of Remark~\ref{weak-singularities2}, 
for every fixed $g<0$, the point $(u, \ell)=(0, 0)$ corresponds to a phase portrait possessing a center 
type singularity). According to \eqref{eq:values-of-g-QES-B} we consider the generic slice given 
by $g=g_4=-1/4$. For this value of the parameter $g$:
\begin{itemize}
	\item we now have the presence of two segments of the black surface (${\cal S}_{6}$); 
	\item the purple straight line (${\cal S}_{4}$) is now drawn as a continuous curve, since it represents 
	a separatrix connection; and
	\item on the yellow segment $3S_{3}$ the corresponding phase portrait possesses a weak focus (of
	order one) and, consequently, this branch of surface (${\cal S}_{3}$) corresponds to a Hopf 
	bifurcation. This means that the phase portrait corresponding to one of the sides of this segment 
	must have a limit cycle; in fact it is in the region $V_{9}$.
\end{itemize}
The corresponding slice is presented in Fig.~\ref{fig:slice-QES-B-04}.

\begin{figure}[h!]
	\centering\includegraphics[width=0.5\textwidth]{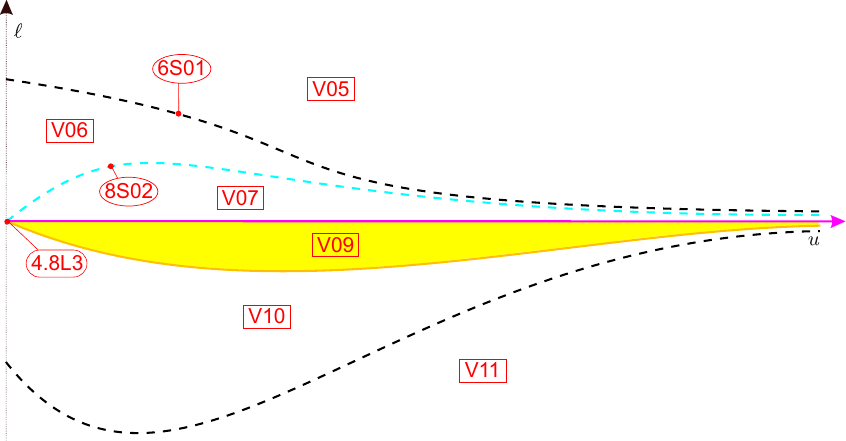} 
	\caption{\small \label{fig:slice-QES-B-05} Singular slice of the parameter space when $g=-1/2$}
\end{figure}

According to Remark~\ref{rmk:another-sv-g-QESB}, we know that for $g<0$ and $\ell>0$ surfaces 
(${\cal S}_{6}$) and (${\cal S}_{8}$) have a common point, for every $u>0$. In fact, this point 
is denoted in Fig.~\ref{fig:slice-QES-B-04} by $6.8L_{1}$. The same remark shows that such an 
intersection point goes to infinity at $g=g_5=-1/2$, and this displacement carries volume region $V_{8}$ 
to infinity. For this singular value of the parameter $g$, the corresponding bifurcation diagram is 
presented in Fig.~\ref{fig:slice-QES-B-05}.

If we consider the generic slice given by $g=g_6=-3/4$ we observe that the intersection point presented 
in Remark~\ref{rmk:another-sv-g-QESB} goes to the complex plane. As there is no other significant 
phenomenon to analyze, we conclude that for the generic value under consideration, the bifurcation diagram 
behaves as the one presented in Fig.~\ref{fig:slice-QES-B-05}.

\begin{figure}[h!]
	\centering\includegraphics[width=0.5\textwidth]{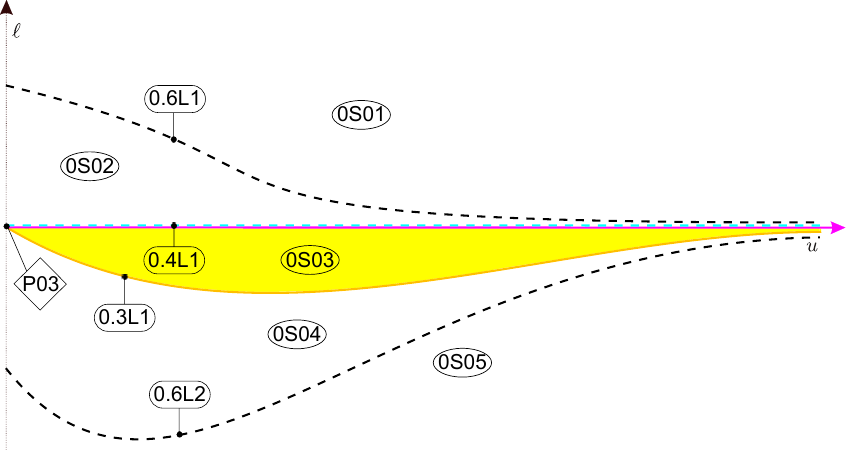} 
	\caption{\small \label{fig:slice-QES-B-07} Singular slice of the parameter space when $g=-1$}
\end{figure}

Now we consider the singular slice $g=g_7=-1$. One may say that this is a quite interesting singular slice, because:
\begin{itemize}
	\item Previously we mentioned that surface (${\cal S}_{0}$), related to a presence of an infinite elliptic--saddle 
	of type $\widehat{\!{1\choose 2}\!\!}\ E-H$, defines the entire plane $g=-1$. As it was pointed out in 
	\cite{Artes-Mota-Rezende-2021c} 	each phase portrait obtained in the study of this slice is topologically 
	equivalent to a phase portrait obtained in a neighborhood of this plane. However, in order to have a coherent 
	bifurcation diagram, this plane must be studied. Here we follow the pattern established in Remark~\ref{rmk-colors}
	and we shall not draw this plane in brown color.
	\item For this value of the parameter $g$, surfaces (${\cal S}_{4}$) and (${\cal S}_{6}$) coincides 
	along $\ell=0$. The remaining parts of the bifurcation diagram behave as in the previous slice.
\end{itemize}
In  Fig.~\ref{fig:slice-QES-B-07} we present the singular slice $g=-1$ completely labeled. In such a 
figure we use the the same pattern as the one applied in Fig.~\ref{fig:slice-QES-A-31} from the previous section.

The next generic slice $g=g_8=-3/2$ deserves a special attention. After passing by an infinite singularity of
type $\widehat{\!{1\choose 2}\!\!}\ E-H$ it is expected to obtain new phase portraits possessing orbits of 
the infinite elliptic--saddle in different positions (when we compare these new phase portrait with the ones 
we had before the bifurcation related to $\widehat{\!{1\choose 2}\!\!}\ E-H$). So, in the slice under 
consideration one may find distinct situations to analyze.

\begin{figure}[h!]
	\centering\includegraphics[width=0.5\textwidth]{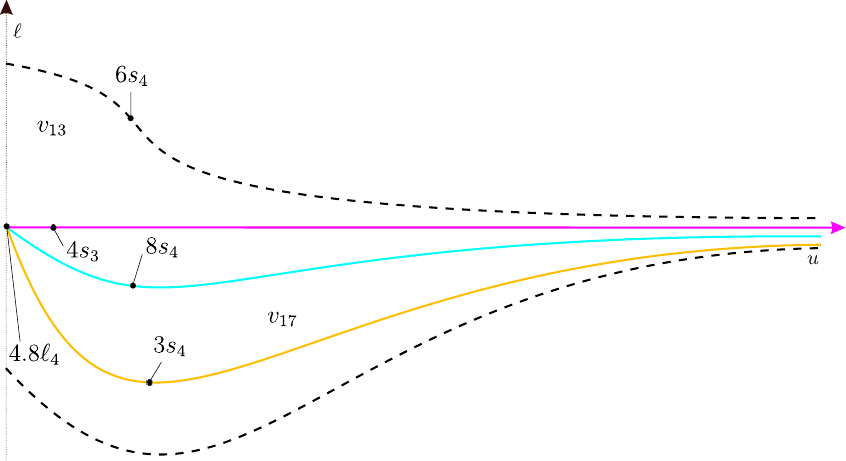} 
	\caption{\small \label{fig:slice-QES-B-08-alg} Generic slice of the parameter space when $g=-3/2$ (only algebraic surfaces)}
\end{figure}

In  Fig.~\ref{fig:slice-QES-B-08-alg} we present such a generic slice, showing only the algebraic surfaces. 
We note the existence of continuous branches of surfaces (${\cal{S}}_{3}$) (in yellow),  (${\cal{S}}_{4}$) 
(in purple), and (${\cal{S}}_{8}$) (in cyan). This means the existence of a weak focus, in the case of surface 
(${\cal{S}}_{3}$), the existence of an algebraic invariant straight line provided by a separatrix connection, 
in the case of surface (${\cal{S}}_{4}$), and the existence of an algebraic invariant parabola formed by 
a separatrix connection, in the case of surface (${\cal{S}}_{8}$). 

We now place for each set of the partition on this slice the local behavior of the flow around the singular points. 
For a specific value of the parameters of each one of the sets in this partition we compute the global phase portrait 
with the numerical program P4 \cite{progP4,Dumortier-Llibre-Artes-2006}.

In this slice we have a partition in two--dimensional unbounded parts. From now on, we use lower--case letters 
provisionally to describe the sets found algebraically in order to do not interfere with the final partition described 
with capital letters. 

For each two--dimensional part we obtain a phase portrait which is coherent with those of all their borders. 
Except for two parts, which are shown in Fig.~\ref{fig:slice-QES-B-08-alg} and named as follows:
\begin{itemize}\label{page:regions-g-8}
	\item $v_{13}$: the region $\{u\ge0, \ell\geq0\}$ bordered by the black curve and infinity;
	\item $v_{17}$: the region bordered by yellow and cyan curves and also by infinity.
\end{itemize}
The study of these parts is important for the coherence of the bifurcation diagram. That is why we have 
decided to present only these parts in the mentioned figures.

\begin{figure}[h!]
	\centering\includegraphics[width=0.7\textwidth]{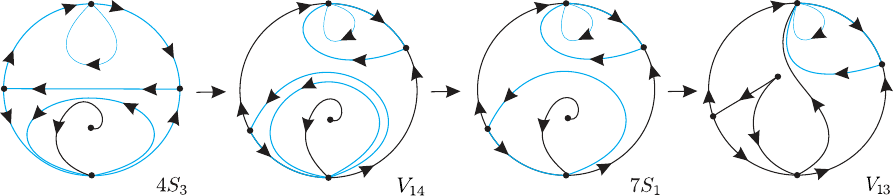} 
	\caption{\small \label{fig:transitionV13-V14} Sequence of phase portraits in parts $V_{13}$ and $V_{14}$ of slice 
		$g=-3/2$ (the labels are according to Fig.~\ref{fig:slice-QES-B-08})}
\end{figure}

We begin with the analysis of part $v_{13}$. The phase portrait in $v_{13}$ near $4s_{3}$ possesses an 
infinite graphic formed by orbits contained in the parabolic sectors of the (infinite) elliptic--saddle. However,
the phase portrait in $v_{13}$ near $6s_{4}$ does not possess such a graphic. Then, there must exist at least 
one element of surface (${\cal{S}}_{7}$) (see $7S_{1}$ in Fig.~\ref{fig:slice-QES-B-08}) dividing part 
$v_{13}$ into two ``new'' parts, $V_{13}$ and $V_{14}$, which represents a bifurcation due to the 
connection between a separatrix of the infinite elliptic--saddle with a separatrix of the infinite saddle 
(see Fig.~\ref{fig:transitionV13-V14} for a sequence of phase portraits in these parts). 

We claim that nonalgebraic surface $7S_{1}$ is unbounded and $4.8\ell_1$ is one of its endpoints. In fact,
numerical verifications indicate the truth of this statement. Note that it is not possible that the starting point of 
this surfaces is on $6s_{4}$, since on black surfaces we have only a $C^\infty$ node--focus bifurcation. On 
the other hand, the endpoint of $7S_{1}$ cannot be on $4s_{3}$ because, in order to have this, first we 
need to break the invariant straight line connecting the opposite infinite saddles. Then, the only possible 
endpoint of surface $7S_{1}$ is $4.8\ell_{1}$, and our claim is proved.

\begin{figure}[h!]
	\centering
	\includegraphics[width=0.7\textwidth]{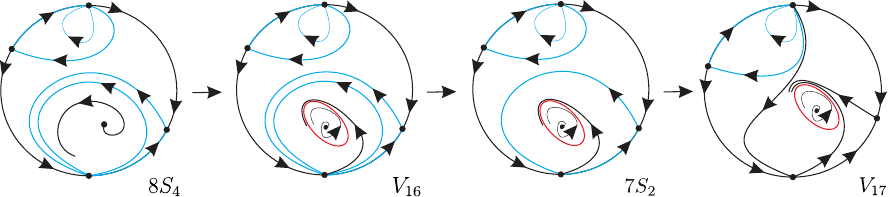}
	\caption{\small \label{fig:transitionV16-V17} Sequence of phase portraits in parts $V_{16}$ and $V_{17}$ of slice $g=-3/2$ (the labels are according to Fig.~\ref{fig:slice-QES-B-08})}
\end{figure}

\begin{figure}[h!]
	\centering
	\includegraphics[width=0.5\textwidth]{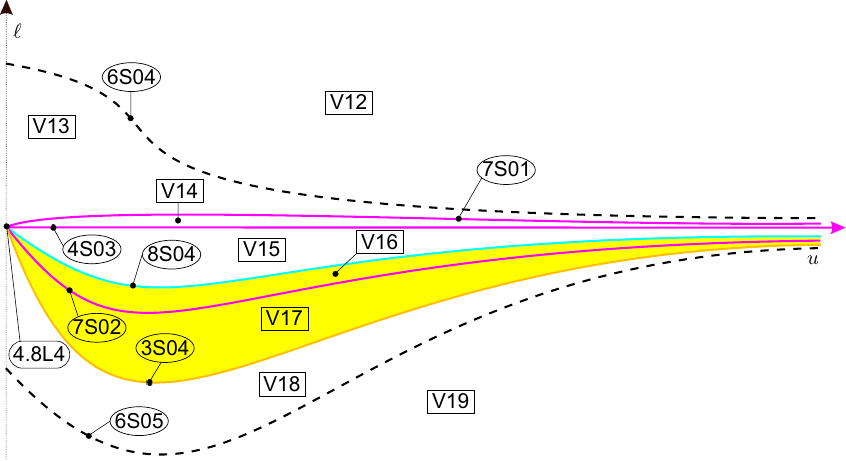}
	\caption{\small \label{fig:slice-QES-B-08} Generic slice of the parameter space when $g=-3/2$}
\end{figure}

Now, we carry out the analysis of part $v_{17}$. We consider the segment $3s_{4}$ in Fig.~\ref{fig:slice-QES-B-08-alg}, 
which is one of the borders of part $v_{17}$. On this segment, the corresponding phase portrait possesses 
a weak focus (of order one) and, consequently, this branch of surface (${\cal{S}}_{3}$) corresponds to a 
Hopf bifurcation. This means that the phase portrait corresponding to one of the sides of this segment must 
have a limit cycle; in fact it is in $v_{17}$. Moreover, the phase portrait in $v_{17}$ near $8s_{4}$ possesses an 
infinite graphic formed by orbits contained in the parabolic sectors of the (infinite) elliptic--saddle. However,
the phase portrait in $v_{17}$ near $3s_{4}$ does not possess such a graphic. Then, there must exist at least 
one element of surface (${\cal{S}}_{7}$) (see $7S_{2}$ in Fig.~\ref{fig:slice-QES-B-08}) dividing part 
$v_{17}$ into two ``new'' parts, $V_{16}$ and $V_{17}$, which represents a bifurcation due to the 
connection between a separatrix of the infinite elliptic--saddle with a separatrix of the infinite saddle 
(see Fig.~\ref{fig:transitionV16-V17} for a sequence of phase portraits in these parts). 

In this paragraph we prove that nonalgebraic surface $7S_{2}$ is unbounded and $4.8\ell_1$ is one of its 
endpoints. Indeed, numerical verifications indicate that this fact is true. Note that if the starting point of 
this surface is any point of $3s_{4}$ then a portion of this subset must not refer to a Hopf bifurcation, 
which contradicts the fact that on $3s_{4}$ we have a weak focus of order one. In addition, the 
endpoint of $7S_{2}$ cannot be on $8s_{4}$ because, in order to have this, first it is necessary 
to break the invariant parabola formed by a separatrix of the infinite elliptic--saddle. So, the only possible 
endpoint of surface $7S_{2}$ is $4.8\ell_{1}$, as we wanted to prove.

The complete bifurcation diagram for this part can be seeing in Fig.~\ref{fig:slice-QES-B-08}.

\begin{figure}[h!]
	\centering
	\includegraphics[width=0.5\textwidth]{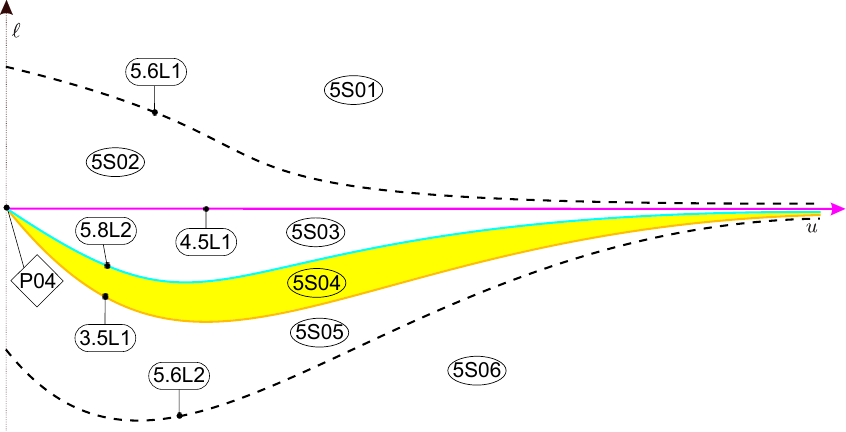}
	\caption{\small \label{fig:slice-QES-B-09} Singular slice of the parameter space when $g=-2$}
\end{figure}

Now we consider the singular slice $g=g_9=-2$. This is another interesting and important singular slice.
\begin{itemize}
	\item Surface ($\mathcal{S}_{5}$)$= g+2$ is related to a coalescence of infinite singular points. 
	Remember that if $\ell\ne0$ the phase portraits obtained in the study of this slice possess at most 
	one pair of infinite singular points and, if $\ell=0$ the corresponding phase portraits have the line 
	at infinity filled up with singularities. Here we follow Remark~\ref{rmk-colors} and we shall not 
	draw the slice $g=-2$ in red color.
	\item By studying the transition among regions and phase portraits from $g=-3/2$ with regions and phase
	portraits from $g=-2$ we observe that $V_{14}$ (respectively $V_{16}$) from slice $g=-3/2$ converges to 
	$4.5L_{1}$ (respectively $5.8L_{2}$) from slice $g=-2$. The correspondence among the remaining regions
	of these slices is clear.
\end{itemize}
In  Fig.~\ref{fig:slice-QES-B-09} we present the slice $g=-2$ completely labeled. In such a 
figure we use the same pattern as the one used in the slices $g=0$ and $g=-1$ in order to present a label for 
each region.

\begin{figure}[h!]
	\centering
	\includegraphics[width=0.5\textwidth]{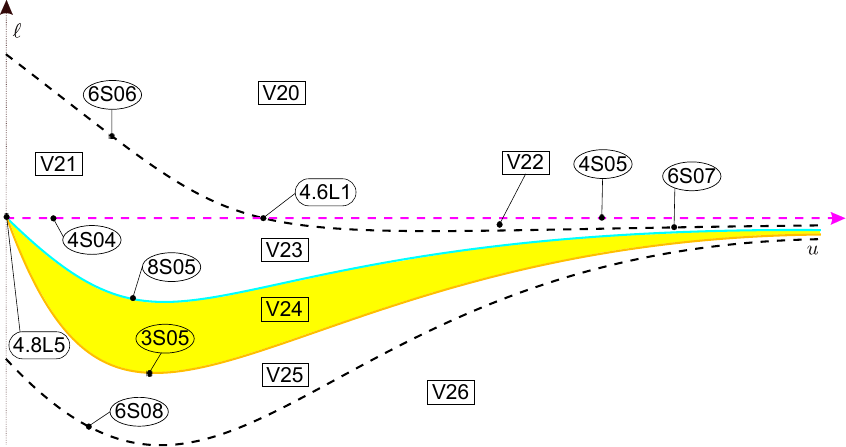}
	\caption{\small \label{fig:slice-QES-B-10} Generic slice of the parameter space when $g=-3$}
\end{figure}

Finally we consider the generic slice $g=g_{10}=-3$. In what follows we present some comments on this
slice.
\begin{itemize}
	\item We observe that due to the nature of the coalescence of infinite singularities on this slice, in the next 
	generic slice $g=g_{10}=-3$ we shall expect to obtain phase portraits with a reduced number of 
	separatrices. In fact, at $g=g_8=-3/2$ we had phase portraits possessing an infinite elliptic--saddle
	and also an infinite saddle. At $g=g_9=-2$ the infinite saddle coalesced with the infinite elliptic--saddle.
	Now, at the generic slice $g=g_{10}=-3$ we have an infinite elliptic--saddle and also an infinite node.
	\item At this value of the parameter $g$ the purple curve (surface ($\mathcal{S}_{4}$)) no longer represents
	a separatrix connection, and this is due to the fact that we do not have an enough number of separatrices in
	order to have an invariant straight line, since we passed by the mentioned coalescence of infinite singularities.
	\item Surfaces ($\mathcal{S}_{4}$) and ($\mathcal{S}_{6}$) have an intersection point along $\ell=0$.
\end{itemize}
The complete bifurcation diagram for this part is presented in Fig.~\ref{fig:slice-QES-B-10}.

Since there is coherence among the generic and singular slices presented before, no more slices are 
needed for the complete coherence of the bifurcation diagram. So, all the values of the parameter 
$g$ in \eqref{eq:values-of-g-QES-B} are sufficient for the coherence of the bifurcation diagram. 
Thus, we can affirm that we have described a complete bifurcation diagram for class $\overline{\QESB}$ 
modulo islands, as we discuss in Sec.~\ref{sec:islands-QESB}.

\subsubsection{Other relevant facts about the bifurcation diagram} \label{sec:islands-QESB}

The bifurcation diagram we have obtained for the class $\overline{\QESB}$ is completely coherent, i.e. in this family, 
by taking any two points in the parameter space and joining them by a continuous curve, along this curve 
the changes in phase portraits that occur when crossing the different bifurcation surfaces we mention can 
be completely explained.

Nevertheless, we cannot be sure that this bifurcation diagram is the complete bifurcation diagram for 
$\overline{\QESB}$ due to the possibility of the existence of ``islands'' inside the parts bordered by unmentioned 
bifurcation surfaces. In case they exist, these ``islands'' would not mean any modification of the nature 
of the singular points. So, on the border of these ``islands'' we could only have bifurcations due to saddle 
connections or multiple limit cycles.

In case there were more bifurcation surfaces, we should still be able to join two representatives of any 
two parts of the 89 parts of $\overline{\QESB}$ found until now with a continuous curve either without crossing 
such a bifurcation surface or, in case the curve crosses it, it must do it an even number of times without 
tangencies, otherwise one must take into account the multiplicity of the tangency, so the total number 
must be even. This is why we call these potential bifurcation surfaces ``\textit{islands}''.

However, we have not found a different phase portrait which could fit in such an island. A potential ``island'' 
would be the set of parameters for which the phase portraits possess a double limit cycle and this ``island'' 
would be inside the parts where $W_{4}<0$ since we have the presence of a focus.

\subsubsection{Completion of the proof of Theorem~\ref{th:main-thm-QES-B}} \label{sec:invariants-QESB} 

In the bifurcation diagram we may have topologically equivalent phase portraits belonging to distinct parts 
of the parameter space. As here we have 89 distinct parts of the parameter space, to help us to identify 
or to distinguish phase portraits, we need to introduce some invariants and we actually choose integer valued, 
character and symbol invariants. Some of them were already used in \cite{Artes-Rezende-Oliveira-2013b} 
and \cite{Artes-Mota-Rezende-2021c}, but we recall them and introduce some needed ones. These invariants 
yield a classification which is easier to grasp.

\begin{definition}\label{def:invariant-I1-QESB}
	We denote by $I_{1}(S)$ a symbol from the set $\{\emptyset, \left[ \times \right], \left[)(\right]\}$ which indicates the following configuration of curves filled up with singularities, respectively: none (nondegenerate systems -- in this case all systems do not contain a curve filled up with singularities), two real straight lines intersecting at a finite point, and an hyperbola. This invariant only makes sense to distinguish the degenerate phase portrait obtained.
\end{definition}

\begin{definition}\label{def:invariant-I2-QESB}
	We denote by $I_{2}(S)$ the sum of the indices of the isolated real finite singular points.
\end{definition}

\begin{definition}\label{def:invariant-I3-QESB}
	We denote by $I_{3}(S)$ the number of real infinite singular points. We note that this number can also be infinite, which is represented by $\infty$.
\end{definition}

\begin{definition}\label{def:invariant-I4-QESB}
	For a given infinite singularity $s$ of a system $S$, let $\l_s$ be the number of global or local separatrices beginning or ending at $s$ and which do not lie on the line at infinity. We have $0\leq\l_s\leq4$. We denote by $I_{4}(S)$ the sequence of all such $\l_s$ when $s$ moves in the set of infinite singular points of the system $S$. We start the sequence at the infinite singular point which receives (or sends) the greatest number of separatrices and take the direction which yields the greatest absolute value, e.g. the values $2110$ and $2011$ for this invariant are symmetrical (and, therefore, they are the same), so we consider $2110$.
\end{definition}

\begin{definition}\label{def:invariant-I5-QESB}
	We denote by $I_{5}(S)$ the number of limit cycles around a foci.
\end{definition}

\begin{definition}\label{def:invariant-I6-QESB}
	We denote by $I_{6}(S)$ an element from the set $\{c, f\}$ indicating the type of the real finite singularity located inside the region bordered by the graphic, which can be either a \textbf{c}enter or a \textbf{f}ocus.
\end{definition}

\begin{definition}\label{def:invariant-I7-QESB}
	We denote by $I_{7}(S)$ a pair $(A,B)$ where $A$ and $B$ represent the number of separatrices arriving or leaving the corresponding parabolic sectors of the singularity $\widehat{\!{1\choose 2}\!\!}\ PHP-E$ at infinity.
\end{definition}

As we have noted previously in Remark~\ref{rem:f-n}, we do not distinguish between phase portraits whose 
only difference is that in one we have a finite node and in the other a focus. Both phase portraits are topologically 
equivalent and they can only be distinguished within the $C^1$ class. In case we may want to distinguish between 
them, a new invariant may easily be introduced.

\begin{theorem} \label{th:QESB-inv}
	Consider the class $\overline{\QESB}$ and all the phase portraits that we have obtained for this family. The values of the affine 
	invariant ${\cal I} = (I_{1}, I_{2}, I_{3}, I_{4}, I_{5}, I_{6}, I_{7})$ given in the diagram from Table \ref{tab:geom-classif-QESB-pt1} 
	yields a partition of these phase portraits of the class $\overline{\QESB}$.
	
	Furthermore, for each value of $\cal I$ in this diagram there corresponds a single phase portrait; i.e. $S$ and $S'$ 
	are such that ${\cal I}(S)={\cal I}(S')$, if and only if $S$ and $S'$ are topologically equivalent.
\end{theorem}

The bifurcation diagram for $\overline{\QESB}$ has 89 parts which produce 27 topologically different phase portraits as 
described in Tables~\ref{tab:geom-classif-QESB-pt1} to \ref{tab:top-equiv-QESB-pt1}. The remaining 62 parts 
do not produce any new phase portrait which was not included in the 27 previous ones. The difference is basically 
the presence of a strong focus instead of a node and vice versa, weak points, and a presence of invariant algebraic 
curves (lines or parabolas) which do not represent a separatrix connection.

The phase portraits having neither limit cycle nor graphic have been denoted surrounded by parenthesis, for example $(V_{20})$; 
the phase portraits having one limit cycle have been denoted surrounded by brackets, for example $[V_{24}]$; the phase portraits 
having one graphic have been denoted surrounded by $\{\ast\}$ and those ones having two or more graphics have been denoted 
surrounded by $\{\!\{\ast\}\!\}$, for example $\{5S_{3}\}$ and $\{\!\{V_{1}\}\!\}$, respectively. Moreover, the phase portraits 
having one limit cycle and more than one graphic have been denoted surrounded by $[\{\{\ast\}\}]$, for example $[\{\{V_{17}\}\}]$.

\begin{proof}[Proof of Theorem~\ref{th:QESB-inv}]
	The above result follows from the results in the previous sections and a careful analysis of the bifurcation diagrams given in 
	Sec.~\ref{subsec:bd-QES-B}, in Figs.~\ref{fig:slice-QES-B-01} to Fig.~\ref{fig:slice-QES-B-10}, the definition of the 
	invariants $I_{j}$ and their explicit values for the corresponding phase portraits.
\end{proof}

We recall some observations regarding the equivalence relations used in this study: the affine and time rescaling, $C^1$ and topological equivalences.

The coarsest one among these three is the topological equivalence and the finest is the affine equivalence. We can have two systems which are topologically equivalent but not $C^1-$equivalent. For example, we could have a system with a finite antisaddle which is a structurally stable node and in another system with a focus, the two systems being topologically equivalent but belonging to distinct $C^1-$equivalence classes, separated by the surface $({\cal S}_{6})$ on which the node turns into a focus.

In Table~\ref{tab:top-equiv-QESB-pt1} we list in the first column 27 parts with all the distinct phase portraits of Fig.~\ref{fig:pp-QES-B}. Corresponding to each part listed in column one we have in each row all parts whose phase portraits are topologically equivalent to the phase portrait appearing in column 1 of the same row.

In the second column we set all the parts whose systems yield topologically equivalent phase portraits to those in the first column, but which may have some algebro--geometric features related to the position of the orbits. In the third column we present all the parts which are topologically equivalent to the ones from the first column having a focus instead of a node.

In the fourth (respectively, fifth; and sixth) column we list all parts whose phase portraits have a node which is at a bifurcation point producing foci close to the node in perturbations, a node--focus to shorten (respectively, a finite weak singular point; and possess an invariant curve (straight line and/or parabola) not yielding a connection of separatrices).

The last column refers to other reasons associated to different geometrical aspects and they are described as follows:
\begin{enumerate}[(1)]
	\vspace{-2mm}
	\item  The phase portraits correspond to symmetric parts of the bifurcation diagram;
	\vspace{-2mm}
	\item  the phase portrait possesses a singularity of type $\widehat{\!{1\choose2}\!\!}\ E-H$ at infinity.
\end{enumerate}

Whenever phase portraits appear in a row in a specific column, the listing is done according to the decreasing dimension of the parts where they appear, always placing the lower dimensions on lower lines.

\subsubsection{Proof of Theorem~\ref{th:main-thm-QES-B}}

The bifurcation diagram described in Sec. \ref{subsec:bd-QES-B}, plus Table~\ref{tab:geom-classif-QESB-pt1} of the geometrical invariants distinguishing the 27 phase portraits, plus Table~\ref{tab:top-equiv-QESB-pt1} giving the equivalences with the remaining phase portraits lead to the proof of Theorem~\ref{th:main-thm-QES-B}.

\begin{table}\caption{\small Geometric classification for the family $\QESB$}\label{tab:geom-classif-QESB-pt1} 
	\begin{center}
		\[
		I_{1}\!=\!
		\left\{
		\begin{array}{ll}
		\left[)(\right] \,\, \left\{\left\{1S_{1}\right\}\right\}, \\  			
		\left[ \times \right] \,\, \left\{\left\{1.1L_{1}\right\}\right\}, \\  			
		\emptyset \,\, \& \,\, I_{2}\!=\! 
		\left\{
		\begin{array}{ll}
		-1 \,\, \left\{\left\{V_{1}\right\}\right\}, \\  			
		1 \,\, \& \,\, I_{3}\!=\!
		\left\{
		\begin{array}{ll}
		1 \,\, \& \,\, I_{4}\!=\!
		\left\{
		\begin{array}{ll}
		20 \,\, \left\{\left\{5.8L_{2}\right\}\right\}, \\  				
		21 \,\, \& \,\, I_{5}\!=\!
		\left\{
		\begin{array}{ll}
		0 \,\, \left\{\left\{5S_{1}\right\}\right\}, \\  
		1 \,\, \left\{\left\{5S_{4}\right\}\right\}, \\  
		\end{array}
		\right. \\
		30 \,\, \left\{\left\{5S_{3}\right\}\right\}, \\  
		\end{array}
		\right. \\
		2 \,\, \& \,\, I_{4}\!=\!
		\left\{
		\begin{array}{ll}
		1010 \,\, \& \,\, I_{5}\!=\! 0\,\, \& \,\, I_{6}\!=\! 
		\left\{
		\begin{array}{ll}
		c \,\, \left\{\left\{4.8L_{3}\right\}\right\}, \\  
		f \,\, \left\{\left\{4S_{2}\right\}\right\}, \\  
		\end{array}
		\right. \\
		1110 \,\, \& \,\, I_{5}\!=\!
		\left\{
		\begin{array}{ll}
		0 \,\, \left\{\left\{V_{5}\right\}\right\}, \\  
		1 \,\, \left\{\left\{V_{9}\right\}\right\}, \\  
		\end{array}
		\right. \\	
		2000 \,\, \& \,\, I_{5}\!=\! 0\,\, \& \,\, I_{6}\!=\! 
		\left\{
		\begin{array}{ll}
		c \,\, \left\{\left\{4.8L_{5}\right\}\right\}, \\  
		f \,\, \left\{\left\{8S_{5}\right\}\right\}, \\  
		\end{array}
		\right. \\	
		2100 \,\, \& \,\, I_{5}\!=\!
		\left\{
		\begin{array}{ll}
		0 \,\, \left\{\left\{V_{20}\right\}\right\}, \\  
		1 \,\, \left\{\left\{V_{24}\right\}\right\}, \\  
		\end{array}
		\right. \\	
		2101 \,\, \left\{\left\{4.8L_{4}\right\}\right\}, \\  				
		2111 \,\, \& \,\, I_{5}\!=\!
		\left\{
		\begin{array}{ll}
		0 \,\, \left\{\left\{7S_{1}\right\}\right\}, \\  
		1 \,\, \left\{\left\{7S_{2}\right\}\right\}, \\  
		\end{array}
		\right. \\	
		2121 \,\, \& \,\, I_{5}\!=\!
		\left\{
		\begin{array}{ll}
		0 \,\, \left\{\left\{V_{12}\right\}\right\}, \\  
		1 \,\, \left\{\left\{V_{17}\right\}\right\}, \\  
		\end{array}
		\right. \\	
		3101 \,\, \left\{\left\{4S_{3}\right\}\right\}, \\  		
		3111 \,\, \left\{\left\{8S_{4}\right\}\right\}, \\  		
		4111 \,\, \& \,\, I_{5}\!=\!
		\left\{
		\begin{array}{ll}
		0 \,\, \& \,\, I_{6}\!=\! f  \,\, \& \,\, I_{7}\!=\!
		\left\{
		\begin{array}{ll}
		(2, 2) \,\, \left\{\left\{V_{14}\right\}\right\}, \\  
		(3, 1) \,\, \left\{\left\{V_{15}\right\}\right\}, \\  
		\end{array}
		\right. \\	
		1 \,\, \left\{\left\{V_{16}\right\}\right\}, \\  
		\end{array}
		\right. \\	
		\end{array}
		\right. \\
		\infty\,\, \& \,\, I_{4}\!=\! 0\,\, \& \,\, I_{5}\!=\! 0\,\, \& \,\, I_{6}\!=\!
		\left\{
		\begin{array}{ll}
		f \,\, \left\{\left\{4.5L_{1}\right\}\right\}, \\  
		c \,\, \left\{\left\{P_{4}\right\}\right\}, \\  
		\end{array}
		\right. \\	
		
		\end{array}
		\right. \\
		\end{array}
		\right.
		\end{array}
		\right.
		\]
	\end{center}
\end{table}

\begin{table}\caption{\small Topological equivalences for the family $\QESB$}\label{tab:top-equiv-QESB-pt1}
{\small	\begin{center}
		\begin{tabular}{ccccccc}
			\hline
			Presented & Identical       & Finite      & Finite       & Finite  &  Possessing       &              \\
			phase     & under           & antisaddle  & antisaddle   & weak    &  invariant curve   &  Other reasons\\
			portrait  & perturbations   & focus       & node--focus  & point   &           (no separatrix)        &              \\
			\hline
			\multirow{1}{*}{$V_{1}$}& $V_{2}$, $V_{3}$, $V_{4}$& & & & & \\  
			& & & & $3S_{1}$, $3S_{2}$ & $4S_{1}$, $8S_{1}$& \\  
			& & & & $3.3L_{1}$  & $4.8L_{1}$, $4.8L_{2}$& \\  
			& & & & &$P_{1}$ & \\  
			\hline
			\multirow{1}{*}{$V_{5}$}& $V_{8}$& $V_{6}$, $V_{7}$ & & & &$V_{10}^{(1)}$, $V_{11}^{(1)}$ \\  
			& & $0S_{2}$, $6S_{1}$ & & &$8S_{3}$ &$0S_{1}^{(2)}$, $0S_{4}^{(1)}$, $0S_{5}^{(2)}$  \\  
			& & $6S_{2}$, $8S_{2}$ & & & & $3S_{3}^{(1)}$, $6S_{3}^{(1)}$ \\  
			& & & $0.6L_{1}$, $6.8L_{1}$ & & & $0.3L_{1}^{(1)}$, $0.6L_{2}^{(1)}$ \\  
			\hline
			\multirow{1}{*}{$V_{9}$}& & & & & & \\  
			& & & & & & $0S_{3}^{(2)}$ \\  
			\hline
			\multirow{1}{*}{$V_{12}$}& &$V_{13}$ & & & & $V_{18}^{(1)}$, $V_{19}^{(1)}$\\  
			& & & $6S_{4}$ & & &$3S_{4}^{(1)}$, $6S_{5}^{(1)}$ \\  
			\hline
			\multirow{1}{*}{$V_{14}$}& & & & & & \\  
			\hline
			\multirow{1}{*}{$V_{15}$}& & & & & & \\  
			\hline
			\multirow{1}{*}{$V_{16}$}& & & & & & \\  
			\hline
			\multirow{1}{*}{$V_{17}$}& & & & & & \\  
			\hline
			\multirow{1}{*}{$V_{20}$}&$V_{22}$ & $V_{21}$, $V_{23}$ & & & &$V_{25}^{(1)}$, $V_{26}^{(1)}$\ \\  
			& &$4S_{4}$ & $6S_{6}$, $6S_{7}$ & &$4S_{5}$ & $3S_{5}^{(1)}$, $6S_{8}^{(1)}$ \\  
			& & &$4.6L_{1}$ & & & \\  
			\hline
			\multirow{1}{*}{$V_{24}$}& & & & & & \\  
			\hline
			\multirow{1}{*}{$1S_{1}$}& $1S_{2}$ & & & & & $1S_{3}^{(1)}$\\  
			&$1.8L_{1}$ & & & & & \\  
			\hline
			\multirow{1}{*}{$4S_{2}$}& & & & & & \\  
			& & & & & $0.4L_{1}$ & \\  
			\hline
			\multirow{1}{*}{$4S_{3}$}& & & & & & \\  
			\hline
			\multirow{1}{*}{$5S_{1}$}& &$5S_{2}$ & & & & $5S_{5}^{(1)}$, $5S_{6}^{(1)}$ \\  
			& & $5.8L_{1}$ & $5.6L_{1}$ & & & $3.5L_{1}^{(1)}$, $5.6L_{2}^{(1)}$ \\  
			\hline
			\multirow{1}{*}{$5S_{3}$}& & & & & & \\  
			\hline
			\multirow{1}{*}{$5S_{4}$}& & & & & & \\  
			\hline
			\multirow{1}{*}{$7S_{1}$}& & & & & & \\  
			\hline
			\multirow{1}{*}{$7S_{2}$}& & & & & & \\  
			\hline
			\multirow{1}{*}{$8S_{4}$}& & & & & & \\  
			\hline
			\multirow{1}{*}{$8S_{5}$}& & & & & & \\  
			\hline
			\multirow{1}{*}{$1.1L_{1}$}& & & & & & \\  
			& $P_{2}$& & & & & \\  
			\hline
			\multirow{1}{*}{$4.5L_{1}$}& & & & & & \\  
			\hline
			\multirow{1}{*}{$4.8L_{3}$}& & & & & & \\  
			& & & & & $P_{3}$& \\  
			\hline
			\multirow{1}{*}{$4.8L_{4}$}& & & & & & \\  
			\hline
			\multirow{1}{*}{$4.8L_{5}$}& & & & & & \\  
			\hline
			\multirow{1}{*}{$5.8L_{2}$}& & & & & & \\  
			\hline
			\multirow{1}{*}{$P_{4}$}& & & & & & \\  
			\hline
		\end{tabular}
	\end{center}}
\end{table}

\subsection{The bifurcation diagram of family $\QESC$}\label{subsec:bd-QES-C}

In this section we present the study of the bifurcation diagram of family $\QESC$, given by 
normal form~\eqref{eq:nf-QES-C}. Note that this family depends on the parameters 
$g\in\mathbb{R}\setminus\{0\}$ (in order to have nondegenerate systems) 
and $\ell\in\mathbb{R}^+\cup\{0\}$ (due to the symmetry we proved before). 
Here we shall consider the bifurcation diagram formed by points with 
Cartesian coordinates $(g,\ell)$ with $\ell\geq0$.

For systems~\eqref{eq:nf-QES-C}, computations show that 
$$\mu_0=\mathbf{D}=\mathbf{P}=0, \quad \mathbf{R}=48g^4x^2,$$
therefore by \cite[Table 5.1]{Artes-Llibre-Schlomiuk-Vulpe-2021a}, for $g\ne0$ systems~\eqref{eq:nf-QES-C} 
possess exactly one real triple finite singular point. 

Now we present the value of the algebraic invariants and T--comitants (with respect to systems~\eqref{eq:nf-QES-C})
which are relevant in our study. Since we have a two--parameter bifurcation diagram, such algebraic tools 
shall give us bifurcation curves.

\medskip

\noindent \textbf{Bifurcation curve in $\mathbb{R}^2$ due to degeneracy of the system}

From the normal form under consideration, calculation show that 
$$\mu_0=0, \quad \mu_1=4g^2x, \quad \mu_2=\mu_3=\mu_4=0.$$
Then by \cite[Lemma 5.2]{Artes-Llibre-Schlomiuk-Vulpe-2021a}, for $g=0$ systems~\eqref{eq:nf-QES-C} 
are reduced to
$$\begin{aligned}
&x^\prime=0,\\
&y^\prime=\ell y+2xy+\ell x^2,
\end{aligned}$$
they are degenerate and therefore we define the bifurcation straight line
$$({\cal L}_{1})\!:g=0.$$

According to \cite[Diagram 12.1]{Artes-Llibre-Schlomiuk-Vulpe-2021a}, for these systems we calculate
$$\eta=0, \quad \widetilde{M}=-32x^2,  \quad \kappa=\widetilde{K}=\widetilde{L}=\kappa_1={K}_1=0,$$
and 
$${L}_2=6\ell^3x^4.$$
As in the case of family $\QESB$, here we also have that ${L}_2=0$ is equivalent to $\ell=0$. So, according to the mentioned reference, 
for $\ell\neq0$ we have a hyperbola filled up with singular points,
and for $\ell=0$ (i.e. at $P_{3}=(g,\ell)=(0,0)$) we have two real straight lines (filled up with singular points) intersecting at a finite point.

\medskip

\noindent \textbf{Bifurcation curves in $\mathbb{R}^2$ due to the presence of invariant algebraic curves}

\medskip

\noindent {\bf (${\cal L}_{4}$)} This curve contains the points of the parameter space in which there appear 
invariant straight lines (see Lemma~\ref{lemma:S4-inv-curves-QES-C}). For systems~\eqref{eq:nf-QES-C} 
we compute the polynomial invariant $B_1$ and we define curve
$$\begin{aligned}
({\cal L}_{4})\!: & \ 8g^6\ell^3=0.
\end{aligned}$$

\medskip

\noindent {\bf (${\cal L}_{8}$)} This curve contains the points of the parameter space in which there appear 
invariant parabolas. According to the conditions stated in Lemma~\ref{lemma:S4-inv-curves-QES-C} we define 
this curve by 
$$\begin{aligned}
({\cal L}_{8})\!: & \ \ell=0.
\end{aligned}$$

We point out that for $g\ne0$, the bifurcation curves (${\cal L}_{4}$) and (${\cal L}_{8}$)
coincide.

\medskip

\noindent \textbf{Bifurcation curve due to multiplicities of infinite singularities}

\medskip

\noindent {\bf (${\cal L}_{5}$)} This is the bifurcation curve due to multiplicity of infinite singular points.  
According to \cite[Lemma 5.5]{Artes-Llibre-Schlomiuk-Vulpe-2021a}, for this family we calculate
$$\eta=0, \quad \widetilde{M}=-8(g-2)^2x^2, \quad C_2=x^2\left[-\ell x+(g-2)y\right],$$
and we observe that along
$$\begin{aligned}
(\mathcal{ L}_{5})\!:& \ g-2=0,
\end{aligned}$$
we have a coalescence of infinite singular points. In addition, due to the mentioned result, along the straight line
$g=2$ the phase portrait corresponding to $\ell=0$ (i.e. the phase portrait corresponding to $P_{1}=(g,\ell)=(2,0)$) 
have the line at infinity filled up with singular points.

\medskip

\noindent \textbf{Bifurcation curve in $\mathbb{R}^2$ due to the infinite elliptic--saddle}

\medskip

\noindent {\bf (${\cal L}_{0}$)} Along the straight line $g=1$ the corresponding phase portraits possess an infinite 
singularity of the type $\widehat{\!{1\choose 2}\!\!}\ E-H$. Due to results on \cite{Artes-Llibre-Schlomiuk-Vulpe-2021a} 
we compute the comitant
$$\widetilde{N}=4 (g-1)x^2$$
and we define 
$$({\cal L}_{0})\!:g-1=0.$$

\medskip

The bifurcation curves listed previously are all algebraic and they, except $({\cal L}_{4})$ and $({\cal L}_{8})$, 
are the bifurcation curves of singularities of systems~\eqref{eq:nf-QES-C} in the parameter space.

Here we shall plot these bifurcation curves in a plane with Cartesian coordinates 
$(g,\ell)$, where the horizontal line is the $g$--axis and $\ell\geq0$. 

\begin{remark} We highlight that since for $g\ne0$ the curve $({\cal L}_{4})$ coincides with $({\cal L}_{8})$, we decided 
	to plot only curve $({\cal L}_{8})$, using the cyan color. In addition, $({\cal L}_{0})$ is drawn in brown, 
	$({\cal L}_{1})$ is drawn in green, and $({\cal L}_{5})$ is drawn in red.
\end{remark}

So, in summary we have the following (distinct) bifurcation curves:
$$\begin{aligned}
(\mathcal{ L}_{0})\!:& \ g-1=0,\\
(\mathcal{ L}_{1})\!:& \ g=0,\\
(\mathcal{ L}_{5})\!:& \ g-2=0,\\
(\mathcal{ L}_{8})\!:& \ \ell=0.\\
\end{aligned}$$
And, as our bifurcation diagram is given by $\{(g,\ell)\in\mathbb{R}^2; \ell\geq0\}$, it is clear
that (in such a set) we have to consider only the curves $g=0, g=1, g=2$, and $\ell=0$, and also the intersection 
among them, i.e. the points $P_{1}=(g,\ell)=(2,0)$, $P_{2}=(g,\ell)=(1,0)$, and $P_{3}=(g,\ell)=(0,0)$.

\begin{figure}[h!]
	\centering
	\includegraphics[width=0.5\textwidth]{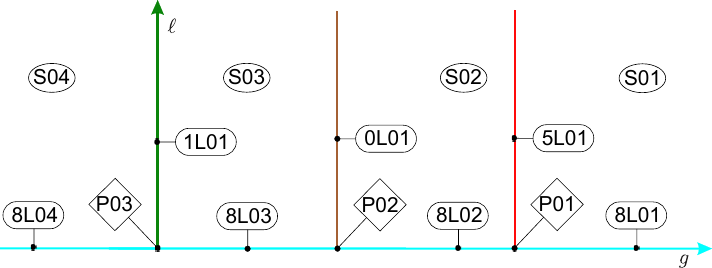} 
	\caption{\small \label{fig:slice-QES-C} Parameter space}
\end{figure}

In  Fig.~\ref{fig:slice-QES-C} we present the bifurcation diagram completely labeled. In such a 
figure we denote an open region by S0i, where i is a number, a bifurcation curve 
$(\mathcal{L}_{j})$ is labeled as jL0k, $k\in\mathbb{N}$, and a point is denoted as in
the previous sections. Moreover, we denote the $\ell$--axis (which represents the degenerate
set) with a dashed and thin black straight line.

From the study of this bifurcation diagram, we obtain phase portraits possessing different types 
of triple finite singular points. In fact, from \cite[Table 6.2]{Artes-Llibre-Schlomiuk-Vulpe-2021a} 
we calculate
$$\widetilde{K}=4gx^2.$$
For nondegenerate systems (i.e. $g\ne0$), this comitant can be positive or negative, depending on the 
sign of the parameter $g$. In what follows we present the different types of triple finite singularities we 
obtained in the study of the bifurcation diagram under consideration.
\begin{enumerate}
	\item If $g>0$ then $\widetilde{K}>0$ and, from the mentioned table we compute 
	$$G_{10}=g^3\ell^3.$$
	Since $g\ne0$, we have that $\text{sign}(G_{10})=\text{sign}(\ell)$ and, from the mentioned table 
	can have two possibilities:
	\begin{itemize}
		\item If $\ell\ne0$ we have a finite semi--elemental triple node $\bar{n}_{(3)}$;
		\item If $\ell=0$ we have a finite nilpotent elliptic--saddle $\widehat{es}_{(3)}$.
	\end{itemize}
	\item If $g<0$ then $\widetilde{K}<0$. Now, from \cite[Diagram 10.2]{Artes-Llibre-Schlomiuk-Vulpe-2021a} 
	we calculate
	$$\kappa=0, \quad \mathcal{F}_1=6g^2\ell.$$
	Since $g\ne0$, we have that $\text{sign}(\mathcal{F}_{1})=\text{sign}(\ell)$ and, from the diagram 
	under consideration again we come across two possibilities:
	\begin{itemize}
		\item If $\ell\ne0$ we have a finite semi--elemental triple saddle $\bar{s}_{(3)}$;
		\item If $\ell=0$ we have a finite nilpotent triple saddle $\widehat{s}_{(3)}$.
	\end{itemize}
\end{enumerate}

By performing the study of this bifurcation diagram we observe that there is coherence among 
all the phase portraits we obtained. Moreover, we point out that in our study we have not found 
any nonalgebraic bifurcation curve and there is no need of it so to complete coherence. So we can affirm that we have described a complete bifurcation 
diagram for class $\overline{\QESC}$ modulo islands, as we discuss in Sec.~\ref{sec:islands-QESC}.

\subsubsection{Other relevant facts about the bifurcation diagram} \label{sec:islands-QESC}

The bifurcation diagram we have obtained for the class $\overline{\QESC}$ is completely coherent, i.e. in this family, 
by taking any two points in the parameter space and joining them by a continuous curve, along this curve 
the changes in phase portraits that occur when crossing the different bifurcation surfaces we mention can 
be completely explained.

Nevertheless, we cannot be sure that this bifurcation diagram is the complete bifurcation diagram for 
$\overline{\QESC}$ due to the possibility of the existence of ``islands'' inside the parts bordered by unmentioned 
bifurcation surfaces. In case they exist, these ``islands'' would not mean any modification of the nature 
of the singular points. So, on the border of these ``islands'' we could only have bifurcations due to saddle 
connections.

In case there were more bifurcation surfaces, we should still be able to join two representatives of any 
two parts of the 14 parts of $\overline{\QESC}$ found until now with a continuous curve either without crossing 
such a bifurcation surface or, in case the curve crosses it, it must do it an even number of times without 
tangencies, otherwise one must take into account the multiplicity of the tangency, so the total number 
must be even. This is why we call these potential bifurcation surfaces ``\textit{islands}''.

However, we have not found a different phase portrait which could fit in such an island. A potential ``island'' 
would be the set of parameters for which the phase portraits have a separatrix connection.

\subsubsection{Completion of the proof of Theorem~\ref{th:main-thm-QES-C}} \label{sec:invariants-QESC} 

In the bifurcation diagram we may have topologically equivalent phase portraits belonging to distinct parts 
of the parameter space. As here we have 14 distinct parts of the parameter space, to help us to identify 
or to distinguish phase portraits, we need to introduce some invariants and we actually choose integer valued, 
character and symbol invariants. Some of them were already used in \cite{Artes-Rezende-Oliveira-2013b} 
and \cite{Artes-Mota-Rezende-2021c}, but we recall them and introduce some needed ones. These invariants 
yield a classification which is easier to grasp.

\begin{definition}\label{def:invariant-I1-QESC}
	We denote by $I_{1}(S)$ a symbol from the set $\{\emptyset, \left[ \times \right], \left[)(\right]\}$ which indicates the following configuration of curves filled up with singularities, respectively: none (nondegenerate systems -- in this case all systems do not contain a curve filled up with singularities), two real straight lines intersecting at a finite point, and an hyperbola. This invariant only makes sense to distinguish the degenerate phase portrait obtained.
\end{definition}

\begin{definition}\label{def:invariant-I2-QESC}
	We denote by $I_{2}(S)$ the sum of the indices of the isolated real finite singular points.
\end{definition}

\begin{definition}\label{def:invariant-I3-QESC}
	We denote by $I_{3}(S)$ the number of real infinite singular points. We note that this number can also be infinite, which is represented by $\infty$.
\end{definition}

\begin{definition}\label{def:invariant-I4-QESC}
	For a given infinite singularity $s$ of a system $S$, let $\l_s$ be the number of global or local separatrices beginning or ending at $s$ and which do not lie on the line at infinity. We have $0\leq\l_s\leq2$. We denote by $I_{4}(S)$ the sequence of all such $\l_s$ when $s$ moves in the set of infinite singular points of the system $S$. We start the sequence at the infinite singular point which receives (or sends) the greatest number of separatrices and take the direction which yields the greatest absolute value, e.g. the values $2100$ and $2001$ for this invariant are symmetrical (and, therefore, they are the same), so we consider $2100$.
\end{definition}

\begin{definition}\label{def:invariant-I5-QESC}
	We denote by $I_{5}(S)$ an element from the set $\{y, n\}$ indicating if the phase portrait has ($y$) or has not ($n$) an infinite elliptic sector.
\end{definition}

\begin{definition}\label{def:invariant-I6-QESC} 
	We denote by $I_{6}(S)$ an element from the set $\{y, n\}$ indicating if the infinite elliptic sector is ($y$) or is not ($n$) bordered by separatrices that connect the finite elliptic--saddle and the infinite multiple point.
\end{definition}

\begin{theorem} \label{th:QESC-inv}
	Consider the class $\overline{\QESC}$ and all the phase portraits that we have obtained for this family. The values of the affine 
	invariant ${\cal I} = (I_{1}, I_{2}, I_{3}, I_{4}, I_{5}, I_{6})$ given in the diagram from Table \ref{tab:geom-classif-QESC-pt1} 
	yields a partition of these phase portraits of the class $\overline{\QESC}$.
	
	Furthermore, for each value of $\cal I$ in this diagram there corresponds a single phase portrait; i.e. $S$ and $S'$ 
	are such that ${\cal I}(S)={\cal I}(S')$, if and only if $S$ and $S'$ are topologically equivalent.
\end{theorem}

The bifurcation diagram for $\overline{\QESC}$ has 14 parts which produce twelve topologically different phase portraits as 
described in Tables~\ref{tab:geom-classif-QESC-pt1} to \ref{tab:top-equiv-QESC-pt1}. The remaining two parts 
do not produce any new phase portrait which was not included in the ten previous ones. The difference is basically 
the presence of invariant algebraic curves (lines or parabolas) which do not represent a separatrix connection or 
a presence of an infinite singularity of type $\widehat{\!{1\choose2}\!\!}\ E-H$.

The phase portraits having no graphics have been denoted surrounded by parenthesis, for example $(S_{1})$ and 
the phase portraits having two or more graphics have been denoted surrounded by $\{\!\{\ast\}\!\}$, for example 
$\{S_{2}\}$.

\begin{proof}[Proof of Theorem~\ref{th:QESC-inv}]
	The above result follows from the results in the previous sections and a careful analysis of the bifurcation diagrams given in 
	Fig.~\ref{fig:slice-QES-C}, the definition of the invariants $I_{j}$ and their explicit values for the corresponding phase portraits.
\end{proof}

In Table~\ref{tab:top-equiv-QESC-pt1} we list in the first column twelve parts with all the distinct phase portraits of Fig.~\ref{fig:pp-QES-C}. Corresponding to each part listed in column one we have in each row all parts whose phase portraits are topologically equivalent to the phase portrait appearing in column 1 of the same row. In the second column we set all the parts whose systems possess an invariant curve (straight line and/or parabola) not yielding a connection of separatrices and in the third column we put the phase portrait possessing a singularity of type $\widehat{\!{1\choose2}\!\!}\ E-H$ at infinity.

Whenever phase portraits appear in a row in a specific column, the listing is done according to the decreasing dimension of the parts where they appear, always placing the lower dimensions on lower lines.

\subsubsection{Proof of Theorem~\ref{th:main-thm-QES-C}}

The bifurcation diagram described in Sec. \ref{subsec:bd-QES-C}, plus Table~\ref{tab:geom-classif-QESC-pt1} of the geometrical invariants distinguishing the ten phase portraits, plus Table~\ref{tab:top-equiv-QESC-pt1} giving the equivalences with the remaining phase portraits lead to the proof of Theorem~\ref{th:main-thm-QES-C}.

\begin{table}\caption{\small Geometric classification for the family $\QESC$}\label{tab:geom-classif-QESC-pt1} 
	\begin{center}
		\[
		I_{1}\!=\!
		\left\{
		\begin{array}{ll}
		\left[)(\right] \,\, \left\{\left\{1L_{1}\right\}\right\}, \\  			
		\left[ \times \right] \,\, \left\{\left\{P_{3}\right\}\right\}, \\  			
		\emptyset \,\, \& \,\, I_{2}\!=\!
		\left\{
		\begin{array}{ll}
		-1 \,\, \left\{\left\{S_{4}\right\}\right\}, \\  			
		1 \,\, \& \,\, I_{3}\!=\!
		\left\{
		\begin{array}{ll}
		1 \,\, \left\{\left\{5L_{1}\right\}\right\}, \\  			
		2 \,\, \& \,\, I_{4}\!=\!
		\left\{
		\begin{array}{ll}
		1010 \,\, \left\{\left\{P_{2}\right\}\right\}, \\  				
		1110 \,\, \left\{\left\{S_{3}\right\}\right\}, \\  				
		2100 \,\, \left\{\left\{S_{1}\right\}\right\}, \\  				
		2101 \,\, \& \,\, I_{5}\!=\!
		\left\{
		\begin{array}{ll}
		n \,\, \left\{\left\{8L_{1}\right\}\right\}, \\  
		y \,\, \& \,\, I_{6}\!=\!
		\left\{
		\begin{array}{ll}
		n \,\, \left\{\left\{8L_{2}\right\}\right\}, \\  
		y \,\, \left\{\left\{8L_{3}\right\}\right\}, \\  
		\end{array}
		\right. \\
		\end{array}
		\right. \\
		2121 \,\, \left\{\left\{S_{2}\right\}\right\}, \\  
		\end{array}
		\right. \\
		\infty\,\, \left\{\left\{P_{1}\right\}\right\}, \\  
		\end{array}
		\right. \\
		\end{array} 
		\right.\\
		\end{array} 
		\right.
		\]
	\end{center}
\end{table}

\begin{table}\caption{\small Topological equivalences for the family $\QESC$}\label{tab:top-equiv-QESC-pt1}
	\begin{center}
		\begin{tabular}{ccccccc}
			\hline
			Presented &   Possessing       &  Possessing           \\
			phase     &   invariant curve   &  $\widehat{\!{1\choose2}\!\!}\ E-H$ \\
			portrait  &     (no separatrix)        &     at infinity         \\
			\hline
			\multirow{1}{*}{$S_{1}$} & & \\  
			\hline
			\multirow{1}{*}{$S_{2}$} & & \\  
			\hline
			\multirow{1}{*}{$S_{3}$} & & \\  
			& &$0L_{1}$ \\  
			\hline
			\multirow{1}{*}{$S_{4}$} & & \\  
			& $8L_{4}$ & \\  
			\hline
			\multirow{1}{*}{$1L_{1}$} & & \\  
			\hline
			\multirow{1}{*}{$5L_{1}$} & & \\  
			\hline
			\multirow{1}{*}{$8L_{1}$} & & \\  
			\hline
			\multirow{1}{*}{$8L_{2}$} & & \\  
			\hline
			\multirow{1}{*}{$8L_{3}$} & & \\  
			\hline
			\multirow{1}{*}{$P_{1}$} & & \\  
			\hline
			\multirow{1}{*}{$P_{2}$} & & \\  
			\hline
			\multirow{1}{*}{$P_{3}$} & & \\  
			\hline
		\end{tabular}
	\end{center}
\end{table}

\bigskip

\noindent\textbf{Acknowledgements.} The first author is partially supported by a MEC/FEDER grant 
number MTM 2016--77278--P and by a CICYT grant number 2017 SGR 1617. The second author 
was partially supported by Conselho Nacional de Desenvolvimento Cient\'ifico e Tecnol\'ogico (CNPq) 
grant number 166449/2020-2. The third author is partially supported by Coordena\c{c}\~ao de 
Aperfei\c{c}oamento de Pessoal de Nivel Superior - Brazil (CAPES) and by Funda\c{c}\~ao de 
Amparo \`a Pesquisa do Estado de S\~ao Paulo (FAPESP) grants 2018/21320-7 and 2019/21181-0.

\bigskip

\appendix
\section{Some incompatibilities in previous classifications}\label{ap:incomp-QES}

It is quite common that by performing the study of a bifurcation diagram that produces some specific 
types of phase portraits, the authors lose one or several phase portraits. This may happen either because 
they do not interpret correctly some of the bifurcation parts or they miss the existence of some nonalgebraic bifurcations. 

In \cite{Artes-Mota-Rezende-2021b} we have decided to start comparing our classification of 
phase portraits with already existing classifications. As we have mentioned in that occasion, we 
plan to do this section in every future work related to classification of phase portraits using normal 
forms. The aim of this study is to detect some incompatibilities in previous papers and also to help 
us look carefully our bifurcation diagram in order to do not lose any phase portrait. Such incompatibilities 
are obtained after we compare all of the phase portraits obtained in our bifurcation diagram with 
phase portraits from some previous papers which possess the same {\it topological configuration 
	of singularities}, according to Def.~1 in \cite{Artes-Llibre-Schlomiuk-Vulpe-2020a}. 

This study also allows the corresponding authors to detect possible mistakes on their works. 
There are some previous papers which are not based on normal forms, but which seek all 
topological realizable phase portraits of a certain codimension 
(see \cite{Artes-Kooij-Llibre-1998,Artes-Llibre-Rezende-2018,Artes-Oliveira-Rezende-2020b,Artes-Mota-Rezende-2021d}). 
We have also crossed results from all the consulted papers with them and no discrepancy has been found. 
Additionally, with this study we are creating a data basis containing all the obtained phase portraits, 
specially containing those phase portraits obtained in our topological studies, in order to create 
an ``encyclopedia'' of phase portraits from quadratic differential systems.

In this paper we are dealing with phase portraits possessing either an infinite nilpotent elliptic--saddle 
or an infinite nilpotent saddle. Regarding the already existing studies related to this paper, in 
\cite{Reyn-Huang-1997} the authors provide a list of phase portraits that have intersection with 
our investigation. We decided to perform a careful analysis of the phase portraits they present and 
also to compare their phase portraits with the ones we obtained.

By doing this study, we have detected some interesting phenomena and also some incompatibilities 
in the mentioned paper. We observe that there are phase portraits in \cite{Reyn-Huang-1997} 
which are topologically equivalent, and this fact allowed us to create sets of topologically equivalent 
phase portraits. In what follows we present such sets. In each set the elements (i.e. phase portraits 
from that paper) are displayed in lines, where in each line we indicate the figure of that paper 
in which the phase portrait appears, followed by the caption of that phase portrait (using the 
notation of the paper under consideration), so one can easily identify all of them in the mentioned
paper.

\begin{center}

$\left\{\!\!\begin{array}{c}
\text{FIGURE 10.1:} \ \gamma<0 \ \text{and} \ \mu=0,\\
\text{FIGURE 10.1:} \ \gamma<0 \ \text{and} \ \mu>0,\\
\text{FIGURE 10.5:} \ \delta>0 \ \text{and} \ \mu=0,\\
\text{FIGURE 10.5:} \ \delta>0 \ \text{and} \ \mu>0,\\
\text{FIGURE 10.5:} \ \delta=0 \ \text{and} \ \mu=0,\\
\text{FIGURE 10.5:} \ \delta=0 \ \text{and} \ \mu>0,\\
\text{FIGURE 11.1:} \ \mu>0 \ \text{and} \ \gamma<-\frac{1}{4\mu},\\
\text{FIGURE 11.5c:} \ \mu>-\frac{1}{4\gamma} \ \text{and} \ \gamma<0\\
\end{array}\!\!\right\}\!,$

\vspace{\baselineskip}

$\left\{\!\!\begin{array}{c}
\text{FIGURE 10.1:} \ \gamma=2 \ \text{and} \ \mu=0,\\
\text{FIGURE 10.2:} \ \gamma=2 \ \text{and} \ \mu=\delta=0\\
\end{array}\!\!\right\}\!,$

\vspace{\baselineskip}

$\left\{\!\!\begin{array}{c}
\text{FIGURE 10.1:} \ \gamma>2 \ \text{and} \ \mu=0,\\
\text{FIGURE 10.2:} \ \gamma>2 \ \text{and} \ \mu=\delta=0\\
\end{array}\!\!\right\}\!,$

\vspace{\baselineskip}

$\left\{\!\!\begin{array}{c}
\text{FIGURE 10.1:} \ 0<\gamma<1 \ \text{and} \ \mu=0,\\
\text{FIGURE 10.4:} \ \mu=\delta=0\\
\end{array}\!\!\right\}\!,$

\vspace{\baselineskip}

$\left\{\!\!\begin{array}{c}
\text{FIGURE 10.1:} \ 1<\gamma<2 \ \text{and} \ \mu=0,\\
\text{FIGURE 10.3:} \ \mu=\delta=0\\
\end{array}\!\!\right\}\!,$

\vspace{\baselineskip}

$\left\{\!\!\begin{array}{c}
\text{FIGURE 10.1:} \ 2<\gamma \ \text{and} \ \mu<0,\\
\text{FIGURE 10.1:} \ \gamma=2 \ \text{and} \ \mu<0,\\
\text{FIGURE 10.2:} \ \kappa=-\infty\\
\end{array}\!\!\right\}\!,$

\vspace{\baselineskip}

$\left\{\!\!\begin{array}{c}
\text{FIGURE 10.2:} \ \kappa=0,\\
\text{FIGURE 10.2:} \ -\infty<\kappa<0,\\
\text{FIGURE 11.1:} \ \gamma>2 \ \text{and} \ \mu<-\frac{1}{4\gamma},\\
\text{FIGURE 11.2:} \ d>d_1(m;g),\\
\text{FIGURE 11.2:} \ d^\varepsilon-1\\
\end{array}\!\!\right\}\!,$

\vspace{\baselineskip}

$\left\{\!\!\begin{array}{c}
\text{FIGURE 10.1:} \ 0<\gamma<1 \ \text{and} \ \mu<0,\\
\text{FIGURE 10.4:} \ \kappa=-\infty\\
\end{array}\!\!\right\}\!,$

\vspace{\baselineskip}

$\left\{\!\!\begin{array}{c}
\text{FIGURE 10.4:} \ \kappa=0,\\
\text{FIGURE 10.4:} \ -\infty<\kappa<0,\\
\text{FIGURE 11.1:} \ 0<\gamma<1 \ \text{and} \ \mu<-\frac{1}{4\gamma},\\
\text{FIGURE 11.4:} \ \delta>\delta_1(\mu;\gamma),\\
\text{FIGURE 11.4:} \ \delta\leq-1\\
\end{array}\!\!\right\}\!,$

\vspace{\baselineskip}

$\left\{\!\!\begin{array}{c}
\text{FIGURE 10.1:} \ 1<\gamma<2 \ \text{and} \ \mu<0,\\
\text{FIGURE 10.3:} \ \kappa=-\infty\\
\end{array}\!\!\right\}\!,$

\vspace{\baselineskip}

$\left\{\!\!\begin{array}{c}
\text{FIGURE 10.3:} \ \kappa=\kappa_5,\\
\text{FIGURE 11.1:} \ \text{u.s.c.},\\
\text{FIGURE 11.3c:} \ \delta=\delta_4(\mu;\gamma),\\
\text{FIGURE 11.3c:} \ \delta=\delta_7(\mu;\gamma)<-1\\
\end{array}\!\!\right\}\!,$

\vspace{\baselineskip}

$\left\{\!\!\begin{array}{c}
	\text{FIGURE 10.3:} \ \kappa=0,\\
	\text{FIGURE 10.3:} \ \kappa_5<\kappa<0,\\
	\text{FIGURE 11.1:} \ 1<\gamma<2 \ \text{and} \ \text{u.s.c.}<\mu<-\frac{1}{4\gamma},\\
	\text{FIGURE 11.3c:} \ \delta>\delta_4(\mu;\gamma),\\
	\text{FIGURE 11.3c:} \ \delta<\delta^{\ast\ast}\\
	\end{array}\!\!\right\}\!,$
	
	\vspace{\baselineskip}

$\left\{\!\!\begin{array}{c}
\text{FIGURE 10.3:} \ -\infty<\kappa<\kappa_5,\\
\text{FIGURE 11.1:} \ 1<\gamma<2 \ \text{and} \ \mu<\text{u.s.c.},\\
\text{FIGURE 11.3c:} \ \delta_5(\mu;\gamma)<\delta<\delta_4(\mu;\gamma),\\
\text{FIGURE 11.3c:} \ \delta_7(\mu;\gamma)<\delta\leq-1\\
\end{array}\!\!\right\}\!,$

\vspace{\baselineskip}

$\left\{\!\!\begin{array}{c}
	\text{FIGURE 11.1:} \ \gamma<0 \ \text{and} \ \mu=0,\\
	\text{FIGURE 11.5c:} \ \mu=0, \gamma<0, \ \text{and} \ \delta>\delta_1(0;\gamma)\\
	\end{array}\!\!\right\}\!,$
	
	\vspace{\baselineskip}

$\left\{\!\!\begin{array}{c}
\text{FIGURE 11.1:} \ \gamma>2 \ \text{and} \ \mu=0,\\
\text{FIGURE 11.2:} \ d=0\\
\end{array}\!\!\right\}\!,$

\vspace{\baselineskip}

$\left\{\!\!\begin{array}{c}
\text{FIGURE 11.1:} \ 0<\gamma<1 \ \text{and} \ \mu=0,\\
\text{FIGURE 11.4:} \ \delta_1^\ast(0;\gamma)<\delta<\delta_1(0;\gamma)\\
\end{array}\!\!\right\}\!,$

\vspace{\baselineskip}

$\left\{\!\!\begin{array}{c}
\text{FIGURE 11.1:} \ 1<\gamma<2 \ \text{and} \ \mu=0,\\
\text{FIGURE 11.3d:} \ \delta=0\\
\end{array}\!\!\right\}\!,$

\vspace{\baselineskip}

$\left\{\!\!\begin{array}{c}
\text{FIGURE 10.1:} \ \gamma<0 \ \text{and} \ \mu<0,\\
\text{FIGURE 10.5:} \ \delta=0 \ \text{and} \ \mu<0\\
\end{array}\!\!\right\}\!,$

\vspace{\baselineskip}

$\left\{\!\!\begin{array}{c}
	\text{FIGURE 10.5:} \ \delta>0 \ \text{and} \ \mu<0,\\
	\text{FIGURE 11.1:} \ \gamma<0 \ \text{and} \ \mu<0,\\
	\text{FIGURE 11.5c:} \ \mu<0, -2<\gamma<0, \ \text{and} \ \delta>\delta_1(\mu;\gamma),\\
	\text{FIGURE 11.5c:} \ \mu<0, -2<\gamma<0, \ \text{and} \ \delta>\delta_3(\mu;\gamma)\\
	\end{array}\!\!\right\}\!,$
	
	\vspace{\baselineskip}

$\left\{\!\!\begin{array}{c}
	\text{FIGURE 11.5c:} \ \mu<0, -2<\gamma<0, \ \text{and} \ \delta=\delta_2(\mu;\gamma),\\
	\text{FIGURE 11.5c:} \ \mu<0, \gamma>-2, \ \text{and} \ \delta=\delta_2(\mu;\gamma)\\
	\end{array}\!\!\right\}\!,$
	
	\vspace{\baselineskip}

$\left\{\!\!\begin{array}{c}
\text{FIGURE 10.1:} \ 2<\gamma \ \text{and} \ \mu>0,\\
\text{FIGURE 10.1:} \ \gamma=2 \ \text{and} \ \mu>0,\\
\text{FIGURE 10.2:} \ \kappa_1<\kappa\leq\infty,\\
\text{FIGURE 11.1:} \ 2<\gamma \ \text{and} \ \mu>0\\
\end{array}\!\!\right\}\!,$

\vspace{\baselineskip}

$\left\{\!\!\begin{array}{c}
\text{FIGURE 10.1:} \ 0<\gamma<1 \ \text{and} \ \mu>0,\\
\text{FIGURE 10.4:} \ \kappa_1<\kappa\leq\infty,\\
\text{FIGURE 11.1:} \ 0<\gamma<1 \ \text{and} \ \mu>0\\
\end{array}\!\!\right\}\!,$ and

\vspace{\baselineskip}

$\left\{\!\!\begin{array}{c}
\text{FIGURE 10.1:} \ 1<\gamma<2 \ \text{and} \ \mu>0,\\
\text{FIGURE 10.3:} \ \kappa_1<\kappa\leq\infty,\\
\text{FIGURE 11.1:} \ 1<\gamma<2 \ \text{and} \ \mu>0\\
\end{array}\!\!\right\}\!.$

\end{center}

We also have a correspondence between the phase portraits from the paper
under consideration and the phase portraits we obtained in our study (the reader
may remember Table~\ref{tab:equiv-pp-famb-famc} of topological equivalence 
between phase portraits from families $\QESB$ and $\QESC$). 
And most important, there is not a single phase portrait in  \cite{Reyn-Huang-1997}  
which is absent in our study. In case that happened and the phase portrait were 
confirmed to exist, we would have a gap in this study.

\begin{table}\caption{\small Correspondence between phase portraits from \cite{Reyn-Huang-1997} and
		phase portraits obtained from the studies of the bifurcation diagrams of families $\QESA$, $\QESB$, and $\QESC$. 
		In the first column we refer to the figures from \cite{Reyn-Huang-1997}, in the second column we list the phase 
		portraits which appear in that figures (using the notation of that paper), and in the third column we indicate the 
		corresponding phase portrait we obtained from the study of families $\QESA$, $\QESB$, or $\QESC$}\label{tab:cor-RH-we-pt1}
	\begin{center}
		{\def\arraystretch{1.5}
			\begin{tabular}{c|cc|c}
				\hline
				\multirow{2}{*}{\bf FIGURE}    & \multicolumn{2}{c|}{\multirow{2}{*}{\bf Phase portrait}} & \textbf{Correspondent in families}   \\
				\textbf{\cite{Reyn-Huang-1997}}		& 																	&										& $\QESA$, $\QESB$, or $\QESC$ \\ 
				\hline
				\multirow{15}{*}{10.1}  &  \multirow{3}{*}{$\gamma>2$}     & \multicolumn{1}{|c|}{$\mu<0$}  &  $4.8L_{3}-\QESB$ \\ \cline{3-4}
				&                                                        &  \multicolumn{1}{|c|}{$\mu=0$}  &  $8L_{3}-\QESC$ \\ \cline{3-4}
				&                                                        &  \multicolumn{1}{|c|}{$\mu>0$}  &  $V_{101}-\QESA$ \\ \cline{2-4}
				&  \multirow{3}{*}{$\gamma=2$}     & \multicolumn{1}{|c|}{$\mu<0$}   &  $4.8L_{3}-\QESB$ \\ \cline{3-4}
				&                                                        &  \multicolumn{1}{|c|}{$\mu=0$}  &  $P_{2}-\QESC$ \\ \cline{3-4}
				&                                                        &  \multicolumn{1}{|c|}{$\mu>0$}  &  $V_{101}-\QESA$ \\ \cline{2-4}
				&  \multirow{3}{*}{$1<\gamma<2$} &  \multicolumn{1}{|c|}{$\mu<0$}  &  $4.8L_{4}-\QESB$ \\ \cline{3-4}
				&                                                        &  \multicolumn{1}{|c|}{$\mu=0$}  &  $8L_{2}-\QESC$ \\ \cline{3-4}
				&                                                        &  \multicolumn{1}{|c|}{$\mu>0$}  &  $V_{188}-\QESA$ \\ \cline{2-4}
				&  \multirow{3}{*}{$0<\gamma<1$} &  \multicolumn{1}{|c|}{$\mu<0$}  &  $4.8L_{5}-\QESB$ \\ \cline{3-4}
				&                                                        &  \multicolumn{1}{|c|}{$\mu=0$}  &  $8L_{1}-\QESC$ \\ \cline{3-4}
				&                                                        &  \multicolumn{1}{|c|}{$\mu>0$}  &  $V_{240}-\QESA$ \\ \cline{2-4}
				&  \multirow{3}{*}{$\gamma<0$}    &  \multicolumn{1}{|c|}{$\mu<0$}  &  $4.8L_{2}-\QESA$ \\ \cline{3-4}
				&                                                        &  \multicolumn{1}{|c|}{$\mu=0$}  &  $V_{1}-\QESB$ \\ \cline{3-4}
				&                                                        &  \multicolumn{1}{|c|}{$\mu>0$}  &  $V_{1}-\QESB$ \\ 
				
				\hline
				\multirow{11}{*}{10.2}&  \multicolumn{2}{c|}{$\kappa=0$}                                 &  $V_{5}-\QESB$ \\ \cline{2-4}
				&   \multicolumn{2}{c|}{$0<\kappa\leq\kappa_3$}            &  $V_{89}-\QESA$ \\ \cline{2-4}
				&   \multicolumn{2}{c|}{$\kappa_3<\kappa<\kappa_2$}  &  $V_{91}-\QESA$ \\ \cline{2-4}
				&   \multicolumn{2}{c|}{$-\infty<\kappa<0$}                     &  $V_{5}-\QESB$ \\ \cline{2-4}
				&   \multicolumn{2}{c|}{$\kappa=-\infty$}                         &  $4.8L_{3}-\QESB$ \\ \cline{2-4}
				&   \multicolumn{2}{c|}{$\kappa=\kappa_2$}                   &  $7S_{7}-\QESA$ \\ \cline{2-4}
				&   \multicolumn{2}{c|}{$\kappa_2<\kappa<\kappa_1$}  &  $V_{94}-\QESA$ \\ \cline{2-4}
				&   \multicolumn{2}{c|}{$\kappa=\kappa_1$}  								&  $4S_{34}-\QESA$ \\ \cline{2-4}
				&   \multicolumn{2}{c|}{$\kappa_1<\kappa\leq\infty$} 		&  $V_{101}-\QESA$ \\ \cline{2-4}
				&  \multirow{2}{*}{$\mu=\delta=0$}	& \multicolumn{1}{|c|}{$\gamma=2$}   &  $P_{2}-\QESC$ \\ \cline{3-4}
				&                                                        			& \multicolumn{1}{|c|}{$\gamma>2$}  &  $8L_{3}-\QESC$ \\
				\hline
			\end{tabular}
		}
	\end{center}
\end{table}

\begin{table}\caption{\small Continuation of Table~\ref{tab:cor-RH-we-pt1}}\label{tab:cor-RH-we-pt2}
	\begin{center}
		{\def\arraystretch{1.5}
			\begin{tabular}{c|cc|c}
				\hline
				\multirow{2}{*}{\bf FIGURE}    & \multicolumn{2}{c|}{\multirow{2}{*}{\bf Phase portrait}} & \textbf{Correspondent in families}   \\
				\textbf{\cite{Reyn-Huang-1997}}							& 																	&										& $\QESA$, $\QESB$, or $\QESC$ \\ 
				\hline
				\multirow{14}{*}{10.3}&  \multicolumn{2}{c|}{$\kappa=0$}                                 &  $V_{12}-\QESB$ \\ \cline{2-4}
				&   \multicolumn{2}{c|}{$0<\kappa\leq\kappa_4$}            &  $V_{168}-\QESA$ \\ \cline{2-4}
				&   \multicolumn{2}{c|}{$\kappa_4<\kappa<\kappa_3$}   &  $V_{170}-\QESA$ \\ \cline{2-4}
				&   \multicolumn{2}{c|}{$\kappa=\kappa_3$}                    &  $7S_{11}-\QESA$ \\ \cline{2-4}
				&   \multicolumn{2}{c|}{$\kappa_3<\kappa<\kappa_2$}   &  $V_{173}-\QESA$ \\ \cline{2-4}
				&   \multicolumn{2}{c|}{$\kappa_5<\kappa<0$}               &  $V_{12}-\QESB$ \\ \cline{2-4}
				&   \multicolumn{2}{c|}{$\kappa=\kappa_5$}                   &  $7S_{1}-\QESB$ \\ \cline{2-4}
				&   \multicolumn{2}{c|}{$-\infty<\kappa<\kappa_5$}  			&  $V_{14}-\QESB$ \\ \cline{2-4}
				&   \multicolumn{2}{c|}{$\kappa=\kappa_2$} 		              &  $8S_{77}-\QESA$ \\ \cline{2-4}
				&   \multicolumn{2}{c|}{$\kappa_2<\kappa<\kappa_1$}  &  $V_{176}-\QESA$ \\ \cline{2-4}
				&   \multicolumn{2}{c|}{$\kappa=\kappa_1$}                  &  $4S_{59}-\QESA$ \\ \cline{2-4}
				&   \multicolumn{2}{c|}{$\kappa=-\infty$}                        &  $4.8L_{4}-\QESB$ \\ \cline{2-4}
				&   \multicolumn{2}{c|}{$\mu=\delta=0$}                         &  $8L_{2}-\QESC$ \\ \cline{2-4}
				&   \multicolumn{2}{c|}{$\kappa_1<\kappa\leq\infty$}      &  $V_{188}-\QESA$ \\ 
				\hline
				\multirow{10}{*}{10.4}&  \multicolumn{2}{c|}{$\kappa=0$}                                 &  $V_{20}-\QESB$ \\ \cline{2-4}
				&   \multicolumn{2}{c|}{$0<\kappa\leq\kappa_3$}            &  $V_{233}-\QESA$ \\ \cline{2-4}
				&   \multicolumn{2}{c|}{$\kappa_3<\kappa<\kappa_2$}   &  $V_{235}-\QESA$ \\ \cline{2-4}
				&   \multicolumn{2}{c|}{$-\infty<\kappa<0$}                     &  $V_{20}-\QESB$ \\ \cline{2-4}
				&   \multicolumn{2}{c|}{$\kappa=-\infty$}                         &  $4.8L_{5}-\QESB$ \\ \cline{2-4}
				&   \multicolumn{2}{c|}{$\mu=\delta=0$}             					   &  $8L_{1}-\QESC$ \\ \cline{2-4}
				&   \multicolumn{2}{c|}{$\kappa=\kappa_2$}                   &  $7S_{15}-\QESA$ \\ \cline{2-4}
				&   \multicolumn{2}{c|}{$\kappa_2<\kappa<\kappa_1$}	&  $V_{238}-\QESA$ \\ \cline{2-4}
				&   \multicolumn{2}{c|}{$\kappa=\kappa_1$} 		              &  $8S_{99}-\QESA$ \\ \cline{2-4}
				&   \multicolumn{2}{c|}{$\kappa_1<\kappa\leq\infty$}      &  $V_{240}-\QESA$ \\ 
				\hline
				\multirow{6}{*}{10.5}  &  \multirow{3}{*}{$\delta>0$}        & \multicolumn{1}{|c|}{$\mu<0$}  &  $V_{1}-\QESA$ \\ \cline{3-4}
				&                                                        &  \multicolumn{1}{|c|}{$\mu=0$}  &  $V_{1}-\QESB$ \\ \cline{3-4}
				&                                                        &  \multicolumn{1}{|c|}{$\mu>0$}  &  $V_{1}-\QESB$ \\ \cline{2-4}
				&  \multirow{3}{*}{$\delta=0$}       & \multicolumn{1}{|c|}{$\mu<0$}   &  $4.8L_{2}-\QESA$ \\ \cline{3-4}
				&                                                        &  \multicolumn{1}{|c|}{$\mu=0$}  &  $V_{1}-\QESB$ \\ \cline{3-4}
				&                                                        &  \multicolumn{1}{|c|}{$\mu>0$}  &  $V_{1}-\QESB$ \\ 
				\hline
			\end{tabular}
		}
	\end{center}
\end{table}

\begin{table}\caption{\small Continuation of Table~\ref{tab:cor-RH-we-pt2}}\label{tab:cor-RH-we-pt3}
	\begin{center}
		{\def\arraystretch{1.5}
			\begin{tabular}{c|cc|c}
				\hline
				\multirow{2}{*}{\bf FIGURE}    & \multicolumn{2}{c|}{\multirow{2}{*}{\bf Phase portrait}} & \textbf{Correspondent in families}   \\
				\textbf{\cite{Reyn-Huang-1997}}								& 																	&										& $\QESA$, $\QESB$, or $\QESC$ \\ 
				\hline
				\multirow{14}{*}{11.1} &  \multirow{3}{*}{$\gamma>2$}      & \multicolumn{1}{|c|}{$\mu<-\frac{1}{4\gamma}$}  		&  $V_{5}-\QESB$ \\ \cline{3-4}
				&                                                          &  \multicolumn{1}{|c|}{$\mu=0$}  																	&  $2.4L_{5}-\QESA$ \\ \cline{3-4}
				&                                                          &  \multicolumn{1}{|c|}{$\mu>0$}  																	&  $V_{101}-\QESA$ \\ \cline{2-4}
				&  \multirow{5}{*}{$1<\gamma<2$}  &  \multicolumn{1}{|c|}{$\mu=\text{u.s.c.}$}                   &  $7S_{1}-\QESB$ \\ \cline{3-4}
				&                                                         &  \multicolumn{1}{|c|}{$\mu<\text{u.s.c.}$}                   &  $V_{14}-\QESB$ \\ \cline{3-4}
				&                                                         &  \multicolumn{1}{|c|}{$\text{u.s.c.}<\mu<-\frac{1}{4\gamma}$}  &  $V_{12}-\QESB$ \\  \cline{3-4}																			
				&                                                          &  \multicolumn{1}{|c|}{$\mu=0$}  																	&  $2.4L_{7}-\QESA$ \\ \cline{3-4}
				&                                                          &  \multicolumn{1}{|c|}{$\mu>0$}  																	&  $V_{188}-\QESA$ \\ \cline{2-4}
				&  \multirow{3}{*}{$0<\gamma<1$}  & \multicolumn{1}{|c|}{$\mu<-\frac{1}{4\gamma}$}  		  &  $V_{20}-\QESB$ \\ \cline{3-4}
				&                                                          &  \multicolumn{1}{|c|}{$\mu=0$}  																	&  $2S_{35}-\QESA$ \\ \cline{3-4}
				&                                                          &  \multicolumn{1}{|c|}{$\mu>0$}  																	&  $V_{240}-\QESA$ \\ \cline{2-4}
				&  \multirow{3}{*}{$\gamma<0$}      & \multicolumn{1}{|c|}{$\mu<0$}  		                                &  $V_{1}-\QESA$ \\ \cline{3-4}
				&                                                          &  \multicolumn{1}{|c|}{$\mu=0$}  																	&  $2S_{6}-\QESA$ \\ \cline{3-4}
				&                                                          &  \multicolumn{1}{|c|}{$\mu>-\frac{1}{4\gamma}$}  		&  $V_{1}-\QESB$ \\ 
				\hline
				\multirow{14}{*}{11.2}&  \multicolumn{2}{c|}{$d>d_1(m;g)$}                              &  $V_{5}-\QESB$ \\ \cline{2-4}
				&   \multicolumn{2}{c|}{$d=d_1(m;g)$}                             &  $4S_{2}-\QESB$ \\ \cline{2-4}
				&   \multicolumn{2}{c|}{$-1<d<d_1(m;g)$}                       &  $V_{9}-\QESB$ \\ \cline{2-4}
				&   \multicolumn{2}{c|}{$d^\varepsilon-1$}                       &  $V_{5}-\QESB$ \\ \cline{2-4}
				&   \multicolumn{2}{c|}{$d_2(0;g)<d<d_{-}(0;g)$}          &  $2S_{18}-\QESA$ \\ \cline{2-4}
				&   \multicolumn{2}{c|}{$0<d<d_2(0;g)$}					               &  $2S_{17}-\QESA$ \\ \cline{2-4}
				&   \multicolumn{2}{c|}{$d_1(0;g)<d<0$}                        &  $2S_{13}-\QESA$ \\ \cline{2-4}
				&   \multicolumn{2}{c|}{$-1<d<d_1(0;g)$}                			&  $2S_{12}-\QESA$ \\ \cline{2-4}
				&   \multicolumn{2}{c|}{$d<-1$} 	                    	              &  $2S_{11}-\QESA$ \\ \cline{2-4}
				&   \multicolumn{2}{c|}{$d^3d_{-}(0;g)=\frac{2}{g}$}   &  $2S_{20}-\QESA$ \\ \cline{2-4}
				&   \multicolumn{2}{c|}{$d=d_2(0;g)$}                            &  $2.7L_{1}-\QESA$ \\ \cline{2-4}
				&   \multicolumn{2}{c|}{$d=0$}                                       &  $2.4L_{5}-\QESA$ \\ \cline{2-4}
				&   \multicolumn{2}{c|}{$d=d_1(0;g)$}                            &  $2.4L_{4}-\QESA$ \\ \cline{2-4}
				&   \multicolumn{2}{c|}{$d=-1$}                                     &  $2.3L_{7}-\QESA$ \\ 
				\hline
			\end{tabular}
		}
	\end{center}
\end{table}

\begin{table}\caption{\small Continuation of Table~\ref{tab:cor-RH-we-pt3}}\label{tab:cor-RH-we-pt4}
	\begin{center}
		{\def\arraystretch{1.5}
			\begin{tabular}{c|cc|c}
				\hline
				\multirow{2}{*}{\bf FIGURE}    & \multicolumn{2}{c|}{\multirow{2}{*}{\bf Phase portrait}} & \textbf{Correspondent in families}   \\
				\textbf{\cite{Reyn-Huang-1997}}								& 																	&										& $\QESA$, $\QESB$, or $\QESC$ \\ 
				\hline
				\multirow{12}{*}{11.3c}&  \multicolumn{2}{c|}{$\delta>\delta_4(\mu;\gamma)$}                                        &  $V_{12}-\QESB$ \\ \cline{2-4}
				&   \multicolumn{2}{c|}{$\delta=\delta_4(\mu;\gamma)$}                                       &  $7S_{1}-\QESB$ \\ \cline{2-4}
				&   \multicolumn{2}{c|}{$\delta_5(\mu;\gamma)<\delta<\delta_4(\mu;\gamma)$}    &  $V_{14}-\QESB$ \\ \cline{2-4}
				&   \multicolumn{2}{c|}{$\delta=\delta_5(\mu;\gamma)$}                                       &  $4S_{3}-\QESB$ \\ \cline{2-4}
				&   \multicolumn{2}{c|}{$\delta_6(\mu;\gamma)<\delta<\delta_5(\mu;\gamma)$}    &  $V_{15}-\QESB$ \\ \cline{2-4}
				&   \multicolumn{2}{c|}{$\delta=\delta_6(\mu;\gamma)$}					                             &  $8S_{4}-\QESB$ \\ \cline{2-4}
				&   \multicolumn{2}{c|}{$\delta^\ast<\delta<\delta_6(\mu;\gamma)$}                     &  $V_{16}-\QESB$ \\ \cline{2-4}
				&   \multicolumn{2}{c|}{$\delta_7(\mu;\gamma)<\delta\leq-1$}                			         &  $V_{14}-\QESB$ \\ \cline{2-4}
				&   \multicolumn{2}{c|}{$\delta=\delta_7(\mu;\gamma)<-1$} 	                    	         &  $7S_{1}-\QESB$ \\ \cline{2-4}
				&   \multicolumn{2}{c|}{$\delta=\delta_7(\mu;\gamma)>-1$}                                 &  $7S_{2}-\QESB$ \\ \cline{2-4}
				&   \multicolumn{2}{c|}{$-1<\delta<\delta_7(\mu;\gamma)$}                                 &  $V_{17}-\QESB$ \\ \cline{2-4}
				&   \multicolumn{2}{c|}{$\delta<\delta^{\ast\ast}$}                                               &  $V_{12}-\QESB$ \\ 
				\hline
				\multirow{14}{*}{11.3d}&  \multicolumn{2}{c|}{$\delta\geq\delta_{-}(0;\gamma)$}                                                 &  $2S_{32}-\QESA$ \\ \cline{2-4}
				&  \multicolumn{2}{c|}{$\delta_{3}(0;\gamma)<\delta<\delta_{-}(0;\gamma)$}                  &  $2S_{30}-\QESA$ \\ \cline{2-4}
				&  \multicolumn{2}{c|}{$\delta=\delta_{3}(0;\gamma)$}                                                     &  $2.7L_{2}-\QESA$ \\ \cline{2-4}
				&  \multicolumn{2}{c|}{$\delta_{2}(0;\gamma)<\delta<\delta_{3}(0;\gamma)$}                  &  $2S_{29}-\QESA$ \\ \cline{2-4}
				&  \multicolumn{2}{c|}{$\delta=\delta_{2}(0;\gamma)$}                                                     &  $2.8L_{9}-\QESA$ \\ \cline{2-4}
				&  \multicolumn{2}{c|}{$0<\delta<\delta_{2}(0;\gamma)$}                                                 &  $2S_{28}-\QESA$ \\ \cline{2-4}
				&  \multicolumn{2}{c|}{$\delta=0$}                                                                                     &  $2.4L_{7}-\QESA$ \\ \cline{2-4}
				&  \multicolumn{2}{c|}{$\delta_{1}^\ast(0;\gamma)<\delta<0$}                                          &  $2S_{26}-\QESA$ \\ \cline{2-4}
				&  \multicolumn{2}{c|}{$\delta=\delta_{1}^\ast(0;\gamma)$}                                              &  $2.4L_{6}-\QESA$ \\ \cline{2-4}
				&  \multicolumn{2}{c|}{$\delta_{2}^\ast(0;\gamma)<\delta<\delta_{1}^\ast(0;\gamma)$}   &  $2S_{25}-\QESA$ \\ \cline{2-4}
				&  \multicolumn{2}{c|}{$\delta=\delta_{2}^\ast(0;\gamma)$}                                              &  $2.8L_{8}-\QESA$ \\ \cline{2-4}
				&  \multicolumn{2}{c|}{$-1<\delta<\delta_{2}^\ast(0;\gamma)$}                                        &  $2S_{24}-\QESA$ \\ \cline{2-4}
				&  \multicolumn{2}{c|}{$\delta=-1$}                                                                                   &  $2.3L_{9}-\QESA$ \\ \cline{2-4}
				&  \multicolumn{2}{c|}{$\delta<-1$}                                                                                   &  $2S_{23}-\QESA$ \\ 
				\hline
			\end{tabular}
		}
	\end{center}
\end{table}

\begin{table}\caption{\small Continuation of Table~\ref{tab:cor-RH-we-pt4}}\label{tab:cor-RH-we-pt5}
	\begin{center}
		{\def\arraystretch{1.45}
			\begin{tabular}{c|cc|c}
				\hline
				\multirow{2}{*}{\bf FIGURE}    & \multicolumn{2}{c|}{\multirow{2}{*}{\bf Phase portrait}} & \textbf{Correspondent in families}   \\
				\textbf{\cite{Reyn-Huang-1997}}								& 																	&										& $\QESA$, $\QESB$, or $\QESC$ \\ 
				\hline
				\multirow{14}{*}{11.4}&  \multicolumn{2}{c|}{$\delta>\delta_1(\mu;\gamma)$}                                       &  $V_{20}-\QESB$ \\ \cline{2-4}
				&   \multicolumn{2}{c|}{$\delta=\delta_1(\mu;\gamma)$}                                      &  $8S_{5}-\QESB$ \\ \cline{2-4}
				&   \multicolumn{2}{c|}{$-1<\delta<\delta_1(\mu;\gamma)$}                                &  $V_{24}-\QESB$ \\ \cline{2-4}
				&   \multicolumn{2}{c|}{$\delta\leq-1$}                                                                &  $V_{20}-\QESB$ \\ \cline{2-4}
				&   \multicolumn{2}{c|}{$\delta\geq\delta_{-}(0;\gamma)=\frac{2}{\gamma}$}    &  $2S_{42}-\QESA$ \\ \cline{2-4}
				&   \multicolumn{2}{c|}{$\delta=\delta_2(0;\gamma)$}					                                &  $2.7L_{3}-\QESA$ \\ \cline{2-4}
				&   \multicolumn{2}{c|}{$\delta=\delta_1(0;\gamma)$}                                         &  $2.8L_{11}-\QESA$ \\ \cline{2-4}
				&   \multicolumn{2}{c|}{$\delta=\delta_1^\ast(0;\gamma)$}                			           &  $2.8L_{10}-\QESA$ \\ \cline{2-4}
				&   \multicolumn{2}{c|}{$\delta=-1$} 	                    	                                         &  $2.3L_{11}-\QESA$ \\ \cline{2-4}
				&   \multicolumn{2}{c|}{$\delta_2(0;\gamma)<\delta<\delta_{-}(0;\gamma)$}      &  $2S_{40}-\QESA$ \\ \cline{2-4}
				&   \multicolumn{2}{c|}{$\delta_1(0;\gamma)<\delta<\delta_{2}(0;\gamma)$}      &  $2S_{39}-\QESA$ \\ \cline{2-4}
				&   \multicolumn{2}{c|}{$\delta_1^\ast(0;\gamma)<\delta<\delta_{1}(0;\gamma)$}& $2S_{35}-\QESA$ \\ \cline{2-4}
				&   \multicolumn{2}{c|}{$-1<\delta<\delta_{1}^\ast(0;\gamma)$}                          & $2S_{34}-\QESA$ \\ \cline{2-4}
				&   \multicolumn{2}{c|}{$\delta<-1$}                                                                    &  $2S_{33}-\QESA$ \\ 
				\hline
				\multirow{19}{*}{11.5c}&                                                                       & \multicolumn{1}{|c|}{$\delta>\delta_1(\mu;\gamma)$}                                   &  $V_{1}-\QESA$ \\ \cline{3-4}
				&                                                                         & \multicolumn{1}{|c|}{$\delta=\delta_1(\mu;\gamma)$}                                   &  $8S_{7}-\QESA$ \\ \cline{3-4}
				&                                                                         & \multicolumn{1}{|c|}{$\delta_2(\mu;\gamma)<\delta<\delta_1(\mu;\gamma)$}&  $V_{12}-\QESA$ \\ \cline{3-4}
				&  \multirow{1}{*}{$\mu<0,$}                           & \multicolumn{1}{|c|}{$\delta=\delta_2(\mu;\gamma)$}                                   &  $7S_{4}-\QESA$ \\ \cline{3-4}
				&  \multirow{1}{*}{$-2<\gamma<0$}                & \multicolumn{1}{|c|}{$-1<\delta<\delta_2(\mu;\gamma)$}                             &  $V_{11}-\QESA$ \\ \cline{3-4}
				&                                                                         & \multicolumn{1}{|c|}{$\delta_3(\mu;\gamma)<\delta\leq-1$}                          &  $V_{9}-\QESA$ \\ \cline{3-4}
				&  																											                  & \multicolumn{1}{|c|}{$\delta=\delta_3(\mu;\gamma)$}                                  &  $4S_{5}-\QESA$ \\ \cline{3-4}
				&                                                                         & \multicolumn{1}{|c|}{$\delta>\delta_3(\mu;\gamma)$}                                  &  $V_{1}-\QESA$ \\ \cline{2-4}
				&   \multicolumn{2}{c|}{$\mu<0, \gamma=-2, \delta=-1$}                   																																													& $3.7L_{1}-\QESA$ \\ \cline{2-4}																
				&  \multirow{2}{*}{$\mu<0, \gamma>-2$}        &  \multicolumn{1}{|c|}{$\delta_2(\mu;\gamma)<\delta<-1$}                            &  $V_{66}-\QESA$ \\ \cline{3-4}
				&                                                                        &  \multicolumn{1}{|c|}{$\delta=\delta_2(\mu;\gamma)$}                                  &  $7S_{4}-\QESA$ \\ \cline{2-4}
				&  \multirow{7}{*}{$\mu=0, \gamma<0$}        &  \multicolumn{1}{|c|}{$\delta>\delta_1(0;\gamma)$}                                      &  $2S_{6}-\QESA$ \\ \cline{3-4}
				&                                                                       &  \multicolumn{1}{|c|}{$\delta=\delta_1(0;\gamma)$}                                      &  $2.8L_{2}-\QESA$ \\ \cline{3-4}
				&                                                                       &  \multicolumn{1}{|c|}{$-1<\delta<\delta_1(0;\gamma)$}                                &  $2S_{5}-\QESA$ \\ \cline{3-4}
				&                                                                       &  \multicolumn{1}{|c|}{$\delta=-1$}                                                                &  $2.3L_{2}-\QESA$ \\ \cline{3-4}
				&                                                                       &  \multicolumn{1}{|c|}{$\delta_3(0;\gamma)<\delta<-1$}                                &  $2S_{4}-\QESA$ \\ \cline{3-4}
				&                                                                       &  \multicolumn{1}{|c|}{$\delta=\delta_3(0;\gamma)$}                                     &  $2.4L_{1}-\QESA$ \\ \cline{3-4}
				&                                                                       &  \multicolumn{1}{|c|}{$\delta>\delta_3(0;\gamma)$}                                     &  $2S_{1}-\QESA$ \\ \cline{2-4}
				&   \multicolumn{2}{c|}{$\mu>-\frac{1}{4\gamma}, \gamma<0$}   																																													& $V_{1}-\QESB$ \\ 
				\hline
			\end{tabular}
		}
	\end{center}
\end{table}

Therefore, as we proved that the phase portraits we obtained are topologically distinct
we conclude that, from the 143 phase portraits from the mentioned paper, the number 
of topologically distinct phase portraits is indeed 94.

From the analysis of the phase portraits we obtained in the closures $\overline{\QESA}$, 
$\overline{\QESB}$, and $\overline{\QESC}$, we observe the existence of 29 phase portraits 
which were not obtained by those authors. One example is our phase portrait $7S_{1}$ in 
$\overline{\QESA}$ which was not found in \cite{Reyn-Huang-1997}.

Another relevant fact we want to add in this section is the following one. In \cite{Artes-Mota-Rezende-2021c} we presented a list of some small prints and incompatibilities found in \cite{Artes-Rezende-Oliveira-2015}. In addition to that list, we point out that in equation (7), corresponding to slices $n_{62}$ up to $n_{69}$, instead of the value $81/40$, the correct is $81/400$. This correction must be made in Figures 89 up to 96 and in Tables 33 up to 37 from that paper.


\begin{thebibliography}{99}

\bibitem{progP4} 
\textsc{Art\'es, J.C., Dumortier, F., Herssens, C., Llibre, J. \& de Maesschalck, P.} 
{Computer program P4 to study phase portraits of planar polynomial differential equations.} 
Available at: \url{http://mat.uab.es/~artes/p4/p4.htm}, 2005.

\bibitem{Artes-Kooij-Llibre-1998}
\textsc{Art\'es, J.C., Kooij, R. \& Llibre, J.}
{Structurally stable quadratic vector fields.} 
\textit{Memoires Amer. Math. Soc.} 
\textbf{134 (639)}, {1998}, {108pp.}

\bibitem{Artes-Llibre-Rezende-2018} 
\textsc{Art\'es, J.C., Llibre, J. \& Rezende, A.C.}
{Structurally unstable quadratic vector fields of codimension one.}
1. ed. 
\textit{Birkh\"auser}. v.1. 2018. 267pp.

\bibitem{Artes-Llibre-Schlomiuk-2006}
\textsc{Art\'es, J.C., Llibre, J. \& Schlomiuk, D.} 
{The geometry of quadratic differential systems with a weak focus of second order. }
\textit{Internat. J. Bifur. Chaos Appl. Sci. Engrg.} 
\textbf{16}, {2006}, {3127--3194.}

\bibitem{Artes-Llibre-Schlomiuk-Vulpe-2020a}
\textsc{Art\'es, J.C., Llibre, J., Schlomiuk, D. \& Vulpe, N.}
{Global topological configurations of singularities for the whole family of quadratic differential systems.} 
\textit{Qual. Theor. Dyn. Syst.}
\textbf{19}, {2020a}, {32pp.}

\bibitem{Artes-Llibre-Schlomiuk-Vulpe-2021a}
\textsc{Art\'es, J.C., Llibre, J., Schlomiuk, D. \& Vulpe, N.}
{Geometric configurations of singularities of planar polynomial differential systems -- A global classification in the quadratic case.} 
1. ed. 
\textit{Birkh\"auser Basel}. v.1. 2021a. 701 pp.

\bibitem{Artes-Llibre-Vulpe-2008}
\textsc{Art\'es, J.C., Llibre, J. \& Vulpe, N.} 
{Singular points of quadratic systems: a complete classification in the coefficient space $\mathbb R^{12}$.}  
\textit{Internat. J. Bifur. Chaos Appl. Sci. Engrg.} 
\textbf{18}, {2008}, {313--362.}

\bibitem{Artes-Mota-Rezende-2021b}
\textsc{Art\'{e}s, J.C., Mota, M.C. \& Rezende, A.C.} 
{Quadratic differential systems with a finite saddle--node and an infinite saddle--node $(1,1)SN$ - (A).} 
\textit{Internat. J. Bifur. Chaos Appl. Sci. Engrg.} 
\textbf{31(2)}, {2021b}, {2150026 -- 24pp.}

\bibitem{Artes-Mota-Rezende-2021c}
\textsc{Art\'{e}s, J.C., Mota, M.C. \& Rezende, A.C.} 
{Quadratic differential systems with a finite saddle--node and an infinite saddle--node $(1,1)SN$ - (B).} 
\textit{Internat. J. Bifur. Chaos Appl. Sci. Engrg.} 
\textbf{31(9)}, {2021c}, {2130026 -- 110pp.}

\bibitem{Artes-Mota-Rezende-2021d}
\textsc{Art\'{e}s, J.C., Mota, M.C. \& Rezende, A.C.} 
{Structurally unstable quadratic vector fields of codimension two: families possessing a finite saddle-node and an infinite saddle-node.} 
\textit{Electron. J. Qual. Theo.}
\textbf{35}, {2021d}, {1--89.}

\bibitem{Artes-Oliveira-Rezende-2020b}
\textsc{Art\'{e}s, J.C., Oliveira, R.D.S. \& Rezende, A.C.} 
{Structurally unstable quadratic vector fields of codimension two: families possessing either a cusp point or two finite saddle--nodes.} 
\textit{J. Dyn. Diff. Equat.} 
{2020b}, {43pp.}

\bibitem{Artes-Rezende-Oliveira-2013b}
\textsc{Art\'{e}s, J.C., Rezende, A.C. \& Oliveira, R.D.S.} 
{Global phase portraits of quadratic polynomial differential systems with a semi--elemental triple node.} 
\textit{Internat. J. Bifur. Chaos Appl. Sci. Engrg.}
\textbf{23}, {2013}, {21pp.}

\bibitem{Artes-Rezende-Oliveira-2014}
\textsc{Art\'{e}s, J.C., Rezende, A.C. \& Oliveira, R.D.S.} 
{The geometry of quadratic polynomial differential systems with a finite and an infinite saddle--node $(A,B)$.} 
\textit{Internat. J. Bifur. Chaos Appl. Sci. Engrg.} 
\textbf{24}, {2014}, {30pp.}

\bibitem{Artes-Rezende-Oliveira-2015}
\textsc{Art\'{e}s, J.C., Rezende, A.C. \& Oliveira, R.D.S.} 
{The geometry of quadratic polynomial differential systems with a finite and an infinite saddle--node $C$.} 
\textit{Internat. J. Bifur. Chaos Appl. Sci. Engrg.} 
\textbf{25}, {2015}, {111pp.}

\bibitem{Dumortier-Llibre-Artes-2006}
\textsc{Dumortier, F., Llibre, J. \& Art\'es, J.C.} 
{Qualitative Theory of Planar Differential Systems.}  
\textit{Universitext, Springer--Verlag, New York--Berlin}. 2006.

\bibitem{Dumortier-Roussarie-Rousseau-1994}
\textsc{Dumortier, F., Roussarie, R. \& Rousseau, C.} 
{Hilbert's 16th problem for quadratic vector fields.} 
\textit{J. Differential Equations.} 
\textbf{110}, {1994}, {66--133.}

\bibitem{Reyn-Huang-1997}
\textsc{Reyn, J.W. \& Huang, X.} 
{Phase portraits of quadratics systems with finite multiplicity three and a degenerate critical point at infinity.} 
\textit{Rocky MT J Math.} 
\textbf{27(3)}, {1997}, {929--978.}

\bibitem{Schlomiuk-Vulpe-2004}
\textsc{Schlomiuk, D. \& Vulpe, N.} 
{Planar quadratic vector fields with invariant lines of total multiplicity at least five.} 
\textit{Qualitative Theory of Dynamical Systems.} 
\textbf{5}, {2004}, {135--194.}

\bibitem{Schlomiuk-Vulpe-2005}
\textsc{Schlomiuk, D. \& Vulpe, N.} 
{Geometry of quadratic differential systems in the neighborhood of the infinity.} 
\textit{J. Differential Equations.} 
\textbf{215}, {2005}, {357--400.}

\bibitem{Vulpe-2011}
\textsc{Vulpe, N.} 
{Characterization of the finite weak singularities of quadratic systems via invariant theory.} 
\textit{Nonlinear Anal.} 
\textbf{74}, {2011}, {6553--6582.}

\end{thebibliography}
\end{document}